\documentclass{amsart}
\usepackage{a4wide,amssymb,color,hyperref}
\usepackage{graphicx} % Required for inserting images

\theoremstyle{plain}

\theoremstyle{definition}

\theoremstyle{remark}

\newtheorem{example}{Example}[section]

\title{Point-transitive and 1-rotational unitals of order 5}
\author{Ivan Hetman, Taras Banakh, Alex Ravsky}
\date{April 2025}

\begin{document}

\begin{abstract} In this paper we introduce enumeration of unitals of order $5$, which are also Steiner systems $S(2,6,126)$, where automorphism group acts transitively and effectively on points or fixes one point.
\end{abstract}

\maketitle

\section{Introduction}
When we started studying design theory, we found that one important objects are called {\bf unitals} and they correspond to $S(2,k+1,k^3+1)$ Steiner system, where $k$ is called the order of an unital. There are lots of examples of unitals of order $3$ and $4$ on websites \cite{Krc} and \cite{NagyWeb}, but it was strange for us that there are only $2$ examples of order $5$ known in \cite{HoCD} and it seems that nodoby tried to generate more. When first new unitals were obtained in \cite{Het} and consequently more in \cite{Het1}, it appeared that there is enough of unitals of order $5$ unknown for some reason. In this paper we present new additional unitals of order $5$ with point-transitive automorphism groups and so-called 1-rotational unitals of order $5$.

In 2023 S.~Stoichev and M.~Gezek \cite{StGez} introduced first major extension with $938$ unitals of order $5$ found as subdesigns in finite projective planes of order $25$. In this paper we are adding additional unitals, but constructed algebraically. Nevertheless, authomorphism groups orders of subdesign unitals are mostly lower than $100$ while automorphism group orders of algebraically generated unitals are multiples of $125$ or $126$. So, the only candidates for intersection are $3$ unitals from Desarguesian projective plane of order $25$. First of them is classical unital, second and third have order of their automorphism group $1000$ and $2000$ respectively, so they could also be present in our lists.

Let's recall the hyperbolic frequency fingerprint introduced in \cite{Het1}. The hyperbolic frequency is calculated as follows. We go through all non-collinear (not on same block) triples $oxy$ and point $p$ on the line (block) $xy$ which is distinct from $x$ and $y$ (so we check all such quadruples). Then we calculate the number of points $u$ on the line $oy$ such that line $pu$ doesn't intersect $ox$. Then we sum frequencies and obtain something like $\{1=31968, 2=534132, 3=2203128, 4=2358972\}$.

According to \cite{GAP4}, there are $16$ groups of order $126$ and $5$ groups of order $125$. We will generate Steiner systems of two types:
\begin{itemize}
\item designs where group of order $126$ acts on $126$ points transitively and effectively;
\item designs where group of order $125$ has two orbits of size $125$ and $1$ on $126$ points. The latter point is the fixed point of the group action.
\end{itemize}
Due to the fact that we are considering non-commutative groups, group action on blocks to obtain full design, will be always from the left. For calculation purposes we needed to transform group from GAP representation to Cayley table. The simplest way we found was to achieve this using LOOPS \cite{LOOPS} package for GAP applying the procedure {\tt CayleyTable(IntoLoop(SmallGroup(126,1)));} where we obtain two-dimensional array which we then convert to 0-based Cayley table by subtracting $1$. Then when we need to obtain all blocks from difference family element, we just use the generated Cayley table.

Before enumeration lists of designs, here is summarized table with number of non-isomorphic unitals corresponding to specific group.
\begin{center}
\begin{tabular}{||c c c c||} 
 \hline
 GAP ID & Structure & Number of unitals & Comments \\ [0.5ex] 
 \hline\hline
 {\tt SmallGroup(125,1)} & $\mathbb Z_{125}$ & 8 & - \\ 
 {\tt SmallGroup(125,2)} & $\mathbb Z_5 \times \mathbb Z_{25}$ & 32 & - \\
 {\tt SmallGroup(125,3)} & $(\mathbb Z_5 \times \mathbb Z_5) \rtimes \mathbb Z_5$ & 20 & Classical unital \\
 {\tt SmallGroup(125,4)} & $\mathbb Z_{25} \rtimes \mathbb Z_5$ & 29 & - \\
 {\tt SmallGroup(125,5)} & $\mathbb Z_5 \times \mathbb Z_5 \times \mathbb Z_5$ & 8 & \cite{HoCD} resolvable unital \\
 \hline\hline
 {\tt SmallGroup(126,1)} & {\tt C7 : C18} & 3 & - \\ 
 {\tt SmallGroup(126,2)} & {\tt C2 x (C7 : C9)} & 33 & - \\
 {\tt SmallGroup(126,3)} & {\tt C7 x D18} & 1 & - \\
 {\tt SmallGroup(126,4)} & {\tt C9 x D14} & 0 & - \\
 {\tt SmallGroup(126,5)} & {\tt D126} & 0 & - \\
 {\tt SmallGroup(126,6)} & $\mathbb Z_{126}$ & 64 & \cite{Het} \\ 
 {\tt SmallGroup(126,7)} & {\tt C3 x (C7 : C6)} & 2 & - \\
 {\tt SmallGroup(126,8)} & {\tt S3 x (C7 : C3)} & 129 & - \\
 {\tt SmallGroup(126,9)} & {\tt C7 : (C3 x S3)} & 0 & - \\
 {\tt SmallGroup(126,10)} & {\tt C6 x (C7 : C3)} & 900 & - \\
 {\tt SmallGroup(126,11)} & {\tt C3 x C3 x D14} & 0 & - \\ 
 {\tt SmallGroup(126,12)} & {\tt C21 x S3} & 74 & - \\
 {\tt SmallGroup(126,13)} & {\tt C3 x D42} & 0 & - \\
 {\tt SmallGroup(126,14)} & {\tt C7 x ((C3 x C3) : C2} & 0 & - \\
 {\tt SmallGroup(126,15)} & {\tt (C3 x C3) : D14} & 0 & - \\
 {\tt SmallGroup(126,16)} & $\mathbb Z_2 \times \mathbb Z_3 \times \mathbb Z_3 \times \mathbb Z_7$ & 8 & \cite{Het1} \\
 \hline
\end{tabular}
\end{center}

\section{1-rotational designs enumeration}
Results will be enumerated in form fingerprint-difference family. Fixed point will be denoted as $\infty$. Two known designs in \cite{HoCD} are marked bold.
\begin{example} $\mathbb Z_{125}$
\begin{enumerate}
\item \{1=25000, 2=580500, 3=3042000, 4=3912500\} [[0, 1, 3, 15, 47, 74], [0, 4, 9, 20, 65, 103], [0, 6, 40, 88, 95, 112], [0, 8, 18, 41, 76, 104], [0, 25, 50, 75, 100, $\infty$]]
\item \{1=36000, 2=615750, 3=2986500, 4=3921750\} [[0, 1, 3, 15, 47, 74], [0, 4, 26, 64, 109, 120], [0, 6, 40, 88, 95, 112], [0, 8, 18, 41, 76, 104], [0, 25, 50, 75, 100, $\infty$]]
\item \{1=49000, 2=708750, 3=2968500, 4=3833750\} [[0, 1, 3, 15, 47, 74], [0, 4, 26, 64, 109, 120], [0, 6, 19, 36, 43, 91], [0, 8, 18, 41, 76, 104], [0, 25, 50, 75, 100, $\infty$]]
\item \{1=54000, 2=662250, 3=3028500, 4=3815250\} [[0, 1, 3, 15, 47, 74], [0, 4, 9, 20, 65, 103], [0, 6, 19, 36, 43, 91], [0, 8, 18, 41, 76, 104], [0, 25, 50, 75, 100, $\infty$]]
\item \{0=1250, 1=42000, 2=635250, 3=2987500, 4=3894000\} [[0, 1, 3, 15, 47, 74], [0, 4, 9, 20, 65, 103], [0, 6, 40, 88, 95, 112], [0, 8, 29, 57, 92, 115], [0, 25, 50, 75, 100, $\infty$]]
\item \{0=1250, 1=45000, 2=658500, 3=3016000, 4=3839250\} [[0, 1, 3, 15, 47, 74], [0, 4, 9, 20, 65, 103], [0, 6, 19, 36, 43, 91], [0, 8, 29, 57, 92, 115], [0, 25, 50, 75, 100, $\infty$]]
\item \{0=1250, 1=51000, 2=627750, 3=3014500, 4=3865500\} [[0, 1, 3, 15, 47, 74], [0, 4, 26, 64, 109, 120], [0, 6, 19, 36, 43, 91], [0, 8, 29, 57, 92, 115], [0, 25, 50, 75, 100, $\infty$]]
\item \{0=2500, 1=47000, 2=682500, 3=3020000, 4=3808000\} [[0, 1, 3, 15, 47, 74], [0, 4, 26, 64, 109, 120], [0, 6, 40, 88, 95, 112], [0, 8, 29, 57, 92, 115], [0, 25, 50, 75, 100, $\infty$]]
\end{enumerate}
\end{example}
\begin{example} $\mathbb Z_5 \times \mathbb Z_{25}$
\begin{enumerate}
\item \{1=20000, 2=578250, 3=3064500, 4=3897250\} [[(0, 0), (0, 1), (0, 3), (1, 0), (2, 13), (3, 6)], [(0, 0), (0, 4), (0, 16), (1, 11), (2, 5), (4, 13)], [(0, 0), (0, 5), (0, 10), (0, 15), (0, 20), $\infty$], [(0, 0), (0, 6), (1, 10), (1, 21), (2, 2), (4, 17)], [(0, 0), (0, 7), (2, 16), (2, 24), (3, 0), (4, 2)]]
\item \{1=21000, 2=624000, 3=3054000, 4=3861000\} [[(0, 0), (0, 1), (0, 3), (1, 0), (2, 2), (4, 9)], [(0, 0), (0, 4), (0, 11), (2, 4), (3, 17), (4, 18)], [(0, 0), (0, 5), (0, 17), (1, 20), (2, 3), (3, 8)], [(0, 0), (0, 6), (0, 16), (1, 12), (2, 21), (4, 2)], [(0, 0), (1, 10), (2, 20), (3, 5), (4, 15), $\infty$]]
\item \{1=26000, 2=609000, 3=3075000, 4=3850000\} [[(0, 0), (0, 1), (0, 3), (0, 18), (1, 1), (2, 18)], [(0, 0), (0, 4), (0, 13), (1, 7), (2, 16), (3, 2)], [(0, 0), (0, 5), (1, 10), (1, 21), (2, 9), (3, 6)], [(0, 0), (0, 6), (1, 20), (2, 13), (3, 17), (4, 19)], [(0, 0), (1, 15), (2, 5), (3, 20), (4, 10), $\infty$]]
\item \{1=28000, 2=589500, 3=3054000, 4=3888500\} [[(0, 0), (0, 1), (0, 3), (1, 0), (1, 19), (2, 3)], [(0, 0), (0, 4), (0, 11), (1, 15), (3, 7), (3, 24)], [(0, 0), (0, 5), (0, 10), (0, 15), (0, 20), $\infty$], [(0, 0), (0, 9), (1, 5), (2, 7), (3, 19), (4, 17)], [(0, 0), (0, 12), (1, 1), (2, 11), (3, 17), (4, 5)]]
\item \{1=31000, 2=607500, 3=3057000, 4=3864500\} [[(0, 0), (0, 1), (0, 3), (1, 0), (2, 13), (3, 6)], [(0, 0), (0, 4), (0, 13), (1, 16), (3, 24), (4, 18)], [(0, 0), (0, 5), (0, 10), (0, 15), (0, 20), $\infty$], [(0, 0), (0, 6), (1, 10), (1, 21), (2, 2), (4, 17)], [(0, 0), (0, 7), (1, 5), (2, 7), (3, 8), (3, 16)]]
\item \{1=32000, 2=582000, 3=3084000, 4=3862000\} [[(0, 0), (0, 1), (0, 3), (1, 0), (1, 19), (2, 3)], [(0, 0), (0, 4), (0, 11), (1, 15), (3, 7), (3, 24)], [(0, 0), (0, 5), (0, 10), (0, 15), (0, 20), $\infty$], [(0, 0), (0, 9), (1, 5), (2, 7), (3, 19), (4, 17)], [(0, 0), (0, 12), (1, 7), (2, 20), (3, 1), (4, 11)]]
\item \{1=32000, 2=601500, 3=3033000, 4=3893500\} [[(0, 0), (0, 1), (0, 3), (0, 9), (1, 0), (3, 19)], [(0, 0), (0, 4), (1, 23), (2, 5), (3, 9), (4, 23)], [(0, 0), (0, 5), (0, 10), (0, 15), (0, 20), $\infty$], [(0, 0), (0, 7), (1, 3), (1, 17), (2, 4), (3, 15)], [(0, 0), (0, 12), (1, 20), (2, 8), (3, 23), (4, 7)]]
\item \{1=32000, 2=648750, 3=3004500, 4=3874750\} [[(0, 0), (0, 1), (0, 3), (1, 0), (1, 19), (2, 3)], [(0, 0), (0, 4), (0, 11), (1, 15), (3, 7), (3, 24)], [(0, 0), (0, 5), (0, 10), (0, 15), (0, 20), $\infty$], [(0, 0), (0, 9), (1, 17), (2, 15), (3, 2), (4, 4)], [(0, 0), (0, 12), (1, 7), (2, 20), (3, 1), (4, 11)]]
\item \{1=33000, 2=651000, 3=3006000, 4=3870000\} [[(0, 0), (0, 1), (0, 3), (1, 1), (3, 19), (4, 13)], [(0, 0), (0, 4), (0, 11), (2, 0), (3, 6), (4, 22)], [(0, 0), (0, 5), (0, 10), (0, 15), (0, 20), $\infty$], [(0, 0), (0, 6), (1, 2), (2, 10), (3, 5), (3, 14)], [(0, 0), (0, 8), (1, 18), (3, 17), (4, 3), (4, 16)]]
\item \{1=34000, 2=599250, 3=3019500, 4=3907250\} [[(0, 0), (0, 1), (0, 3), (1, 0), (2, 2), (4, 9)], [(0, 0), (0, 4), (0, 11), (2, 4), (3, 17), (4, 18)], [(0, 0), (0, 5), (0, 13), (2, 22), (3, 2), (4, 10)], [(0, 0), (0, 6), (0, 16), (1, 12), (2, 21), (4, 2)], [(0, 0), (1, 10), (2, 20), (3, 5), (4, 15), $\infty$]]
\item \{1=35000, 2=678750, 3=2992500, 4=3853750\} [[(0, 0), (0, 1), (0, 3), (1, 1), (3, 19), (4, 13)], [(0, 0), (0, 4), (0, 11), (2, 0), (3, 6), (4, 22)], [(0, 0), (0, 5), (0, 10), (0, 15), (0, 20), $\infty$], [(0, 0), (0, 6), (1, 2), (2, 10), (3, 5), (3, 14)], [(0, 0), (0, 8), (1, 5), (1, 17), (2, 16), (4, 15)]]
\item \{1=37000, 2=620250, 3=3139500, 4=3763250\} [[(0, 0), (0, 1), (0, 3), (1, 1), (1, 17), (3, 4)], [(0, 0), (0, 4), (1, 11), (1, 24), (2, 5), (4, 20)], [(0, 0), (0, 5), (0, 10), (0, 15), (0, 20), $\infty$], [(0, 0), (0, 6), (1, 18), (2, 8), (3, 11), (3, 18)], [(0, 0), (0, 8), (1, 4), (2, 6), (2, 17), (3, 14)]]
\item \{1=39000, 2=588000, 3=3021000, 4=3912000\} [[(0, 0), (0, 1), (0, 3), (1, 0), (2, 13), (3, 6)], [(0, 0), (0, 4), (0, 13), (1, 16), (3, 24), (4, 18)], [(0, 0), (0, 5), (0, 10), (0, 15), (0, 20), $\infty$], [(0, 0), (0, 6), (1, 10), (1, 21), (2, 2), (4, 17)], [(0, 0), (0, 7), (2, 16), (2, 24), (3, 0), (4, 2)]]
\item \{1=40000, 2=635250, 3=3049500, 4=3835250\} [[(0, 0), (0, 1), (0, 3), (0, 9), (1, 0), (3, 19)], [(0, 0), (0, 4), (1, 6), (2, 20), (3, 24), (4, 6)], [(0, 0), (0, 5), (0, 10), (0, 15), (0, 20), $\infty$], [(0, 0), (0, 7), (2, 17), (3, 3), (4, 4), (4, 15)], [(0, 0), (0, 12), (1, 5), (2, 14), (3, 4), (4, 17)]]
\item \{1=41000, 2=740250, 3=3028500, 4=3750250\} [[(0, 0), (0, 1), (0, 3), (1, 1), (3, 19), (4, 13)], [(0, 0), (0, 4), (0, 11), (2, 0), (3, 6), (4, 22)], [(0, 0), (0, 5), (0, 10), (0, 15), (0, 20), $\infty$], [(0, 0), (0, 6), (2, 1), (2, 17), (3, 21), (4, 4)], [(0, 0), (0, 8), (1, 5), (1, 17), (2, 16), (4, 15)]]
\item \{1=44000, 2=672000, 3=3009000, 4=3835000\} [[(0, 0), (0, 1), (0, 3), (1, 0), (2, 2), (4, 9)], [(0, 0), (0, 4), (0, 11), (2, 4), (3, 17), (4, 18)], [(0, 0), (0, 5), (0, 17), (1, 20), (2, 3), (3, 8)], [(0, 0), (0, 6), (0, 15), (1, 4), (3, 10), (4, 19)], [(0, 0), (1, 10), (2, 20), (3, 5), (4, 15), $\infty$]]
\item \{1=46000, 2=620250, 3=3016500, 4=3877250\} [[(0, 0), (0, 1), (0, 3), (0, 9), (1, 0), (3, 19)], [(0, 0), (0, 4), (1, 6), (2, 20), (3, 24), (4, 6)], [(0, 0), (0, 5), (0, 10), (0, 15), (0, 20), $\infty$], [(0, 0), (0, 7), (1, 3), (1, 17), (2, 4), (3, 15)], [(0, 0), (0, 12), (1, 5), (2, 14), (3, 4), (4, 17)]]
\item \{1=47000, 2=613500, 3=2997000, 4=3902500\} [[(0, 0), (0, 1), (0, 3), (0, 9), (1, 0), (3, 19)], [(0, 0), (0, 4), (1, 6), (2, 20), (3, 24), (4, 6)], [(0, 0), (0, 5), (0, 10), (0, 15), (0, 20), $\infty$], [(0, 0), (0, 7), (1, 3), (1, 17), (2, 4), (3, 15)], [(0, 0), (0, 12), (1, 20), (2, 8), (3, 23), (4, 7)]]
\item \{1=56000, 2=690750, 3=3037500, 4=3775750\} [[(0, 0), (0, 1), (0, 3), (1, 1), (1, 17), (3, 4)], [(0, 0), (0, 4), (1, 9), (3, 24), (4, 5), (4, 18)], [(0, 0), (0, 5), (0, 10), (0, 15), (0, 20), $\infty$], [(0, 0), (0, 6), (2, 13), (2, 20), (3, 23), (4, 13)], [(0, 0), (0, 8), (2, 19), (3, 2), (3, 16), (4, 4)]]
\item \{1=57000, 2=683250, 3=2995500, 4=3824250\} [[(0, 0), (0, 1), (0, 3), (1, 1), (1, 17), (3, 4)], [(0, 0), (0, 4), (1, 9), (3, 24), (4, 5), (4, 18)], [(0, 0), (0, 5), (0, 10), (0, 15), (0, 20), $\infty$], [(0, 0), (0, 6), (1, 18), (2, 8), (3, 11), (3, 18)], [(0, 0), (0, 8), (2, 19), (3, 2), (3, 16), (4, 4)]]
\item \{0=1250, 1=25000, 2=608250, 3=3029500, 4=3896000\} [[(0, 0), (0, 1), (0, 3), (1, 0), (2, 13), (3, 6)], [(0, 0), (0, 4), (0, 13), (1, 16), (3, 24), (4, 18)], [(0, 0), (0, 5), (0, 10), (0, 15), (0, 20), $\infty$], [(0, 0), (0, 6), (1, 14), (3, 4), (4, 10), (4, 21)], [(0, 0), (0, 7), (1, 5), (2, 7), (3, 8), (3, 16)]]
\item \{0=1250, 1=26000, 2=627750, 3=3038500, 4=3866500\} [[(0, 0), (0, 1), (0, 3), (0, 18), (1, 1), (2, 18)], [(0, 0), (0, 4), (0, 13), (1, 7), (2, 16), (3, 2)], [(0, 0), (0, 5), (2, 24), (3, 21), (4, 9), (4, 20)], [(0, 0), (0, 6), (1, 20), (2, 13), (3, 17), (4, 19)], [(0, 0), (1, 15), (2, 5), (3, 20), (4, 10), $\infty$]]
\item \{0=1250, 1=31000, 2=638250, 3=3080500, 4=3809000\} [[(0, 0), (0, 1), (0, 3), (1, 1), (1, 9), (4, 8)], [(0, 0), (0, 4), (2, 0), (2, 7), (2, 12), (3, 14)], [(0, 0), (0, 6), (0, 16), (1, 10), (2, 22), (3, 2)], [(0, 0), (0, 11), (1, 24), (2, 20), (3, 23), (4, 14)], [(0, 0), (1, 15), (2, 5), (3, 20), (4, 10), $\infty$]]
\item \{0=1250, 1=33000, 2=650250, 3=3080500, 4=3795000\} [[(0, 0), (0, 1), (0, 3), (1, 1), (1, 9), (4, 8)], [(0, 0), (0, 4), (2, 15), (3, 4), (3, 17), (3, 22)], [(0, 0), (0, 6), (0, 15), (2, 4), (3, 9), (4, 21)], [(0, 0), (0, 11), (1, 22), (2, 13), (3, 16), (4, 12)], [(0, 0), (1, 15), (2, 5), (3, 20), (4, 10), $\infty$]]
\item \{0=1250, 1=35000, 2=588750, 3=3086500, 4=3848500\} [[(0, 0), (0, 1), (0, 3), (1, 1), (1, 9), (4, 8)], [(0, 0), (0, 4), (2, 0), (2, 7), (2, 12), (3, 14)], [(0, 0), (0, 6), (0, 15), (2, 4), (3, 9), (4, 21)], [(0, 0), (0, 11), (1, 24), (2, 20), (3, 23), (4, 14)], [(0, 0), (1, 15), (2, 5), (3, 20), (4, 10), $\infty$]]
\item \{0=1250, 1=40000, 2=615000, 3=3019000, 4=3884750\} [[(0, 0), (0, 1), (0, 3), (1, 0), (1, 19), (2, 3)], [(0, 0), (0, 4), (0, 11), (1, 15), (3, 7), (3, 24)], [(0, 0), (0, 5), (0, 10), (0, 15), (0, 20), $\infty$], [(0, 0), (0, 9), (1, 17), (2, 15), (3, 2), (4, 4)], [(0, 0), (0, 12), (1, 1), (2, 11), (3, 17), (4, 5)]]
\item \{0=1250, 1=46000, 2=616500, 3=3025000, 4=3871250\} [[(0, 0), (0, 1), (0, 3), (1, 1), (1, 9), (4, 8)], [(0, 0), (0, 4), (2, 0), (2, 7), (2, 12), (3, 14)], [(0, 0), (0, 6), (0, 15), (2, 4), (3, 9), (4, 21)], [(0, 0), (0, 11), (1, 22), (2, 13), (3, 16), (4, 12)], [(0, 0), (1, 15), (2, 5), (3, 20), (4, 10), $\infty$]]
\item \{0=1250, 1=47000, 2=678000, 3=2944000, 4=3889750\} [[(0, 0), (0, 1), (0, 3), (1, 0), (2, 2), (4, 9)], [(0, 0), (0, 4), (0, 11), (2, 4), (3, 17), (4, 18)], [(0, 0), (0, 5), (0, 13), (2, 22), (3, 2), (4, 10)], [(0, 0), (0, 6), (0, 15), (1, 4), (3, 10), (4, 19)], [(0, 0), (1, 10), (2, 20), (3, 5), (4, 15), $\infty$]]
\item \{0=1250, 1=48000, 2=641250, 3=3041500, 4=3828000\} [[(0, 0), (0, 1), (0, 3), (1, 1), (1, 17), (3, 4)], [(0, 0), (0, 4), (1, 9), (3, 24), (4, 5), (4, 18)], [(0, 0), (0, 5), (0, 10), (0, 15), (0, 20), $\infty$], [(0, 0), (0, 6), (1, 18), (2, 8), (3, 11), (3, 18)], [(0, 0), (0, 8), (1, 4), (2, 6), (2, 17), (3, 14)]]
\item \{0=1250, 1=51000, 2=704250, 3=3011500, 4=3792000\} [[(0, 0), (0, 1), (0, 3), (1, 1), (3, 19), (4, 13)], [(0, 0), (0, 4), (0, 11), (2, 0), (3, 6), (4, 22)], [(0, 0), (0, 5), (0, 10), (0, 15), (0, 20), $\infty$], [(0, 0), (0, 6), (2, 1), (2, 17), (3, 21), (4, 4)], [(0, 0), (0, 8), (1, 18), (3, 17), (4, 3), (4, 16)]]
\item \{0=3750, 1=34000, 2=600750, 3=3001500, 4=3920000\} [[(0, 0), (0, 1), (0, 3), (0, 18), (1, 1), (2, 18)], [(0, 0), (0, 4), (0, 13), (1, 7), (2, 16), (3, 2)], [(0, 0), (0, 5), (1, 10), (1, 21), (2, 9), (3, 6)], [(0, 0), (0, 6), (1, 12), (2, 14), (3, 18), (4, 11)], [(0, 0), (1, 15), (2, 5), (3, 20), (4, 10), $\infty$]]
\item \{0=3750, 1=41000, 2=651000, 3=2946000, 4=3918250\} [[(0, 0), (0, 1), (0, 3), (0, 18), (1, 1), (2, 18)], [(0, 0), (0, 4), (0, 16), (2, 2), (3, 13), (4, 22)], [(0, 0), (0, 5), (1, 10), (1, 21), (2, 9), (3, 6)], [(0, 0), (0, 6), (1, 12), (2, 14), (3, 18), (4, 11)], [(0, 0), (1, 15), (2, 5), (3, 20), (4, 10), $\infty$]]
\end{enumerate}
\end{example}

\begin{example} $(\mathbb Z_5 \times \mathbb Z_5) \rtimes \mathbb Z_5$, where the generator $1$ of  the group $\mathbb Z_5$ acts on generators $(0,1)$ and $(1,0)$ of the group $\mathbb Z_5 \times \mathbb Z_5$ mapping them to $(0,1)$ and $(1,1)$, respectively.
\begin{enumerate}
\item {\bf \{4=7560000\} [[((0, 0), 0), ((0, 0), 1), ((1, 0), 0), ((1, 1), 2), ((2, 1), 1), ((2, 2), 2)], [((0, 0), 0), ((0, 0), 2), ((2, 1), 3), ((2, 4), 1), ((3, 2), 0), ((4, 2), 3)], [((0, 0), 0), ((0, 1), 0), ((0, 2), 0), ((0, 3), 0), ((0, 4), 0), $\infty$], [((0, 0), 0), ((0, 1), 3), ((2, 0), 2), ((2, 3), 3), ((4, 1), 2), ((4, 3), 0)], [((0, 0), 0), ((0, 1), 4), ((1, 1), 3), ((1, 4), 2), ((2, 0), 3), ((4, 4), 4)]]}
\item \{2=684000, 3=3054000, 4=3822000\} [[((0, 0), 0), ((0, 0), 1), ((0, 0), 2), ((0, 0), 3), ((0, 0), 4), $\infty$], [((0, 0), 0), ((0, 1), 0), ((1, 0), 0), ((1, 2), 1), ((4, 0), 2), ((4, 3), 3)], [((0, 0), 0), ((0, 1), 1), ((2, 2), 0), ((2, 4), 0), ((4, 1), 2), ((4, 2), 1)], [((0, 0), 0), ((0, 1), 2), ((1, 0), 4), ((1, 4), 1), ((2, 0), 2), ((4, 3), 0)], [((0, 0), 0), ((0, 2), 2), ((2, 0), 3), ((2, 3), 0), ((3, 4), 0), ((4, 3), 4)]]
\item \{1=12000, 2=456000, 3=3060000, 4=4032000\} [[((0, 0), 0), ((0, 0), 1), ((1, 0), 0), ((1, 1), 2), ((2, 0), 2), ((2, 4), 1)], [((0, 0), 0), ((0, 0), 2), ((2, 2), 1), ((2, 4), 3), ((3, 1), 0), ((4, 0), 3)], [((0, 0), 0), ((0, 1), 0), ((0, 2), 0), ((0, 3), 0), ((0, 4), 0), $\infty$], [((0, 0), 0), ((0, 1), 3), ((2, 0), 3), ((2, 2), 2), ((4, 2), 2), ((4, 4), 0)], [((0, 0), 0), ((0, 1), 4), ((1, 0), 3), ((1, 3), 2), ((2, 1), 3), ((4, 3), 4)]]
\item \{1=20000, 2=591000, 3=3024000, 4=3925000\} [[((0, 0), 0), ((0, 0), 1), ((0, 0), 2), ((0, 0), 3), ((0, 0), 4), $\infty$], [((0, 0), 0), ((0, 1), 0), ((1, 0), 0), ((1, 2), 4), ((4, 0), 1), ((4, 3), 0)], [((0, 0), 0), ((0, 1), 1), ((1, 0), 3), ((1, 1), 2), ((3, 0), 4), ((3, 2), 4)], [((0, 0), 0), ((0, 1), 2), ((1, 1), 3), ((2, 1), 0), ((3, 0), 3), ((3, 1), 1)], [((0, 0), 0), ((0, 2), 2), ((1, 0), 2), ((1, 3), 4), ((2, 4), 0), ((4, 0), 4)]]
\item \{1=24000, 2=591000, 3=2994000, 4=3951000\} [[((0, 0), 0), ((0, 0), 1), ((0, 0), 2), ((0, 0), 3), ((0, 0), 4), $\infty$], [((0, 0), 0), ((0, 1), 0), ((1, 0), 0), ((1, 1), 2), ((4, 0), 3), ((4, 4), 0)], [((0, 0), 0), ((0, 1), 1), ((1, 0), 1), ((1, 4), 2), ((2, 2), 0), ((4, 0), 2)], [((0, 0), 0), ((0, 2), 0), ((2, 0), 1), ((2, 2), 3), ((3, 0), 2), ((3, 3), 4)], [((0, 0), 0), ((0, 2), 1), ((2, 1), 2), ((2, 4), 3), ((3, 1), 1), ((4, 3), 2)]]
\item \{1=27000, 2=588750, 3=3013500, 4=3930750\} [[((0, 0), 0), ((0, 0), 1), ((0, 0), 2), ((0, 0), 3), ((0, 0), 4), $\infty$], [((0, 0), 0), ((0, 1), 0), ((0, 2), 1), ((1, 0), 0), ((3, 0), 3), ((4, 4), 0)], [((0, 0), 0), ((0, 1), 2), ((1, 0), 3), ((2, 2), 4), ((2, 3), 2), ((3, 4), 1)], [((0, 0), 0), ((0, 1), 4), ((0, 3), 1), ((1, 3), 4), ((3, 2), 1), ((4, 3), 3)], [((0, 0), 0), ((0, 2), 0), ((2, 0), 4), ((2, 3), 1), ((3, 0), 0), ((4, 2), 2)]]
\item \{1=27000, 2=593250, 3=3028500, 4=3911250\} [[((0, 0), 0), ((0, 0), 1), ((0, 0), 2), ((0, 0), 3), ((0, 0), 4), $\infty$], [((0, 0), 0), ((0, 1), 0), ((1, 0), 0), ((1, 1), 2), ((3, 0), 4), ((4, 3), 0)], [((0, 0), 0), ((0, 1), 1), ((0, 3), 2), ((1, 3), 1), ((2, 4), 4), ((4, 4), 3)], [((0, 0), 0), ((0, 1), 3), ((1, 1), 1), ((1, 4), 4), ((3, 1), 0), ((3, 2), 4)], [((0, 0), 0), ((0, 2), 0), ((1, 2), 4), ((2, 3), 0), ((4, 0), 3), ((4, 2), 2)]]
\item \{1=32000, 2=585000, 3=2976000, 4=3967000\} [[((0, 0), 0), ((0, 0), 1), ((1, 0), 0), ((1, 1), 2), ((2, 0), 2), ((2, 4), 1)], [((0, 0), 0), ((0, 0), 2), ((2, 1), 4), ((2, 4), 0), ((4, 2), 2), ((4, 4), 4)], [((0, 0), 0), ((0, 1), 0), ((0, 2), 0), ((0, 3), 0), ((0, 4), 0), $\infty$], [((0, 0), 0), ((0, 1), 3), ((1, 2), 2), ((2, 3), 0), ((3, 0), 4), ((3, 2), 2)], [((0, 0), 0), ((0, 1), 4), ((1, 0), 3), ((1, 3), 2), ((2, 1), 3), ((4, 3), 4)]]
\item \{1=35000, 2=674250, 3=3007500, 4=3843250\} [[((0, 0), 0), ((0, 0), 1), ((0, 0), 2), ((0, 0), 3), ((0, 0), 4), $\infty$], [((0, 0), 0), ((0, 1), 0), ((0, 2), 1), ((1, 0), 0), ((1, 1), 2), ((2, 3), 3)], [((0, 0), 0), ((0, 1), 3), ((0, 3), 2), ((1, 2), 3), ((2, 1), 2), ((3, 1), 4)], [((0, 0), 0), ((0, 1), 4), ((1, 4), 3), ((2, 1), 4), ((3, 0), 2), ((3, 2), 4)], [((0, 0), 0), ((0, 2), 0), ((1, 0), 4), ((2, 0), 0), ((3, 1), 3), ((4, 0), 3)]]
\item \{1=38000, 2=564750, 3=3007500, 4=3949750\} [[((0, 0), 0), ((0, 0), 1), ((0, 0), 2), ((0, 0), 3), ((0, 0), 4), $\infty$], [((0, 0), 0), ((0, 1), 0), ((1, 0), 0), ((1, 1), 3), ((1, 2), 2), ((4, 0), 4)], [((0, 0), 0), ((0, 1), 1), ((2, 0), 0), ((2, 2), 0), ((4, 3), 0), ((4, 4), 2)], [((0, 0), 0), ((0, 2), 1), ((1, 1), 0), ((3, 1), 4), ((3, 3), 3), ((4, 3), 1)], [((0, 0), 0), ((0, 2), 3), ((1, 4), 4), ((2, 3), 1), ((3, 1), 0), ((4, 0), 2)]]
\item \{1=40000, 2=669000, 3=3003000, 4=3848000\} [[((0, 0), 0), ((0, 0), 1), ((0, 0), 2), ((0, 0), 3), ((0, 0), 4), $\infty$], [((0, 0), 0), ((0, 1), 0), ((1, 0), 0), ((1, 2), 1), ((4, 0), 2), ((4, 3), 3)], [((0, 0), 0), ((0, 1), 1), ((1, 2), 0), ((1, 3), 4), ((2, 4), 2), ((4, 2), 0)], [((0, 0), 0), ((0, 1), 2), ((1, 0), 1), ((1, 1), 4), ((3, 1), 1), ((3, 4), 1)], [((0, 0), 0), ((0, 2), 2), ((2, 0), 3), ((2, 3), 0), ((3, 4), 0), ((4, 3), 4)]]
\item \{1=60000, 2=657000, 3=2931000, 4=3912000\} [[((0, 0), 0), ((0, 0), 1), ((0, 0), 2), ((0, 0), 3), ((0, 0), 4), $\infty$], [((0, 0), 0), ((0, 1), 0), ((1, 0), 0), ((1, 2), 3), ((4, 0), 1), ((4, 2), 3)], [((0, 0), 0), ((0, 1), 1), ((1, 2), 2), ((1, 3), 1), ((2, 1), 1), ((4, 4), 3)], [((0, 0), 0), ((0, 1), 2), ((1, 0), 3), ((1, 1), 1), ((3, 0), 2), ((3, 2), 2)], [((0, 0), 0), ((0, 2), 1), ((2, 2), 3), ((2, 4), 2), ((3, 4), 0), ((4, 3), 1)]]
\item \{1=82000, 2=693000, 3=3006000, 4=3779000\} [[((0, 0), 0), ((0, 0), 1), ((0, 0), 2), ((0, 0), 3), ((0, 0), 4), $\infty$], [((0, 0), 0), ((0, 1), 0), ((1, 0), 0), ((1, 2), 1), ((4, 0), 2), ((4, 3), 3)], [((0, 0), 0), ((0, 1), 1), ((2, 0), 2), ((2, 2), 2), ((4, 2), 2), ((4, 3), 1)], [((0, 0), 0), ((0, 1), 2), ((1, 3), 1), ((1, 4), 4), ((2, 0), 1), ((4, 2), 4)], [((0, 0), 0), ((0, 2), 2), ((2, 0), 0), ((2, 2), 3), ((3, 3), 0), ((4, 1), 4)]]
\item \{1=144000, 2=630000, 3=3000000, 4=3786000\} [[((0, 0), 0), ((0, 0), 1), ((0, 0), 2), ((0, 0), 3), ((0, 0), 4), $\infty$], [((0, 0), 0), ((0, 1), 0), ((1, 0), 0), ((1, 2), 3), ((4, 0), 2), ((4, 3), 0)], [((0, 0), 0), ((0, 1), 1), ((1, 1), 0), ((1, 2), 4), ((2, 4), 1), ((4, 1), 4)], [((0, 0), 0), ((0, 1), 2), ((1, 0), 4), ((1, 4), 1), ((3, 1), 3), ((3, 3), 3)], [((0, 0), 0), ((0, 2), 1), ((2, 1), 4), ((2, 3), 3), ((3, 1), 0), ((4, 4), 4)]]
\item \{0=1250, 1=20000, 2=571500, 3=2986000, 4=3981250\} [[((0, 0), 0), ((0, 0), 1), ((0, 0), 2), ((0, 0), 3), ((0, 0), 4), $\infty$], [((0, 0), 0), ((0, 1), 0), ((0, 2), 1), ((2, 0), 0), ((3, 1), 4), ((4, 3), 3)], [((0, 0), 0), ((0, 1), 2), ((1, 0), 0), ((1, 3), 1), ((2, 1), 2), ((2, 4), 4)], [((0, 0), 0), ((0, 1), 3), ((1, 3), 3), ((3, 2), 2), ((3, 4), 2), ((4, 0), 3)], [((0, 0), 0), ((0, 1), 4), ((1, 2), 4), ((2, 2), 4), ((2, 4), 1), ((3, 0), 2)]]
\item \{0=1250, 1=24000, 2=636750, 3=3005500, 4=3892500\} [[((0, 0), 0), ((0, 0), 1), ((0, 0), 2), ((0, 0), 3), ((0, 0), 4), $\infty$], [((0, 0), 0), ((0, 1), 0), ((1, 0), 0), ((1, 1), 1), ((1, 2), 3), ((3, 2), 3)], [((0, 0), 0), ((0, 1), 3), ((1, 1), 0), ((2, 3), 1), ((3, 2), 2), ((4, 2), 1)], [((0, 0), 0), ((0, 1), 4), ((0, 3), 1), ((1, 4), 3), ((2, 1), 0), ((2, 3), 0)], [((0, 0), 0), ((0, 2), 1), ((1, 3), 0), ((2, 4), 4), ((3, 3), 2), ((4, 3), 0)]]
\item \{0=1250, 1=27000, 2=566250, 3=2969500, 4=3996000\} [[((0, 0), 0), ((0, 0), 1), ((0, 0), 2), ((0, 0), 3), ((0, 0), 4), $\infty$], [((0, 0), 0), ((0, 1), 0), ((0, 2), 1), ((1, 0), 0), ((3, 3), 1), ((4, 1), 4)], [((0, 0), 0), ((0, 1), 2), ((1, 1), 4), ((2, 0), 3), ((3, 1), 3), ((4, 3), 4)], [((0, 0), 0), ((0, 1), 3), ((0, 3), 1), ((1, 0), 3), ((2, 3), 0), ((4, 2), 2)], [((0, 0), 0), ((0, 1), 4), ((2, 0), 2), ((2, 2), 2), ((4, 0), 2), ((4, 3), 0)]]
\item \{0=1250, 1=34000, 2=592500, 3=3031000, 4=3901250\} [[((0, 0), 0), ((0, 0), 1), ((0, 0), 2), ((0, 0), 3), ((0, 0), 4), $\infty$], [((0, 0), 0), ((0, 1), 0), ((0, 3), 1), ((1, 0), 0), ((3, 0), 3), ((3, 4), 0)], [((0, 0), 0), ((0, 1), 1), ((2, 1), 1), ((3, 2), 4), ((4, 0), 2), ((4, 2), 4)], [((0, 0), 0), ((0, 1), 2), ((1, 0), 1), ((1, 3), 1), ((2, 3), 3), ((4, 2), 2)], [((0, 0), 0), ((0, 1), 4), ((1, 1), 1), ((1, 3), 4), ((3, 3), 1), ((4, 1), 2)]]
\item \{0=1250, 1=35000, 2=669000, 3=3034000, 4=3820750\} [[((0, 0), 0), ((0, 0), 1), ((0, 0), 2), ((0, 0), 3), ((0, 0), 4), $\infty$], [((0, 0), 0), ((0, 1), 0), ((0, 2), 1), ((1, 0), 0), ((1, 2), 0), ((3, 4), 2)], [((0, 0), 0), ((0, 1), 2), ((1, 4), 3), ((2, 3), 2), ((3, 2), 4), ((4, 1), 2)], [((0, 0), 0), ((0, 1), 3), ((0, 2), 2), ((3, 1), 2), ((3, 4), 3), ((4, 2), 4)], [((0, 0), 0), ((0, 2), 3), ((1, 3), 1), ((2, 1), 2), ((3, 4), 4), ((4, 1), 1)]]
\item \{0=2500, 1=31000, 2=545250, 3=2997500, 4=3983750\} [[((0, 0), 0), ((0, 0), 1), ((0, 0), 2), ((0, 0), 3), ((0, 0), 4), $\infty$], [((0, 0), 0), ((0, 1), 0), ((0, 2), 1), ((1, 0), 0), ((3, 4), 4), ((4, 2), 3)], [((0, 0), 0), ((0, 1), 2), ((0, 3), 2), ((1, 4), 3), ((2, 4), 0), ((3, 2), 4)], [((0, 0), 0), ((0, 1), 3), ((1, 3), 0), ((2, 0), 0), ((3, 1), 4), ((3, 3), 3)], [((0, 0), 0), ((0, 1), 4), ((1, 0), 3), ((1, 3), 1), ((2, 0), 4), ((4, 1), 4)]]
\end{enumerate}
\end{example}

\begin{example} $\mathbb Z_{25} \rtimes \mathbb Z_5$, where the generator $1$ of $\mathbb Z_5$ acts on the generator $1$ of $\mathbb Z_{25}$ mapping it to $6$.
\begin{enumerate}
\item \{1=26000, 2=501750, 3=2977500, 4=4054750\} [[(0, 0), (0, 1), (1, 0), (3, 3), (6, 4), (9, 4)], [(0, 0), (0, 2), (7, 1), (9, 0), (19, 1), (24, 3)], [(0, 0), (1, 2), (8, 0), (10, 4), (20, 2), (22, 0)], [(0, 0), (2, 2), (5, 1), (16, 4), (18, 4), (19, 0)], [(0, 0), (5, 0), (10, 0), (15, 0), (20, 0), $\infty$]]
\item \{1=26000, 2=588750, 3=2971500, 4=3973750\} [[(0, 0), (0, 1), (0, 2), (0, 3), (0, 4), $\infty$], [(0, 0), (1, 0), (2, 1), (3, 3), (7, 2), (17, 0)], [(0, 0), (1, 2), (10, 4), (11, 4), (13, 0), (22, 2)], [(0, 0), (1, 3), (18, 2), (19, 3), (21, 4), (24, 3)], [(0, 0), (1, 4), (5, 3), (10, 0), (17, 4), (18, 3)]]
\item \{1=27000, 2=610500, 3=3003000, 4=3919500\} [[(0, 0), (0, 1), (0, 2), (0, 3), (0, 4), $\infty$], [(0, 0), (1, 0), (2, 1), (18, 4), (19, 3), (23, 4)], [(0, 0), (1, 2), (6, 1), (10, 2), (15, 0), (18, 0)], [(0, 0), (1, 3), (3, 4), (12, 4), (13, 2), (23, 3)], [(0, 0), (2, 0), (5, 2), (11, 0), (21, 4), (23, 2)]]
\item \{1=27000, 2=615750, 3=3022500, 4=3894750\} [[(0, 0), (0, 1), (0, 2), (0, 3), (0, 4), $\infty$], [(0, 0), (1, 0), (3, 0), (4, 3), (10, 4), (11, 0)], [(0, 0), (1, 1), (5, 2), (6, 0), (15, 3), (23, 1)], [(0, 0), (1, 4), (4, 1), (9, 2), (17, 3), (18, 0)], [(0, 0), (2, 3), (4, 2), (7, 1), (8, 3), (22, 3)]]
\item \{1=28000, 2=576750, 3=2989500, 4=3965750\} [[(0, 0), (0, 1), (0, 2), (0, 3), (0, 4), $\infty$], [(0, 0), (1, 0), (2, 1), (7, 2), (14, 1), (24, 4)], [(0, 0), (1, 2), (5, 0), (10, 4), (16, 1), (24, 1)], [(0, 0), (1, 3), (8, 0), (13, 2), (15, 3), (16, 3)], [(0, 0), (2, 3), (13, 4), (16, 0), (19, 3), (22, 0)]]
\item \{1=29000, 2=542250, 3=2983500, 4=4005250\} [[(0, 0), (0, 1), (0, 2), (0, 3), (0, 4), $\infty$], [(0, 0), (1, 0), (2, 1), (4, 0), (8, 3), (18, 2)], [(0, 0), (1, 4), (6, 4), (20, 4), (21, 2), (23, 3)], [(0, 0), (2, 0), (10, 0), (13, 2), (15, 4), (21, 3)], [(0, 0), (2, 2), (6, 1), (7, 0), (16, 4), (19, 0)]]
\item \{1=29000, 2=583500, 3=3000000, 4=3947500\} [[(0, 0), (0, 1), (0, 2), (0, 3), (0, 4), $\infty$], [(0, 0), (1, 0), (2, 1), (5, 1), (17, 1), (21, 3)], [(0, 0), (1, 2), (2, 2), (17, 3), (18, 2), (20, 0)], [(0, 0), (1, 3), (2, 4), (20, 2), (21, 4), (23, 3)], [(0, 0), (2, 3), (4, 2), (10, 2), (12, 1), (14, 3)]]
\item \{1=29000, 2=608250, 3=2980500, 4=3942250\} [[(0, 0), (0, 1), (0, 2), (0, 3), (0, 4), $\infty$], [(0, 0), (1, 0), (2, 1), (5, 0), (7, 0), (15, 2)], [(0, 0), (1, 2), (3, 0), (5, 3), (16, 1), (20, 3)], [(0, 0), (1, 4), (4, 2), (7, 1), (8, 0), (20, 4)], [(0, 0), (3, 3), (6, 1), (7, 2), (9, 1), (16, 0)]]
\item \{1=31000, 2=583500, 3=2997000, 4=3948500\} [[(0, 0), (0, 1), (1, 0), (6, 4), (8, 0), (21, 3)], [(0, 0), (0, 2), (3, 4), (9, 3), (12, 0), (14, 4)], [(0, 0), (1, 1), (8, 1), (13, 4), (19, 4), (23, 3)], [(0, 0), (2, 0), (6, 2), (14, 3), (21, 0), (23, 1)], [(0, 0), (5, 0), (10, 0), (15, 0), (20, 0), $\infty$]]
\item \{1=31000, 2=597000, 3=2991000, 4=3941000\} [[(0, 0), (0, 1), (0, 2), (0, 3), (0, 4), $\infty$], [(0, 0), (1, 0), (2, 4), (16, 3), (19, 0), (21, 3)], [(0, 0), (1, 1), (6, 2), (10, 3), (11, 1), (13, 3)], [(0, 0), (1, 2), (8, 1), (22, 4), (23, 2), (24, 1)], [(0, 0), (3, 0), (9, 1), (10, 1), (12, 0), (17, 3)]]
\item \{1=32000, 2=576750, 3=3028500, 4=3922750\} [[(0, 0), (0, 1), (0, 2), (0, 3), (0, 4), $\infty$], [(0, 0), (1, 0), (2, 1), (6, 2), (21, 2), (24, 4)], [(0, 0), (1, 2), (5, 1), (7, 2), (14, 2), (24, 3)], [(0, 0), (1, 4), (3, 1), (7, 0), (16, 2), (18, 3)], [(0, 0), (2, 0), (7, 4), (16, 0), (19, 2), (22, 0)]]
\item \{1=32000, 2=613500, 3=2946000, 4=3968500\} [[(0, 0), (0, 1), (0, 2), (0, 3), (0, 4), $\infty$], [(0, 0), (1, 0), (2, 1), (9, 0), (12, 2), (14, 0)], [(0, 0), (1, 2), (4, 0), (8, 1), (10, 3), (22, 0)], [(0, 0), (1, 3), (6, 2), (12, 1), (21, 2), (22, 3)], [(0, 0), (2, 0), (13, 4), (17, 1), (18, 2), (20, 3)]]
\item \{1=33000, 2=656250, 3=3010500, 4=3860250\} [[(0, 0), (0, 1), (0, 2), (0, 3), (0, 4), $\infty$], [(0, 0), (1, 0), (3, 0), (4, 2), (20, 4), (21, 3)], [(0, 0), (1, 1), (6, 1), (10, 1), (16, 4), (23, 4)], [(0, 0), (1, 4), (4, 0), (11, 4), (12, 2), (17, 1)], [(0, 0), (2, 1), (16, 1), (17, 2), (19, 3), (24, 1)]]
\item \{1=34000, 2=607500, 3=2961000, 4=3957500\} [[(0, 0), (0, 1), (0, 2), (0, 3), (0, 4), $\infty$], [(0, 0), (1, 0), (2, 1), (3, 4), (6, 0), (13, 2)], [(0, 0), (1, 2), (2, 2), (3, 0), (5, 3), (10, 0)], [(0, 0), (1, 3), (11, 2), (18, 2), (22, 1), (23, 3)], [(0, 0), (2, 0), (6, 2), (10, 2), (11, 1), (12, 1)]]
\item \{1=34000, 2=693000, 3=2967000, 4=3866000\} [[(0, 0), (0, 1), (1, 0), (2, 3), (13, 4), (20, 3)], [(0, 0), (0, 2), (6, 2), (7, 0), (17, 1), (20, 1)], [(0, 0), (1, 1), (7, 3), (11, 0), (13, 2), (17, 0)], [(0, 0), (3, 0), (11, 2), (12, 2), (13, 3), (14, 1)], [(0, 0), (5, 0), (10, 0), (15, 0), (20, 0), $\infty$]]
\item \{1=35000, 2=598500, 3=2982000, 4=3944500\} [[(0, 0), (0, 1), (0, 2), (0, 3), (0, 4), $\infty$], [(0, 0), (1, 0), (2, 1), (4, 0), (9, 0), (12, 2)], [(0, 0), (1, 2), (3, 1), (5, 3), (15, 3), (19, 4)], [(0, 0), (1, 3), (7, 1), (15, 1), (18, 4), (19, 0)], [(0, 0), (2, 0), (13, 0), (18, 2), (22, 1), (23, 3)]]
\item \{1=36000, 2=583500, 3=3000000, 4=3940500\} [[(0, 0), (0, 1), (0, 2), (0, 3), (0, 4), $\infty$], [(0, 0), (1, 0), (2, 1), (5, 0), (7, 4), (16, 1)], [(0, 0), (1, 2), (3, 4), (5, 2), (15, 4), (17, 4)], [(0, 0), (1, 3), (5, 1), (8, 3), (12, 4), (22, 4)], [(0, 0), (2, 2), (14, 3), (15, 2), (18, 0), (23, 4)]]
\item \{1=36000, 2=624750, 3=3034500, 4=3864750\} [[(0, 0), (0, 1), (0, 2), (0, 3), (0, 4), $\infty$], [(0, 0), (1, 0), (2, 3), (5, 3), (6, 1), (20, 1)], [(0, 0), (1, 1), (12, 0), (14, 0), (19, 2), (21, 1)], [(0, 0), (1, 2), (9, 0), (13, 1), (16, 3), (22, 4)], [(0, 0), (2, 1), (3, 0), (10, 0), (12, 2), (18, 2)]]
\item \{1=39000, 2=608250, 3=2962500, 4=3950250\} [[(0, 0), (0, 1), (0, 2), (0, 3), (0, 4), $\infty$], [(0, 0), (1, 0), (2, 1), (7, 4), (16, 0), (19, 2)], [(0, 0), (1, 2), (6, 1), (7, 2), (10, 2), (14, 2)], [(0, 0), (1, 4), (2, 2), (3, 3), (12, 1), (21, 4)], [(0, 0), (2, 4), (15, 4), (17, 1), (22, 2), (24, 2)]]
\item \{1=40000, 2=606750, 3=3016500, 4=3896750\} [[(0, 0), (0, 1), (0, 2), (0, 3), (0, 4), $\infty$], [(0, 0), (1, 0), (2, 2), (7, 4), (13, 3), (15, 2)], [(0, 0), (1, 4), (3, 0), (5, 0), (16, 0), (19, 3)], [(0, 0), (2, 3), (3, 1), (6, 0), (20, 4), (21, 0)], [(0, 0), (2, 4), (5, 3), (11, 2), (17, 3), (20, 1)]]
\item \{1=41000, 2=597000, 3=3003000, 4=3919000\} [[(0, 0), (0, 1), (0, 2), (0, 3), (0, 4), $\infty$], [(0, 0), (1, 0), (2, 1), (3, 4), (16, 0), (24, 1)], [(0, 0), (1, 3), (2, 0), (4, 3), (9, 2), (19, 1)], [(0, 0), (2, 2), (9, 4), (14, 4), (15, 4), (18, 1)], [(0, 0), (3, 0), (7, 0), (12, 1), (17, 3), (23, 2)]]
\item \{1=43000, 2=615750, 3=3037500, 4=3863750\} [[(0, 0), (0, 1), (0, 2), (0, 3), (0, 4), $\infty$], [(0, 0), (1, 0), (2, 1), (8, 3), (11, 3), (22, 1)], [(0, 0), (1, 4), (10, 4), (12, 4), (20, 4), (24, 1)], [(0, 0), (2, 2), (3, 4), (6, 0), (16, 1), (23, 1)], [(0, 0), (2, 4), (5, 4), (8, 1), (13, 3), (17, 2)]]
\item \{1=44000, 2=636000, 3=2982000, 4=3898000\} [[(0, 0), (0, 1), (0, 2), (0, 3), (0, 4), $\infty$], [(0, 0), (1, 0), (2, 1), (9, 0), (12, 2), (20, 1)], [(0, 0), (1, 2), (4, 2), (18, 3), (19, 3), (24, 4)], [(0, 0), (1, 3), (7, 2), (12, 4), (13, 2), (23, 2)], [(0, 0), (3, 3), (5, 0), (14, 1), (15, 3), (17, 1)]]
\item \{1=45000, 2=623250, 3=3043500, 4=3848250\} [[(0, 0), (0, 1), (0, 2), (0, 3), (0, 4), $\infty$], [(0, 0), (1, 0), (3, 1), (7, 2), (18, 1), (22, 3)], [(0, 0), (1, 1), (3, 0), (19, 4), (21, 0), (23, 3)], [(0, 0), (1, 2), (4, 2), (5, 3), (11, 0), (16, 0)], [(0, 0), (1, 3), (3, 4), (18, 3), (19, 0), (21, 4)]]
\item \{0=1250, 1=27000, 2=561000, 3=2986000, 4=3984750\} [[(0, 0), (0, 1), (0, 2), (0, 3), (0, 4), $\infty$], [(0, 0), (1, 0), (2, 1), (12, 2), (17, 1), (20, 3)], [(0, 0), (1, 2), (6, 2), (7, 0), (9, 0), (16, 4)], [(0, 0), (1, 3), (2, 2), (5, 3), (13, 1), (22, 1)], [(0, 0), (3, 0), (6, 4), (15, 1), (22, 4), (24, 4)]]
\item \{0=1250, 1=27000, 2=580500, 3=2941000, 4=4010250\} [[(0, 0), (0, 1), (0, 2), (0, 3), (0, 4), $\infty$], [(0, 0), (1, 0), (2, 1), (3, 0), (20, 3), (21, 0)], [(0, 0), (1, 2), (5, 4), (8, 1), (11, 4), (17, 0)], [(0, 0), (1, 4), (4, 2), (8, 3), (10, 0), (14, 3)], [(0, 0), (2, 4), (5, 3), (15, 2), (22, 3), (24, 2)]]
\item \{0=1250, 1=33000, 2=591000, 3=2983000, 4=3951750\} [[(0, 0), (0, 1), (0, 2), (0, 3), (0, 4), $\infty$], [(0, 0), (1, 0), (2, 3), (6, 0), (14, 0), (17, 4)], [(0, 0), (1, 1), (4, 0), (5, 4), (13, 2), (19, 4)], [(0, 0), (2, 0), (9, 4), (12, 2), (18, 1), (20, 2)], [(0, 0), (2, 4), (6, 2), (10, 3), (15, 0), (22, 3)]]
\item \{0=1250, 1=43000, 2=636750, 3=3017500, 4=3861500\} [[(0, 0), (0, 1), (0, 2), (0, 3), (0, 4), $\infty$], [(0, 0), (1, 0), (2, 1), (4, 1), (17, 3), (20, 0)], [(0, 0), (1, 2), (2, 4), (11, 1), (20, 2), (22, 2)], [(0, 0), (1, 3), (2, 0), (7, 4), (12, 0), (24, 1)], [(0, 0), (1, 4), (3, 2), (9, 0), (13, 3), (15, 3)]]
\item \{0=2500, 1=35000, 2=597000, 3=3014000, 4=3911500\} [[(0, 0), (0, 1), (0, 2), (0, 3), (0, 4), $\infty$], [(0, 0), (1, 0), (2, 3), (3, 4), (15, 4), (19, 3)], [(0, 0), (1, 1), (6, 3), (8, 4), (11, 4), (15, 1)], [(0, 0), (1, 2), (7, 4), (17, 4), (22, 2), (23, 2)], [(0, 0), (1, 4), (3, 1), (4, 1), (6, 4), (8, 0)]]
\end{enumerate}
\end{example}

\begin{example} $\mathbb Z_5 \times \mathbb Z_5 \times \mathbb Z_5$
\begin{enumerate}
\item \{2=366000, 3=2688000, 4=4506000\} [[(0, 0, 0), (0, 0, 1), (0, 0, 2), (0, 0, 3), (0, 0, 4), $\infty$], [(0, 0, 0), (0, 1, 0), (1, 0, 0), (1, 1, 1), (2, 3, 1), (4, 3, 2)], [(0, 0, 0), (0, 1, 2), (1, 2, 1), (1, 4, 2), (2, 0, 0), (2, 1, 4)], [(0, 0, 0), (0, 1, 3), (2, 2, 0), (2, 4, 0), (3, 2, 2), (3, 4, 4)], [(0, 0, 0), (0, 2, 3), (1, 1, 2), (2, 1, 0), (3, 0, 1), (3, 2, 0)]]
\item \{1=39000, 2=783000, 3=2931000, 4=3807000\} [[(0, 0, 0), (0, 0, 1), (0, 0, 2), (0, 0, 3), (0, 0, 4), $\infty$], [(0, 0, 0), (0, 1, 0), (0, 2, 1), (1, 0, 0), (1, 1, 2), (3, 2, 3)], [(0, 0, 0), (0, 1, 3), (1, 0, 1), (1, 3, 1), (3, 1, 1), (4, 3, 0)], [(0, 0, 0), (0, 1, 4), (0, 3, 3), (2, 0, 0), (3, 2, 1), (4, 2, 0)], [(0, 0, 0), (0, 2, 3), (1, 1, 0), (2, 1, 3), (3, 2, 4), (4, 4, 1)]]
\item \{1=54000, 2=760500, 3=2964000, 4=3781500\} [[(0, 0, 0), (0, 0, 1), (0, 0, 2), (0, 0, 3), (0, 0, 4), $\infty$], [(0, 0, 0), (0, 1, 0), (0, 2, 1), (1, 0, 0), (1, 1, 2), (3, 2, 3)], [(0, 0, 0), (0, 1, 3), (1, 0, 1), (1, 3, 1), (3, 1, 1), (4, 3, 0)], [(0, 0, 0), (0, 1, 4), (0, 3, 1), (1, 4, 4), (2, 4, 3), (3, 1, 4)], [(0, 0, 0), (0, 2, 3), (1, 1, 0), (2, 1, 3), (3, 2, 4), (4, 4, 1)]]
\item \{1=62000, 2=631500, 3=2901000, 4=3965500\} [[(0, 0, 0), (0, 0, 1), (0, 0, 2), (0, 0, 3), (0, 0, 4), $\infty$], [(0, 0, 0), (0, 1, 0), (1, 0, 0), (1, 1, 1), (2, 3, 1), (4, 3, 2)], [(0, 0, 0), (0, 1, 2), (1, 2, 1), (1, 4, 2), (2, 0, 0), (2, 1, 4)], [(0, 0, 0), (0, 1, 3), (2, 2, 0), (2, 4, 0), (3, 2, 2), (3, 4, 4)], [(0, 0, 0), (0, 2, 3), (2, 0, 3), (2, 2, 2), (3, 1, 3), (4, 1, 1)]]
\item {\bf \{1=76000, 2=699000, 3=2934000, 4=3851000\} [[(0, 0, 0), (0, 0, 1), (0, 0, 2), (0, 0, 3), (0, 0, 4), $\infty$], [(0, 0, 0), (0, 1, 0), (1, 0, 0), (1, 1, 1), (2, 3, 1), (4, 3, 2)], [(0, 0, 0), (0, 1, 2), (1, 2, 1), (1, 4, 2), (2, 0, 0), (2, 1, 4)], [(0, 0, 0), (0, 1, 3), (2, 2, 4), (2, 4, 1), (3, 2, 3), (3, 4, 3)], [(0, 0, 0), (0, 2, 3), (2, 0, 3), (2, 2, 2), (3, 1, 3), (4, 1, 1)]]}
\item \{1=80000, 2=690000, 3=2712000, 4=4078000\} [[(0, 0, 0), (0, 0, 1), (0, 0, 2), (0, 0, 3), (0, 0, 4), $\infty$], [(0, 0, 0), (0, 1, 0), (1, 0, 0), (1, 1, 1), (2, 3, 1), (4, 3, 2)], [(0, 0, 0), (0, 1, 2), (3, 0, 3), (3, 1, 2), (4, 2, 0), (4, 4, 1)], [(0, 0, 0), (0, 1, 3), (2, 2, 4), (2, 4, 1), (3, 2, 3), (3, 4, 3)], [(0, 0, 0), (0, 2, 3), (1, 1, 2), (2, 1, 0), (3, 0, 1), (3, 2, 0)]]
\item \{0=1250, 1=47000, 2=705750, 3=2975500, 4=3830500\} [[(0, 0, 0), (0, 0, 1), (0, 0, 2), (0, 0, 3), (0, 0, 4), $\infty$], [(0, 0, 0), (0, 1, 0), (0, 2, 1), (1, 0, 0), (1, 1, 2), (3, 2, 3)], [(0, 0, 0), (0, 1, 3), (1, 0, 1), (1, 3, 1), (3, 1, 1), (4, 3, 0)], [(0, 0, 0), (0, 1, 4), (0, 3, 3), (2, 0, 0), (3, 2, 1), (4, 2, 0)], [(0, 0, 0), (0, 2, 3), (1, 3, 2), (2, 0, 4), (3, 1, 0), (4, 1, 3)]]
\item \{0=1250, 1=52000, 2=712500, 3=3004000, 4=3790250\} [[(0, 0, 0), (0, 0, 1), (0, 0, 2), (0, 0, 3), (0, 0, 4), $\infty$], [(0, 0, 0), (0, 1, 0), (0, 2, 1), (1, 0, 0), (1, 1, 2), (3, 2, 3)], [(0, 0, 0), (0, 1, 3), (1, 0, 1), (1, 3, 1), (3, 1, 1), (4, 3, 0)], [(0, 0, 0), (0, 1, 4), (0, 3, 1), (1, 4, 4), (2, 4, 3), (3, 1, 4)], [(0, 0, 0), (0, 2, 3), (1, 3, 2), (2, 0, 4), (3, 1, 0), (4, 1, 3)]]
\end{enumerate}
\end{example}

\section{Transitive designs enumeration}
{\tt SmallGroup(126,6)} is $\mathbb Z_{126}$ and {\tt SmallGroup(126,16)} is $\mathbb Z_2 \times \mathbb Z_3 \times \mathbb Z_3 \times \mathbb Z_7$, so we don't include them and just point to \cite{Het} and \cite{Het1}. Results will be enumerated in form fingerprint-difference family.

Despite there are two designs that share same fingerprints $\{1=30240, 2=594972, 3=2954952, 4=3979836\}$ and $\{1=33264, 2=565488, 3=3024000, 4=3937248\}$, computer calculations show that they are non-isomorphic.

\begin{example} {\tt SmallGroup(126,1) = C7 : C18}
\begin{enumerate}
\item \{1=36288, 2=608580, 3=3010392, 4=3904740\} [[0, 1, 3, 6, 11, 17], [0, 2, 10, 81, 87, 93], [0, 5, 7, 39, 41, 115], [0, 8, 40, 76, 104, 125], [0, 9, 32, 90, 107, 124], [0, 19, 25, 46, 69, 91]]
\item \{1=36288, 2=692496, 3=3014928, 4=3816288\} [[0, 1, 3, 6, 11, 17], [0, 2, 4, 42, 78, 121], [0, 7, 72, 96, 112, 117], [0, 8, 31, 39, 113, 124], [0, 9, 35, 99, 102, 107], [0, 15, 69, 73, 85, 100]]
\item \{1=43344, 2=642600, 3=2981664, 4=3892392\} [[0, 1, 3, 6, 11, 17], [0, 2, 10, 31, 37, 61], [0, 5, 24, 25, 35, 87], [0, 7, 55, 79, 99, 105], [0, 9, 26, 29, 75, 100], [0, 12, 15, 56, 77, 102]]
\end{enumerate}
\end{example}

\begin{example} {\tt SmallGroup(126,2) = C2 x (C7 : C9)}
\begin{enumerate}
\item \{1=22176, 2=584388, 3=2983176, 4=3970260\} [[0, 1, 3, 6, 11, 17], [0, 2, 7, 23, 42, 71], [0, 4, 21, 72, 78, 85], [0, 5, 32, 53, 56, 76], [0, 9, 77, 83, 117, 124]]
\item \{1=26208, 2=644112, 3=2975616, 4=3914064\} [[0, 1, 3, 6, 11, 17], [0, 2, 4, 16, 60, 99], [0, 5, 26, 68, 75, 105], [0, 8, 43, 59, 64, 88], [0, 9, 37, 62, 97, 100]]
\item \{1=28224, 2=607068, 3=2992248, 4=3932460\} [[0, 1, 3, 6, 11, 17], [0, 2, 7, 14, 35, 64], [0, 4, 15, 21, 88, 117], [0, 8, 44, 54, 73, 84], [0, 9, 87, 90, 99, 119]]
\item \{1=29232, 2=600264, 3=2993760, 4=3936744\} [[0, 1, 3, 6, 11, 17], [0, 2, 4, 23, 72, 99], [0, 5, 58, 68, 90, 112], [0, 8, 26, 62, 69, 83], [0, 14, 48, 63, 82, 120]]
\item \{1=30240, 2=576072, 3=3054240, 4=3899448\} [[0, 1, 3, 6, 11, 17], [0, 2, 4, 21, 88, 108], [0, 7, 46, 94, 99, 104], [0, 8, 33, 55, 103, 124], [0, 14, 73, 92, 122, 123]]
\item \{1=30240, 2=634284, 3=3007368, 4=3888108\} [[0, 1, 3, 6, 11, 17], [0, 2, 4, 21, 101, 111], [0, 5, 24, 27, 78, 85], [0, 7, 57, 79, 81, 115], [0, 8, 39, 67, 92, 125]]
\item \{1=31248, 2=567756, 3=2995272, 4=3965724\} [[0, 1, 3, 6, 11, 17], [0, 2, 4, 52, 64, 73], [0, 5, 26, 57, 110, 121], [0, 8, 27, 46, 61, 107], [0, 10, 70, 71, 83, 95]]
\item \{1=31248, 2=588924, 3=2977128, 4=3962700\} [[0, 1, 3, 6, 11, 17], [0, 2, 4, 41, 42, 66], [0, 5, 26, 61, 73, 125], [0, 7, 15, 76, 81, 104], [0, 8, 68, 95, 98, 102]]
\item \{1=31248, 2=637308, 3=3049704, 4=3841740\} [[0, 1, 3, 6, 11, 17], [0, 2, 4, 23, 42, 82], [0, 5, 18, 43, 68, 106], [0, 7, 10, 35, 74, 96], [0, 8, 63, 85, 99, 104]]
\item \{1=32256, 2=560952, 3=3011904, 4=3954888\} [[0, 1, 3, 6, 11, 17], [0, 2, 4, 16, 67, 68], [0, 5, 51, 70, 113, 124], [0, 8, 48, 62, 91, 122], [0, 9, 34, 38, 59, 73]]
\item \{1=34272, 2=669816, 3=2948400, 4=3907512\} [[0, 1, 3, 6, 11, 17], [0, 2, 4, 51, 55, 108], [0, 5, 56, 66, 89, 115], [0, 9, 10, 68, 83, 94], [0, 12, 31, 98, 104, 110]]
\item \{1=36288, 2=629748, 3=3007368, 4=3886596\} [[0, 1, 3, 6, 11, 17], [0, 2, 4, 19, 38, 70], [0, 8, 16, 78, 111, 119], [0, 9, 47, 108, 109, 121], [0, 14, 33, 59, 72, 104]]
\item \{1=37296, 2=582876, 3=3022488, 4=3917340\} [[0, 1, 3, 6, 11, 17], [0, 2, 4, 32, 60, 90], [0, 7, 37, 40, 62, 84], [0, 8, 52, 68, 98, 102], [0, 25, 34, 44, 65, 113]]
\item \{1=37296, 2=613872, 3=2984688, 4=3924144\} [[0, 1, 3, 6, 11, 17], [0, 2, 7, 10, 41, 96], [0, 5, 18, 43, 84, 122], [0, 12, 32, 70, 93, 114], [0, 15, 77, 79, 99, 120]]
\item \{1=37296, 2=622944, 3=3014928, 4=3884832\} [[0, 1, 3, 6, 11, 17], [0, 2, 4, 29, 32, 101], [0, 7, 26, 75, 97, 119], [0, 10, 12, 76, 86, 99], [0, 13, 31, 46, 104, 121]]
\item \{1=37296, 2=636552, 3=2993760, 4=3892392\} [[0, 1, 3, 6, 11, 17], [0, 2, 4, 16, 34, 99], [0, 5, 12, 13, 61, 73], [0, 8, 41, 53, 63, 107], [0, 9, 69, 75, 102, 118]]
\item \{1=38304, 2=588168, 3=3027024, 4=3906504\} [[0, 1, 3, 6, 11, 17], [0, 2, 4, 37, 50, 87], [0, 5, 90, 96, 107, 108], [0, 7, 22, 34, 86, 93], [0, 13, 38, 83, 101, 118]]
\item \{1=38304, 2=645624, 3=2984688, 4=3891384\} [[0, 1, 3, 6, 11, 17], [0, 2, 4, 35, 40, 65], [0, 5, 85, 87, 108, 121], [0, 7, 69, 90, 102, 118], [0, 8, 33, 43, 97, 124]]
\item \{1=39312, 2=575316, 3=3043656, 4=3901716\} [[0, 1, 3, 6, 11, 17], [0, 2, 4, 14, 39, 119], [0, 7, 38, 96, 103, 108], [0, 8, 31, 95, 109, 113], [0, 10, 18, 35, 71, 94]]
\item \{1=40320, 2=612360, 3=2969568, 4=3937752\} [[0, 1, 3, 6, 11, 17], [0, 2, 4, 21, 32, 78], [0, 7, 23, 53, 83, 117], [0, 9, 33, 62, 81, 95], [0, 14, 34, 57, 58, 86]]
\item \{1=41328, 2=628236, 3=2998296, 4=3892140\} [[0, 1, 3, 6, 11, 17], [0, 2, 4, 23, 69, 124], [0, 5, 26, 55, 61, 116], [0, 7, 14, 43, 53, 94], [0, 9, 38, 48, 52, 79]]
\item \{1=42336, 2=615384, 3=3027024, 4=3875256\} [[0, 1, 3, 6, 11, 17], [0, 2, 4, 34, 68, 90], [0, 5, 57, 61, 93, 111], [0, 7, 15, 27, 29, 107], [0, 14, 41, 78, 97, 101]]
\item \{1=43344, 2=695520, 3=3024000, 4=3797136\} [[0, 1, 3, 6, 11, 17], [0, 2, 4, 26, 41, 122], [0, 5, 47, 48, 83, 108], [0, 7, 30, 42, 88, 92], [0, 8, 16, 58, 82, 106]]
\item \{1=45360, 2=625968, 3=3008880, 4=3879792\} [[0, 1, 3, 6, 11, 17], [0, 2, 4, 52, 73, 122], [0, 5, 33, 62, 81, 90], [0, 8, 61, 103, 107, 111], [0, 9, 75, 102, 113, 120]]
\item \{0=1260, 1=28224, 2=613872, 3=3006864, 4=3909780\} [[0, 1, 3, 6, 11, 17], [0, 2, 4, 23, 35, 81], [0, 5, 75, 78, 90, 94], [0, 7, 22, 74, 87, 117], [0, 25, 31, 56, 108, 113]]
\item \{0=1260, 1=29232, 2=635796, 3=3020472, 4=3873240\} [[0, 1, 3, 6, 11, 17], [0, 2, 4, 16, 82, 123], [0, 5, 32, 51, 104, 105], [0, 8, 24, 86, 90, 116], [0, 9, 29, 56, 79, 84]]
\item \{0=1260, 1=31248, 2=628992, 3=3012912, 4=3885588\} [[0, 1, 3, 6, 11, 17], [0, 2, 4, 21, 68, 110], [0, 7, 39, 71, 85, 105], [0, 8, 84, 88, 116, 118], [0, 13, 55, 98, 101, 108]]
\item \{0=1260, 1=33264, 2=594216, 3=3012912, 4=3918348\} [[0, 1, 3, 6, 11, 17], [0, 2, 4, 16, 83, 122], [0, 5, 19, 90, 113, 114], [0, 8, 23, 53, 57, 120], [0, 9, 18, 60, 75, 105]]
\item \{0=1260, 1=34272, 2=574560, 3=2982672, 4=3967236\} [[0, 1, 3, 6, 11, 17], [0, 2, 4, 26, 64, 97], [0, 5, 19, 39, 46, 74], [0, 7, 16, 103, 104, 109], [0, 8, 30, 88, 90, 93]]
\item \{0=1260, 1=38304, 2=664524, 3=2993256, 4=3862656\} [[0, 1, 3, 6, 11, 17], [0, 2, 4, 35, 112, 115], [0, 7, 39, 85, 99, 105], [0, 8, 15, 48, 51, 86], [0, 13, 55, 107, 108, 117]]
\item \{0=1260, 1=44352, 2=610092, 3=2993256, 4=3911040\} [[0, 1, 3, 6, 11, 17], [0, 2, 4, 35, 112, 119], [0, 7, 42, 81, 86, 111], [0, 8, 38, 78, 85, 90], [0, 14, 67, 73, 93, 95]]
\item \{0=2520, 1=28224, 2=656964, 3=2979144, 4=3893148\} [[0, 1, 3, 6, 11, 17], [0, 2, 4, 23, 72, 104], [0, 5, 61, 101, 102, 114], [0, 7, 40, 71, 76, 84], [0, 10, 58, 70, 74, 86]]
\item \{0=2520, 1=35280, 2=569268, 3=2997288, 4=3955644\} [[0, 1, 3, 6, 11, 17], [0, 2, 4, 16, 42, 122], [0, 5, 28, 33, 72, 88], [0, 8, 13, 61, 75, 83], [0, 12, 32, 39, 71, 86]]
\end{enumerate}
\end{example}

\begin{example} {\tt SmallGroup(126,3) = C7 x D18}
\begin{enumerate}
\item \{0=1260, 1=35280, 2=638820, 3=2950920, 4=3933720\} [[0, 1, 4, 7, 13, 19], [0, 2, 3, 61, 96, 98], [0, 5, 17, 24, 55, 86], [0, 6, 39, 45, 50, 95], [0, 9, 10, 49, 109, 113], [0, 11, 33, 66, 91, 110]]
\end{enumerate}
\end{example}

\begin{example} There are no designs for {\tt SmallGroup(126,4) = C9 x D14} and {\tt SmallGroup(126,5) = D126} and {\tt SmallGroup(126,9) = C7 : (C3 x S3)} and {\tt SmallGroup(126,11) = C3 x C3 x D14} and {\tt SmallGroup(126,13) = C3 x D42} and {\tt SmallGroup(126,14) = C7 x ((C3 x C3) : C2)} and {\tt SmallGroup(126,15) = (C3 x C3) : D14}
\end{example}

\begin{example} {\tt SmallGroup(126,7) = C3 x (C7 : C6)}
\begin{enumerate}
\item \{0=1260, 1=33264, 2=596484, 3=2993256, 4=3935736\} [[0, 1, 2, 5, 8, 14], [0, 3, 7, 11, 18, 33], [0, 4, 15, 24, 29, 99], [0, 6, 44, 96, 101, 116], [0, 9, 36, 62, 94, 118], [0, 10, 61, 90, 93, 119], [0, 12, 19, 28, 52, 73], [0, 20, 22, 103, 105, 125], [0, 48, 64, 89, 91, 109], [0, 49, 54, 58, 71, 112]]
\item \{1=49392, 2=595728, 3=3002832, 4=3912048\} [[0, 1, 2, 5, 8, 14], [0, 3, 7, 11, 18, 33], [0, 4, 20, 93, 109, 119], [0, 6, 62, 77, 117, 122], [0, 9, 36, 43, 56, 75], [0, 10, 50, 53, 92, 114], [0, 12, 15, 47, 79, 86], [0, 19, 30, 40, 67, 120], [0, 23, 32, 69, 113, 118], [0, 35, 63, 72, 85, 110]]
\end{enumerate}
\end{example}

\begin{example} {\tt SmallGroup(126,8) = S3 x (C7 : C3)}
\begin{enumerate}
\item \{1=24192, 2=574560, 3=3011904, 4=3949344\} [[0, 1, 3, 6, 11, 17], [0, 2, 7, 20, 24, 111], [0, 4, 21, 60, 88, 121], [0, 10, 59, 83, 117, 122], [0, 16, 18, 78, 100, 125]]
\item \{1=25200, 2=567000, 3=3020976, 4=3946824\} [[0, 1, 3, 6, 11, 17], [0, 2, 4, 23, 71, 76], [0, 5, 42, 47, 78, 94], [0, 7, 26, 41, 66, 119], [0, 10, 70, 97, 101, 122]]
\item \{1=26208, 2=573048, 3=3040128, 4=3920616\} [[0, 1, 9, 15, 25, 33], [0, 2, 8, 78, 93, 106], [0, 3, 16, 19, 49, 96], [0, 4, 20, 57, 75, 119], [0, 5, 32, 58, 90, 100], [0, 7, 44, 48, 85, 108]]
\item \{1=26208, 2=573804, 3=2993256, 4=3966732\} [[0, 1, 2, 5, 8, 14], [0, 3, 4, 22, 52, 56], [0, 6, 74, 95, 99, 116], [0, 7, 54, 64, 94, 123], [0, 12, 19, 41, 104, 110], [0, 13, 32, 46, 98, 125], [0, 17, 21, 63, 77, 115], [0, 31, 40, 100, 117, 124]]
\item \{1=27216, 2=529200, 3=2987712, 4=4015872\} [[0, 1, 9, 15, 25, 33], [0, 2, 8, 44, 59, 74], [0, 3, 35, 67, 85, 96], [0, 7, 34, 45, 75, 107], [0, 10, 49, 62, 100, 115], [0, 20, 24, 60, 106, 110]]
\item \{1=27216, 2=550368, 3=3024000, 4=3958416\} [[0, 1, 3, 6, 11, 17], [0, 2, 7, 26, 63, 123], [0, 4, 9, 35, 76, 104], [0, 10, 61, 73, 74, 85], [0, 12, 43, 59, 83, 116]]
\item \{1=27216, 2=557172, 3=3033576, 4=3942036\} [[0, 1, 2, 5, 8, 14], [0, 3, 7, 29, 71, 90], [0, 4, 30, 31, 50, 93], [0, 6, 24, 49, 108, 122], [0, 9, 62, 67, 81, 83], [0, 10, 12, 20, 22, 54], [0, 17, 59, 92, 101, 121], [0, 21, 61, 75, 106, 110]]
\item \{1=27216, 2=563976, 3=3035088, 4=3933720\} [[0, 1, 2, 5, 8, 14], [0, 3, 7, 10, 51, 86], [0, 4, 76, 94, 118, 124], [0, 6, 24, 49, 108, 122], [0, 9, 52, 67, 77, 117], [0, 12, 50, 79, 81, 113], [0, 13, 20, 93, 107, 120], [0, 17, 74, 95, 112, 123]]
\item \{1=27216, 2=619920, 3=3014928, 4=3897936\} [[0, 1, 3, 6, 11, 17], [0, 2, 7, 26, 60, 124], [0, 4, 47, 62, 66, 113], [0, 5, 20, 38, 59, 99], [0, 9, 25, 31, 52, 119]]
\item \{1=27216, 2=631260, 3=2985192, 4=3916332\} [[0, 1, 2, 5, 8, 14], [0, 3, 4, 36, 89, 100], [0, 6, 74, 95, 99, 116], [0, 10, 26, 64, 79, 121], [0, 12, 34, 50, 65, 70], [0, 17, 24, 66, 108, 125], [0, 21, 61, 75, 106, 110], [0, 30, 35, 46, 76, 118]]
\item \{1=28224, 2=542052, 3=3030552, 4=3959172\} [[0, 1, 2, 5, 8, 14], [0, 3, 29, 36, 64, 82], [0, 6, 24, 49, 108, 122], [0, 7, 38, 59, 116, 125], [0, 13, 50, 52, 79, 101], [0, 17, 74, 95, 112, 123], [0, 21, 61, 75, 106, 110], [0, 27, 37, 39, 91, 109]]
\item \{1=28224, 2=584388, 3=2971080, 4=3976308\} [[0, 1, 3, 6, 11, 17], [0, 2, 7, 31, 35, 114], [0, 4, 64, 88, 97, 121], [0, 5, 24, 36, 90, 104], [0, 12, 22, 37, 68, 96]]
\item \{1=28224, 2=612360, 3=2975616, 4=3943800\} [[0, 1, 2, 5, 8, 14], [0, 3, 4, 55, 96, 102], [0, 6, 74, 95, 99, 116], [0, 7, 28, 40, 84, 124], [0, 9, 33, 69, 83, 101], [0, 10, 17, 48, 56, 98], [0, 12, 44, 52, 68, 110], [0, 19, 31, 91, 120, 125]]
\item \{1=29232, 2=539784, 3=3061296, 4=3929688\} [[0, 1, 12, 18, 60, 69], [0, 2, 4, 52, 114, 125], [0, 5, 10, 41, 46, 104], [0, 9, 24, 79, 91, 108], [0, 15, 43, 65, 67, 106], [0, 20, 22, 95, 97, 124]]
\item \{1=29232, 2=560196, 3=2964024, 4=4006548\} [[0, 1, 9, 15, 39, 48], [0, 2, 8, 26, 40, 55], [0, 3, 78, 115, 121, 125], [0, 5, 23, 35, 50, 120], [0, 13, 37, 57, 93, 119], [0, 16, 25, 46, 95, 104]]
\item \{1=29232, 2=569268, 3=3049704, 4=3911796\} [[0, 1, 3, 6, 11, 17], [0, 2, 7, 19, 46, 73], [0, 5, 61, 76, 88, 108], [0, 9, 10, 40, 59, 63], [0, 18, 55, 94, 96, 113]]
\item \{1=29232, 2=581364, 3=2973096, 4=3976308\} [[0, 1, 2, 5, 8, 14], [0, 3, 4, 19, 79, 92], [0, 6, 10, 31, 56, 82], [0, 9, 16, 18, 36, 63], [0, 12, 49, 72, 116, 119], [0, 15, 51, 53, 64, 73], [0, 17, 74, 95, 112, 123], [0, 23, 69, 85, 90, 114]]
\item \{1=29232, 2=600264, 3=2982672, 4=3947832\} [[0, 1, 3, 6, 11, 17], [0, 2, 7, 8, 16, 29], [0, 4, 28, 43, 104, 123], [0, 9, 47, 93, 105, 114], [0, 10, 68, 102, 110, 113], [0, 13, 22, 37, 73, 109], [0, 21, 25, 72, 77, 120]]
\item \{1=29232, 2=608580, 3=2989224, 4=3932964\} [[0, 1, 3, 6, 11, 17], [0, 2, 4, 16, 20, 93], [0, 5, 13, 23, 96, 113], [0, 10, 67, 68, 70, 111], [0, 27, 38, 48, 51, 106]]
\item \{1=29232, 2=611604, 3=2997288, 4=3921876\} [[0, 1, 2, 5, 8, 14], [0, 3, 10, 15, 68, 123], [0, 4, 54, 72, 94, 110], [0, 6, 38, 41, 64, 67], [0, 7, 16, 27, 93, 99], [0, 9, 46, 82, 108, 121], [0, 12, 55, 59, 75, 90], [0, 17, 21, 63, 77, 115]]
\item \{1=29232, 2=628236, 3=2992248, 4=3910284\} [[0, 1, 9, 15, 25, 33], [0, 2, 21, 53, 99, 121], [0, 3, 7, 71, 82, 107], [0, 4, 14, 47, 60, 66], [0, 20, 40, 43, 74, 95], [0, 24, 73, 76, 96, 118]]
\item \{1=29232, 2=644112, 3=3006864, 4=3879792\} [[0, 1, 2, 5, 8, 14], [0, 3, 7, 36, 69, 114], [0, 4, 28, 29, 75, 115], [0, 6, 74, 95, 99, 116], [0, 10, 27, 50, 85, 123], [0, 13, 23, 106, 110, 124], [0, 17, 59, 92, 101, 121], [0, 21, 58, 78, 91, 120]]
\item \{1=29232, 2=659232, 3=2995776, 4=3875760\} [[0, 1, 2, 5, 8, 14], [0, 3, 15, 16, 46, 103], [0, 4, 51, 78, 90, 114], [0, 6, 59, 85, 92, 113], [0, 7, 36, 53, 79, 95], [0, 17, 24, 66, 108, 125], [0, 19, 20, 22, 49, 81], [0, 26, 39, 41, 73, 118]]
\item \{1=30240, 2=555660, 3=2998296, 4=3975804\} [[0, 1, 2, 5, 8, 14], [0, 3, 7, 16, 35, 96], [0, 4, 54, 72, 94, 110], [0, 6, 59, 85, 92, 113], [0, 9, 38, 39, 73, 111], [0, 10, 17, 48, 56, 98], [0, 15, 50, 77, 81, 97], [0, 23, 29, 33, 93, 100]]
\item \{1=30240, 2=594972, 3=2954952, 4=3979836\} [[0, 1, 2, 5, 8, 14], [0, 3, 7, 15, 34, 52], [0, 6, 74, 95, 99, 116], [0, 9, 36, 59, 78, 91], [0, 13, 22, 47, 97, 114], [0, 16, 35, 46, 84, 90], [0, 17, 21, 63, 77, 115], [0, 33, 37, 64, 103, 109]]
\item \{1=30240, 2=610848, 3=3002832, 4=3916080\} [[0, 1, 2, 5, 8, 14], [0, 3, 4, 10, 84, 93], [0, 6, 38, 41, 64, 67], [0, 7, 20, 42, 71, 92], [0, 9, 36, 68, 70, 112], [0, 12, 43, 99, 103, 108], [0, 17, 21, 63, 77, 115], [0, 23, 35, 90, 101, 122]]
\item \{1=30240, 2=617652, 3=3039624, 4=3872484\} [[0, 1, 3, 6, 11, 17], [0, 2, 7, 31, 41, 66], [0, 4, 45, 83, 93, 121], [0, 9, 36, 79, 92, 109], [0, 10, 50, 64, 102, 118], [0, 21, 23, 43, 77, 91], [0, 39, 68, 69, 73, 122]]
\item \{1=31248, 2=578340, 3=2989224, 4=3961188\} [[0, 1, 3, 6, 11, 17], [0, 2, 4, 29, 56, 94], [0, 5, 40, 63, 95, 111], [0, 9, 35, 98, 99, 122], [0, 15, 51, 59, 67, 110]]
\item \{1=31248, 2=591948, 3=2944872, 4=3991932\} [[0, 1, 2, 5, 8, 14], [0, 3, 4, 10, 67, 94], [0, 6, 24, 49, 108, 122], [0, 7, 48, 52, 57, 71], [0, 9, 36, 83, 85, 87], [0, 12, 16, 62, 100, 106], [0, 15, 35, 63, 93, 107], [0, 17, 59, 92, 101, 121]]
\item \{1=31248, 2=594216, 3=3030048, 4=3904488\} [[0, 1, 3, 6, 11, 17], [0, 2, 7, 20, 76, 98], [0, 4, 38, 73, 102, 107], [0, 5, 25, 47, 69, 77], [0, 10, 15, 85, 96, 110]]
\item \{1=31248, 2=595728, 3=3015936, 4=3917088\} [[0, 1, 2, 5, 8, 14], [0, 3, 4, 15, 34, 109], [0, 6, 74, 95, 99, 116], [0, 7, 62, 66, 85, 100], [0, 9, 36, 59, 78, 91], [0, 10, 17, 48, 56, 98], [0, 12, 60, 64, 73, 122], [0, 47, 52, 57, 115, 119]]
\item \{1=31248, 2=649404, 3=3004344, 4=3875004\} [[0, 1, 2, 5, 8, 14], [0, 3, 4, 48, 104, 109], [0, 6, 38, 41, 64, 67], [0, 7, 10, 59, 107, 120], [0, 9, 18, 40, 53, 95], [0, 15, 23, 47, 90, 117], [0, 17, 24, 66, 108, 125], [0, 34, 39, 73, 84, 100]]
\item \{1=32256, 2=591948, 3=2977128, 4=3958668\} [[0, 1, 3, 6, 11, 17], [0, 2, 7, 10, 113, 125], [0, 5, 32, 42, 45, 87], [0, 9, 16, 61, 76, 120], [0, 19, 47, 108, 116, 117]]
\item \{1=32256, 2=595728, 3=3065328, 4=3866688\} [[0, 1, 2, 5, 8, 14], [0, 3, 7, 42, 47, 125], [0, 4, 22, 91, 107, 120], [0, 6, 74, 95, 99, 116], [0, 9, 28, 35, 38, 68], [0, 10, 17, 48, 56, 98], [0, 12, 46, 73, 121, 122], [0, 23, 27, 76, 90, 124]]
\item \{1=32256, 2=599508, 3=3051720, 4=3876516\} [[0, 1, 2, 5, 8, 14], [0, 3, 7, 63, 81, 99], [0, 6, 59, 85, 92, 113], [0, 9, 26, 36, 38, 57], [0, 10, 17, 48, 56, 98], [0, 12, 32, 34, 62, 77], [0, 16, 66, 72, 95, 124], [0, 23, 51, 67, 90, 100]]
\item \{1=32256, 2=625212, 3=2982168, 4=3920364\} [[0, 1, 2, 5, 8, 14], [0, 3, 4, 16, 84, 88], [0, 6, 74, 95, 99, 116], [0, 7, 38, 63, 73, 87], [0, 10, 26, 28, 35, 113], [0, 17, 59, 92, 101, 121], [0, 21, 58, 78, 91, 120], [0, 23, 27, 76, 90, 124]]
\item \{1=32256, 2=629748, 3=2992248, 4=3905748\} [[0, 1, 3, 6, 11, 17], [0, 2, 4, 16, 59, 73], [0, 5, 23, 25, 75, 121], [0, 9, 19, 95, 104, 118], [0, 13, 45, 70, 79, 85]]
\item \{1=33264, 2=565488, 3=3024000, 4=3937248\} [[0, 1, 2, 5, 8, 14], [0, 3, 4, 28, 88, 120], [0, 6, 24, 49, 108, 122], [0, 7, 53, 95, 102, 117], [0, 9, 19, 39, 75, 83], [0, 16, 37, 69, 99, 109], [0, 17, 38, 41, 81, 84], [0, 23, 59, 60, 90, 119]]
\item \{1=33264, 2=568512, 3=3015936, 4=3942288\} [[0, 1, 2, 5, 8, 14], [0, 3, 10, 15, 84, 115], [0, 4, 9, 23, 36, 54], [0, 6, 24, 49, 108, 122], [0, 7, 21, 72, 88, 99], [0, 12, 38, 85, 98, 100], [0, 17, 74, 95, 112, 123], [0, 27, 55, 57, 109, 119]]
\item \{1=33264, 2=572292, 3=2955960, 4=3998484\} [[0, 1, 2, 5, 8, 14], [0, 3, 16, 21, 22, 99], [0, 4, 27, 97, 119, 121], [0, 6, 59, 85, 92, 113], [0, 12, 35, 76, 86, 96], [0, 15, 58, 72, 87, 114], [0, 17, 38, 41, 81, 84], [0, 19, 48, 106, 110, 122]]
\item \{1=33264, 2=579852, 3=3016440, 4=3930444\} [[0, 1, 3, 6, 11, 17], [0, 2, 4, 29, 32, 79], [0, 7, 30, 54, 76, 96], [0, 10, 12, 15, 102, 125], [0, 20, 57, 90, 100, 124]]
\item \{1=33264, 2=587412, 3=3026520, 4=3912804\} [[0, 1, 9, 15, 25, 33], [0, 2, 22, 49, 80, 84], [0, 5, 62, 72, 95, 101], [0, 10, 59, 83, 105, 120], [0, 21, 61, 75, 106, 110], [0, 30, 35, 46, 76, 118]]
\item \{1=33264, 2=629748, 3=3043656, 4=3853332\} [[0, 1, 2, 5, 8, 14], [0, 3, 4, 21, 71, 106], [0, 6, 24, 49, 108, 122], [0, 7, 23, 66, 90, 99], [0, 10, 13, 43, 63, 105], [0, 15, 29, 40, 98, 100], [0, 17, 38, 41, 81, 84], [0, 32, 46, 55, 75, 117]]
\item \{1=33264, 2=635796, 3=2963016, 4=3927924\} [[0, 1, 2, 5, 8, 14], [0, 3, 4, 21, 81, 116], [0, 6, 10, 31, 56, 82], [0, 7, 44, 73, 103, 111], [0, 13, 16, 52, 96, 110], [0, 15, 58, 72, 87, 114], [0, 17, 24, 66, 108, 125], [0, 20, 22, 41, 43, 90]]
\item \{1=34272, 2=590436, 3=2976120, 4=3959172\} [[0, 1, 2, 5, 8, 14], [0, 3, 7, 23, 51, 96], [0, 4, 57, 82, 114, 120], [0, 6, 59, 85, 92, 113], [0, 9, 13, 36, 40, 72], [0, 12, 15, 34, 64, 84], [0, 17, 74, 95, 112, 123], [0, 30, 35, 46, 76, 118]]
\item \{1=34272, 2=604800, 3=3007872, 4=3913056\} [[0, 1, 2, 5, 8, 14], [0, 3, 4, 10, 69, 93], [0, 6, 24, 49, 108, 122], [0, 7, 55, 75, 85, 98], [0, 9, 31, 44, 59, 103], [0, 12, 22, 99, 109, 114], [0, 17, 21, 63, 77, 115], [0, 23, 34, 50, 90, 121]]
\item \{1=34272, 2=640332, 3=2951928, 4=3933468\} [[0, 1, 3, 6, 11, 17], [0, 2, 7, 8, 16, 29], [0, 4, 9, 106, 110, 118], [0, 5, 59, 69, 78, 79], [0, 10, 28, 40, 93, 99], [0, 13, 37, 70, 75, 82], [0, 24, 48, 108, 117, 121]]
\item \{1=35280, 2=590436, 3=3037608, 4=3896676\} [[0, 1, 2, 5, 8, 14], [0, 3, 4, 50, 89, 96], [0, 6, 74, 95, 99, 116], [0, 9, 32, 55, 62, 113], [0, 10, 61, 72, 93, 107], [0, 12, 28, 82, 104, 123], [0, 17, 59, 92, 101, 121], [0, 29, 70, 84, 91, 120]]
\item \{1=35280, 2=595728, 3=3002832, 4=3926160\} [[0, 1, 3, 6, 11, 17], [0, 2, 4, 35, 40, 47], [0, 5, 18, 48, 51, 88], [0, 7, 10, 42, 90, 101], [0, 13, 22, 63, 84, 113]]
\item \{1=35280, 2=601020, 3=2974104, 4=3949596\} [[0, 1, 2, 5, 8, 14], [0, 3, 21, 36, 100, 124], [0, 6, 10, 31, 56, 82], [0, 7, 93, 102, 107, 125], [0, 13, 55, 65, 88, 92], [0, 15, 29, 69, 76, 118], [0, 17, 38, 41, 81, 84], [0, 23, 24, 63, 112, 120]]
\item \{1=35280, 2=616140, 3=2995272, 4=3913308\} [[0, 1, 9, 15, 39, 48], [0, 2, 8, 13, 24, 38], [0, 3, 16, 57, 68, 115], [0, 4, 23, 71, 112, 122], [0, 5, 18, 61, 94, 125], [0, 10, 26, 56, 58, 107]]
\item \{1=35280, 2=656208, 3=3003840, 4=3864672\} [[0, 1, 2, 5, 8, 14], [0, 3, 4, 10, 84, 105], [0, 6, 24, 49, 108, 122], [0, 7, 46, 78, 113, 125], [0, 9, 35, 41, 77, 124], [0, 15, 57, 71, 112, 119], [0, 17, 59, 92, 101, 121], [0, 23, 45, 70, 86, 90]]
\item \{1=36288, 2=588168, 3=3062304, 4=3873240\} [[0, 1, 2, 5, 8, 14], [0, 3, 7, 21, 73, 91], [0, 6, 59, 85, 92, 113], [0, 9, 19, 62, 69, 109], [0, 12, 46, 61, 101, 108], [0, 17, 38, 41, 81, 84], [0, 23, 45, 70, 86, 90], [0, 31, 40, 100, 117, 124]]
\item \{1=36288, 2=591192, 3=2965536, 4=3966984\} [[0, 1, 9, 15, 39, 48], [0, 2, 35, 47, 64, 99], [0, 4, 54, 78, 105, 121], [0, 5, 46, 69, 93, 107], [0, 10, 29, 66, 81, 123], [0, 20, 22, 38, 40, 42]]
\item \{1=36288, 2=591192, 3=3007872, 4=3924648\} [[0, 1, 2, 5, 8, 14], [0, 3, 10, 15, 25, 92], [0, 4, 27, 47, 79, 117], [0, 6, 24, 49, 108, 122], [0, 9, 36, 53, 99, 101], [0, 16, 54, 58, 70, 107], [0, 17, 21, 63, 77, 115], [0, 19, 23, 84, 90, 113]]
\item \{1=36288, 2=594972, 3=2952936, 4=3975804\} [[0, 1, 3, 6, 11, 17], [0, 2, 4, 29, 82, 106], [0, 5, 25, 39, 56, 124], [0, 12, 53, 79, 85, 116], [0, 15, 20, 51, 71, 100]]
\item \{1=36288, 2=599508, 3=3019464, 4=3904740\} [[0, 1, 2, 5, 8, 14], [0, 3, 7, 35, 38, 111], [0, 4, 21, 39, 45, 116], [0, 6, 24, 49, 108, 122], [0, 9, 36, 43, 56, 75], [0, 15, 32, 60, 63, 86], [0, 17, 59, 92, 101, 121], [0, 20, 22, 34, 47, 66]]
\item \{1=36288, 2=626724, 3=2966040, 4=3930948\} [[0, 1, 2, 5, 8, 14], [0, 3, 4, 35, 68, 123], [0, 6, 24, 49, 108, 122], [0, 7, 34, 41, 90, 117], [0, 9, 36, 44, 76, 106], [0, 10, 17, 48, 56, 98], [0, 13, 30, 82, 87, 120], [0, 37, 40, 78, 92, 109]]
\item \{1=36288, 2=645624, 3=3006864, 4=3871224\} [[0, 1, 2, 5, 8, 14], [0, 3, 16, 28, 96, 113], [0, 4, 39, 57, 73, 91], [0, 6, 24, 49, 108, 122], [0, 7, 9, 27, 50, 87], [0, 12, 35, 55, 67, 105], [0, 17, 74, 95, 112, 123], [0, 25, 37, 38, 109, 110]]
\item \{1=37296, 2=591192, 3=2948400, 4=3983112\} [[0, 1, 2, 5, 8, 14], [0, 3, 4, 19, 73, 109], [0, 6, 74, 95, 99, 116], [0, 9, 26, 47, 62, 125], [0, 12, 27, 87, 108, 113], [0, 17, 59, 92, 101, 121], [0, 29, 70, 84, 91, 120], [0, 31, 51, 106, 110, 114]]
\item \{1=37296, 2=617652, 3=2993256, 4=3911796\} [[0, 1, 2, 5, 8, 14], [0, 3, 4, 19, 50, 115], [0, 6, 10, 31, 56, 82], [0, 9, 36, 77, 96, 107], [0, 12, 32, 103, 104, 118], [0, 13, 52, 75, 117, 121], [0, 17, 24, 66, 108, 125], [0, 23, 59, 60, 90, 119]]
\item \{1=37296, 2=619920, 3=2987712, 4=3915072\} [[0, 1, 3, 6, 11, 17], [0, 2, 7, 19, 49, 104], [0, 5, 27, 64, 83, 107], [0, 10, 25, 45, 90, 95], [0, 12, 28, 77, 81, 119]]
\item \{1=37296, 2=629748, 3=3007368, 4=3885588\} [[0, 1, 2, 5, 8, 14], [0, 3, 4, 21, 43, 65], [0, 6, 38, 41, 64, 67], [0, 7, 76, 115, 118, 122], [0, 9, 13, 70, 93, 125], [0, 10, 17, 48, 56, 98], [0, 12, 59, 73, 87, 105], [0, 15, 33, 100, 106, 110]]
\item \{1=38304, 2=594216, 3=2973600, 4=3953880\} [[0, 1, 2, 5, 8, 14], [0, 3, 7, 19, 70, 122], [0, 6, 10, 31, 56, 82], [0, 9, 21, 24, 51, 100], [0, 12, 32, 90, 97, 116], [0, 17, 38, 41, 81, 84], [0, 27, 40, 72, 94, 124], [0, 33, 37, 64, 103, 109]]
\item \{1=38304, 2=666036, 3=2969064, 4=3886596\} [[0, 1, 9, 15, 39, 48], [0, 2, 33, 35, 81, 101], [0, 3, 19, 41, 43, 116], [0, 4, 5, 82, 113, 114], [0, 20, 22, 77, 79, 118], [0, 23, 51, 67, 90, 100]]
\item \{1=39312, 2=619164, 3=2957976, 4=3943548\} [[0, 1, 2, 5, 8, 14], [0, 3, 15, 16, 46, 122], [0, 4, 55, 73, 75, 93], [0, 6, 38, 41, 64, 67], [0, 7, 52, 71, 85, 111], [0, 9, 32, 34, 36, 81], [0, 17, 59, 92, 101, 121], [0, 18, 58, 61, 112, 120]]
\item \{1=39312, 2=621432, 3=3002832, 4=3896424\} [[0, 1, 2, 5, 8, 14], [0, 3, 16, 21, 78, 99], [0, 4, 22, 91, 107, 120], [0, 6, 59, 85, 92, 113], [0, 7, 35, 36, 51, 86], [0, 10, 15, 19, 67, 115], [0, 17, 24, 66, 108, 125], [0, 32, 46, 55, 75, 117]]
\item \{1=39312, 2=626724, 3=3004344, 4=3889620\} [[0, 1, 2, 5, 8, 14], [0, 3, 4, 16, 57, 125], [0, 6, 59, 85, 92, 113], [0, 7, 35, 54, 88, 102], [0, 9, 36, 43, 56, 75], [0, 12, 19, 37, 79, 99], [0, 17, 21, 63, 77, 115], [0, 23, 51, 67, 90, 100]]
\item \{1=39312, 2=655452, 3=2986200, 4=3879036\} [[0, 1, 3, 6, 11, 17], [0, 2, 4, 19, 52, 98], [0, 5, 36, 53, 68, 117], [0, 10, 73, 76, 83, 122], [0, 13, 37, 75, 86, 97]]
\item \{1=39312, 2=658476, 3=3005352, 4=3856860\} [[0, 1, 2, 5, 8, 14], [0, 3, 7, 46, 78, 109], [0, 4, 58, 72, 76, 90], [0, 6, 10, 31, 56, 82], [0, 9, 42, 53, 68, 96], [0, 12, 13, 24, 65, 85], [0, 17, 74, 95, 112, 123], [0, 19, 39, 67, 73, 121]]
\item \{1=39312, 2=692496, 3=3005856, 4=3822336\} [[0, 1, 3, 6, 11, 17], [0, 2, 4, 16, 37, 67], [0, 5, 35, 42, 47, 110], [0, 9, 62, 75, 86, 112], [0, 12, 49, 61, 93, 94]]
\item \{1=40320, 2=603288, 3=3046176, 4=3870216\} [[0, 1, 9, 15, 54, 63], [0, 2, 8, 44, 59, 74], [0, 3, 24, 29, 33, 69], [0, 4, 32, 52, 79, 85], [0, 10, 40, 107, 111, 124], [0, 18, 71, 109, 110, 112]]
\item \{1=40320, 2=607068, 3=3019464, 4=3893148\} [[0, 1, 2, 5, 8, 14], [0, 3, 4, 31, 35, 63], [0, 6, 59, 85, 92, 113], [0, 7, 10, 15, 100, 111], [0, 9, 36, 41, 60, 73], [0, 13, 49, 78, 101, 117], [0, 17, 24, 66, 108, 125], [0, 20, 27, 57, 91, 120]]
\item \{1=40320, 2=619920, 3=2936304, 4=3963456\} [[0, 1, 3, 6, 11, 17], [0, 2, 7, 34, 47, 115], [0, 4, 14, 81, 93, 124], [0, 9, 10, 46, 61, 76], [0, 13, 28, 71, 91, 100]]
\item \{1=40320, 2=627480, 3=2967552, 4=3924648\} [[0, 1, 2, 5, 8, 14], [0, 3, 16, 36, 85, 99], [0, 4, 49, 72, 84, 105], [0, 6, 21, 46, 77, 102], [0, 7, 27, 50, 51, 96], [0, 13, 39, 73, 75, 107], [0, 17, 24, 66, 108, 125], [0, 23, 37, 60, 74, 109]]
\item \{1=40320, 2=640332, 3=2999304, 4=3880044\} [[0, 1, 9, 15, 39, 48], [0, 2, 8, 79, 94, 107], [0, 3, 29, 87, 88, 91], [0, 4, 62, 68, 116, 125], [0, 5, 19, 43, 83, 118], [0, 21, 61, 75, 106, 110]]
\item \{1=41328, 2=575316, 3=3025512, 4=3917844\} [[0, 1, 3, 6, 11, 17], [0, 2, 4, 19, 94, 122], [0, 5, 35, 39, 74, 115], [0, 9, 44, 55, 59, 102], [0, 12, 45, 83, 88, 121]]
\item \{1=41328, 2=608580, 3=3025512, 4=3884580\} [[0, 1, 2, 5, 8, 14], [0, 3, 4, 28, 69, 119], [0, 6, 21, 46, 77, 102], [0, 7, 27, 34, 67, 112], [0, 10, 93, 99, 114, 120], [0, 16, 70, 76, 81, 118], [0, 17, 59, 92, 101, 121], [0, 37, 41, 78, 106, 110]]
\item \{1=41328, 2=625212, 3=2981160, 4=3912300\} [[0, 1, 2, 5, 8, 14], [0, 3, 4, 35, 40, 81], [0, 6, 59, 85, 92, 113], [0, 7, 10, 45, 58, 93], [0, 13, 28, 38, 78, 104], [0, 17, 74, 95, 112, 123], [0, 27, 37, 39, 91, 109], [0, 32, 46, 55, 75, 117]]
\item \{1=42336, 2=575316, 3=3010392, 4=3931956\} [[0, 1, 2, 5, 8, 14], [0, 3, 4, 15, 68, 105], [0, 6, 10, 31, 56, 82], [0, 7, 61, 90, 101, 116], [0, 9, 21, 26, 73, 77], [0, 12, 55, 67, 69, 104], [0, 17, 74, 95, 112, 123], [0, 27, 40, 72, 94, 124]]
\item \{1=42336, 2=628992, 3=3009888, 4=3878784\} [[0, 1, 9, 15, 54, 63], [0, 2, 7, 26, 70, 87], [0, 4, 22, 91, 107, 120], [0, 5, 12, 13, 82, 122], [0, 10, 35, 56, 68, 114], [0, 24, 48, 108, 117, 121], [0, 30, 57, 83, 88, 119], [0, 39, 42, 74, 95, 124]]
\item \{1=42336, 2=629748, 3=2973096, 4=3914820\} [[0, 1, 9, 15, 39, 48], [0, 2, 3, 44, 81, 86], [0, 4, 34, 65, 79, 118], [0, 5, 59, 106, 115, 124], [0, 10, 52, 56, 85, 122], [0, 23, 24, 97, 107, 108]]
\item \{1=42336, 2=653184, 3=2953440, 4=3911040\} [[0, 1, 9, 15, 25, 33], [0, 2, 8, 42, 57, 72], [0, 3, 32, 78, 85, 88], [0, 4, 35, 46, 66, 95], [0, 12, 22, 63, 84, 123], [0, 13, 23, 106, 110, 124]]
\item \{1=42336, 2=663012, 3=2940840, 4=3913812\} [[0, 1, 3, 6, 11, 17], [0, 2, 4, 26, 35, 122], [0, 5, 29, 44, 69, 112], [0, 7, 31, 46, 47, 125], [0, 9, 28, 38, 104, 113]]
\item \{1=43344, 2=576828, 3=3022488, 4=3917340\} [[0, 1, 12, 18, 60, 69], [0, 2, 4, 64, 68, 117], [0, 3, 5, 13, 63, 112], [0, 7, 58, 72, 82, 101], [0, 9, 25, 91, 92, 120], [0, 10, 46, 56, 88, 116], [0, 16, 70, 76, 81, 118], [0, 24, 43, 55, 94, 108]]
\item \{1=43344, 2=617652, 3=3043656, 4=3855348\} [[0, 1, 2, 5, 8, 14], [0, 3, 4, 19, 46, 96], [0, 6, 59, 85, 92, 113], [0, 9, 37, 62, 64, 68], [0, 12, 21, 41, 61, 98], [0, 17, 74, 95, 112, 123], [0, 18, 54, 81, 94, 103], [0, 30, 53, 58, 72, 122]]
\item \{1=43344, 2=632016, 3=3025008, 4=3859632\} [[0, 1, 9, 15, 39, 48], [0, 2, 8, 96, 109, 118], [0, 3, 18, 44, 75, 103], [0, 4, 23, 79, 113, 123], [0, 5, 25, 32, 35, 64], [0, 31, 34, 74, 95, 122]]
\item \{1=43344, 2=643356, 3=3017448, 4=3855852\} [[0, 1, 9, 15, 54, 63], [0, 2, 8, 42, 57, 72], [0, 3, 12, 50, 51, 123], [0, 5, 32, 35, 76, 120], [0, 18, 22, 81, 95, 125], [0, 37, 40, 78, 92, 109]]
\item \{1=44352, 2=647136, 3=3003840, 4=3864672\} [[0, 1, 2, 5, 8, 14], [0, 3, 7, 46, 64, 69], [0, 6, 74, 95, 99, 116], [0, 9, 28, 49, 77, 118], [0, 10, 22, 59, 83, 109], [0, 17, 24, 66, 108, 125], [0, 29, 39, 73, 88, 123], [0, 31, 40, 100, 117, 124]]
\item \{1=45360, 2=570780, 3=2995272, 4=3948588\} [[0, 1, 3, 6, 11, 17], [0, 2, 7, 20, 57, 77], [0, 4, 16, 38, 54, 99], [0, 10, 49, 101, 102, 107], [0, 19, 64, 95, 104, 115]]
\item \{1=45360, 2=606312, 3=3012912, 4=3895416\} [[0, 1, 2, 5, 8, 14], [0, 3, 16, 21, 32, 42], [0, 6, 74, 95, 99, 116], [0, 7, 43, 53, 102, 113], [0, 9, 36, 50, 52, 98], [0, 10, 25, 60, 66, 71], [0, 12, 55, 59, 75, 90], [0, 17, 38, 41, 81, 84]]
\item \{1=45360, 2=638820, 3=3018456, 4=3857364\} [[0, 1, 3, 6, 11, 17], [0, 2, 4, 16, 105, 125], [0, 5, 29, 31, 62, 76], [0, 9, 36, 71, 113, 115], [0, 20, 24, 60, 106, 110], [0, 25, 59, 92, 109, 124], [0, 26, 44, 45, 99, 117]]
\item \{1=46368, 2=595728, 3=3005856, 4=3912048\} [[0, 1, 3, 6, 11, 17], [0, 2, 4, 16, 46, 104], [0, 5, 36, 81, 114, 123], [0, 10, 15, 41, 109, 118], [0, 20, 53, 66, 67, 115]]
\item \{1=47376, 2=617652, 3=3031560, 4=3863412\} [[0, 1, 9, 15, 39, 48], [0, 2, 46, 101, 118, 119], [0, 3, 32, 64, 105, 116], [0, 4, 27, 50, 77, 96], [0, 5, 58, 70, 72, 125], [0, 23, 63, 88, 90, 103]]
\item \{1=50400, 2=615384, 3=2969568, 4=3924648\} [[0, 1, 3, 6, 11, 17], [0, 2, 7, 8, 16, 29], [0, 4, 20, 27, 42, 82], [0, 5, 46, 51, 61, 115], [0, 15, 41, 43, 53, 91], [0, 21, 33, 77, 81, 123], [0, 30, 57, 83, 88, 119]]
\item \{1=52416, 2=639576, 3=2992752, 4=3875256\} [[0, 1, 2, 5, 8, 14], [0, 3, 4, 35, 105, 110], [0, 6, 38, 41, 64, 67], [0, 7, 58, 72, 82, 101], [0, 9, 12, 63, 79, 122], [0, 10, 54, 60, 91, 94], [0, 13, 32, 55, 93, 113], [0, 17, 24, 66, 108, 125]]
\item \{0=1260, 1=25200, 2=591192, 3=3003840, 4=3938508\} [[0, 1, 3, 6, 11, 17], [0, 2, 7, 20, 38, 96], [0, 4, 50, 51, 92, 105], [0, 5, 35, 90, 91, 119], [0, 12, 32, 71, 72, 124]]
\item \{0=1260, 1=25200, 2=592704, 3=3022992, 4=3917844\} [[0, 1, 2, 5, 8, 14], [0, 3, 7, 40, 52, 111], [0, 4, 23, 29, 62, 112], [0, 6, 59, 85, 92, 113], [0, 9, 19, 36, 64, 66], [0, 10, 17, 48, 56, 98], [0, 16, 22, 73, 79, 100], [0, 21, 58, 78, 91, 120]]
\item \{0=1260, 1=28224, 2=582120, 3=2975616, 4=3972780\} [[0, 1, 9, 15, 54, 63], [0, 2, 8, 13, 24, 38], [0, 3, 45, 70, 88, 111], [0, 10, 42, 56, 75, 118], [0, 18, 53, 59, 110, 112], [0, 27, 48, 68, 82, 91]]
\item \{0=1260, 1=29232, 2=577584, 3=3015936, 4=3935988\} [[0, 1, 3, 6, 11, 17], [0, 2, 7, 21, 71, 91], [0, 4, 36, 55, 57, 125], [0, 5, 12, 39, 60, 114], [0, 13, 40, 87, 88, 99]]
\item \{0=1260, 1=31248, 2=582120, 3=3037104, 4=3908268\} [[0, 1, 2, 5, 8, 14], [0, 3, 7, 29, 56, 87], [0, 4, 9, 106, 110, 118], [0, 6, 24, 49, 108, 122], [0, 12, 28, 81, 93, 99], [0, 13, 55, 75, 91, 109], [0, 16, 58, 70, 85, 97], [0, 17, 59, 92, 101, 121]]
\item \{0=1260, 1=31248, 2=594972, 3=2993256, 4=3939264\} [[0, 1, 2, 5, 8, 14], [0, 3, 4, 22, 86, 101], [0, 6, 38, 41, 64, 67], [0, 7, 56, 71, 76, 124], [0, 13, 78, 108, 113, 117], [0, 15, 35, 63, 93, 107], [0, 17, 74, 95, 112, 123], [0, 21, 39, 73, 96, 110]]
\item \{0=1260, 1=32256, 2=592704, 3=3033072, 4=3900708\} [[0, 1, 2, 5, 8, 14], [0, 3, 7, 16, 52, 71], [0, 4, 9, 37, 112, 116], [0, 6, 24, 49, 108, 122], [0, 12, 13, 32, 115, 124], [0, 17, 59, 92, 101, 121], [0, 20, 22, 41, 43, 90], [0, 39, 46, 73, 103, 104]]
\item \{0=1260, 1=32256, 2=597240, 3=2999808, 4=3929436\} [[0, 1, 9, 15, 25, 33], [0, 2, 18, 96, 111, 122], [0, 3, 44, 66, 95, 101], [0, 7, 54, 64, 94, 123], [0, 10, 34, 45, 75, 90], [0, 20, 24, 60, 106, 110]]
\item \{0=1260, 1=32256, 2=652428, 3=2988216, 4=3885840\} [[0, 1, 2, 5, 8, 14], [0, 3, 4, 16, 69, 79], [0, 6, 59, 85, 92, 113], [0, 7, 26, 39, 94, 100], [0, 9, 90, 95, 117, 121], [0, 17, 24, 66, 108, 125], [0, 20, 22, 44, 91, 93], [0, 47, 52, 57, 115, 119]]
\item \{0=1260, 1=33264, 2=604800, 3=3018960, 4=3901716\} [[0, 1, 3, 6, 11, 17], [0, 2, 7, 24, 52, 123], [0, 4, 38, 54, 66, 71], [0, 10, 42, 72, 97, 116], [0, 15, 90, 95, 101, 124]]
\item \{0=1260, 1=33264, 2=622944, 3=3031056, 4=3871476\} [[0, 1, 3, 6, 11, 17], [0, 2, 7, 34, 108, 109], [0, 4, 9, 35, 75, 76], [0, 5, 62, 66, 87, 119], [0, 15, 18, 59, 63, 104]]
\item \{0=1260, 1=33264, 2=638064, 3=2976624, 4=3910788\} [[0, 1, 2, 5, 8, 14], [0, 3, 7, 49, 82, 113], [0, 4, 16, 59, 75, 110], [0, 6, 38, 41, 64, 67], [0, 15, 58, 72, 87, 114], [0, 17, 21, 63, 77, 115], [0, 19, 45, 46, 79, 112], [0, 33, 47, 54, 94, 116]]
\item \{0=1260, 1=34272, 2=601776, 3=3007872, 4=3914820\} [[0, 1, 2, 5, 8, 14], [0, 3, 4, 23, 28, 100], [0, 6, 21, 46, 77, 102], [0, 7, 27, 37, 47, 124], [0, 9, 36, 61, 74, 93], [0, 12, 34, 62, 98, 119], [0, 16, 51, 91, 120, 122], [0, 17, 24, 66, 108, 125]]
\item \{0=1260, 1=34272, 2=631260, 3=2984184, 4=3909024\} [[0, 1, 3, 6, 11, 17], [0, 2, 4, 40, 47, 98], [0, 5, 29, 52, 71, 112], [0, 7, 20, 35, 91, 111], [0, 16, 37, 58, 85, 105]]
\item \{0=1260, 1=34272, 2=648648, 3=3080448, 4=3795372\} [[0, 1, 2, 5, 8, 14], [0, 3, 4, 65, 96, 117], [0, 6, 74, 95, 99, 116], [0, 9, 36, 50, 52, 98], [0, 10, 54, 62, 115, 122], [0, 12, 30, 61, 76, 90], [0, 17, 59, 92, 101, 121], [0, 34, 39, 73, 84, 100]]
\item \{0=1260, 1=36288, 2=560196, 3=2963016, 4=3999240\} [[0, 1, 3, 6, 11, 17], [0, 2, 4, 20, 32, 101], [0, 7, 35, 81, 85, 86], [0, 10, 58, 88, 93, 99], [0, 15, 27, 51, 105, 119]]
\item \{0=1260, 1=36288, 2=598752, 3=2965536, 4=3958164\} [[0, 1, 2, 5, 8, 14], [0, 3, 7, 16, 69, 103], [0, 4, 93, 107, 109, 119], [0, 6, 24, 49, 108, 122], [0, 10, 28, 72, 76, 125], [0, 12, 47, 54, 96, 114], [0, 17, 21, 63, 77, 115], [0, 34, 39, 73, 84, 100]]
\item \{0=1260, 1=36288, 2=620676, 3=3005352, 4=3896424\} [[0, 1, 2, 5, 8, 14], [0, 3, 4, 64, 104, 124], [0, 6, 10, 31, 56, 82], [0, 7, 52, 59, 105, 115], [0, 12, 45, 48, 68, 121], [0, 17, 74, 95, 112, 123], [0, 20, 22, 34, 47, 66], [0, 24, 37, 61, 109, 118]]
\item \{0=1260, 1=37296, 2=592704, 3=3040128, 4=3888612\} [[0, 1, 3, 6, 11, 17], [0, 2, 4, 32, 67, 125], [0, 7, 15, 46, 65, 118], [0, 13, 50, 61, 93, 113], [0, 18, 45, 87, 88, 101]]
\item \{0=1260, 1=37296, 2=599508, 3=3023496, 4=3898440\} [[0, 1, 2, 5, 8, 14], [0, 3, 7, 29, 42, 87], [0, 4, 9, 23, 36, 54], [0, 6, 24, 49, 108, 122], [0, 10, 56, 63, 104, 123], [0, 12, 32, 43, 47, 84], [0, 17, 59, 92, 101, 121], [0, 19, 31, 91, 120, 125], [0, 20, 22, 27, 73, 75], [0, 37, 48, 52, 83, 109]]
\item \{0=1260, 1=37296, 2=604800, 3=3019968, 4=3896676\} [[0, 1, 2, 5, 8, 14], [0, 3, 7, 36, 65, 88], [0, 4, 64, 71, 92, 112], [0, 6, 74, 95, 99, 116], [0, 10, 40, 107, 111, 124], [0, 13, 46, 94, 100, 125], [0, 17, 38, 41, 81, 84], [0, 27, 55, 57, 109, 119]]
\item \{0=1260, 1=37296, 2=619164, 3=3014424, 4=3887856\} [[0, 1, 2, 5, 8, 14], [0, 3, 22, 29, 65, 99], [0, 4, 39, 49, 52, 71], [0, 6, 59, 85, 92, 113], [0, 10, 57, 97, 118, 119], [0, 17, 24, 66, 108, 125], [0, 18, 30, 61, 77, 101], [0, 23, 37, 90, 95, 96]]
\item \{0=1260, 1=37296, 2=623700, 3=3019464, 4=3878280\} [[0, 1, 9, 15, 39, 48], [0, 2, 33, 44, 47, 49], [0, 3, 42, 68, 103, 104], [0, 4, 59, 77, 92, 108], [0, 10, 54, 60, 91, 94], [0, 16, 46, 97, 112, 118]]
\item \{0=1260, 1=38304, 2=608580, 3=2999304, 4=3912552\} [[0, 1, 3, 6, 11, 17], [0, 2, 4, 16, 41, 91], [0, 5, 39, 61, 62, 101], [0, 12, 86, 88, 110, 113], [0, 15, 35, 66, 72, 104]]
\item \{0=1260, 1=39312, 2=579096, 3=2988720, 4=3951612\} [[0, 1, 3, 6, 11, 17], [0, 2, 7, 37, 58, 124], [0, 4, 16, 38, 49, 110], [0, 10, 20, 51, 66, 97], [0, 13, 52, 63, 96, 112]]
\item \{0=1260, 1=39312, 2=641844, 3=3004344, 4=3873240\} [[0, 1, 9, 15, 39, 48], [0, 2, 49, 60, 66, 101], [0, 4, 16, 95, 100, 115], [0, 10, 18, 35, 37, 85], [0, 20, 22, 44, 91, 93], [0, 21, 57, 58, 72, 97]]
\item \{0=1260, 1=40320, 2=635796, 3=2941848, 4=3940776\} [[0, 1, 3, 6, 11, 17], [0, 2, 7, 31, 35, 79], [0, 4, 15, 96, 98, 123], [0, 5, 44, 69, 81, 91], [0, 9, 36, 47, 49, 51], [0, 10, 20, 56, 61, 94], [0, 12, 38, 41, 73, 76]]
\item \{0=1260, 1=42336, 2=638820, 3=2995272, 4=3882312\} [[0, 1, 3, 6, 11, 17], [0, 2, 4, 8, 10, 21], [0, 5, 12, 70, 96, 98], [0, 9, 36, 44, 76, 106], [0, 13, 81, 93, 122, 123], [0, 15, 23, 47, 90, 117], [0, 16, 35, 40, 41, 104]]
\item \{0=1260, 1=44352, 2=604044, 3=2958984, 4=3951360\} [[0, 1, 2, 5, 8, 14], [0, 3, 4, 15, 94, 96], [0, 6, 74, 95, 99, 116], [0, 7, 56, 87, 108, 110], [0, 9, 19, 52, 62, 91], [0, 16, 40, 105, 121, 124], [0, 17, 38, 41, 81, 84], [0, 20, 22, 103, 117, 122]]
\item \{0=1260, 1=44352, 2=634284, 3=2972088, 4=3908016\} [[0, 1, 3, 6, 11, 17], [0, 2, 4, 32, 38, 61], [0, 7, 57, 90, 117, 123], [0, 12, 13, 16, 104, 118], [0, 18, 50, 55, 88, 115]]
\item \{0=1260, 1=44352, 2=683424, 3=2975616, 4=3855348\} [[0, 1, 2, 5, 8, 14], [0, 3, 7, 52, 56, 92], [0, 4, 9, 31, 97, 101], [0, 6, 21, 46, 77, 102], [0, 12, 41, 61, 115, 124], [0, 16, 45, 58, 72, 88], [0, 17, 74, 95, 112, 123], [0, 19, 48, 106, 110, 122]]
\item \{0=1260, 1=53424, 2=636552, 3=3012912, 4=3855852\} [[0, 1, 3, 6, 11, 17], [0, 2, 4, 16, 40, 118], [0, 5, 19, 20, 61, 93], [0, 10, 43, 66, 85, 97], [0, 12, 39, 98, 104, 121]]
\item \{0=2520, 1=31248, 2=576072, 3=2947392, 4=4002768\} [[0, 1, 3, 6, 11, 17], [0, 2, 7, 37, 58, 85], [0, 4, 27, 75, 79, 81], [0, 5, 12, 38, 83, 105], [0, 10, 21, 25, 69, 100]]
\end{enumerate}
\end{example}

\begin{example} {\tt SmallGroup(126,10) = C6 x (C7 : C3)}
\begin{enumerate}
\item \{1=18144, 2=621432, 3=3026016, 4=3894408\} [[0, 1, 2, 5, 8, 14], [0, 3, 7, 60, 63, 109], [0, 4, 45, 57, 76, 111], [0, 6, 27, 48, 75, 108], [0, 10, 30, 34, 35, 41]]
\item \{1=19152, 2=584388, 3=2957976, 4=3998484\} [[0, 1, 2, 5, 8, 14], [0, 3, 7, 41, 64, 66], [0, 4, 10, 48, 88, 106], [0, 12, 31, 44, 76, 116], [0, 21, 25, 43, 73, 121]]
\item \{1=19152, 2=611604, 3=2969064, 4=3960180\} [[0, 1, 2, 5, 8, 14], [0, 3, 7, 45, 72, 119], [0, 4, 31, 46, 47, 77], [0, 6, 29, 66, 82, 113], [0, 9, 32, 44, 67, 109]]
\item \{1=20160, 2=513324, 3=2952936, 4=4073580\} [[0, 1, 2, 5, 8, 14], [0, 3, 4, 84, 107, 121], [0, 7, 28, 42, 62, 68], [0, 9, 41, 66, 72, 76], [0, 10, 20, 50, 83, 120]]
\item \{1=20160, 2=579096, 3=2971584, 4=3989160\} [[0, 1, 2, 5, 8, 14], [0, 3, 10, 24, 69, 76], [0, 6, 41, 86, 93, 120], [0, 7, 9, 31, 36, 63], [0, 12, 23, 27, 50, 83], [0, 16, 28, 58, 72, 104], [0, 20, 22, 53, 99, 101]]
\item \{1=20160, 2=588924, 3=2981160, 4=3969756\} [[0, 1, 2, 5, 8, 14], [0, 3, 7, 35, 101, 114], [0, 4, 28, 59, 61, 72], [0, 6, 31, 66, 96, 123], [0, 9, 18, 19, 74, 84]]
\item \{1=21168, 2=526932, 3=2997288, 4=4014612\} [[0, 1, 2, 5, 8, 14], [0, 3, 4, 48, 77, 101], [0, 6, 13, 32, 45, 65], [0, 7, 28, 61, 107, 111], [0, 9, 24, 79, 91, 108], [0, 33, 54, 64, 94, 103], [0, 35, 40, 81, 102, 124]]
\item \{1=21168, 2=542052, 3=2971080, 4=4025700\} [[0, 1, 2, 5, 8, 14], [0, 3, 4, 21, 82, 119], [0, 6, 39, 41, 100, 112], [0, 7, 9, 27, 72, 113], [0, 10, 13, 18, 87, 93]]
\item \{1=21168, 2=542808, 3=2992752, 4=4003272\} [[0, 1, 2, 5, 8, 14], [0, 3, 10, 31, 81, 95], [0, 4, 42, 46, 70, 103], [0, 7, 44, 75, 77, 110], [0, 18, 71, 102, 107, 112]]
\item \{1=21168, 2=557928, 3=2985696, 4=3995208\} [[0, 1, 2, 5, 8, 14], [0, 3, 4, 27, 67, 73], [0, 6, 51, 66, 74, 78], [0, 7, 37, 59, 102, 111], [0, 20, 35, 83, 115, 125]]
\item \{1=21168, 2=560196, 3=2991240, 4=3987396\} [[0, 1, 2, 5, 8, 14], [0, 3, 4, 29, 39, 55], [0, 6, 31, 52, 96, 106], [0, 7, 15, 79, 81, 113], [0, 9, 13, 65, 71, 111]]
\item \{1=21168, 2=573804, 3=2977128, 4=3987900\} [[0, 1, 2, 5, 8, 14], [0, 3, 4, 10, 93, 115], [0, 6, 47, 88, 91, 111], [0, 9, 29, 71, 75, 92], [0, 15, 27, 81, 96, 108]]
\item \{1=21168, 2=600264, 3=2966544, 4=3972024\} [[0, 1, 2, 5, 8, 14], [0, 3, 7, 60, 85, 112], [0, 4, 49, 55, 95, 122], [0, 6, 26, 67, 102, 117], [0, 10, 34, 54, 73, 106]]
\item \{1=21168, 2=618408, 3=2991744, 4=3928680\} [[0, 1, 2, 5, 8, 14], [0, 3, 4, 36, 73, 116], [0, 6, 46, 55, 99, 114], [0, 12, 43, 51, 67, 108], [0, 15, 53, 82, 93, 118]]
\item \{1=22176, 2=534492, 3=3001320, 4=4002012\} [[0, 1, 2, 5, 8, 14], [0, 3, 7, 29, 49, 70], [0, 4, 23, 79, 81, 91], [0, 6, 34, 72, 73, 118], [0, 12, 52, 68, 94, 112]]
\item \{1=22176, 2=559440, 3=3051216, 4=3927168\} [[0, 1, 3, 6, 11, 17], [0, 2, 7, 61, 90, 107], [0, 4, 15, 73, 81, 84], [0, 5, 21, 36, 44, 122], [0, 13, 53, 104, 112, 114]]
\item \{1=22176, 2=563220, 3=3037608, 4=3936996\} [[0, 1, 2, 5, 8, 14], [0, 3, 7, 36, 55, 118], [0, 4, 16, 54, 82, 111], [0, 6, 46, 52, 58, 74], [0, 18, 19, 24, 61, 124]]
\item \{1=22176, 2=576828, 3=2951928, 4=4009068\} [[0, 1, 2, 5, 8, 14], [0, 3, 7, 45, 104, 112], [0, 4, 65, 67, 100, 101], [0, 6, 50, 57, 60, 115], [0, 10, 25, 58, 74, 102]]
\item \{1=22176, 2=577584, 3=2978640, 4=3981600\} [[0, 1, 2, 5, 8, 14], [0, 3, 4, 36, 64, 118], [0, 6, 35, 77, 112, 122], [0, 7, 29, 56, 72, 104], [0, 15, 27, 39, 41, 125]]
\item \{1=22176, 2=582120, 3=2969568, 4=3986136\} [[0, 1, 2, 5, 8, 14], [0, 3, 16, 32, 87, 96], [0, 6, 59, 72, 108, 123], [0, 10, 50, 76, 86, 97], [0, 12, 22, 83, 93, 119]]
\item \{1=22176, 2=586656, 3=2947392, 4=4003776\} [[0, 1, 2, 5, 8, 14], [0, 3, 4, 38, 93, 112], [0, 6, 12, 53, 77, 98], [0, 9, 61, 86, 105, 124], [0, 10, 27, 49, 79, 85]]
\item \{1=22176, 2=593460, 3=2990232, 4=3954132\} [[0, 1, 2, 5, 8, 14], [0, 3, 7, 59, 78, 90], [0, 4, 91, 101, 117, 118], [0, 6, 13, 67, 71, 86], [0, 10, 45, 85, 97, 98]]
\item \{1=22176, 2=600264, 3=2979648, 4=3957912\} [[0, 1, 2, 5, 8, 14], [0, 3, 16, 32, 91, 109], [0, 6, 36, 72, 92, 96], [0, 7, 23, 27, 34, 107], [0, 18, 61, 86, 108, 123]]
\item \{1=22176, 2=614628, 3=2953944, 4=3969252\} [[0, 1, 2, 5, 8, 14], [0, 3, 7, 10, 35, 58], [0, 4, 9, 86, 93, 97], [0, 6, 21, 22, 75, 99], [0, 12, 48, 68, 120, 125]]
\item \{1=22176, 2=625212, 3=2947896, 4=3964716\} [[0, 1, 2, 5, 8, 14], [0, 3, 10, 84, 106, 108], [0, 4, 28, 63, 78, 107], [0, 6, 25, 50, 54, 118], [0, 7, 35, 41, 49, 79]]
\item \{1=23184, 2=562464, 3=3032064, 4=3942288\} [[0, 1, 2, 5, 8, 14], [0, 3, 4, 43, 77, 115], [0, 7, 55, 62, 90, 105], [0, 9, 36, 68, 87, 99], [0, 10, 29, 50, 60, 119], [0, 16, 70, 81, 93, 107], [0, 31, 34, 74, 95, 122]]
\item \{1=23184, 2=578340, 3=3030552, 4=3927924\} [[0, 1, 2, 5, 8, 14], [0, 3, 4, 38, 66, 68], [0, 6, 44, 75, 89, 110], [0, 9, 34, 48, 53, 123], [0, 10, 51, 56, 84, 121], [0, 15, 29, 69, 93, 107], [0, 21, 23, 43, 77, 91]]
\item \{1=23184, 2=581364, 3=2984184, 4=3971268\} [[0, 1, 2, 5, 8, 14], [0, 3, 4, 63, 68, 106], [0, 6, 27, 48, 57, 67], [0, 10, 21, 30, 94, 124], [0, 12, 15, 43, 71, 73]]
\item \{1=23184, 2=586656, 3=2941344, 4=4008816\} [[0, 1, 2, 5, 8, 14], [0, 3, 4, 46, 73, 98], [0, 7, 47, 67, 102, 113], [0, 9, 51, 85, 103, 107], [0, 12, 32, 44, 94, 104]]
\item \{1=23184, 2=586656, 3=2987712, 4=3962448\} [[0, 1, 2, 5, 8, 14], [0, 3, 7, 29, 34, 125], [0, 4, 37, 38, 78, 113], [0, 6, 32, 70, 100, 106], [0, 13, 15, 24, 65, 121]]
\item \{1=23184, 2=592704, 3=2994768, 4=3949344\} [[0, 1, 2, 5, 8, 14], [0, 3, 4, 63, 68, 99], [0, 6, 24, 42, 70, 92], [0, 7, 29, 49, 52, 84], [0, 13, 15, 33, 53, 65]]
\item \{1=23184, 2=621432, 3=2993760, 4=3921624\} [[0, 1, 2, 5, 8, 14], [0, 3, 10, 82, 84, 122], [0, 4, 15, 61, 76, 95], [0, 6, 30, 33, 34, 120], [0, 19, 24, 55, 96, 123]]
\item \{1=23184, 2=627480, 3=3008880, 4=3900456\} [[0, 1, 3, 6, 11, 17], [0, 2, 4, 20, 52, 125], [0, 5, 26, 79, 91, 95], [0, 12, 31, 62, 68, 94], [0, 15, 19, 24, 97, 112]]
\item \{1=23184, 2=643356, 3=3014424, 4=3879036\} [[0, 1, 2, 5, 8, 14], [0, 3, 7, 63, 82, 103], [0, 4, 36, 69, 91, 104], [0, 6, 56, 72, 107, 117], [0, 15, 40, 71, 112, 124], [0, 21, 41, 44, 74, 95], [0, 35, 39, 73, 101, 122]]
\item \{1=24192, 2=536004, 3=3009384, 4=3990420\} [[0, 1, 2, 5, 8, 14], [0, 3, 7, 48, 66, 85], [0, 6, 47, 70, 99, 115], [0, 10, 28, 34, 64, 97], [0, 12, 29, 52, 67, 117]]
\item \{1=24192, 2=545832, 3=2999808, 4=3990168\} [[0, 1, 2, 5, 8, 14], [0, 3, 7, 69, 103, 114], [0, 4, 50, 55, 111, 116], [0, 6, 24, 34, 45, 61], [0, 18, 29, 30, 31, 88]]
\item \{1=24192, 2=557172, 3=2961000, 4=4017636\} [[0, 1, 2, 5, 8, 14], [0, 3, 7, 30, 47, 102], [0, 4, 17, 65, 79, 103], [0, 9, 10, 69, 74, 76], [0, 12, 22, 52, 88, 90]]
\item \{1=24192, 2=562464, 3=3018960, 4=3954384\} [[0, 1, 2, 5, 8, 14], [0, 3, 4, 51, 64, 76], [0, 6, 53, 79, 106, 108], [0, 7, 20, 22, 46, 48], [0, 10, 21, 56, 62, 95], [0, 12, 15, 44, 75, 98], [0, 32, 37, 70, 109, 121]]
\item \{1=24192, 2=565488, 3=3016944, 4=3953376\} [[0, 1, 2, 5, 8, 14], [0, 3, 21, 30, 57, 73], [0, 4, 47, 59, 109, 113], [0, 6, 33, 61, 83, 104], [0, 7, 9, 29, 90, 114]]
\item \{1=24192, 2=567000, 3=2949408, 4=4019400\} [[0, 1, 2, 5, 8, 14], [0, 3, 4, 31, 86, 118], [0, 6, 90, 104, 113, 119], [0, 9, 12, 37, 71, 73], [0, 18, 44, 55, 63, 106]]
\item \{1=24192, 2=567756, 3=2989224, 4=3978828\} [[0, 1, 2, 5, 8, 14], [0, 3, 10, 50, 51, 86], [0, 4, 17, 43, 93, 102], [0, 7, 16, 87, 99, 114], [0, 12, 45, 76, 77, 98]]
\item \{1=24192, 2=585144, 3=2973600, 4=3977064\} [[0, 1, 2, 5, 8, 14], [0, 3, 4, 10, 26, 115], [0, 6, 19, 27, 49, 65], [0, 9, 32, 34, 75, 122], [0, 15, 21, 76, 106, 117]]
\item \{1=24192, 2=592704, 3=2984688, 4=3958416\} [[0, 1, 2, 5, 8, 14], [0, 3, 7, 10, 40, 121], [0, 4, 21, 38, 77, 95], [0, 6, 26, 86, 87, 88], [0, 9, 27, 75, 81, 99], [0, 16, 18, 20, 22, 63], [0, 28, 32, 69, 106, 110]]
\item \{1=24192, 2=595728, 3=2993760, 4=3946320\} [[0, 1, 2, 5, 8, 14], [0, 3, 4, 31, 39, 79], [0, 7, 64, 66, 110, 114], [0, 9, 50, 55, 95, 102], [0, 10, 34, 45, 49, 59]]
\item \{1=24192, 2=595728, 3=3005856, 4=3934224\} [[0, 1, 2, 5, 8, 14], [0, 3, 10, 29, 36, 87], [0, 4, 30, 38, 101, 114], [0, 6, 19, 48, 84, 116], [0, 7, 43, 51, 65, 81]]
\item \{1=24192, 2=601020, 3=2980152, 4=3954636\} [[0, 1, 2, 5, 8, 14], [0, 3, 4, 16, 123, 125], [0, 6, 12, 47, 52, 119], [0, 9, 24, 79, 91, 108], [0, 13, 22, 37, 73, 109], [0, 15, 35, 63, 76, 118], [0, 18, 46, 60, 100, 124]]
\item \{1=24192, 2=604800, 3=3001824, 4=3929184\} [[0, 1, 2, 5, 8, 14], [0, 3, 7, 30, 48, 86], [0, 6, 62, 72, 100, 102], [0, 9, 19, 43, 74, 97], [0, 12, 55, 106, 122, 123]]
\item \{1=24192, 2=609336, 3=3024000, 4=3902472\} [[0, 1, 2, 5, 8, 14], [0, 3, 7, 46, 68, 84], [0, 4, 31, 43, 97, 117], [0, 10, 61, 72, 93, 107], [0, 13, 38, 40, 83, 123], [0, 15, 35, 63, 76, 118], [0, 48, 51, 74, 95, 125]]
\item \{1=25200, 2=535248, 3=3024000, 4=3975552\} [[0, 1, 3, 6, 11, 17], [0, 2, 4, 20, 78, 101], [0, 5, 41, 81, 111, 112], [0, 7, 69, 85, 95, 117], [0, 10, 97, 103, 104, 108]]
\item \{1=25200, 2=548100, 3=2997288, 4=3989412\} [[0, 1, 2, 5, 8, 14], [0, 3, 7, 38, 60, 90], [0, 4, 30, 76, 101, 117], [0, 9, 35, 37, 57, 102], [0, 12, 22, 77, 83, 114]]
\item \{1=25200, 2=575316, 3=3005352, 4=3954132\} [[0, 1, 2, 5, 8, 14], [0, 3, 4, 47, 63, 117], [0, 6, 13, 61, 72, 95], [0, 7, 28, 49, 98, 120], [0, 10, 42, 77, 83, 87]]
\item \{1=25200, 2=575316, 3=3028536, 4=3930948\} [[0, 1, 2, 5, 8, 14], [0, 3, 4, 58, 68, 102], [0, 6, 39, 62, 91, 124], [0, 7, 22, 50, 53, 116], [0, 15, 21, 51, 70, 106]]
\item \{1=25200, 2=576072, 3=2998800, 4=3959928\} [[0, 1, 2, 5, 8, 14], [0, 3, 10, 40, 82, 104], [0, 4, 22, 81, 92, 115], [0, 6, 7, 53, 101, 112], [0, 13, 28, 33, 47, 61]]
\item \{1=25200, 2=581364, 3=2968056, 4=3985380\} [[0, 1, 2, 5, 8, 14], [0, 3, 4, 27, 33, 85], [0, 7, 15, 29, 56, 66], [0, 9, 36, 70, 83, 101], [0, 12, 28, 35, 50, 115], [0, 21, 31, 52, 77, 99], [0, 37, 41, 78, 106, 110]]
\item \{1=25200, 2=584388, 3=2968056, 4=3982356\} [[0, 1, 2, 5, 8, 14], [0, 3, 10, 63, 66, 72], [0, 4, 21, 40, 45, 103], [0, 6, 42, 43, 65, 88], [0, 12, 32, 54, 96, 118]]
\item \{1=25200, 2=588168, 3=3026016, 4=3920616\} [[0, 1, 2, 5, 8, 14], [0, 3, 4, 58, 82, 112], [0, 6, 26, 54, 61, 75], [0, 7, 56, 64, 103, 113], [0, 18, 31, 50, 78, 120]]
\item \{1=25200, 2=590436, 3=2998296, 4=3946068\} [[0, 1, 2, 5, 8, 14], [0, 3, 7, 46, 90, 93], [0, 4, 43, 45, 69, 94], [0, 9, 44, 56, 68, 111], [0, 13, 48, 71, 83, 97]]
\item \{1=25200, 2=596484, 3=2977128, 4=3961188\} [[0, 1, 2, 5, 8, 14], [0, 3, 7, 38, 51, 119], [0, 4, 16, 39, 95, 109], [0, 9, 18, 29, 36, 48], [0, 10, 56, 69, 101, 125], [0, 13, 28, 40, 93, 98], [0, 15, 35, 63, 76, 118]]
\item \{1=25200, 2=597996, 3=3016440, 4=3920364\} [[0, 1, 2, 5, 8, 14], [0, 3, 7, 48, 66, 85], [0, 6, 41, 77, 83, 119], [0, 10, 26, 74, 97, 101], [0, 12, 52, 63, 117, 118]]
\item \{1=25200, 2=599508, 3=2964024, 4=3971268\} [[0, 1, 2, 5, 8, 14], [0, 3, 4, 16, 75, 86], [0, 6, 19, 40, 43, 121], [0, 7, 52, 67, 100, 101], [0, 10, 61, 69, 95, 112]]
\item \{1=25200, 2=619920, 3=3002832, 4=3912048\} [[0, 1, 2, 5, 8, 14], [0, 3, 7, 68, 90, 106], [0, 6, 13, 46, 81, 121], [0, 9, 36, 50, 69, 82], [0, 10, 20, 58, 72, 96], [0, 16, 19, 74, 95, 114], [0, 18, 37, 78, 79, 98]]
\item \{1=25200, 2=624456, 3=2963520, 4=3946824\} [[0, 1, 2, 5, 8, 14], [0, 3, 4, 45, 60, 110], [0, 6, 23, 26, 66, 94], [0, 7, 68, 107, 111, 114], [0, 9, 51, 67, 83, 87]]
\item \{1=25200, 2=647892, 3=2963016, 4=3923892\} [[0, 1, 2, 5, 8, 14], [0, 3, 7, 69, 90, 93], [0, 4, 22, 71, 117, 124], [0, 6, 51, 55, 91, 115], [0, 9, 12, 44, 59, 112]]
\item \{1=26208, 2=531468, 3=3052728, 4=3949596\} [[0, 1, 2, 5, 8, 14], [0, 3, 21, 38, 69, 114], [0, 6, 22, 24, 97, 113], [0, 9, 42, 78, 119, 121], [0, 10, 30, 67, 73, 116]]
\item \{1=26208, 2=535248, 3=2981664, 4=4016880\} [[0, 1, 2, 5, 8, 14], [0, 3, 16, 22, 107, 115], [0, 4, 51, 64, 98, 119], [0, 6, 52, 58, 99, 104], [0, 7, 21, 23, 56, 106]]
\item \{1=26208, 2=539784, 3=2956464, 4=4037544\} [[0, 1, 2, 5, 8, 14], [0, 3, 4, 16, 81, 106], [0, 6, 39, 50, 76, 102], [0, 10, 33, 40, 43, 59], [0, 13, 32, 72, 95, 113]]
\item \{1=26208, 2=546588, 3=3022488, 4=3964716\} [[0, 1, 2, 5, 8, 14], [0, 3, 7, 50, 87, 112], [0, 4, 28, 73, 95, 96], [0, 6, 19, 24, 77, 101], [0, 12, 13, 31, 75, 114]]
\item \{1=26208, 2=554904, 3=2965536, 4=4013352\} [[0, 1, 2, 5, 8, 14], [0, 3, 10, 24, 78, 116], [0, 6, 37, 38, 68, 114], [0, 7, 55, 91, 110, 117], [0, 12, 19, 22, 79, 98]]
\item \{1=26208, 2=559440, 3=2984688, 4=3989664\} [[0, 1, 3, 6, 11, 17], [0, 2, 4, 26, 50, 64], [0, 5, 25, 32, 39, 109], [0, 7, 54, 61, 66, 107], [0, 9, 24, 79, 91, 108], [0, 19, 23, 67, 90, 121], [0, 28, 51, 76, 103, 118]]
\item \{1=26208, 2=563220, 3=2965032, 4=4005540\} [[0, 1, 2, 5, 8, 14], [0, 3, 10, 60, 72, 114], [0, 4, 6, 62, 63, 97], [0, 9, 36, 52, 65, 84], [0, 15, 53, 86, 108, 121], [0, 19, 48, 91, 120, 122], [0, 21, 30, 51, 77, 98]]
\item \{1=26208, 2=563976, 3=3020976, 4=3948840\} [[0, 1, 2, 5, 8, 14], [0, 3, 4, 36, 63, 117], [0, 6, 13, 43, 56, 57], [0, 15, 21, 50, 74, 87], [0, 16, 42, 67, 94, 108]]
\item \{1=26208, 2=564732, 3=2998296, 4=3970764\} [[0, 1, 2, 5, 8, 14], [0, 3, 21, 38, 60, 87], [0, 6, 26, 74, 108, 121], [0, 9, 10, 29, 64, 79], [0, 20, 47, 57, 100, 115]]
\item \{1=26208, 2=572292, 3=3011400, 4=3950100\} [[0, 1, 2, 5, 8, 14], [0, 3, 4, 95, 114, 117], [0, 6, 45, 54, 83, 96], [0, 9, 12, 53, 76, 90], [0, 10, 19, 75, 109, 122]]
\item \{1=26208, 2=577584, 3=2990736, 4=3965472\} [[0, 1, 2, 5, 8, 14], [0, 3, 10, 72, 74, 86], [0, 4, 69, 102, 107, 122], [0, 6, 42, 48, 52, 85], [0, 9, 49, 62, 103, 111]]
\item \{1=26208, 2=579096, 3=3050208, 4=3904488\} [[0, 1, 2, 5, 8, 14], [0, 3, 16, 36, 60, 112], [0, 4, 33, 34, 55, 67], [0, 7, 21, 27, 79, 109], [0, 12, 66, 76, 100, 110]]
\item \{1=26208, 2=588168, 3=2986704, 4=3958920\} [[0, 1, 2, 5, 8, 14], [0, 3, 21, 33, 46, 125], [0, 4, 20, 68, 87, 95], [0, 7, 45, 60, 65, 81], [0, 9, 35, 37, 57, 120]]
\item \{1=26208, 2=588924, 3=2962008, 4=3982860\} [[0, 1, 2, 5, 8, 14], [0, 3, 4, 26, 50, 93], [0, 6, 21, 43, 67, 94], [0, 7, 48, 58, 69, 110], [0, 10, 46, 49, 83, 87]]
\item \{1=26208, 2=588924, 3=2966040, 4=3978828\} [[0, 1, 2, 5, 8, 14], [0, 3, 7, 55, 77, 78], [0, 4, 59, 71, 117, 122], [0, 6, 39, 93, 95, 125], [0, 10, 45, 81, 91, 111]]
\item \{1=26208, 2=596484, 3=2950920, 4=3986388\} [[0, 1, 2, 5, 8, 14], [0, 3, 16, 42, 45, 108], [0, 4, 29, 31, 34, 88], [0, 6, 7, 50, 61, 118], [0, 12, 43, 94, 106, 125]]
\item \{1=26208, 2=599508, 3=2992248, 4=3942036\} [[0, 1, 2, 5, 8, 14], [0, 3, 10, 81, 85, 103], [0, 4, 18, 22, 51, 95], [0, 7, 46, 70, 82, 114], [0, 9, 35, 66, 72, 118]]
\item \{1=26208, 2=604044, 3=2971080, 4=3958668\} [[0, 1, 2, 5, 8, 14], [0, 3, 4, 35, 60, 66], [0, 7, 15, 92, 119, 120], [0, 10, 41, 51, 117, 124], [0, 21, 31, 52, 77, 99], [0, 23, 37, 90, 95, 96], [0, 32, 76, 82, 104, 118]]
\item \{1=26208, 2=604044, 3=2985192, 4=3944556\} [[0, 1, 2, 5, 8, 14], [0, 3, 4, 98, 107, 121], [0, 7, 16, 40, 54, 115], [0, 9, 47, 58, 77, 108], [0, 10, 73, 81, 106, 114]]
\item \{1=26208, 2=609336, 3=2989728, 4=3934728\} [[0, 1, 2, 5, 8, 14], [0, 3, 21, 64, 99, 111], [0, 4, 19, 38, 60, 76], [0, 6, 31, 74, 78, 119], [0, 9, 27, 46, 67, 104]]
\item \{1=26208, 2=611604, 3=2992248, 4=3929940\} [[0, 1, 2, 5, 8, 14], [0, 3, 4, 82, 121, 124], [0, 6, 37, 61, 103, 107], [0, 9, 21, 50, 51, 125], [0, 12, 19, 43, 75, 94]]
\item \{1=26208, 2=614628, 3=2997288, 4=3921876\} [[0, 1, 2, 5, 8, 14], [0, 3, 4, 55, 71, 96], [0, 7, 66, 103, 108, 113], [0, 10, 30, 77, 87, 119], [0, 12, 54, 84, 100, 115]]
\item \{1=26208, 2=614628, 3=3023496, 4=3895668\} [[0, 1, 2, 5, 8, 14], [0, 3, 4, 27, 93, 125], [0, 6, 45, 69, 79, 90], [0, 7, 10, 55, 111, 121], [0, 16, 33, 74, 98, 123]]
\item \{1=26208, 2=615384, 3=2943360, 4=3975048\} [[0, 1, 2, 5, 8, 14], [0, 3, 16, 21, 107, 124], [0, 4, 71, 116, 117, 122], [0, 6, 36, 44, 63, 123], [0, 13, 47, 69, 79, 95]]
\item \{1=26208, 2=628236, 3=3062808, 4=3842748\} [[0, 1, 2, 5, 8, 14], [0, 3, 7, 64, 87, 111], [0, 4, 20, 27, 67, 99], [0, 6, 83, 90, 101, 116], [0, 10, 45, 92, 96, 109]]
\item \{1=26208, 2=633528, 3=2989728, 4=3910536\} [[0, 1, 2, 5, 8, 14], [0, 3, 7, 37, 69, 108], [0, 4, 34, 38, 71, 75], [0, 6, 22, 82, 93, 125], [0, 9, 88, 94, 106, 120]]
\item \{1=26208, 2=642600, 3=3021984, 4=3869208\} [[0, 1, 2, 5, 8, 14], [0, 3, 7, 41, 51, 66], [0, 4, 9, 106, 110, 118], [0, 6, 39, 43, 68, 84], [0, 10, 13, 74, 95, 108], [0, 15, 32, 61, 86, 90], [0, 47, 57, 69, 102, 119]]
\item \{1=26208, 2=652428, 3=2940840, 4=3940524\} [[0, 1, 2, 5, 8, 14], [0, 3, 4, 19, 84, 116], [0, 6, 16, 30, 52, 58], [0, 10, 21, 49, 96, 108], [0, 15, 33, 53, 115, 118]]
\item \{1=27216, 2=552636, 3=3022488, 4=3957660\} [[0, 1, 2, 5, 8, 14], [0, 3, 22, 39, 63, 124], [0, 6, 10, 51, 64, 100], [0, 7, 42, 62, 111, 113], [0, 18, 86, 93, 107, 125], [0, 21, 31, 52, 77, 99], [0, 28, 32, 69, 106, 110]]
\item \{1=27216, 2=560952, 3=2981664, 4=3990168\} [[0, 1, 2, 5, 8, 14], [0, 3, 7, 25, 48, 115], [0, 4, 40, 85, 122, 123], [0, 10, 13, 67, 86, 87], [0, 18, 30, 100, 106, 110], [0, 20, 22, 41, 43, 90], [0, 24, 60, 72, 108, 109]]
\item \{1=27216, 2=564732, 3=3001320, 4=3966732\} [[0, 1, 2, 5, 8, 14], [0, 3, 7, 104, 113, 115], [0, 4, 17, 42, 65, 99], [0, 9, 26, 29, 62, 74], [0, 10, 69, 72, 77, 122]]
\item \{1=27216, 2=573048, 3=3002832, 4=3956904\} [[0, 1, 2, 5, 8, 14], [0, 3, 7, 51, 86, 122], [0, 4, 15, 26, 62, 91], [0, 6, 36, 49, 78, 110], [0, 16, 73, 84, 95, 112]]
\item \{1=27216, 2=573048, 3=3005856, 4=3953880\} [[0, 1, 2, 5, 8, 14], [0, 3, 4, 55, 60, 100], [0, 6, 44, 73, 88, 95], [0, 7, 46, 49, 101, 108], [0, 10, 15, 35, 40, 98]]
\item \{1=27216, 2=579096, 3=2996784, 4=3956904\} [[0, 1, 2, 5, 8, 14], [0, 3, 4, 23, 29, 43], [0, 6, 32, 33, 58, 77], [0, 7, 60, 64, 76, 91], [0, 13, 22, 85, 92, 104]]
\item \{1=27216, 2=581364, 3=2979144, 4=3972276\} [[0, 1, 2, 5, 8, 14], [0, 3, 10, 25, 77, 82], [0, 4, 55, 69, 102, 108], [0, 6, 60, 66, 100, 104], [0, 7, 26, 35, 59, 121]]
\item \{1=27216, 2=581364, 3=3019464, 4=3931956\} [[0, 1, 2, 5, 8, 14], [0, 3, 7, 72, 92, 108], [0, 6, 23, 74, 96, 115], [0, 9, 44, 66, 81, 99], [0, 10, 12, 30, 41, 103]]
\item \{1=27216, 2=583632, 3=2989728, 4=3959424\} [[0, 1, 2, 5, 8, 14], [0, 3, 10, 30, 40, 96], [0, 4, 38, 64, 104, 114], [0, 7, 16, 71, 107, 112], [0, 13, 55, 60, 102, 125]]
\item \{1=27216, 2=583632, 3=3012912, 4=3936240\} [[0, 1, 2, 5, 8, 14], [0, 3, 4, 77, 81, 119], [0, 6, 36, 83, 98, 104], [0, 7, 20, 47, 71, 90], [0, 13, 16, 40, 96, 122]]
\item \{1=27216, 2=586656, 3=3031056, 4=3915072\} [[0, 1, 2, 5, 8, 14], [0, 3, 25, 51, 100, 118], [0, 4, 36, 43, 48, 79], [0, 10, 27, 72, 105, 113], [0, 15, 39, 47, 73, 117], [0, 16, 19, 74, 95, 114], [0, 31, 35, 54, 94, 112]]
\item \{1=27216, 2=589680, 3=2996784, 4=3946320\} [[0, 1, 2, 5, 8, 14], [0, 3, 10, 64, 115, 122], [0, 4, 59, 77, 92, 108], [0, 6, 33, 42, 110, 117], [0, 9, 19, 88, 93, 125], [0, 16, 54, 69, 94, 99], [0, 23, 37, 90, 95, 96]]
\item \{1=27216, 2=589680, 3=3012912, 4=3930192\} [[0, 1, 2, 5, 8, 14], [0, 3, 10, 66, 91, 92], [0, 4, 9, 84, 93, 112], [0, 6, 26, 69, 75, 109], [0, 12, 16, 35, 79, 113]]
\item \{1=27216, 2=594216, 3=3003840, 4=3934728\} [[0, 1, 2, 5, 8, 14], [0, 3, 7, 38, 52, 106], [0, 4, 17, 59, 69, 114], [0, 9, 43, 54, 95, 98], [0, 13, 30, 60, 99, 124]]
\item \{1=27216, 2=594972, 3=2990232, 4=3947580\} [[0, 1, 2, 5, 8, 14], [0, 3, 7, 58, 95, 116], [0, 4, 9, 49, 50, 109], [0, 10, 34, 63, 77, 122], [0, 12, 62, 86, 90, 93]]
\item \{1=27216, 2=598752, 3=2998800, 4=3935232\} [[0, 1, 2, 5, 8, 14], [0, 3, 10, 54, 58, 121], [0, 4, 16, 55, 82, 92], [0, 6, 25, 30, 44, 111], [0, 19, 31, 106, 110, 125], [0, 20, 22, 100, 102, 104], [0, 24, 69, 81, 108, 115]]
\item \{1=27216, 2=601020, 3=3015432, 4=3916332\} [[0, 1, 2, 5, 8, 14], [0, 3, 4, 73, 99, 122], [0, 6, 51, 55, 78, 86], [0, 12, 21, 76, 81, 92], [0, 15, 34, 60, 72, 91]]
\item \{1=27216, 2=605556, 3=2967048, 4=3960180\} [[0, 1, 2, 5, 8, 14], [0, 3, 16, 49, 63, 116], [0, 4, 35, 56, 78, 99], [0, 6, 7, 45, 55, 121], [0, 9, 21, 26, 73, 77], [0, 29, 70, 84, 106, 110], [0, 31, 34, 74, 95, 122]]
\item \{1=27216, 2=611604, 3=2998296, 4=3922884\} [[0, 1, 2, 5, 8, 14], [0, 3, 7, 25, 48, 64], [0, 4, 30, 35, 65, 125], [0, 9, 26, 40, 44, 117], [0, 10, 63, 87, 109, 119]]
\item \{1=27216, 2=616896, 3=3019968, 4=3895920\} [[0, 1, 2, 5, 8, 14], [0, 3, 7, 40, 94, 123], [0, 4, 30, 50, 111, 122], [0, 6, 35, 59, 62, 90], [0, 20, 39, 51, 97, 118]]
\item \{1=27216, 2=618408, 3=2987712, 4=3926664\} [[0, 1, 3, 6, 11, 17], [0, 2, 4, 16, 26, 84], [0, 5, 56, 61, 81, 92], [0, 9, 85, 100, 105, 120], [0, 19, 20, 40, 78, 110]]
\item \{1=27216, 2=625212, 3=3018456, 4=3889116\} [[0, 1, 2, 5, 8, 14], [0, 3, 7, 56, 78, 96], [0, 6, 41, 77, 83, 93], [0, 9, 49, 72, 88, 118], [0, 15, 39, 47, 73, 117], [0, 20, 22, 53, 99, 101], [0, 24, 30, 104, 108, 112]]
\item \{1=27216, 2=634284, 3=2965032, 4=3933468\} [[0, 1, 2, 5, 8, 14], [0, 3, 4, 41, 100, 115], [0, 6, 24, 69, 106, 113], [0, 7, 23, 57, 76, 96], [0, 13, 53, 65, 72, 107]]
\item \{1=28224, 2=531468, 3=2984184, 4=4016124\} [[0, 1, 2, 5, 8, 14], [0, 3, 16, 60, 74, 100], [0, 4, 33, 62, 97, 118], [0, 9, 37, 53, 90, 103], [0, 15, 24, 88, 99, 108], [0, 19, 39, 73, 84, 113], [0, 20, 22, 50, 52, 98]]
\item \{1=28224, 2=544320, 3=3007872, 4=3979584\} [[0, 1, 2, 5, 8, 14], [0, 3, 16, 30, 66, 100], [0, 4, 39, 43, 88, 106], [0, 6, 35, 42, 72, 86], [0, 7, 22, 50, 81, 103]]
\item \{1=28224, 2=548856, 3=2968560, 4=4014360\} [[0, 1, 2, 5, 8, 14], [0, 3, 16, 45, 46, 99], [0, 4, 19, 26, 56, 125], [0, 9, 53, 60, 84, 94], [0, 20, 23, 38, 41, 79], [0, 31, 34, 74, 95, 122], [0, 32, 55, 63, 75, 105]]
\item \{1=28224, 2=563976, 3=3009888, 4=3957912\} [[0, 1, 2, 5, 8, 14], [0, 3, 7, 25, 52, 103], [0, 4, 37, 86, 106, 124], [0, 10, 45, 67, 70, 112], [0, 16, 41, 60, 74, 98]]
\item \{1=28224, 2=569268, 3=3030552, 4=3931956\} [[0, 1, 2, 5, 8, 14], [0, 3, 4, 21, 69, 107], [0, 6, 24, 76, 91, 98], [0, 7, 44, 50, 55, 111], [0, 13, 64, 86, 90, 121]]
\item \{1=28224, 2=570024, 3=3016944, 4=3944808\} [[0, 1, 2, 5, 8, 14], [0, 3, 4, 52, 54, 59], [0, 6, 36, 72, 89, 103], [0, 15, 40, 50, 55, 118], [0, 16, 29, 98, 100, 119]]
\item \{1=28224, 2=575316, 3=2996280, 4=3960180\} [[0, 1, 2, 5, 8, 14], [0, 3, 21, 72, 93, 112], [0, 4, 22, 84, 108, 121], [0, 6, 33, 59, 65, 104], [0, 7, 52, 73, 88, 124]]
\item \{1=28224, 2=579096, 3=3026016, 4=3926664\} [[0, 1, 2, 5, 8, 14], [0, 3, 10, 37, 64, 95], [0, 4, 33, 38, 78, 82], [0, 9, 42, 70, 86, 100], [0, 12, 52, 65, 112, 114]]
\item \{1=28224, 2=585144, 3=3034080, 4=3912552\} [[0, 1, 2, 5, 8, 14], [0, 3, 7, 10, 59, 102], [0, 4, 67, 96, 101, 112], [0, 6, 35, 57, 79, 90], [0, 9, 25, 40, 55, 115]]
\item \{1=28224, 2=585900, 3=3014424, 4=3931452\} [[0, 1, 2, 5, 8, 14], [0, 3, 7, 16, 68, 76], [0, 4, 36, 43, 57, 88], [0, 10, 28, 29, 64, 117], [0, 21, 31, 50, 60, 94]]
\item \{1=28224, 2=586656, 3=2988720, 4=3956400\} [[0, 1, 2, 5, 8, 14], [0, 3, 16, 69, 92, 125], [0, 4, 56, 71, 79, 123], [0, 7, 22, 26, 95, 109], [0, 9, 32, 73, 78, 86]]
\item \{1=28224, 2=589680, 3=2996784, 4=3945312\} [[0, 1, 2, 5, 8, 14], [0, 3, 4, 70, 105, 118], [0, 6, 36, 57, 92, 110], [0, 10, 58, 100, 112, 117], [0, 15, 27, 37, 74, 114]]
\item \{1=28224, 2=591948, 3=3067848, 4=3871980\} [[0, 1, 2, 5, 8, 14], [0, 3, 29, 47, 82, 101], [0, 6, 41, 91, 95, 122], [0, 7, 30, 35, 97, 104], [0, 10, 56, 60, 93, 124], [0, 18, 76, 114, 115, 118], [0, 20, 22, 68, 70, 112]]
\item \{1=28224, 2=597996, 3=2970072, 4=3963708\} [[0, 1, 2, 5, 8, 14], [0, 3, 4, 103, 104, 108], [0, 6, 12, 74, 79, 115], [0, 10, 26, 62, 90, 101], [0, 18, 22, 64, 119, 121]]
\item \{1=28224, 2=600264, 3=2994768, 4=3936744\} [[0, 1, 2, 5, 8, 14], [0, 3, 4, 49, 81, 117], [0, 6, 46, 73, 106, 124], [0, 7, 53, 57, 95, 118], [0, 9, 42, 50, 90, 92]]
\item \{1=28224, 2=602532, 3=3004344, 4=3924900\} [[0, 1, 2, 5, 8, 14], [0, 3, 29, 66, 82, 101], [0, 4, 88, 99, 103, 115], [0, 6, 7, 61, 67, 104], [0, 10, 21, 27, 74, 77], [0, 15, 19, 54, 94, 98], [0, 33, 37, 47, 109, 116]]
\item \{1=28224, 2=607068, 3=3012408, 4=3912300\} [[0, 1, 2, 5, 8, 14], [0, 3, 7, 36, 51, 121], [0, 4, 30, 41, 60, 84], [0, 12, 43, 59, 72, 102], [0, 15, 35, 66, 81, 108]]
\item \{1=28224, 2=611604, 3=2988216, 4=3931956\} [[0, 1, 2, 5, 8, 14], [0, 3, 4, 19, 55, 94], [0, 6, 42, 82, 111, 119], [0, 9, 18, 78, 101, 121], [0, 13, 28, 84, 100, 116]]
\item \{1=28224, 2=616896, 3=3003840, 4=3911040\} [[0, 1, 2, 5, 8, 14], [0, 3, 7, 65, 84, 102], [0, 4, 10, 18, 87, 116], [0, 9, 38, 49, 93, 94], [0, 12, 21, 32, 103, 124]]
\item \{1=28224, 2=619920, 3=2991744, 4=3920112\} [[0, 1, 2, 5, 8, 14], [0, 3, 4, 19, 49, 77], [0, 6, 39, 58, 60, 101], [0, 9, 78, 100, 111, 124], [0, 20, 27, 81, 85, 123]]
\item \{1=28224, 2=627480, 3=2962512, 4=3941784\} [[0, 1, 2, 5, 8, 14], [0, 3, 10, 24, 29, 58], [0, 6, 40, 70, 77, 112], [0, 7, 43, 49, 51, 121], [0, 13, 52, 71, 84, 98]]
\item \{1=28224, 2=628236, 3=2954952, 4=3948588\} [[0, 1, 2, 5, 8, 14], [0, 3, 4, 33, 48, 68], [0, 7, 34, 87, 91, 93], [0, 9, 27, 54, 76, 118], [0, 10, 37, 56, 79, 110], [0, 13, 24, 64, 71, 116], [0, 20, 22, 86, 88, 121]]
\item \{1=28224, 2=628236, 3=3002328, 4=3901212\} [[0, 1, 2, 5, 8, 14], [0, 3, 10, 100, 119, 125], [0, 4, 22, 26, 38, 70], [0, 7, 35, 36, 64, 104], [0, 16, 21, 33, 74, 102]]
\item \{1=28224, 2=628992, 3=3004848, 4=3897936\} [[0, 1, 2, 5, 8, 14], [0, 3, 4, 18, 48, 88], [0, 9, 52, 67, 90, 96], [0, 10, 33, 91, 109, 118], [0, 13, 20, 38, 66, 103]]
\item \{1=28224, 2=633528, 3=2997792, 4=3900456\} [[0, 1, 2, 5, 8, 14], [0, 3, 4, 37, 45, 79], [0, 6, 26, 52, 66, 99], [0, 9, 53, 86, 101, 106], [0, 10, 43, 83, 84, 90]]
\item \{1=28224, 2=651672, 3=2952432, 4=3927672\} [[0, 1, 2, 5, 8, 14], [0, 3, 10, 37, 63, 119], [0, 4, 47, 70, 79, 84], [0, 12, 55, 65, 94, 113], [0, 18, 19, 88, 112, 120]]
\item \{1=28224, 2=656208, 3=3017952, 4=3857616\} [[0, 1, 2, 5, 8, 14], [0, 3, 7, 45, 55, 98], [0, 4, 20, 35, 79, 86], [0, 6, 69, 106, 108, 109], [0, 9, 62, 68, 91, 111]]
\item \{1=28224, 2=669816, 3=2986704, 4=3875256\} [[0, 1, 2, 5, 8, 14], [0, 3, 4, 56, 85, 96], [0, 6, 29, 34, 38, 111], [0, 15, 40, 53, 55, 92], [0, 18, 37, 88, 112, 123]]
\item \{1=29232, 2=517104, 3=2999808, 4=4013856\} [[0, 1, 2, 5, 8, 14], [0, 3, 4, 37, 52, 58], [0, 6, 13, 69, 99, 119], [0, 7, 15, 73, 112, 115], [0, 12, 60, 81, 110, 125]]
\item \{1=29232, 2=542052, 3=2978136, 4=4010580\} [[0, 1, 2, 5, 8, 14], [0, 3, 16, 113, 123, 125], [0, 4, 30, 89, 90, 110], [0, 6, 41, 44, 102, 103], [0, 9, 32, 60, 75, 109]]
\item \{1=29232, 2=551124, 3=2992248, 4=3987396\} [[0, 1, 2, 5, 8, 14], [0, 3, 7, 57, 63, 90], [0, 4, 17, 64, 88, 108], [0, 9, 75, 94, 112, 120], [0, 16, 34, 77, 81, 114]]
\item \{1=29232, 2=558684, 3=2995272, 4=3976812\} [[0, 1, 2, 5, 8, 14], [0, 3, 7, 82, 86, 87], [0, 4, 42, 67, 96, 107], [0, 6, 45, 62, 72, 120], [0, 10, 21, 49, 100, 122]]
\item \{1=29232, 2=558684, 3=3006360, 4=3965724\} [[0, 1, 3, 6, 11, 17], [0, 2, 4, 19, 23, 118], [0, 5, 32, 42, 115, 125], [0, 9, 43, 51, 97, 103], [0, 15, 29, 69, 93, 107], [0, 20, 22, 53, 99, 101], [0, 31, 34, 74, 95, 122]]
\item \{1=29232, 2=559440, 3=2971584, 4=3999744\} [[0, 1, 2, 5, 8, 14], [0, 3, 4, 71, 77, 90], [0, 6, 10, 51, 89, 118], [0, 9, 33, 52, 74, 108], [0, 15, 27, 68, 99, 117]]
\item \{1=29232, 2=560952, 3=2987712, 4=3982104\} [[0, 1, 2, 5, 8, 14], [0, 3, 10, 29, 54, 125], [0, 4, 36, 76, 92, 99], [0, 7, 41, 95, 104, 115], [0, 12, 48, 53, 84, 85]]
\item \{1=29232, 2=565488, 3=3013920, 4=3951360\} [[0, 1, 2, 5, 8, 14], [0, 3, 4, 48, 104, 119], [0, 6, 35, 67, 74, 89], [0, 9, 16, 70, 91, 117], [0, 10, 19, 37, 73, 124]]
\item \{1=29232, 2=571536, 3=3047184, 4=3912048\} [[0, 1, 2, 5, 8, 14], [0, 3, 21, 25, 65, 93], [0, 4, 23, 41, 48, 115], [0, 7, 35, 59, 61, 72], [0, 15, 33, 68, 74, 119]]
\item \{1=29232, 2=573048, 3=2981664, 4=3976056\} [[0, 1, 2, 5, 8, 14], [0, 3, 7, 65, 81, 91], [0, 4, 22, 61, 73, 87], [0, 9, 10, 13, 71, 99], [0, 15, 34, 66, 115, 124]]
\item \{1=29232, 2=573804, 3=3013416, 4=3943548\} [[0, 1, 2, 5, 8, 14], [0, 3, 4, 52, 93, 125], [0, 6, 48, 106, 111, 120], [0, 7, 15, 84, 88, 115], [0, 16, 22, 26, 68, 83]]
\item \{1=29232, 2=576828, 3=2982168, 4=3971772\} [[0, 1, 2, 5, 8, 14], [0, 3, 16, 45, 113, 115], [0, 4, 33, 39, 56, 79], [0, 7, 29, 34, 97, 114], [0, 18, 70, 93, 112, 124]]
\item \{1=29232, 2=579852, 3=2983176, 4=3967740\} [[0, 1, 2, 5, 8, 14], [0, 3, 7, 23, 35, 65], [0, 4, 99, 104, 106, 119], [0, 12, 54, 96, 107, 115], [0, 18, 22, 37, 88, 95]]
\item \{1=29232, 2=581364, 3=3004344, 4=3945060\} [[0, 1, 2, 5, 8, 14], [0, 3, 4, 52, 58, 101], [0, 6, 46, 79, 82, 114], [0, 7, 31, 78, 92, 99], [0, 9, 74, 76, 111, 125]]
\item \{1=29232, 2=582876, 3=3023496, 4=3924396\} [[0, 1, 2, 5, 8, 14], [0, 3, 4, 27, 34, 108], [0, 6, 16, 45, 91, 113], [0, 9, 32, 43, 63, 106], [0, 12, 23, 29, 115, 122]]
\item \{1=29232, 2=583632, 3=2972592, 4=3974544\} [[0, 1, 2, 5, 8, 14], [0, 3, 4, 26, 29, 75], [0, 6, 35, 48, 114, 120], [0, 7, 37, 65, 105, 107], [0, 10, 52, 101, 117, 121]]
\item \{1=29232, 2=585144, 3=3000816, 4=3944808\} [[0, 1, 2, 5, 8, 14], [0, 3, 10, 74, 112, 118], [0, 4, 30, 47, 49, 92], [0, 9, 19, 65, 70, 104], [0, 13, 20, 33, 87, 90]]
\item \{1=29232, 2=585900, 3=3009384, 4=3935484\} [[0, 1, 2, 5, 8, 14], [0, 3, 7, 19, 46, 47], [0, 6, 68, 84, 119, 123], [0, 9, 36, 43, 56, 75], [0, 13, 52, 83, 96, 109], [0, 15, 38, 41, 70, 122], [0, 23, 35, 85, 90, 125]]
\item \{1=29232, 2=588168, 3=2993760, 4=3948840\} [[0, 1, 2, 5, 8, 14], [0, 3, 29, 66, 113, 120], [0, 4, 48, 65, 107, 112], [0, 6, 19, 52, 55, 78], [0, 7, 58, 72, 85, 98], [0, 10, 21, 56, 62, 95], [0, 20, 22, 92, 94, 96]]
\item \{1=29232, 2=588924, 3=2968056, 4=3973788\} [[0, 1, 2, 5, 8, 14], [0, 3, 7, 48, 94, 107], [0, 4, 42, 85, 112, 123], [0, 6, 13, 37, 53, 110], [0, 9, 39, 50, 56, 76]]
\item \{1=29232, 2=588924, 3=2995272, 4=3946572\} [[0, 1, 3, 6, 11, 17], [0, 2, 4, 20, 61, 123], [0, 5, 19, 64, 89, 94], [0, 10, 46, 95, 111, 114], [0, 25, 34, 53, 55, 125]]
\item \{1=29232, 2=588924, 3=3023496, 4=3918348\} [[0, 1, 3, 6, 11, 17], [0, 2, 7, 10, 70, 96], [0, 5, 19, 45, 48, 62], [0, 9, 64, 66, 98, 125], [0, 16, 57, 105, 119, 121], [0, 23, 50, 51, 90, 113], [0, 39, 42, 58, 73, 92]]
\item \{1=29232, 2=589680, 3=2943360, 4=3997728\} [[0, 1, 2, 5, 8, 14], [0, 3, 7, 39, 64, 86], [0, 4, 17, 56, 96, 101], [0, 9, 32, 42, 68, 111], [0, 20, 60, 76, 81, 124]]
\item \{1=29232, 2=590436, 3=3013416, 4=3926916\} [[0, 1, 2, 5, 8, 14], [0, 3, 4, 47, 60, 110], [0, 6, 12, 32, 92, 109], [0, 9, 29, 49, 53, 73], [0, 10, 52, 62, 105, 120]]
\item \{1=29232, 2=591192, 3=2990736, 4=3948840\} [[0, 1, 2, 5, 8, 14], [0, 3, 16, 50, 78, 120], [0, 4, 38, 58, 89, 111], [0, 9, 10, 18, 88, 112], [0, 24, 28, 87, 103, 105]]
\item \{1=29232, 2=599508, 3=2975112, 4=3956148\} [[0, 1, 2, 5, 8, 14], [0, 3, 4, 33, 41, 118], [0, 7, 52, 59, 63, 71], [0, 9, 12, 47, 88, 125], [0, 13, 57, 95, 104, 113]]
\item \{1=29232, 2=601020, 3=3006360, 4=3923388\} [[0, 1, 2, 5, 8, 14], [0, 3, 16, 30, 83, 100], [0, 6, 50, 51, 52, 86], [0, 7, 9, 43, 98, 119], [0, 18, 32, 59, 79, 115]]
\item \{1=29232, 2=602532, 3=2954952, 4=3973284\} [[0, 1, 2, 5, 8, 14], [0, 3, 4, 45, 69, 100], [0, 6, 24, 34, 61, 84], [0, 7, 55, 75, 82, 101], [0, 9, 12, 74, 95, 107], [0, 10, 63, 83, 117, 120], [0, 23, 27, 76, 90, 124]]
\item \{1=29232, 2=602532, 3=3020472, 4=3907764\} [[0, 1, 2, 5, 8, 14], [0, 3, 25, 29, 83, 101], [0, 4, 34, 47, 79, 89], [0, 6, 38, 48, 77, 125], [0, 10, 26, 91, 120, 124], [0, 13, 22, 37, 73, 109], [0, 16, 51, 106, 110, 122]]
\item \{1=29232, 2=606312, 3=2963520, 4=3960936\} [[0, 1, 2, 5, 8, 14], [0, 3, 7, 56, 108, 118], [0, 4, 39, 61, 101, 103], [0, 6, 25, 45, 72, 94], [0, 13, 57, 70, 99, 116]]
\item \{1=29232, 2=609336, 3=3011904, 4=3909528\} [[0, 1, 3, 6, 11, 17], [0, 2, 4, 31, 99, 119], [0, 5, 46, 51, 52, 125], [0, 7, 16, 86, 104, 118], [0, 12, 45, 60, 70, 81]]
\item \{1=29232, 2=620676, 3=3008376, 4=3901716\} [[0, 1, 2, 5, 8, 14], [0, 3, 29, 48, 90, 125], [0, 4, 46, 57, 72, 108], [0, 6, 44, 68, 89, 112], [0, 9, 10, 70, 105, 109]]
\item \{1=29232, 2=621432, 3=2997792, 4=3911544\} [[0, 1, 2, 5, 8, 14], [0, 3, 4, 39, 108, 125], [0, 7, 53, 56, 69, 110], [0, 9, 46, 59, 119, 123], [0, 13, 61, 93, 112, 117]]
\item \{1=29232, 2=625212, 3=3028536, 4=3877020\} [[0, 1, 2, 5, 8, 14], [0, 3, 10, 30, 37, 108], [0, 4, 60, 72, 104, 105], [0, 6, 41, 68, 103, 119], [0, 12, 22, 47, 75, 113]]
\item \{1=29232, 2=650916, 3=2975112, 4=3904740\} [[0, 1, 2, 5, 8, 14], [0, 3, 7, 63, 104, 122], [0, 4, 37, 58, 97, 103], [0, 6, 16, 22, 52, 101], [0, 10, 71, 95, 99, 121]]
\item \{1=29232, 2=655452, 3=2979144, 4=3896172\} [[0, 1, 2, 5, 8, 14], [0, 3, 21, 30, 78, 114], [0, 4, 39, 86, 116, 122], [0, 6, 51, 52, 65, 68], [0, 7, 22, 35, 42, 124]]
\item \{1=30240, 2=535248, 3=3000816, 4=3993696\} [[0, 1, 2, 5, 8, 14], [0, 3, 16, 21, 73, 84], [0, 4, 19, 46, 108, 121], [0, 6, 30, 58, 66, 116], [0, 15, 55, 59, 107, 109]]
\item \{1=30240, 2=552636, 3=3053736, 4=3923388\} [[0, 1, 2, 5, 8, 14], [0, 3, 4, 10, 76, 89], [0, 6, 27, 46, 93, 104], [0, 9, 68, 83, 98, 105], [0, 13, 61, 63, 84, 106]]
\item \{1=30240, 2=553392, 3=2988720, 4=3987648\} [[0, 1, 2, 5, 8, 14], [0, 3, 16, 45, 74, 96], [0, 4, 6, 35, 38, 115], [0, 7, 49, 85, 111, 122], [0, 10, 72, 73, 100, 119]]
\item \{1=30240, 2=557172, 3=2971080, 4=4001508\} [[0, 1, 2, 5, 8, 14], [0, 3, 7, 47, 99, 122], [0, 4, 21, 36, 85, 100], [0, 6, 10, 42, 64, 115], [0, 18, 70, 71, 98, 101]]
\item \{1=30240, 2=561708, 3=2958984, 4=4009068\} [[0, 1, 2, 5, 8, 14], [0, 3, 4, 106, 108, 117], [0, 6, 23, 58, 63, 69], [0, 7, 28, 53, 88, 94], [0, 15, 42, 75, 99, 119]]
\item \{1=30240, 2=562464, 3=3009888, 4=3957408\} [[0, 1, 2, 5, 8, 14], [0, 3, 4, 18, 43, 48], [0, 9, 55, 62, 65, 109], [0, 10, 28, 61, 85, 92], [0, 12, 29, 88, 107, 115]]
\item \{1=30240, 2=573048, 3=3018960, 4=3937752\} [[0, 1, 2, 5, 8, 14], [0, 3, 10, 57, 103, 121], [0, 4, 28, 33, 43, 59], [0, 7, 44, 93, 106, 111], [0, 9, 18, 62, 63, 100]]
\item \{1=30240, 2=576828, 3=3044664, 4=3908268\} [[0, 1, 2, 5, 8, 14], [0, 3, 4, 60, 99, 119], [0, 6, 38, 52, 85, 98], [0, 7, 22, 56, 63, 109], [0, 12, 23, 51, 53, 122]]
\item \{1=30240, 2=579852, 3=2980152, 4=3969756\} [[0, 1, 2, 5, 8, 14], [0, 3, 10, 42, 78, 86], [0, 4, 65, 88, 101, 117], [0, 6, 58, 68, 97, 125], [0, 19, 30, 81, 84, 92]]
\item \{1=30240, 2=581364, 3=2978136, 4=3970260\} [[0, 1, 2, 5, 8, 14], [0, 3, 16, 72, 108, 115], [0, 4, 6, 31, 71, 85], [0, 7, 37, 64, 109, 123], [0, 9, 32, 39, 104, 120], [0, 15, 38, 41, 70, 122], [0, 20, 22, 77, 79, 118]]
\item \{1=30240, 2=583632, 3=3007872, 4=3938256\} [[0, 1, 2, 5, 8, 14], [0, 3, 4, 35, 86, 112], [0, 6, 49, 102, 104, 111], [0, 9, 59, 64, 88, 118], [0, 13, 22, 61, 83, 90]]
\item \{1=30240, 2=588168, 3=2993760, 4=3947832\} [[0, 1, 2, 5, 8, 14], [0, 3, 10, 67, 103, 108], [0, 4, 18, 27, 96, 112], [0, 7, 55, 75, 82, 101], [0, 9, 12, 54, 84, 121], [0, 25, 37, 38, 109, 110], [0, 29, 58, 66, 72, 105]]
\item \{1=30240, 2=594972, 3=2954952, 4=3979836\} [[0, 1, 2, 5, 8, 14], [0, 3, 4, 98, 105, 113], [0, 6, 43, 81, 85, 109], [0, 9, 36, 52, 65, 84], [0, 10, 21, 27, 74, 77], [0, 12, 32, 69, 96, 122], [0, 13, 39, 91, 119, 120]]
\item \{1=30240, 2=594972, 3=2967048, 4=3967740\} [[0, 1, 2, 5, 8, 14], [0, 3, 4, 92, 106, 115], [0, 6, 10, 67, 96, 107], [0, 7, 44, 45, 73, 108], [0, 13, 32, 52, 79, 112]]
\item \{1=30240, 2=596484, 3=2965032, 4=3968244\} [[0, 1, 3, 6, 11, 17], [0, 2, 4, 16, 40, 85], [0, 5, 13, 73, 106, 119], [0, 10, 52, 68, 72, 108], [0, 15, 25, 34, 54, 114]]
\item \{1=30240, 2=600264, 3=2993760, 4=3935736\} [[0, 1, 3, 6, 11, 17], [0, 2, 4, 20, 67, 97], [0, 5, 26, 47, 61, 120], [0, 7, 43, 71, 90, 103], [0, 10, 35, 37, 66, 111]]
\item \{1=30240, 2=625212, 3=3003336, 4=3901212\} [[0, 1, 2, 5, 8, 14], [0, 3, 7, 49, 84, 107], [0, 4, 15, 41, 61, 124], [0, 6, 42, 48, 96, 118], [0, 13, 28, 31, 71, 117]]
\item \{1=30240, 2=626724, 3=3043656, 4=3859380\} [[0, 1, 3, 6, 11, 17], [0, 2, 4, 20, 72, 101], [0, 5, 12, 29, 55, 102], [0, 7, 53, 91, 98, 114], [0, 13, 15, 60, 104, 106]]
\item \{1=30240, 2=629748, 3=3021480, 4=3878532\} [[0, 1, 2, 5, 8, 14], [0, 3, 7, 60, 64, 101], [0, 4, 16, 36, 71, 72], [0, 6, 70, 81, 99, 117], [0, 10, 37, 40, 87, 114]]
\item \{1=30240, 2=638064, 3=2980656, 4=3911040\} [[0, 1, 2, 5, 8, 14], [0, 3, 7, 10, 94, 107], [0, 4, 34, 35, 91, 110], [0, 6, 27, 30, 64, 96], [0, 12, 43, 47, 101, 124]]
\item \{1=30240, 2=638064, 3=3017952, 4=3873744\} [[0, 1, 2, 5, 8, 14], [0, 3, 10, 54, 114, 123], [0, 4, 40, 70, 77, 102], [0, 7, 58, 71, 93, 117], [0, 13, 30, 79, 95, 118]]
\item \{1=30240, 2=668304, 3=3013920, 4=3847536\} [[0, 1, 2, 5, 8, 14], [0, 3, 25, 65, 115, 121], [0, 4, 93, 107, 109, 119], [0, 6, 45, 69, 72, 106], [0, 7, 26, 87, 91, 100], [0, 23, 27, 76, 90, 124], [0, 48, 51, 74, 95, 125]]
\item \{1=30240, 2=678132, 3=2979144, 4=3872484\} [[0, 1, 2, 5, 8, 14], [0, 3, 7, 73, 108, 120], [0, 4, 49, 84, 102, 117], [0, 6, 53, 85, 116, 119], [0, 9, 26, 68, 70, 115]]
\item \{1=31248, 2=558684, 3=2984184, 4=3985884\} [[0, 1, 2, 5, 8, 14], [0, 3, 7, 78, 91, 109], [0, 4, 17, 85, 106, 108], [0, 10, 35, 71, 87, 123], [0, 16, 51, 90, 116, 124]]
\item \{1=31248, 2=559440, 3=2993760, 4=3975552\} [[0, 1, 2, 5, 8, 14], [0, 3, 7, 69, 85, 125], [0, 4, 16, 41, 73, 96], [0, 6, 45, 55, 59, 112], [0, 12, 13, 32, 95, 115]]
\item \{1=31248, 2=565488, 3=3020976, 4=3942288\} [[0, 1, 2, 5, 8, 14], [0, 3, 16, 46, 93, 112], [0, 4, 9, 39, 89, 97], [0, 10, 60, 62, 83, 117], [0, 15, 34, 63, 74, 88]]
\item \{1=31248, 2=566244, 3=2987208, 4=3975300\} [[0, 1, 2, 5, 8, 14], [0, 3, 4, 43, 91, 98], [0, 6, 32, 63, 72, 76], [0, 9, 12, 47, 52, 94], [0, 10, 60, 75, 116, 117]]
\item \{1=31248, 2=566244, 3=2991240, 4=3971268\} [[0, 1, 2, 5, 8, 14], [0, 3, 7, 39, 108, 115], [0, 4, 57, 60, 96, 100], [0, 6, 44, 92, 117, 120], [0, 12, 82, 103, 104, 107]]
\item \{1=31248, 2=571536, 3=3014928, 4=3942288\} [[0, 1, 2, 5, 8, 14], [0, 3, 7, 40, 46, 58], [0, 6, 34, 75, 79, 106], [0, 10, 43, 59, 94, 121], [0, 16, 66, 70, 76, 96]]
\item \{1=31248, 2=571536, 3=3021984, 4=3935232\} [[0, 1, 2, 5, 8, 14], [0, 3, 4, 35, 36, 59], [0, 6, 13, 77, 99, 124], [0, 7, 26, 38, 106, 113], [0, 12, 32, 44, 102, 125]]
\item \{1=31248, 2=581364, 3=2987208, 4=3960180\} [[0, 1, 2, 5, 8, 14], [0, 3, 29, 36, 47, 106], [0, 6, 48, 89, 99, 114], [0, 13, 24, 45, 61, 71], [0, 16, 25, 72, 85, 100]]
\item \{1=31248, 2=582120, 3=2950416, 4=3996216\} [[0, 1, 2, 5, 8, 14], [0, 3, 21, 63, 113, 122], [0, 7, 64, 103, 106, 125], [0, 9, 27, 29, 55, 67], [0, 10, 33, 40, 41, 79]]
\item \{1=31248, 2=582120, 3=3041136, 4=3905496\} [[0, 1, 2, 5, 8, 14], [0, 3, 4, 40, 98, 105], [0, 6, 81, 84, 109, 121], [0, 10, 12, 70, 90, 110], [0, 13, 52, 65, 96, 122]]
\item \{1=31248, 2=582876, 3=3002328, 4=3943548\} [[0, 1, 2, 5, 8, 14], [0, 3, 10, 40, 54, 121], [0, 4, 16, 55, 86, 106], [0, 6, 38, 69, 79, 125], [0, 7, 22, 59, 95, 102]]
\item \{1=31248, 2=584388, 3=2985192, 4=3959172\} [[0, 1, 2, 5, 8, 14], [0, 3, 4, 75, 86, 95], [0, 6, 40, 61, 65, 111], [0, 9, 49, 62, 83, 109], [0, 15, 24, 38, 67, 121]]
\item \{1=31248, 2=586656, 3=2960496, 4=3981600\} [[0, 1, 2, 5, 8, 14], [0, 3, 22, 29, 81, 125], [0, 4, 17, 23, 92, 93], [0, 7, 30, 42, 46, 87], [0, 10, 41, 51, 111, 123]]
\item \{1=31248, 2=588924, 3=2954952, 4=3984876\} [[0, 1, 2, 5, 8, 14], [0, 3, 7, 38, 93, 119], [0, 4, 52, 67, 71, 91], [0, 6, 10, 21, 70, 96], [0, 12, 39, 43, 46, 47]]
\item \{1=31248, 2=591948, 3=3002328, 4=3934476\} [[0, 1, 2, 5, 8, 14], [0, 3, 16, 32, 45, 116], [0, 6, 46, 69, 96, 113], [0, 7, 38, 62, 110, 111], [0, 18, 47, 54, 94, 123], [0, 20, 22, 29, 31, 33], [0, 24, 52, 64, 102, 108]]
\item \{1=31248, 2=597996, 3=2991240, 4=3939516\} [[0, 1, 2, 5, 8, 14], [0, 3, 16, 37, 93, 103], [0, 4, 42, 49, 59, 75], [0, 6, 26, 56, 64, 114], [0, 7, 57, 61, 63, 120]]
\item \{1=31248, 2=598752, 3=2984688, 4=3945312\} [[0, 1, 2, 5, 8, 14], [0, 3, 21, 45, 63, 98], [0, 7, 61, 69, 93, 99], [0, 9, 81, 83, 85, 120], [0, 10, 37, 43, 106, 109], [0, 12, 38, 41, 73, 76], [0, 16, 19, 74, 95, 114]]
\item \{1=31248, 2=600264, 3=2992752, 4=3935736\} [[0, 1, 2, 5, 8, 14], [0, 3, 10, 46, 110, 114], [0, 4, 15, 17, 59, 90], [0, 7, 37, 64, 109, 123], [0, 9, 36, 70, 83, 101], [0, 12, 19, 61, 62, 120], [0, 21, 30, 51, 77, 98]]
\item \{1=31248, 2=601776, 3=3039120, 4=3887856\} [[0, 1, 2, 5, 8, 14], [0, 3, 10, 63, 72, 114], [0, 4, 17, 26, 61, 86], [0, 9, 27, 66, 91, 98], [0, 15, 40, 71, 112, 124], [0, 16, 21, 35, 77, 83], [0, 25, 37, 38, 109, 110]]
\item \{1=31248, 2=602532, 3=2997288, 4=3928932\} [[0, 1, 2, 5, 8, 14], [0, 3, 4, 60, 86, 105], [0, 6, 48, 92, 109, 119], [0, 9, 26, 68, 103, 106], [0, 10, 21, 27, 74, 77], [0, 15, 19, 54, 94, 98], [0, 16, 34, 91, 120, 125]]
\item \{1=31248, 2=606312, 3=3002832, 4=3919608\} [[0, 1, 2, 5, 8, 14], [0, 3, 16, 46, 100, 117], [0, 4, 60, 64, 89, 111], [0, 9, 34, 36, 47, 66], [0, 10, 28, 56, 70, 102], [0, 13, 24, 98, 103, 107], [0, 20, 22, 95, 97, 124]]
\item \{1=31248, 2=611604, 3=2971080, 4=3946068\} [[0, 1, 2, 5, 8, 14], [0, 3, 21, 72, 101, 119], [0, 4, 28, 87, 102, 117], [0, 6, 40, 56, 64, 71], [0, 7, 37, 49, 53, 88]]
\item \{1=31248, 2=616896, 3=3019968, 4=3891888\} [[0, 1, 2, 5, 8, 14], [0, 3, 7, 22, 57, 76], [0, 4, 18, 50, 55, 77], [0, 10, 54, 82, 83, 122], [0, 12, 23, 29, 53, 109]]
\item \{1=31248, 2=617652, 3=2973096, 4=3938004\} [[0, 1, 2, 5, 8, 14], [0, 3, 7, 69, 85, 102], [0, 6, 57, 59, 74, 98], [0, 10, 22, 25, 35, 93], [0, 13, 15, 65, 84, 94]]
\item \{1=31248, 2=617652, 3=2977128, 4=3933972\} [[0, 1, 3, 6, 11, 17], [0, 2, 4, 16, 20, 42], [0, 5, 24, 29, 104, 120], [0, 12, 44, 66, 101, 110], [0, 23, 33, 74, 96, 115]]
\item \{1=31248, 2=620676, 3=2942856, 4=3965220\} [[0, 1, 2, 5, 8, 14], [0, 3, 7, 69, 98, 110], [0, 4, 56, 71, 78, 120], [0, 6, 13, 32, 45, 67], [0, 9, 46, 47, 99, 104]]
\item \{1=31248, 2=622188, 3=2981160, 4=3925404\} [[0, 1, 2, 5, 8, 14], [0, 3, 16, 33, 64, 100], [0, 4, 30, 62, 97, 120], [0, 6, 31, 43, 71, 122], [0, 15, 69, 77, 81, 119]]
\item \{1=31248, 2=625968, 3=2970576, 4=3932208\} [[0, 1, 2, 5, 8, 14], [0, 3, 16, 33, 60, 64], [0, 4, 15, 66, 75, 123], [0, 6, 13, 31, 37, 114], [0, 12, 59, 92, 94, 119], [0, 21, 23, 79, 90, 120], [0, 24, 69, 81, 108, 115]]
\item \{1=31248, 2=628992, 3=2993760, 4=3906000\} [[0, 1, 2, 5, 8, 14], [0, 3, 21, 99, 100, 114], [0, 6, 12, 30, 69, 98], [0, 13, 33, 48, 64, 86], [0, 15, 24, 26, 90, 124]]
\item \{1=31248, 2=650916, 3=2982168, 4=3895668\} [[0, 1, 2, 5, 8, 14], [0, 3, 16, 29, 112, 121], [0, 6, 31, 77, 104, 114], [0, 7, 26, 53, 68, 106], [0, 9, 90, 91, 118, 122]]
\item \{1=31248, 2=652428, 3=3001320, 4=3875004\} [[0, 1, 2, 5, 8, 14], [0, 3, 4, 75, 105, 107], [0, 6, 23, 52, 95, 122], [0, 9, 18, 67, 101, 112], [0, 13, 45, 79, 82, 118]]
\item \{1=31248, 2=663012, 3=2989224, 4=3876516\} [[0, 1, 3, 6, 11, 17], [0, 2, 4, 43, 64, 68], [0, 5, 58, 66, 71, 98], [0, 7, 50, 82, 86, 123], [0, 12, 24, 28, 51, 104]]
\item \{1=32256, 2=529200, 3=3047184, 4=3951360\} [[0, 1, 2, 5, 8, 14], [0, 3, 4, 41, 52, 122], [0, 6, 33, 69, 99, 109], [0, 7, 23, 88, 91, 121], [0, 16, 19, 46, 98, 125]]
\item \{1=32256, 2=561708, 3=3029544, 4=3936492\} [[0, 1, 2, 5, 8, 14], [0, 3, 7, 56, 72, 92], [0, 4, 9, 49, 82, 120], [0, 6, 51, 102, 104, 106], [0, 18, 37, 61, 101, 112]]
\item \{1=32256, 2=563976, 3=3030048, 4=3933720\} [[0, 1, 2, 5, 8, 14], [0, 3, 4, 10, 109, 121], [0, 6, 26, 35, 49, 65], [0, 7, 21, 31, 57, 93], [0, 12, 45, 69, 79, 116]]
\item \{1=32256, 2=571536, 3=3002832, 4=3953376\} [[0, 1, 2, 5, 8, 14], [0, 3, 16, 30, 100, 118], [0, 4, 35, 61, 88, 114], [0, 10, 63, 66, 110, 125], [0, 20, 22, 95, 97, 124], [0, 23, 46, 90, 103, 104], [0, 27, 55, 57, 109, 119]]
\item \{1=32256, 2=575316, 3=3035592, 4=3916836\} [[0, 1, 2, 5, 8, 14], [0, 3, 10, 37, 41, 123], [0, 4, 33, 58, 60, 113], [0, 7, 36, 51, 53, 95], [0, 12, 24, 28, 102, 120]]
\item \{1=32256, 2=582876, 3=3033576, 4=3911292\} [[0, 1, 2, 5, 8, 14], [0, 3, 4, 29, 37, 69], [0, 7, 42, 53, 55, 67], [0, 9, 41, 43, 57, 65], [0, 13, 45, 95, 115, 118]]
\item \{1=32256, 2=585144, 3=2976624, 4=3965976\} [[0, 1, 2, 5, 8, 14], [0, 3, 7, 72, 91, 92], [0, 4, 63, 87, 106, 111], [0, 6, 22, 38, 57, 82], [0, 16, 27, 73, 88, 101]]
\item \{1=32256, 2=586656, 3=3002832, 4=3938256\} [[0, 1, 2, 5, 8, 14], [0, 3, 21, 51, 78, 92], [0, 4, 16, 37, 67, 113], [0, 6, 33, 55, 82, 117], [0, 15, 41, 63, 79, 107]]
\item \{1=32256, 2=588924, 3=2984184, 4=3954636\} [[0, 1, 2, 5, 8, 14], [0, 3, 10, 30, 109, 124], [0, 4, 48, 50, 51, 85], [0, 6, 52, 76, 81, 97], [0, 9, 13, 24, 64, 104]]
\item \{1=32256, 2=588924, 3=3033576, 4=3905244\} [[0, 1, 2, 5, 8, 14], [0, 3, 7, 63, 65, 76], [0, 4, 9, 42, 48, 71], [0, 6, 32, 67, 69, 116], [0, 10, 29, 62, 73, 84]]
\item \{1=32256, 2=590436, 3=2970072, 4=3967236\} [[0, 1, 2, 5, 8, 14], [0, 3, 4, 19, 23, 48], [0, 9, 12, 62, 67, 103], [0, 20, 26, 57, 68, 73], [0, 29, 32, 72, 109, 122]]
\item \{1=32256, 2=594216, 3=2940336, 4=3993192\} [[0, 1, 2, 5, 8, 14], [0, 3, 4, 39, 41, 68], [0, 6, 36, 62, 66, 116], [0, 12, 51, 59, 70, 118], [0, 16, 20, 33, 71, 86]]
\item \{1=32256, 2=594216, 3=2977632, 4=3955896\} [[0, 1, 2, 5, 8, 14], [0, 3, 21, 57, 64, 113], [0, 4, 23, 26, 67, 91], [0, 6, 13, 60, 63, 116], [0, 7, 42, 104, 122, 123]]
\item \{1=32256, 2=594216, 3=2997792, 4=3935736\} [[0, 1, 2, 5, 8, 14], [0, 3, 16, 32, 36, 96], [0, 6, 44, 75, 117, 118], [0, 7, 30, 35, 78, 106], [0, 10, 12, 93, 108, 123]]
\item \{1=32256, 2=597996, 3=3000312, 4=3929436\} [[0, 1, 2, 5, 8, 14], [0, 3, 4, 68, 70, 82], [0, 6, 22, 37, 91, 108], [0, 7, 36, 47, 57, 94], [0, 10, 33, 88, 95, 111]]
\item \{1=32256, 2=599508, 3=3028536, 4=3899700\} [[0, 1, 2, 5, 8, 14], [0, 3, 4, 19, 43, 84], [0, 6, 50, 100, 110, 122], [0, 10, 56, 63, 104, 123], [0, 13, 37, 57, 93, 119], [0, 15, 29, 32, 72, 113], [0, 31, 35, 54, 94, 112]]
\item \{1=32256, 2=601020, 3=2993256, 4=3933468\} [[0, 1, 2, 5, 8, 14], [0, 3, 16, 50, 96, 106], [0, 4, 33, 34, 88, 103], [0, 7, 23, 28, 85, 114], [0, 9, 19, 36, 49, 81], [0, 10, 61, 72, 93, 107], [0, 21, 31, 52, 77, 99]]
\item \{1=32256, 2=602532, 3=2995272, 4=3929940\} [[0, 1, 2, 5, 8, 14], [0, 3, 4, 26, 101, 106], [0, 6, 49, 61, 75, 95], [0, 7, 22, 54, 82, 100], [0, 19, 38, 92, 93, 98]]
\item \{1=32256, 2=604800, 3=2955456, 4=3967488\} [[0, 1, 2, 5, 8, 14], [0, 3, 29, 30, 48, 122], [0, 6, 47, 60, 70, 71], [0, 7, 54, 81, 94, 116], [0, 12, 64, 101, 117, 121], [0, 13, 20, 22, 55, 57], [0, 16, 21, 35, 77, 83]]
\item \{1=32256, 2=606312, 3=2980656, 4=3940776\} [[0, 1, 2, 5, 8, 14], [0, 3, 4, 98, 101, 103], [0, 6, 23, 53, 74, 79], [0, 7, 9, 84, 91, 94], [0, 12, 61, 69, 93, 106]]
\item \{1=32256, 2=606312, 3=3004848, 4=3916584\} [[0, 1, 2, 5, 8, 14], [0, 3, 7, 23, 25, 125], [0, 4, 38, 61, 102, 113], [0, 13, 53, 82, 103, 118], [0, 15, 29, 51, 70, 117]]
\item \{1=32256, 2=608580, 3=2972088, 4=3947076\} [[0, 1, 2, 5, 8, 14], [0, 3, 21, 42, 82, 122], [0, 4, 18, 29, 71, 75], [0, 7, 47, 101, 102, 123], [0, 9, 33, 37, 61, 68]]
\item \{1=32256, 2=609336, 3=3001824, 4=3916584\} [[0, 1, 2, 5, 8, 14], [0, 3, 4, 23, 99, 115], [0, 6, 24, 37, 42, 113], [0, 7, 27, 40, 75, 77], [0, 13, 33, 88, 120, 125]]
\item \{1=32256, 2=610092, 3=3064824, 4=3852828\} [[0, 1, 2, 5, 8, 14], [0, 3, 16, 25, 60, 75], [0, 4, 37, 70, 89, 90], [0, 6, 64, 94, 122, 123], [0, 9, 35, 36, 67, 98], [0, 20, 22, 100, 102, 104], [0, 21, 28, 77, 85, 105]]
\item \{1=32256, 2=612360, 3=2981664, 4=3933720\} [[0, 1, 2, 5, 8, 14], [0, 3, 4, 56, 84, 123], [0, 6, 34, 36, 61, 110], [0, 7, 70, 76, 95, 96], [0, 12, 59, 92, 94, 119], [0, 15, 37, 51, 65, 109], [0, 20, 22, 29, 31, 33]]
\item \{1=32256, 2=619164, 3=2981160, 4=3927420\} [[0, 1, 2, 5, 8, 14], [0, 3, 4, 34, 103, 114], [0, 6, 38, 65, 94, 124], [0, 9, 60, 77, 85, 92], [0, 12, 48, 79, 84, 95]]
\item \{1=32256, 2=625212, 3=3003336, 4=3899196\} [[0, 1, 2, 5, 8, 14], [0, 3, 7, 40, 47, 60], [0, 4, 26, 58, 84, 97], [0, 9, 35, 38, 62, 113], [0, 16, 34, 91, 120, 125], [0, 18, 30, 100, 106, 110], [0, 32, 37, 70, 109, 121]]
\item \{1=32256, 2=625968, 3=2974608, 4=3927168\} [[0, 1, 2, 5, 8, 14], [0, 3, 4, 38, 48, 116], [0, 6, 37, 78, 88, 107], [0, 7, 34, 49, 61, 63], [0, 10, 30, 56, 81, 117], [0, 35, 57, 64, 115, 119], [0, 39, 54, 73, 79, 95]]
\item \{1=32256, 2=643356, 3=2998296, 4=3886092\} [[0, 1, 2, 5, 8, 14], [0, 3, 7, 81, 102, 114], [0, 4, 36, 50, 103, 107], [0, 6, 42, 85, 88, 117], [0, 19, 46, 57, 101, 119], [0, 21, 31, 52, 77, 99], [0, 39, 43, 54, 74, 94]]
\item \{1=32256, 2=645624, 3=2987712, 4=3894408\} [[0, 1, 2, 5, 8, 14], [0, 3, 7, 19, 77, 124], [0, 6, 49, 52, 78, 123], [0, 10, 20, 58, 72, 96], [0, 12, 64, 83, 97, 120], [0, 28, 39, 53, 73, 115], [0, 31, 34, 74, 95, 122]]
\item \{1=32256, 2=647892, 3=2992248, 4=3887604\} [[0, 1, 2, 5, 8, 14], [0, 3, 4, 16, 19, 43], [0, 6, 46, 87, 93, 125], [0, 9, 53, 64, 98, 101], [0, 10, 56, 63, 104, 123], [0, 13, 23, 58, 90, 118], [0, 15, 33, 91, 100, 120]]
\item \{1=32256, 2=666036, 3=3009384, 4=3852324\} [[0, 1, 2, 5, 8, 14], [0, 3, 7, 58, 87, 116], [0, 4, 23, 91, 114, 117], [0, 6, 30, 57, 82, 121], [0, 12, 13, 15, 79, 92]]
\item \{1=33264, 2=514836, 3=3027528, 4=3984372\} [[0, 1, 2, 5, 8, 14], [0, 3, 21, 82, 90, 113], [0, 4, 41, 46, 61, 81], [0, 7, 31, 93, 102, 117], [0, 9, 25, 91, 92, 120], [0, 15, 40, 71, 112, 124], [0, 24, 35, 47, 86, 108]]
\item \{1=33264, 2=515592, 3=2966544, 4=4044600\} [[0, 1, 2, 5, 8, 14], [0, 3, 25, 30, 39, 115], [0, 4, 38, 61, 67, 109], [0, 6, 10, 27, 51, 124], [0, 16, 57, 71, 91, 116]]
\item \{1=33264, 2=536760, 3=3000816, 4=3989160\} [[0, 1, 2, 5, 8, 14], [0, 3, 16, 72, 96, 112], [0, 4, 6, 54, 67, 78], [0, 7, 21, 47, 77, 103], [0, 9, 36, 86, 104, 113], [0, 12, 13, 23, 85, 117], [0, 20, 39, 44, 73, 109]]
\item \{1=33264, 2=559440, 3=2987712, 4=3979584\} [[0, 1, 2, 5, 8, 14], [0, 3, 4, 37, 79, 103], [0, 6, 27, 54, 82, 125], [0, 7, 31, 34, 93, 112], [0, 16, 22, 72, 74, 121]]
\item \{1=33264, 2=560952, 3=2973600, 4=3992184\} [[0, 1, 2, 5, 8, 14], [0, 3, 4, 37, 79, 103], [0, 6, 16, 51, 59, 94], [0, 7, 9, 56, 58, 66], [0, 18, 24, 46, 78, 119]]
\item \{1=33264, 2=561708, 3=2989224, 4=3975804\} [[0, 1, 2, 5, 8, 14], [0, 3, 4, 31, 56, 103], [0, 6, 27, 33, 93, 124], [0, 7, 9, 81, 84, 91], [0, 13, 22, 58, 71, 78]]
\item \{1=33264, 2=564732, 3=2958984, 4=4003020\} [[0, 1, 2, 5, 8, 14], [0, 3, 7, 31, 39, 77], [0, 4, 40, 43, 84, 97], [0, 10, 22, 34, 70, 119], [0, 15, 19, 24, 55, 122]]
\item \{1=33264, 2=565488, 3=3024000, 4=3937248\} [[0, 1, 2, 5, 8, 14], [0, 3, 4, 29, 43, 94], [0, 6, 50, 60, 104, 113], [0, 7, 30, 35, 56, 64], [0, 12, 48, 65, 76, 117]]
\item \{1=33264, 2=568512, 3=2969568, 4=3988656\} [[0, 1, 2, 5, 8, 14], [0, 3, 4, 37, 117, 124], [0, 6, 13, 41, 43, 122], [0, 7, 23, 47, 52, 94], [0, 9, 49, 63, 69, 100]]
\item \{1=33264, 2=570780, 3=2962008, 4=3993948\} [[0, 1, 3, 6, 11, 17], [0, 2, 7, 10, 46, 49], [0, 5, 33, 50, 53, 85], [0, 13, 73, 87, 101, 118], [0, 15, 21, 40, 42, 123]]
\item \{1=33264, 2=581364, 3=3012408, 4=3932964\} [[0, 1, 2, 5, 8, 14], [0, 3, 4, 18, 48, 115], [0, 9, 34, 54, 96, 108], [0, 10, 68, 95, 109, 113], [0, 16, 21, 60, 75, 112]]
\item \{1=33264, 2=583632, 3=3001824, 4=3941280\} [[0, 1, 2, 5, 8, 14], [0, 3, 4, 55, 96, 98], [0, 6, 60, 63, 84, 112], [0, 7, 29, 34, 62, 66], [0, 10, 15, 81, 88, 106]]
\item \{1=33264, 2=583632, 3=3048192, 4=3894912\} [[0, 1, 3, 6, 11, 17], [0, 2, 4, 16, 41, 110], [0, 5, 26, 36, 66, 111], [0, 13, 52, 75, 90, 99], [0, 15, 29, 69, 93, 107], [0, 20, 22, 32, 34, 81], [0, 31, 35, 54, 94, 112]]
\item \{1=33264, 2=584388, 3=2943864, 4=3998484\} [[0, 1, 2, 5, 8, 14], [0, 3, 10, 87, 91, 104], [0, 6, 35, 41, 49, 89], [0, 9, 47, 54, 59, 109], [0, 13, 99, 101, 103, 119]]
\item \{1=33264, 2=585144, 3=2988720, 4=3952872\} [[0, 1, 2, 5, 8, 14], [0, 3, 7, 29, 45, 98], [0, 4, 71, 78, 92, 110], [0, 6, 13, 50, 64, 120], [0, 9, 10, 52, 96, 99]]
\item \{1=33264, 2=587412, 3=2973096, 4=3966228\} [[0, 1, 2, 5, 8, 14], [0, 3, 7, 37, 91, 111], [0, 4, 68, 79, 93, 116], [0, 9, 13, 56, 81, 122], [0, 15, 40, 53, 64, 92]]
\item \{1=33264, 2=588924, 3=2967048, 4=3970764\} [[0, 1, 2, 5, 8, 14], [0, 3, 16, 63, 81, 117], [0, 6, 29, 71, 106, 124], [0, 7, 23, 44, 56, 60], [0, 25, 47, 49, 79, 95]]
\item \{1=33264, 2=589680, 3=2987712, 4=3949344\} [[0, 1, 2, 5, 8, 14], [0, 3, 4, 57, 109, 116], [0, 6, 26, 61, 107, 123], [0, 7, 58, 70, 73, 77], [0, 12, 38, 45, 51, 67]]
\item \{1=33264, 2=591948, 3=2965032, 4=3969756\} [[0, 1, 2, 5, 8, 14], [0, 3, 10, 36, 67, 93], [0, 4, 26, 35, 64, 99], [0, 6, 27, 85, 97, 108], [0, 15, 32, 51, 82, 124]]
\item \{1=33264, 2=592704, 3=2985696, 4=3948336\} [[0, 1, 2, 5, 8, 14], [0, 3, 16, 63, 116, 121], [0, 4, 47, 61, 71, 118], [0, 6, 41, 66, 88, 89], [0, 20, 22, 32, 34, 81], [0, 21, 31, 52, 77, 99], [0, 23, 26, 74, 95, 119]]
\item \{1=33264, 2=593460, 3=2971080, 4=3962196\} [[0, 1, 2, 5, 8, 14], [0, 3, 21, 25, 42, 114], [0, 6, 53, 56, 66, 110], [0, 7, 16, 72, 75, 108], [0, 13, 23, 33, 104, 119]]
\item \{1=33264, 2=595728, 3=2940336, 4=3990672\} [[0, 1, 2, 5, 8, 14], [0, 3, 4, 19, 48, 92], [0, 6, 42, 45, 65, 94], [0, 9, 47, 57, 95, 108], [0, 13, 55, 58, 69, 117]]
\item \{1=33264, 2=595728, 3=2993760, 4=3937248\} [[0, 1, 3, 6, 11, 17], [0, 2, 4, 16, 70, 96], [0, 5, 53, 73, 85, 116], [0, 9, 34, 88, 109, 124], [0, 10, 15, 56, 64, 105], [0, 12, 55, 59, 75, 90], [0, 19, 48, 91, 120, 122]]
\item \{1=33264, 2=598752, 3=3079440, 4=3848544\} [[0, 1, 2, 5, 8, 14], [0, 3, 25, 48, 90, 96], [0, 4, 61, 68, 87, 110], [0, 7, 15, 62, 120, 124], [0, 18, 59, 92, 102, 122], [0, 19, 28, 93, 107, 123], [0, 20, 22, 47, 49, 51]]
\item \{1=33264, 2=600264, 3=3040128, 4=3886344\} [[0, 1, 2, 5, 8, 14], [0, 3, 10, 45, 102, 125], [0, 6, 53, 64, 97, 116], [0, 9, 42, 48, 71, 111], [0, 13, 24, 96, 98, 107]]
\item \{1=33264, 2=607824, 3=3011904, 4=3907008\} [[0, 1, 2, 5, 8, 14], [0, 3, 4, 67, 91, 101], [0, 6, 13, 61, 104, 109], [0, 7, 23, 49, 90, 112], [0, 10, 27, 56, 59, 108], [0, 15, 65, 78, 94, 119], [0, 33, 58, 68, 72, 82]]
\item \{1=33264, 2=610092, 3=2973096, 4=3943548\} [[0, 1, 2, 5, 8, 14], [0, 3, 4, 49, 72, 113], [0, 6, 19, 23, 37, 123], [0, 7, 20, 43, 73, 108], [0, 9, 40, 53, 96, 103]]
\item \{1=33264, 2=612360, 3=3008880, 4=3905496\} [[0, 1, 2, 5, 8, 14], [0, 3, 7, 55, 103, 111], [0, 4, 15, 43, 71, 122], [0, 6, 30, 49, 70, 90], [0, 13, 57, 72, 107, 117]]
\item \{1=33264, 2=613116, 3=2969064, 4=3944556\} [[0, 1, 2, 5, 8, 14], [0, 3, 7, 22, 64, 94], [0, 4, 75, 110, 117, 123], [0, 12, 28, 85, 90, 116], [0, 15, 27, 40, 41, 109]]
\item \{1=33264, 2=615384, 3=2979648, 4=3931704\} [[0, 1, 2, 5, 8, 14], [0, 3, 25, 29, 48, 122], [0, 4, 34, 75, 89, 108], [0, 6, 35, 88, 101, 125], [0, 10, 15, 59, 70, 91]]
\item \{1=33264, 2=615384, 3=3012912, 4=3898440\} [[0, 1, 2, 5, 8, 14], [0, 3, 16, 45, 49, 100], [0, 4, 36, 78, 101, 120], [0, 6, 7, 31, 94, 96], [0, 13, 18, 22, 70, 124]]
\item \{1=33264, 2=615384, 3=3016944, 4=3894408\} [[0, 1, 2, 5, 8, 14], [0, 3, 4, 36, 47, 123], [0, 6, 33, 46, 56, 70], [0, 7, 29, 82, 84, 103], [0, 13, 53, 55, 85, 117]]
\item \{1=33264, 2=617652, 3=2980152, 4=3928932\} [[0, 1, 3, 6, 11, 17], [0, 2, 4, 19, 28, 60], [0, 5, 38, 74, 76, 117], [0, 9, 78, 85, 113, 125], [0, 10, 26, 27, 72, 104]]
\item \{1=33264, 2=617652, 3=3044664, 4=3864420\} [[0, 1, 2, 5, 8, 14], [0, 3, 7, 50, 85, 122], [0, 4, 35, 54, 99, 108], [0, 6, 25, 48, 62, 116], [0, 15, 43, 59, 72, 118]]
\item \{1=33264, 2=619164, 3=2967048, 4=3940524\} [[0, 1, 2, 5, 8, 14], [0, 3, 7, 81, 115, 123], [0, 6, 43, 48, 71, 103], [0, 9, 36, 86, 104, 113], [0, 10, 52, 56, 85, 122], [0, 12, 24, 39, 73, 106], [0, 13, 65, 72, 97, 114]]
\item \{1=33264, 2=620676, 3=2994264, 4=3911796\} [[0, 1, 2, 5, 8, 14], [0, 3, 22, 25, 65, 109], [0, 4, 52, 75, 78, 115], [0, 6, 27, 40, 82, 96], [0, 7, 47, 57, 90, 123]]
\item \{1=33264, 2=621432, 3=3012912, 4=3892392\} [[0, 1, 2, 5, 8, 14], [0, 3, 4, 18, 77, 98], [0, 9, 36, 53, 85, 112], [0, 10, 29, 40, 94, 104], [0, 13, 15, 79, 87, 93], [0, 16, 19, 74, 95, 114], [0, 24, 73, 76, 96, 118]]
\item \{1=33264, 2=626724, 3=3010392, 4=3889620\} [[0, 1, 2, 5, 8, 14], [0, 3, 4, 26, 39, 98], [0, 6, 55, 63, 66, 74], [0, 7, 38, 57, 78, 123], [0, 9, 19, 32, 103, 104]]
\item \{1=33264, 2=628236, 3=2963016, 4=3935484\} [[0, 1, 2, 5, 8, 14], [0, 3, 21, 60, 83, 113], [0, 4, 17, 38, 85, 89], [0, 9, 16, 62, 72, 109], [0, 10, 27, 63, 88, 120]]
\item \{1=33264, 2=628236, 3=2978136, 4=3920364\} [[0, 1, 2, 5, 8, 14], [0, 3, 4, 60, 94, 113], [0, 6, 48, 93, 111, 123], [0, 7, 9, 84, 86, 98], [0, 10, 27, 70, 112, 119]]
\item \{1=33264, 2=634284, 3=3008376, 4=3884076\} [[0, 1, 2, 5, 8, 14], [0, 3, 7, 21, 93, 110], [0, 4, 34, 71, 86, 124], [0, 6, 23, 57, 63, 78], [0, 13, 46, 88, 98, 99]]
\item \{1=33264, 2=639576, 3=2983680, 4=3903480\} [[0, 1, 2, 5, 8, 14], [0, 3, 21, 47, 63, 124], [0, 4, 27, 28, 41, 97], [0, 9, 36, 105, 114, 123], [0, 15, 39, 76, 106, 122], [0, 16, 34, 91, 120, 125], [0, 20, 22, 68, 70, 112]]
\item \{1=33264, 2=651672, 3=2974608, 4=3900456\} [[0, 1, 2, 5, 8, 14], [0, 3, 4, 46, 67, 105], [0, 6, 27, 66, 79, 118], [0, 10, 28, 35, 61, 90], [0, 12, 43, 72, 83, 94]]
\item \{1=33264, 2=672084, 3=3000312, 4=3854340\} [[0, 1, 2, 5, 8, 14], [0, 3, 10, 50, 78, 110], [0, 4, 18, 30, 82, 111], [0, 7, 35, 41, 98, 123], [0, 23, 25, 63, 119, 121]]
\item \{1=33264, 2=673596, 3=2982168, 4=3870972\} [[0, 1, 2, 5, 8, 14], [0, 3, 7, 25, 39, 55], [0, 4, 9, 68, 78, 82], [0, 10, 45, 58, 86, 110], [0, 13, 30, 57, 104, 121]]
\item \{1=33264, 2=689472, 3=3024000, 4=3813264\} [[0, 1, 2, 5, 8, 14], [0, 3, 10, 50, 78, 96], [0, 6, 48, 63, 119, 121], [0, 7, 46, 92, 97, 122], [0, 12, 21, 31, 76, 112]]
\item \{1=34272, 2=557928, 3=2987712, 4=3980088\} [[0, 1, 2, 5, 8, 14], [0, 3, 4, 39, 52, 94], [0, 6, 12, 27, 32, 93], [0, 9, 61, 102, 106, 113], [0, 10, 18, 41, 64, 112]]
\item \{1=34272, 2=563976, 3=2987712, 4=3974040\} [[0, 1, 2, 5, 8, 14], [0, 3, 16, 29, 32, 119], [0, 6, 30, 35, 42, 84], [0, 7, 38, 52, 120, 124], [0, 12, 19, 75, 97, 106]]
\item \{1=34272, 2=564732, 3=3009384, 4=3951612\} [[0, 1, 2, 5, 8, 14], [0, 3, 7, 25, 30, 78], [0, 4, 10, 54, 108, 123], [0, 9, 49, 72, 100, 117], [0, 13, 20, 23, 82, 121]]
\item \{1=34272, 2=566244, 3=2958984, 4=4000500\} [[0, 1, 2, 5, 8, 14], [0, 3, 4, 40, 92, 113], [0, 7, 44, 50, 63, 88], [0, 12, 37, 45, 69, 71], [0, 13, 67, 90, 107, 123]]
\item \{1=34272, 2=572292, 3=2995272, 4=3958164\} [[0, 1, 2, 5, 8, 14], [0, 3, 10, 25, 84, 92], [0, 4, 6, 19, 48, 70], [0, 13, 20, 24, 63, 103], [0, 28, 61, 85, 98, 120]]
\item \{1=34272, 2=579096, 3=3005856, 4=3940776\} [[0, 1, 2, 5, 8, 14], [0, 3, 4, 49, 59, 81], [0, 6, 19, 37, 74, 98], [0, 7, 9, 38, 62, 78], [0, 21, 33, 87, 110, 124]]
\item \{1=34272, 2=585144, 3=2998800, 4=3941784\} [[0, 1, 2, 5, 8, 14], [0, 3, 10, 29, 64, 117], [0, 4, 87, 113, 118, 119], [0, 9, 44, 49, 59, 108], [0, 13, 37, 96, 99, 114]]
\item \{1=34272, 2=585900, 3=3009384, 4=3930444\} [[0, 1, 2, 5, 8, 14], [0, 3, 10, 30, 40, 100], [0, 4, 35, 63, 98, 122], [0, 6, 70, 77, 81, 119], [0, 7, 41, 42, 43, 92]]
\item \{1=34272, 2=588168, 3=2961504, 4=3976056\} [[0, 1, 2, 5, 8, 14], [0, 3, 4, 52, 65, 100], [0, 6, 12, 19, 75, 122], [0, 10, 20, 58, 72, 96], [0, 15, 44, 69, 101, 119], [0, 37, 41, 78, 106, 110], [0, 48, 51, 74, 95, 125]]
\item \{1=34272, 2=588168, 3=3020976, 4=3916584\} [[0, 1, 2, 5, 8, 14], [0, 3, 4, 31, 40, 72], [0, 6, 33, 35, 93, 123], [0, 10, 34, 56, 67, 113], [0, 13, 28, 38, 79, 82], [0, 37, 49, 87, 100, 109], [0, 47, 57, 69, 102, 119]]
\item \{1=34272, 2=589680, 3=3032064, 4=3903984\} [[0, 1, 2, 5, 8, 14], [0, 3, 10, 41, 85, 124], [0, 4, 6, 48, 96, 118], [0, 7, 53, 69, 91, 103], [0, 15, 25, 37, 98, 122]]
\item \{1=34272, 2=590436, 3=2977128, 4=3958164\} [[0, 1, 3, 6, 11, 17], [0, 2, 4, 20, 55, 65], [0, 5, 25, 81, 85, 99], [0, 7, 29, 53, 91, 111], [0, 10, 46, 49, 98, 103]]
\item \{1=34272, 2=592704, 3=2972592, 4=3960432\} [[0, 1, 2, 5, 8, 14], [0, 3, 21, 22, 46, 87], [0, 4, 27, 39, 97, 116], [0, 7, 37, 64, 109, 123], [0, 10, 19, 50, 56, 100], [0, 12, 32, 55, 60, 117], [0, 24, 31, 105, 108, 113]]
\item \{1=34272, 2=594216, 3=2925216, 4=4006296\} [[0, 1, 2, 5, 8, 14], [0, 3, 16, 22, 32, 109], [0, 6, 33, 35, 98, 110], [0, 7, 38, 92, 115, 124], [0, 10, 51, 56, 84, 121], [0, 19, 28, 93, 107, 123], [0, 21, 23, 79, 90, 120]]
\item \{1=34272, 2=597240, 3=3015936, 4=3912552\} [[0, 1, 2, 5, 8, 14], [0, 3, 4, 47, 58, 104], [0, 6, 10, 37, 77, 93], [0, 7, 49, 65, 79, 101], [0, 13, 18, 71, 92, 123]]
\item \{1=34272, 2=597996, 3=3034584, 4=3893148\} [[0, 1, 2, 5, 8, 14], [0, 3, 21, 60, 81, 108], [0, 4, 10, 49, 57, 71], [0, 7, 36, 47, 102, 107], [0, 12, 38, 68, 79, 112]]
\item \{1=34272, 2=598752, 3=3021984, 4=3904992\} [[0, 1, 2, 5, 8, 14], [0, 3, 4, 52, 86, 92], [0, 6, 12, 41, 51, 95], [0, 9, 19, 40, 83, 109], [0, 16, 57, 105, 119, 121], [0, 20, 24, 60, 106, 110], [0, 23, 68, 69, 90, 122]]
\item \{1=34272, 2=600264, 3=2989728, 4=3935736\} [[0, 1, 2, 5, 8, 14], [0, 3, 7, 30, 67, 93], [0, 4, 16, 69, 100, 104], [0, 6, 71, 77, 106, 120], [0, 9, 24, 43, 58, 63]]
\item \{1=34272, 2=603288, 3=3045168, 4=3877272\} [[0, 1, 2, 5, 8, 14], [0, 3, 4, 31, 37, 119], [0, 6, 47, 78, 99, 104], [0, 9, 36, 71, 102, 121], [0, 10, 22, 64, 86, 93], [0, 23, 52, 90, 101, 114], [0, 33, 59, 92, 115, 125]]
\item \{1=34272, 2=606312, 3=3012912, 4=3906504\} [[0, 1, 2, 5, 8, 14], [0, 3, 10, 42, 111, 116], [0, 4, 15, 22, 90, 99], [0, 7, 35, 64, 108, 121], [0, 16, 46, 60, 71, 74]]
\item \{1=34272, 2=609336, 3=2990736, 4=3925656\} [[0, 1, 2, 5, 8, 14], [0, 3, 4, 40, 60, 116], [0, 6, 21, 45, 71, 98], [0, 7, 37, 64, 109, 123], [0, 10, 34, 95, 111, 114], [0, 29, 67, 87, 91, 120], [0, 33, 59, 92, 115, 125]]
\item \{1=34272, 2=610092, 3=2976120, 4=3939516\} [[0, 1, 2, 5, 8, 14], [0, 3, 25, 54, 75, 83], [0, 4, 27, 46, 73, 98], [0, 7, 43, 45, 92, 105], [0, 16, 69, 100, 101, 106]]
\item \{1=34272, 2=610092, 3=3000312, 4=3915324\} [[0, 1, 2, 5, 8, 14], [0, 3, 16, 51, 72, 108], [0, 4, 6, 26, 27, 116], [0, 7, 21, 39, 90, 94], [0, 18, 59, 92, 102, 122], [0, 20, 22, 50, 52, 98], [0, 35, 57, 64, 115, 119]]
\item \{1=34272, 2=610848, 3=2955456, 4=3959424\} [[0, 1, 2, 5, 8, 14], [0, 3, 10, 50, 54, 69], [0, 4, 18, 21, 61, 81], [0, 7, 20, 38, 92, 111], [0, 9, 83, 97, 114, 115]]
\item \{1=34272, 2=611604, 3=3004344, 4=3909780\} [[0, 1, 2, 5, 8, 14], [0, 3, 4, 72, 88, 101], [0, 6, 32, 62, 77, 98], [0, 7, 38, 92, 105, 109], [0, 9, 43, 87, 93, 95]]
\item \{1=34272, 2=615384, 3=2966544, 4=3943800\} [[0, 1, 2, 5, 8, 14], [0, 3, 4, 39, 86, 118], [0, 6, 13, 36, 81, 92], [0, 7, 21, 31, 44, 104], [0, 10, 49, 61, 111, 114]]
\item \{1=34272, 2=616140, 3=3001320, 4=3908268\} [[0, 1, 2, 5, 8, 14], [0, 3, 16, 54, 106, 108], [0, 4, 18, 78, 98, 102], [0, 7, 10, 62, 72, 84], [0, 12, 19, 55, 69, 95]]
\item \{1=34272, 2=617652, 3=2975112, 4=3932964\} [[0, 1, 2, 5, 8, 14], [0, 3, 7, 16, 85, 102], [0, 6, 22, 42, 91, 96], [0, 9, 19, 25, 47, 94], [0, 10, 35, 46, 70, 106]]
\item \{1=34272, 2=622188, 3=2983176, 4=3920364\} [[0, 1, 2, 5, 8, 14], [0, 3, 10, 21, 67, 124], [0, 4, 35, 70, 98, 105], [0, 6, 7, 65, 72, 120], [0, 9, 29, 41, 42, 100]]
\item \{1=34272, 2=624456, 3=2972592, 4=3928680\} [[0, 1, 2, 5, 8, 14], [0, 3, 7, 10, 34, 77], [0, 4, 43, 69, 94, 100], [0, 9, 24, 52, 81, 112], [0, 18, 61, 75, 76, 92]]
\item \{1=34272, 2=628236, 3=2978136, 4=3919356\} [[0, 1, 2, 5, 8, 14], [0, 3, 7, 23, 69, 103], [0, 4, 17, 20, 77, 96], [0, 12, 28, 35, 82, 101], [0, 13, 47, 58, 95, 98]]
\item \{1=34272, 2=631260, 3=3004344, 4=3890124\} [[0, 1, 2, 5, 8, 14], [0, 3, 7, 16, 37, 59], [0, 4, 23, 60, 62, 81], [0, 6, 50, 64, 101, 124], [0, 10, 22, 49, 104, 115]]
\item \{1=34272, 2=637308, 3=2992248, 4=3896172\} [[0, 1, 2, 5, 8, 14], [0, 3, 4, 33, 46, 75], [0, 7, 50, 82, 87, 111], [0, 9, 35, 74, 88, 101], [0, 13, 30, 38, 78, 108]]
\item \{1=34272, 2=638064, 3=2983680, 4=3903984\} [[0, 1, 2, 5, 8, 14], [0, 3, 4, 29, 41, 96], [0, 6, 39, 93, 115, 124], [0, 12, 24, 69, 90, 105], [0, 13, 87, 99, 116, 118]]
\item \{1=34272, 2=641088, 3=2984688, 4=3899952\} [[0, 1, 2, 5, 8, 14], [0, 3, 10, 24, 78, 116], [0, 6, 22, 38, 104, 109], [0, 7, 44, 52, 64, 76], [0, 12, 54, 55, 98, 123]]
\item \{1=34272, 2=648648, 3=2984688, 4=3892392\} [[0, 1, 3, 6, 11, 17], [0, 2, 4, 26, 68, 73], [0, 5, 48, 61, 69, 122], [0, 7, 43, 47, 71, 120], [0, 13, 37, 52, 102, 124]]
\item \{1=34272, 2=660744, 3=3010896, 4=3854088\} [[0, 1, 2, 5, 8, 14], [0, 3, 7, 22, 23, 110], [0, 4, 67, 70, 82, 107], [0, 6, 42, 94, 96, 114], [0, 15, 50, 76, 78, 91]]
\item \{1=34272, 2=669060, 3=3052728, 4=3803940\} [[0, 1, 3, 6, 11, 17], [0, 2, 4, 20, 63, 90], [0, 5, 56, 61, 108, 115], [0, 7, 35, 79, 86, 120], [0, 12, 59, 64, 65, 118]]
\item \{1=34272, 2=675108, 3=2941848, 4=3908772\} [[0, 1, 2, 5, 8, 14], [0, 3, 7, 23, 72, 107], [0, 4, 16, 26, 39, 88], [0, 13, 47, 78, 95, 117], [0, 15, 77, 81, 97, 125]]
\item \{1=34272, 2=680400, 3=2998800, 4=3846528\} [[0, 1, 2, 5, 8, 14], [0, 3, 21, 64, 96, 113], [0, 4, 20, 56, 76, 104], [0, 6, 86, 92, 97, 120], [0, 12, 43, 46, 66, 88]]
\item \{1=35280, 2=568512, 3=3007872, 4=3948336\} [[0, 1, 2, 5, 8, 14], [0, 3, 4, 19, 43, 84], [0, 6, 33, 54, 73, 76], [0, 10, 21, 56, 62, 95], [0, 15, 35, 99, 101, 124], [0, 20, 22, 50, 52, 98], [0, 28, 32, 69, 106, 110]]
\item \{1=35280, 2=573804, 3=3026520, 4=3924396\} [[0, 1, 2, 5, 8, 14], [0, 3, 4, 10, 79, 102], [0, 6, 33, 43, 50, 125], [0, 7, 40, 71, 72, 100], [0, 9, 65, 66, 91, 96]]
\item \{1=35280, 2=579852, 3=3046680, 4=3898188\} [[0, 1, 2, 5, 8, 14], [0, 3, 10, 58, 74, 125], [0, 4, 17, 46, 79, 89], [0, 9, 62, 78, 84, 98], [0, 16, 68, 70, 82, 119]]
\item \{1=35280, 2=580608, 3=3005856, 4=3938256\} [[0, 1, 2, 5, 8, 14], [0, 3, 4, 76, 117, 125], [0, 6, 26, 29, 52, 98], [0, 9, 38, 41, 61, 119], [0, 12, 43, 50, 109, 121], [0, 16, 21, 35, 77, 83], [0, 20, 22, 68, 70, 112]]
\item \{1=35280, 2=582120, 3=3022992, 4=3919608\} [[0, 1, 2, 5, 8, 14], [0, 3, 7, 48, 55, 99], [0, 4, 43, 54, 114, 124], [0, 9, 27, 61, 112, 117], [0, 10, 42, 68, 108, 115]]
\item \{1=35280, 2=588924, 3=2992248, 4=3943548\} [[0, 1, 2, 5, 8, 14], [0, 3, 25, 29, 65, 87], [0, 4, 10, 54, 76, 114], [0, 6, 22, 71, 75, 90], [0, 7, 21, 64, 105, 125]]
\item \{1=35280, 2=588924, 3=3000312, 4=3935484\} [[0, 1, 2, 5, 8, 14], [0, 3, 21, 65, 87, 108], [0, 4, 36, 47, 78, 102], [0, 6, 10, 64, 83, 111], [0, 7, 29, 51, 116, 117]]
\item \{1=35280, 2=589680, 3=2966544, 4=3968496\} [[0, 1, 2, 5, 8, 14], [0, 3, 4, 43, 48, 59], [0, 6, 19, 67, 77, 101], [0, 7, 50, 95, 119, 121], [0, 9, 40, 60, 86, 122]]
\item \{1=35280, 2=590436, 3=2996280, 4=3938004\} [[0, 1, 2, 5, 8, 14], [0, 3, 29, 33, 81, 99], [0, 4, 68, 87, 89, 110], [0, 6, 10, 22, 73, 93], [0, 13, 28, 38, 40, 116]]
\item \{1=35280, 2=591948, 3=2995272, 4=3937500\} [[0, 1, 2, 5, 8, 14], [0, 3, 4, 36, 52, 93], [0, 6, 12, 47, 51, 103], [0, 10, 33, 71, 77, 98], [0, 13, 45, 66, 73, 78]]
\item \{1=35280, 2=594216, 3=3000816, 4=3929688\} [[0, 1, 2, 5, 8, 14], [0, 3, 21, 38, 107, 123], [0, 6, 32, 33, 39, 104], [0, 7, 48, 65, 76, 113], [0, 10, 50, 87, 88, 95]]
\item \{1=35280, 2=597240, 3=3000816, 4=3926664\} [[0, 1, 2, 5, 8, 14], [0, 3, 7, 48, 96, 110], [0, 4, 62, 78, 103, 118], [0, 6, 25, 41, 46, 119], [0, 15, 24, 88, 99, 108], [0, 16, 70, 81, 93, 107], [0, 35, 39, 73, 101, 122]]
\item \{1=35280, 2=597996, 3=2964024, 4=3962700\} [[0, 1, 2, 5, 8, 14], [0, 3, 16, 64, 99, 125], [0, 4, 20, 63, 87, 113], [0, 6, 76, 86, 97, 119], [0, 12, 35, 54, 79, 110]]
\item \{1=35280, 2=598752, 3=2969568, 4=3956400\} [[0, 1, 2, 5, 8, 14], [0, 3, 10, 42, 60, 82], [0, 6, 24, 30, 64, 99], [0, 7, 46, 71, 92, 97], [0, 12, 21, 22, 70, 103]]
\item \{1=35280, 2=601776, 3=3008880, 4=3914064\} [[0, 1, 2, 5, 8, 14], [0, 3, 7, 59, 63, 96], [0, 4, 17, 29, 35, 116], [0, 12, 24, 45, 69, 101], [0, 13, 72, 97, 123, 124]]
\item \{1=35280, 2=602532, 3=3041640, 4=3880548\} [[0, 1, 2, 5, 8, 14], [0, 3, 4, 21, 63, 65], [0, 7, 53, 68, 75, 91], [0, 9, 13, 79, 87, 124], [0, 10, 34, 69, 73, 93]]
\item \{1=35280, 2=613872, 3=2965536, 4=3945312\} [[0, 1, 2, 5, 8, 14], [0, 3, 21, 33, 47, 100], [0, 4, 31, 73, 106, 108], [0, 6, 22, 78, 79, 82], [0, 13, 46, 75, 87, 90]]
\item \{1=35280, 2=613872, 3=2979648, 4=3931200\} [[0, 1, 2, 5, 8, 14], [0, 3, 4, 31, 84, 89], [0, 6, 30, 87, 93, 114], [0, 10, 20, 88, 101, 118], [0, 12, 43, 73, 83, 119]]
\item \{1=35280, 2=616140, 3=3016440, 4=3892140\} [[0, 1, 3, 6, 11, 17], [0, 2, 4, 31, 91, 94], [0, 5, 52, 55, 57, 89], [0, 9, 59, 69, 71, 95], [0, 23, 30, 84, 106, 107]]
\item \{1=35280, 2=616896, 3=2993760, 4=3914064\} [[0, 1, 2, 5, 8, 14], [0, 3, 7, 41, 66, 104], [0, 4, 36, 59, 87, 94], [0, 6, 26, 61, 77, 114], [0, 10, 28, 57, 79, 82]]
\item \{1=35280, 2=616896, 3=3009888, 4=3897936\} [[0, 1, 2, 5, 8, 14], [0, 3, 10, 69, 75, 109], [0, 4, 15, 50, 76, 77], [0, 6, 12, 24, 90, 113], [0, 9, 44, 64, 112, 123]]
\item \{1=35280, 2=617652, 3=2981160, 4=3925908\} [[0, 1, 2, 5, 8, 14], [0, 3, 21, 36, 48, 92], [0, 4, 10, 19, 123, 124], [0, 6, 63, 64, 78, 108], [0, 13, 18, 50, 61, 106]]
\item \{1=35280, 2=618408, 3=3004848, 4=3901464\} [[0, 1, 2, 5, 8, 14], [0, 3, 4, 29, 56, 112], [0, 6, 16, 34, 69, 113], [0, 7, 22, 52, 75, 81], [0, 18, 35, 68, 106, 122]]
\item \{1=35280, 2=622188, 3=3015432, 4=3887100\} [[0, 1, 2, 5, 8, 14], [0, 3, 4, 77, 85, 106], [0, 6, 29, 50, 70, 90], [0, 7, 10, 58, 93, 104], [0, 13, 52, 71, 84, 92]]
\item \{1=35280, 2=622944, 3=2985696, 4=3916080\} [[0, 1, 2, 5, 8, 14], [0, 3, 4, 58, 105, 120], [0, 6, 60, 62, 75, 99], [0, 9, 36, 71, 102, 121], [0, 10, 35, 56, 68, 114], [0, 12, 27, 28, 48, 76], [0, 15, 39, 47, 73, 117]]
\item \{1=35280, 2=633528, 3=2979648, 4=3911544\} [[0, 1, 2, 5, 8, 14], [0, 3, 10, 30, 55, 106], [0, 4, 46, 65, 87, 102], [0, 6, 34, 61, 89, 116], [0, 25, 50, 72, 82, 103]]
\item \{1=35280, 2=633528, 3=2991744, 4=3899448\} [[0, 1, 2, 5, 8, 14], [0, 3, 4, 10, 81, 93], [0, 6, 49, 62, 78, 110], [0, 7, 45, 52, 71, 77], [0, 9, 36, 88, 100, 115], [0, 18, 59, 92, 102, 122], [0, 20, 22, 29, 31, 33]]
\item \{1=35280, 2=635040, 3=2981664, 4=3908016\} [[0, 1, 2, 5, 8, 14], [0, 3, 16, 87, 101, 125], [0, 4, 26, 27, 47, 90], [0, 6, 53, 56, 79, 120], [0, 7, 28, 57, 84, 119], [0, 15, 24, 88, 99, 108], [0, 20, 22, 32, 34, 81]]
\item \{1=35280, 2=638820, 3=3009384, 4=3876516\} [[0, 1, 2, 5, 8, 14], [0, 3, 25, 47, 65, 120], [0, 6, 22, 32, 70, 115], [0, 7, 10, 68, 75, 109], [0, 9, 36, 77, 96, 107], [0, 23, 29, 88, 90, 123], [0, 28, 49, 52, 74, 95]]
\item \{1=35280, 2=639576, 3=2986704, 4=3898440\} [[0, 1, 2, 5, 8, 14], [0, 3, 4, 21, 58, 63], [0, 7, 30, 68, 79, 84], [0, 9, 34, 91, 98, 124], [0, 10, 18, 52, 76, 103]]
\item \{1=35280, 2=641088, 3=2966544, 4=3917088\} [[0, 1, 3, 6, 11, 17], [0, 2, 7, 26, 60, 108], [0, 4, 15, 62, 91, 115], [0, 5, 35, 40, 78, 106], [0, 16, 21, 52, 72, 114]]
\item \{1=35280, 2=650916, 3=2943864, 4=3929940\} [[0, 1, 2, 5, 8, 14], [0, 3, 7, 16, 35, 72], [0, 4, 26, 36, 77, 122], [0, 6, 60, 101, 104, 114], [0, 15, 27, 78, 109, 125]]
\item \{1=35280, 2=654696, 3=3028032, 4=3841992\} [[0, 1, 2, 5, 8, 14], [0, 3, 10, 78, 84, 120], [0, 4, 46, 59, 66, 94], [0, 6, 48, 57, 101, 117], [0, 13, 55, 65, 102, 121]]
\item \{1=35280, 2=669060, 3=2942856, 4=3912804\} [[0, 1, 2, 5, 8, 14], [0, 3, 36, 48, 66, 121], [0, 6, 42, 77, 100, 110], [0, 7, 15, 51, 84, 87], [0, 16, 59, 67, 69, 117]]
\item \{1=35280, 2=672084, 3=2967048, 4=3885588\} [[0, 1, 2, 5, 8, 14], [0, 3, 4, 21, 88, 116], [0, 6, 23, 47, 104, 122], [0, 7, 44, 52, 98, 120], [0, 10, 13, 43, 99, 118]]
\item \{1=36288, 2=570024, 3=3015936, 4=3937752\} [[0, 1, 2, 5, 8, 14], [0, 3, 4, 72, 117, 123], [0, 6, 27, 67, 86, 108], [0, 7, 31, 46, 73, 87], [0, 9, 51, 59, 94, 96]]
\item \{1=36288, 2=577584, 3=3002832, 4=3943296\} [[0, 1, 2, 5, 8, 14], [0, 3, 7, 65, 104, 123], [0, 4, 70, 71, 95, 121], [0, 6, 37, 84, 90, 118], [0, 9, 44, 56, 81, 111]]
\item \{1=36288, 2=579852, 3=2981160, 4=3962700\} [[0, 1, 2, 5, 8, 14], [0, 3, 16, 75, 83, 114], [0, 4, 37, 77, 90, 123], [0, 6, 35, 70, 79, 121], [0, 10, 25, 76, 81, 115]]
\item \{1=36288, 2=580608, 3=3006864, 4=3936240\} [[0, 1, 2, 5, 8, 14], [0, 3, 10, 41, 51, 96], [0, 4, 21, 47, 76, 115], [0, 6, 24, 33, 60, 65], [0, 25, 29, 44, 49, 120]]
\item \{1=36288, 2=582120, 3=3047184, 4=3894408\} [[0, 1, 2, 5, 8, 14], [0, 3, 4, 48, 69, 71], [0, 6, 12, 38, 57, 107], [0, 10, 46, 59, 61, 122], [0, 16, 18, 29, 72, 84]]
\item \{1=36288, 2=584388, 3=3029544, 4=3909780\} [[0, 1, 2, 5, 8, 14], [0, 3, 21, 38, 83, 104], [0, 6, 7, 16, 69, 76], [0, 9, 36, 105, 114, 123], [0, 10, 30, 56, 81, 117], [0, 12, 15, 96, 102, 106], [0, 27, 55, 57, 109, 119]]
\item \{1=36288, 2=588924, 3=2999304, 4=3935484\} [[0, 1, 2, 5, 8, 14], [0, 3, 21, 78, 84, 100], [0, 4, 19, 87, 112, 120], [0, 6, 13, 33, 36, 109], [0, 10, 34, 72, 74, 118]]
\item \{1=36288, 2=588924, 3=3043656, 4=3891132\} [[0, 1, 2, 5, 8, 14], [0, 3, 4, 36, 93, 125], [0, 6, 10, 84, 87, 123], [0, 12, 28, 55, 62, 100], [0, 18, 35, 38, 94, 101]]
\item \{1=36288, 2=590436, 3=2979144, 4=3954132\} [[0, 1, 2, 5, 8, 14], [0, 3, 4, 50, 63, 105], [0, 6, 19, 70, 100, 110], [0, 9, 21, 55, 87, 91], [0, 10, 20, 34, 60, 82]]
\item \{1=36288, 2=591192, 3=3033072, 4=3899448\} [[0, 1, 2, 5, 8, 14], [0, 3, 7, 51, 86, 119], [0, 4, 47, 49, 62, 82], [0, 6, 39, 53, 79, 102], [0, 9, 36, 70, 83, 101], [0, 10, 26, 91, 120, 124], [0, 20, 22, 59, 61, 106]]
\item \{1=36288, 2=593460, 3=3048696, 4=3881556\} [[0, 1, 2, 5, 8, 14], [0, 3, 4, 35, 60, 64], [0, 7, 38, 41, 65, 68], [0, 9, 12, 58, 63, 66], [0, 13, 40, 73, 115, 118], [0, 18, 37, 81, 103, 109], [0, 20, 22, 86, 88, 121]]
\item \{1=36288, 2=594972, 3=3012408, 4=3916332\} [[0, 1, 2, 5, 8, 14], [0, 3, 4, 66, 115, 118], [0, 6, 22, 38, 59, 83], [0, 7, 35, 42, 48, 70], [0, 21, 31, 52, 77, 99], [0, 23, 37, 90, 95, 96], [0, 29, 39, 73, 104, 116]]
\item \{1=36288, 2=597240, 3=3011904, 4=3914568\} [[0, 1, 2, 5, 8, 14], [0, 3, 4, 54, 58, 66], [0, 6, 34, 71, 77, 99], [0, 7, 44, 79, 103, 112], [0, 15, 24, 43, 117, 118]]
\item \{1=36288, 2=598752, 3=2995776, 4=3929184\} [[0, 1, 2, 5, 8, 14], [0, 3, 4, 37, 47, 58], [0, 6, 16, 42, 69, 110], [0, 9, 63, 74, 76, 117], [0, 10, 22, 35, 71, 94]]
\item \{1=36288, 2=607068, 3=2955960, 4=3960684\} [[0, 1, 2, 5, 8, 14], [0, 3, 7, 64, 78, 95], [0, 4, 20, 31, 82, 104], [0, 6, 75, 77, 79, 113], [0, 12, 28, 38, 97, 108]]
\item \{1=36288, 2=607068, 3=2967048, 4=3949596\} [[0, 1, 2, 5, 8, 14], [0, 3, 10, 50, 113, 123], [0, 4, 58, 72, 76, 90], [0, 6, 12, 70, 117, 119], [0, 9, 24, 36, 42, 55], [0, 13, 63, 66, 86, 107], [0, 21, 28, 77, 85, 105]]
\item \{1=36288, 2=610848, 3=2984688, 4=3928176\} [[0, 1, 2, 5, 8, 14], [0, 3, 4, 49, 75, 123], [0, 6, 33, 36, 94, 107], [0, 7, 27, 48, 67, 104], [0, 10, 13, 41, 56, 92], [0, 31, 37, 69, 83, 109], [0, 40, 47, 52, 115, 124]]
\item \{1=36288, 2=611604, 3=2982168, 4=3929940\} [[0, 1, 2, 5, 8, 14], [0, 3, 7, 10, 41, 113], [0, 4, 34, 48, 93, 95], [0, 6, 19, 36, 40, 69], [0, 12, 21, 22, 44, 121]]
\item \{1=36288, 2=612360, 3=2955456, 4=3955896\} [[0, 1, 2, 5, 8, 14], [0, 3, 4, 22, 56, 101], [0, 6, 16, 68, 89, 120], [0, 13, 32, 64, 78, 117], [0, 21, 28, 96, 106, 123]]
\item \{1=36288, 2=613872, 3=2934288, 4=3975552\} [[0, 1, 2, 5, 8, 14], [0, 3, 4, 69, 81, 112], [0, 6, 21, 36, 47, 73], [0, 10, 35, 75, 91, 98], [0, 19, 38, 57, 92, 110]]
\item \{1=36288, 2=625212, 3=2947896, 4=3950604\} [[0, 1, 2, 5, 8, 14], [0, 3, 4, 82, 88, 116], [0, 6, 29, 36, 43, 72], [0, 7, 59, 90, 107, 124], [0, 10, 34, 63, 87, 100]]
\item \{1=36288, 2=625212, 3=3022488, 4=3876012\} [[0, 1, 2, 5, 8, 14], [0, 3, 4, 57, 82, 117], [0, 6, 26, 41, 53, 74], [0, 7, 10, 28, 79, 115], [0, 12, 44, 76, 81, 110]]
\item \{1=36288, 2=626724, 3=3011400, 4=3885588\} [[0, 1, 2, 5, 8, 14], [0, 3, 16, 37, 75, 87], [0, 4, 33, 35, 38, 76], [0, 7, 46, 95, 108, 121], [0, 12, 13, 79, 84, 106]]
\item \{1=36288, 2=626724, 3=3023496, 4=3873492\} [[0, 1, 2, 5, 8, 14], [0, 3, 29, 51, 100, 122], [0, 4, 17, 35, 38, 125], [0, 7, 23, 75, 106, 113], [0, 9, 13, 46, 60, 74]]
\item \{1=36288, 2=629748, 3=2944872, 4=3949092\} [[0, 1, 2, 5, 8, 14], [0, 3, 7, 81, 112, 122], [0, 4, 21, 58, 69, 111], [0, 6, 22, 44, 92, 98], [0, 16, 40, 41, 46, 83]]
\item \{1=36288, 2=632016, 3=2982672, 4=3909024\} [[0, 1, 2, 5, 8, 14], [0, 3, 4, 41, 88, 106], [0, 6, 34, 39, 49, 58], [0, 12, 45, 90, 96, 116], [0, 16, 27, 76, 107, 115]]
\item \{1=36288, 2=632772, 3=3002328, 4=3888612\} [[0, 1, 2, 5, 8, 14], [0, 3, 25, 47, 69, 93], [0, 4, 33, 58, 96, 115], [0, 7, 10, 28, 95, 112], [0, 15, 18, 65, 88, 121]]
\item \{1=36288, 2=635796, 3=2957976, 4=3929940\} [[0, 1, 2, 5, 8, 14], [0, 3, 4, 65, 101, 120], [0, 6, 22, 34, 82, 87], [0, 7, 37, 92, 116, 125], [0, 9, 21, 55, 100, 119]]
\item \{1=36288, 2=643356, 3=3013416, 4=3866940\} [[0, 1, 2, 5, 8, 14], [0, 3, 4, 27, 29, 93], [0, 6, 38, 96, 99, 114], [0, 7, 36, 46, 75, 78], [0, 12, 21, 22, 71, 116]]
\item \{1=36288, 2=646380, 3=2983176, 4=3894156\} [[0, 1, 2, 5, 8, 14], [0, 3, 16, 57, 72, 73], [0, 4, 22, 47, 110, 121], [0, 6, 24, 42, 71, 78], [0, 7, 44, 59, 112, 123]]
\item \{1=36288, 2=646380, 3=2998296, 4=3879036\} [[0, 1, 2, 5, 8, 14], [0, 3, 10, 29, 74, 113], [0, 4, 22, 43, 62, 115], [0, 6, 7, 85, 111, 125], [0, 9, 36, 50, 69, 82], [0, 15, 55, 75, 87, 114], [0, 21, 30, 51, 77, 98]]
\item \{1=36288, 2=648648, 3=3000816, 4=3874248\} [[0, 1, 2, 5, 8, 14], [0, 3, 4, 27, 68, 69], [0, 7, 26, 45, 84, 107], [0, 9, 38, 42, 85, 98], [0, 12, 21, 24, 65, 120]]
\item \{1=36288, 2=649404, 3=2997288, 4=3877020\} [[0, 1, 2, 5, 8, 14], [0, 3, 7, 10, 40, 63], [0, 4, 42, 70, 78, 98], [0, 9, 26, 36, 38, 57], [0, 16, 27, 58, 82, 107], [0, 18, 21, 64, 77, 116], [0, 28, 35, 66, 91, 120]]
\item \{1=36288, 2=650916, 3=3016440, 4=3856356\} [[0, 1, 2, 5, 8, 14], [0, 3, 10, 72, 74, 125], [0, 4, 60, 70, 84, 107], [0, 6, 62, 67, 97, 112], [0, 18, 23, 71, 76, 119]]
\item \{1=36288, 2=652428, 3=2983176, 4=3888108\} [[0, 1, 2, 5, 8, 14], [0, 3, 7, 19, 102, 103], [0, 6, 30, 32, 61, 100], [0, 12, 27, 66, 77, 101], [0, 16, 70, 81, 93, 107], [0, 34, 39, 50, 73, 121], [0, 37, 41, 78, 106, 110]]
\item \{1=36288, 2=659232, 3=2955456, 4=3909024\} [[0, 1, 2, 5, 8, 14], [0, 3, 4, 40, 81, 92], [0, 6, 31, 57, 89, 115], [0, 13, 47, 78, 106, 114], [0, 15, 37, 83, 120, 122]]
\item \{1=36288, 2=677376, 3=2989728, 4=3856608\} [[0, 1, 2, 5, 8, 14], [0, 3, 7, 38, 48, 122], [0, 4, 35, 50, 60, 106], [0, 10, 18, 21, 22, 109], [0, 25, 52, 76, 83, 96]]
\item \{1=37296, 2=553392, 3=3012912, 4=3956400\} [[0, 1, 2, 5, 8, 14], [0, 3, 21, 25, 64, 73], [0, 4, 6, 57, 70, 95], [0, 9, 32, 42, 62, 113], [0, 15, 37, 58, 82, 104]]
\item \{1=37296, 2=560952, 3=3025008, 4=3936744\} [[0, 1, 2, 5, 8, 14], [0, 3, 16, 49, 72, 100], [0, 4, 56, 78, 95, 112], [0, 6, 34, 79, 101, 103], [0, 7, 22, 67, 93, 97]]
\item \{1=37296, 2=563976, 3=2987712, 4=3971016\} [[0, 1, 2, 5, 8, 14], [0, 3, 21, 45, 60, 125], [0, 6, 10, 40, 63, 88], [0, 7, 23, 41, 44, 101], [0, 19, 67, 69, 98, 100]]
\item \{1=37296, 2=569268, 3=2978136, 4=3975300\} [[0, 1, 2, 5, 8, 14], [0, 3, 4, 37, 81, 112], [0, 6, 12, 53, 56, 59], [0, 9, 26, 54, 106, 122], [0, 16, 79, 95, 102, 125]]
\item \{1=37296, 2=571536, 3=3029040, 4=3922128\} [[0, 1, 2, 5, 8, 14], [0, 3, 4, 35, 47, 73], [0, 6, 66, 89, 110, 115], [0, 9, 26, 62, 64, 121], [0, 15, 37, 51, 65, 109], [0, 21, 23, 79, 90, 120], [0, 24, 30, 104, 108, 112]]
\item \{1=37296, 2=573048, 3=2976624, 4=3973032\} [[0, 1, 2, 5, 8, 14], [0, 3, 4, 71, 95, 124], [0, 6, 46, 56, 68, 78], [0, 7, 99, 110, 117, 122], [0, 12, 30, 61, 76, 86]]
\item \{1=37296, 2=577584, 3=3046176, 4=3898944\} [[0, 1, 2, 5, 8, 14], [0, 3, 7, 48, 49, 54], [0, 4, 34, 86, 88, 113], [0, 6, 19, 22, 53, 65], [0, 10, 52, 58, 93, 118]]
\item \{1=37296, 2=580608, 3=2990736, 4=3951360\} [[0, 1, 2, 5, 8, 14], [0, 3, 4, 10, 30, 101], [0, 6, 49, 63, 96, 107], [0, 7, 21, 42, 58, 62], [0, 12, 48, 75, 86, 92]]
\item \{1=37296, 2=585144, 3=3030048, 4=3907512\} [[0, 1, 2, 5, 8, 14], [0, 3, 10, 86, 96, 123], [0, 4, 21, 45, 120, 122], [0, 6, 23, 27, 67, 87], [0, 7, 63, 68, 103, 113]]
\item \{1=37296, 2=588168, 3=3000816, 4=3933720\} [[0, 1, 2, 5, 8, 14], [0, 3, 4, 35, 46, 72], [0, 6, 13, 65, 103, 111], [0, 7, 63, 68, 82, 121], [0, 9, 41, 81, 92, 104]]
\item \{1=37296, 2=590436, 3=3053736, 4=3878532\} [[0, 1, 2, 5, 8, 14], [0, 3, 4, 19, 62, 67], [0, 6, 30, 86, 99, 124], [0, 10, 56, 69, 101, 125], [0, 12, 38, 41, 73, 76], [0, 16, 57, 105, 119, 121], [0, 18, 63, 78, 100, 116]]
\item \{1=37296, 2=591192, 3=3053232, 4=3878280\} [[0, 1, 2, 5, 8, 14], [0, 3, 4, 41, 46, 96], [0, 6, 35, 39, 102, 120], [0, 7, 71, 85, 95, 113], [0, 12, 22, 69, 72, 82]]
\item \{1=37296, 2=594972, 3=2945880, 4=3981852\} [[0, 1, 2, 5, 8, 14], [0, 3, 4, 82, 86, 91], [0, 6, 68, 78, 88, 106], [0, 10, 28, 35, 42, 77], [0, 12, 52, 92, 96, 109]]
\item \{1=37296, 2=595728, 3=2974608, 4=3952368\} [[0, 1, 2, 5, 8, 14], [0, 3, 4, 73, 118, 125], [0, 6, 13, 46, 79, 113], [0, 9, 48, 61, 65, 108], [0, 12, 69, 109, 116, 124]]
\item \{1=37296, 2=598752, 3=2996784, 4=3927168\} [[0, 1, 2, 5, 8, 14], [0, 3, 10, 46, 91, 101], [0, 4, 26, 62, 65, 109], [0, 6, 58, 97, 104, 123], [0, 9, 59, 88, 102, 118]]
\item \{1=37296, 2=598752, 3=3022992, 4=3900960\} [[0, 1, 2, 5, 8, 14], [0, 3, 7, 99, 104, 115], [0, 4, 41, 92, 97, 122], [0, 6, 26, 60, 112, 125], [0, 9, 27, 50, 98, 124]]
\item \{1=37296, 2=598752, 3=3024000, 4=3899952\} [[0, 1, 2, 5, 8, 14], [0, 3, 10, 30, 48, 96], [0, 7, 22, 44, 109, 114], [0, 9, 55, 59, 111, 121], [0, 16, 51, 58, 66, 95]]
\item \{1=37296, 2=610092, 3=2933784, 4=3978828\} [[0, 1, 2, 5, 8, 14], [0, 3, 7, 78, 84, 123], [0, 4, 17, 26, 59, 82], [0, 9, 46, 104, 106, 119], [0, 10, 28, 41, 91, 97]]
\item \{1=37296, 2=610848, 3=2993760, 4=3918096\} [[0, 1, 3, 6, 11, 17], [0, 2, 7, 28, 55, 116], [0, 4, 21, 36, 62, 69], [0, 5, 26, 66, 67, 121], [0, 16, 22, 52, 59, 91]]
\item \{1=37296, 2=613116, 3=2982168, 4=3927420\} [[0, 1, 2, 5, 8, 14], [0, 3, 7, 82, 108, 111], [0, 4, 37, 66, 77, 102], [0, 9, 35, 42, 96, 106], [0, 13, 84, 90, 103, 121]]
\item \{1=37296, 2=613116, 3=2986200, 4=3923388\} [[0, 1, 2, 5, 8, 14], [0, 3, 10, 21, 31, 119], [0, 4, 35, 70, 105, 112], [0, 7, 9, 64, 86, 90], [0, 20, 48, 85, 91, 92]]
\item \{1=37296, 2=615384, 3=2998800, 4=3908520\} [[0, 1, 2, 5, 8, 14], [0, 3, 7, 87, 110, 124], [0, 4, 49, 108, 117, 118], [0, 6, 13, 16, 37, 71], [0, 9, 42, 43, 58, 114]]
\item \{1=37296, 2=616140, 3=2997288, 4=3909276\} [[0, 1, 2, 5, 8, 14], [0, 3, 4, 27, 94, 121], [0, 6, 51, 56, 103, 122], [0, 7, 30, 42, 59, 75], [0, 16, 22, 34, 41, 112]]
\item \{1=37296, 2=618408, 3=2998800, 4=3905496\} [[0, 1, 2, 5, 8, 14], [0, 3, 21, 25, 123, 124], [0, 4, 17, 58, 71, 107], [0, 7, 28, 56, 87, 92], [0, 9, 44, 55, 63, 83]]
\item \{1=37296, 2=619164, 3=3009384, 4=3894156\} [[0, 1, 2, 5, 8, 14], [0, 3, 4, 29, 107, 122], [0, 6, 52, 62, 81, 95], [0, 9, 16, 60, 77, 82], [0, 10, 40, 46, 72, 119]]
\item \{1=37296, 2=619920, 3=2990736, 4=3912048\} [[0, 1, 2, 5, 8, 14], [0, 3, 4, 23, 43, 117], [0, 6, 37, 59, 79, 125], [0, 7, 16, 77, 99, 122], [0, 12, 48, 54, 85, 100]]
\item \{1=37296, 2=619920, 3=3022992, 4=3879792\} [[0, 1, 2, 5, 8, 14], [0, 3, 4, 81, 90, 99], [0, 6, 50, 52, 86, 100], [0, 10, 20, 82, 92, 112], [0, 15, 43, 57, 101, 125]]
\item \{1=37296, 2=622188, 3=2996280, 4=3904236\} [[0, 1, 2, 5, 8, 14], [0, 3, 4, 69, 86, 109], [0, 7, 40, 56, 100, 101], [0, 9, 53, 60, 90, 96], [0, 12, 24, 27, 28, 98]]
\item \{1=37296, 2=626724, 3=2965032, 4=3930948\} [[0, 1, 2, 5, 8, 14], [0, 3, 25, 54, 65, 121], [0, 4, 6, 70, 72, 108], [0, 7, 20, 22, 46, 48], [0, 9, 19, 40, 88, 107], [0, 15, 39, 47, 73, 117], [0, 23, 26, 75, 90, 92]]
\item \{1=37296, 2=632772, 3=3012408, 4=3877524\} [[0, 1, 2, 5, 8, 14], [0, 3, 21, 56, 78, 83], [0, 4, 76, 97, 104, 114], [0, 6, 36, 49, 59, 71], [0, 12, 41, 66, 72, 123]]
\item \{1=37296, 2=633528, 3=2895984, 4=3993192\} [[0, 1, 2, 5, 8, 14], [0, 3, 4, 57, 92, 115], [0, 6, 27, 50, 74, 96], [0, 7, 30, 34, 41, 110], [0, 9, 19, 60, 83, 94]]
\item \{1=37296, 2=637308, 3=2926728, 4=3958668\} [[0, 1, 2, 5, 8, 14], [0, 3, 16, 49, 63, 67], [0, 6, 13, 52, 92, 97], [0, 9, 43, 54, 109, 122], [0, 19, 28, 93, 107, 123], [0, 20, 22, 68, 70, 112], [0, 24, 31, 105, 108, 113]]
\item \{1=37296, 2=638064, 3=3014928, 4=3869712\} [[0, 1, 2, 5, 8, 14], [0, 3, 16, 72, 93, 109], [0, 6, 12, 88, 110, 125], [0, 9, 73, 108, 112, 113], [0, 15, 63, 67, 78, 115]]
\item \{1=37296, 2=638820, 3=3002328, 4=3881556\} [[0, 1, 2, 5, 8, 14], [0, 3, 4, 27, 86, 90], [0, 6, 41, 49, 68, 96], [0, 10, 54, 85, 109, 112], [0, 12, 59, 69, 81, 106]]
\item \{1=37296, 2=642600, 3=2959488, 4=3920616\} [[0, 1, 2, 5, 8, 14], [0, 3, 16, 32, 81, 123], [0, 6, 49, 70, 89, 94], [0, 9, 83, 90, 119, 121], [0, 13, 39, 73, 75, 107], [0, 21, 30, 51, 77, 98], [0, 40, 47, 52, 115, 124]]
\item \{1=37296, 2=642600, 3=3002832, 4=3877272\} [[0, 1, 2, 5, 8, 14], [0, 3, 7, 37, 65, 122], [0, 4, 26, 49, 62, 93], [0, 6, 10, 54, 100, 110], [0, 18, 22, 38, 84, 95]]
\item \{1=37296, 2=644112, 3=2999808, 4=3878784\} [[0, 1, 2, 5, 8, 14], [0, 3, 4, 19, 50, 76], [0, 6, 62, 63, 82, 121], [0, 9, 38, 75, 100, 122], [0, 20, 22, 29, 31, 33], [0, 21, 39, 60, 77, 106], [0, 23, 52, 90, 101, 114]]
\item \{1=37296, 2=647892, 3=3009384, 4=3865428\} [[0, 1, 2, 5, 8, 14], [0, 3, 10, 40, 41, 121], [0, 4, 34, 36, 49, 79], [0, 6, 26, 52, 58, 112], [0, 16, 28, 70, 110, 115]]
\item \{1=37296, 2=650916, 3=2999304, 4=3872484\} [[0, 1, 2, 5, 8, 14], [0, 3, 4, 18, 50, 112], [0, 9, 35, 40, 59, 88], [0, 10, 12, 94, 104, 107], [0, 19, 55, 95, 119, 121]]
\item \{1=37296, 2=653940, 3=2919672, 4=3949092\} [[0, 1, 2, 5, 8, 14], [0, 3, 7, 54, 57, 115], [0, 4, 81, 117, 120, 122], [0, 6, 31, 87, 88, 106], [0, 9, 26, 39, 58, 70]]
\item \{1=37296, 2=681156, 3=2979144, 4=3862404\} [[0, 1, 2, 5, 8, 14], [0, 3, 4, 46, 58, 66], [0, 6, 22, 50, 70, 92], [0, 7, 68, 84, 117, 120], [0, 9, 19, 40, 61, 109]]
\item \{1=37296, 2=687960, 3=2966544, 4=3868200\} [[0, 1, 2, 5, 8, 14], [0, 3, 4, 22, 48, 51], [0, 6, 19, 43, 65, 125], [0, 10, 62, 83, 113, 120], [0, 15, 29, 69, 93, 107], [0, 24, 52, 64, 102, 108], [0, 27, 40, 72, 94, 124]]
\item \{1=38304, 2=555660, 3=2983176, 4=3982860\} [[0, 1, 2, 5, 8, 14], [0, 3, 4, 16, 96, 106], [0, 6, 37, 57, 72, 99], [0, 7, 38, 98, 108, 121], [0, 10, 45, 61, 67, 118]]
\item \{1=38304, 2=579852, 3=3033576, 4=3908268\} [[0, 1, 2, 5, 8, 14], [0, 3, 10, 54, 67, 121], [0, 4, 28, 70, 97, 112], [0, 6, 13, 24, 111, 119], [0, 7, 9, 26, 58, 99]]
\item \{1=38304, 2=584388, 3=2988216, 4=3949092\} [[0, 1, 2, 5, 8, 14], [0, 3, 22, 33, 39, 109], [0, 6, 19, 61, 107, 124], [0, 7, 75, 78, 102, 106], [0, 16, 52, 67, 74, 118]]
\item \{1=38304, 2=593460, 3=2953944, 4=3974292\} [[0, 1, 2, 5, 8, 14], [0, 3, 16, 51, 81, 103], [0, 4, 6, 19, 94, 119], [0, 9, 27, 37, 97, 112], [0, 12, 30, 70, 96, 122]]
\item \{1=38304, 2=593460, 3=2998296, 4=3929940\} [[0, 1, 2, 5, 8, 14], [0, 3, 25, 36, 72, 100], [0, 4, 15, 39, 50, 89], [0, 6, 67, 70, 86, 118], [0, 23, 33, 35, 54, 74]]
\item \{1=38304, 2=598752, 3=2989728, 4=3933216\} [[0, 1, 2, 5, 8, 14], [0, 3, 4, 38, 70, 116], [0, 6, 12, 39, 54, 90], [0, 9, 32, 75, 81, 117], [0, 16, 34, 61, 68, 123]]
\item \{1=38304, 2=598752, 3=3020976, 4=3901968\} [[0, 1, 2, 5, 8, 14], [0, 3, 7, 25, 84, 123], [0, 4, 35, 58, 104, 118], [0, 9, 70, 75, 91, 117], [0, 10, 21, 41, 61, 98]]
\item \{1=38304, 2=599508, 3=3010392, 4=3911796\} [[0, 1, 2, 5, 8, 14], [0, 3, 4, 29, 101, 110], [0, 6, 22, 32, 44, 67], [0, 7, 39, 42, 62, 77], [0, 10, 30, 82, 96, 107]]
\item \{1=38304, 2=607068, 3=3013416, 4=3901212\} [[0, 1, 2, 5, 8, 14], [0, 3, 16, 46, 122, 123], [0, 4, 33, 71, 120, 124], [0, 7, 30, 49, 103, 121], [0, 12, 69, 81, 117, 119]]
\item \{1=38304, 2=607068, 3=3017448, 4=3897180\} [[0, 1, 2, 5, 8, 14], [0, 3, 4, 21, 81, 112], [0, 6, 47, 48, 79, 125], [0, 7, 38, 41, 65, 68], [0, 10, 26, 91, 120, 124], [0, 12, 46, 61, 106, 117], [0, 35, 39, 73, 101, 122]]
\item \{1=38304, 2=607824, 3=2976624, 4=3937248\} [[0, 1, 2, 5, 8, 14], [0, 3, 7, 21, 78, 108], [0, 4, 22, 46, 49, 75], [0, 6, 26, 74, 93, 104], [0, 10, 25, 28, 86, 114]]
\item \{1=38304, 2=610848, 3=2993760, 4=3917088\} [[0, 1, 2, 5, 8, 14], [0, 3, 4, 37, 63, 88], [0, 7, 56, 61, 92, 119], [0, 9, 35, 36, 67, 98], [0, 15, 38, 104, 106, 108], [0, 16, 40, 113, 117, 124], [0, 20, 22, 103, 105, 125]]
\item \{1=38304, 2=610848, 3=3013920, 4=3896928\} [[0, 1, 2, 5, 8, 14], [0, 3, 7, 19, 58, 111], [0, 6, 27, 57, 76, 99], [0, 9, 40, 64, 78, 107], [0, 16, 22, 46, 83, 96]]
\item \{1=38304, 2=611604, 3=2995272, 4=3914820\} [[0, 1, 2, 5, 8, 14], [0, 3, 21, 65, 83, 113], [0, 6, 41, 63, 119, 121], [0, 7, 40, 51, 87, 90], [0, 10, 30, 56, 81, 117], [0, 20, 22, 74, 76, 78], [0, 39, 42, 58, 73, 92]]
\item \{1=38304, 2=611604, 3=3014424, 4=3895668\} [[0, 1, 2, 5, 8, 14], [0, 3, 7, 51, 114, 121], [0, 4, 43, 46, 73, 122], [0, 6, 45, 52, 83, 109], [0, 9, 36, 50, 69, 82], [0, 10, 30, 56, 81, 117], [0, 21, 23, 79, 90, 120]]
\item \{1=38304, 2=613116, 3=3004344, 4=3904236\} [[0, 1, 2, 5, 8, 14], [0, 3, 21, 51, 56, 107], [0, 4, 20, 33, 46, 103], [0, 9, 50, 55, 85, 119], [0, 18, 41, 63, 78, 97]]
\item \{1=38304, 2=613872, 3=2976624, 4=3931200\} [[0, 1, 2, 5, 8, 14], [0, 3, 7, 19, 50, 125], [0, 6, 23, 25, 46, 87], [0, 9, 98, 101, 114, 121], [0, 15, 29, 82, 100, 104]]
\item \{1=38304, 2=613872, 3=3016944, 4=3890880\} [[0, 1, 2, 5, 8, 14], [0, 3, 10, 66, 68, 107], [0, 4, 21, 57, 61, 104], [0, 6, 36, 63, 88, 120], [0, 12, 15, 16, 27, 94]]
\item \{1=38304, 2=616896, 3=2997792, 4=3907008\} [[0, 1, 2, 5, 8, 14], [0, 3, 4, 21, 31, 79], [0, 6, 50, 85, 110, 124], [0, 7, 22, 49, 101, 104], [0, 15, 33, 64, 67, 93]]
\item \{1=38304, 2=616896, 3=3003840, 4=3900960\} [[0, 1, 2, 5, 8, 14], [0, 3, 21, 33, 90, 92], [0, 4, 10, 37, 52, 121], [0, 6, 7, 35, 75, 108], [0, 12, 31, 72, 112, 119]]
\item \{1=38304, 2=621432, 3=3006864, 4=3893400\} [[0, 1, 2, 5, 8, 14], [0, 3, 10, 101, 103, 119], [0, 4, 33, 42, 78, 82], [0, 9, 44, 74, 86, 90], [0, 15, 55, 75, 87, 114], [0, 18, 21, 64, 77, 116], [0, 20, 22, 47, 49, 51]]
\item \{1=38304, 2=622188, 3=2977128, 4=3922380\} [[0, 1, 2, 5, 8, 14], [0, 3, 7, 19, 73, 118], [0, 6, 40, 47, 87, 117], [0, 9, 12, 70, 75, 92], [0, 15, 34, 67, 91, 108]]
\item \{1=38304, 2=626724, 3=2973096, 4=3921876\} [[0, 1, 2, 5, 8, 14], [0, 3, 10, 21, 57, 123], [0, 4, 60, 98, 107, 114], [0, 6, 41, 44, 45, 89], [0, 13, 75, 83, 90, 122]]
\item \{1=38304, 2=635040, 3=2963520, 4=3923136\} [[0, 1, 2, 5, 8, 14], [0, 3, 4, 10, 59, 115], [0, 6, 12, 88, 114, 123], [0, 15, 43, 50, 77, 112], [0, 18, 30, 86, 96, 98]]
\item \{1=38304, 2=638064, 3=3015936, 4=3867696\} [[0, 1, 2, 5, 8, 14], [0, 3, 7, 40, 51, 113], [0, 4, 27, 48, 55, 125], [0, 6, 29, 43, 64, 87], [0, 10, 21, 42, 50, 73]]
\item \{1=38304, 2=639576, 3=3010896, 4=3871224\} [[0, 1, 2, 5, 8, 14], [0, 3, 7, 21, 25, 76], [0, 4, 35, 56, 104, 116], [0, 13, 48, 90, 96, 102], [0, 18, 44, 50, 109, 124]]
\item \{1=38304, 2=641088, 3=2996784, 4=3883824\} [[0, 1, 2, 5, 8, 14], [0, 3, 7, 59, 95, 124], [0, 4, 36, 43, 121, 122], [0, 6, 31, 44, 106, 115], [0, 12, 24, 35, 38, 104]]
\item \{1=38304, 2=647892, 3=2974104, 4=3899700\} [[0, 1, 2, 5, 8, 14], [0, 3, 10, 36, 48, 83], [0, 4, 38, 43, 79, 125], [0, 7, 34, 44, 93, 108], [0, 16, 57, 63, 71, 97]]
\item \{1=38304, 2=659988, 3=2946888, 4=3914820\} [[0, 1, 3, 6, 11, 17], [0, 2, 4, 28, 114, 123], [0, 5, 42, 76, 86, 107], [0, 7, 16, 43, 68, 88], [0, 9, 26, 56, 84, 120]]
\item \{1=38304, 2=664524, 3=2998296, 4=3858876\} [[0, 1, 2, 5, 8, 14], [0, 3, 4, 16, 33, 68], [0, 7, 35, 36, 88, 109], [0, 10, 87, 102, 111, 112], [0, 20, 22, 77, 79, 118], [0, 24, 43, 55, 94, 108], [0, 34, 39, 50, 73, 121]]
\item \{1=38304, 2=681156, 3=2969064, 4=3871476\} [[0, 1, 2, 5, 8, 14], [0, 3, 29, 87, 98, 101], [0, 4, 10, 33, 68, 94], [0, 7, 37, 51, 90, 104], [0, 12, 28, 39, 55, 120]]
\item \{1=38304, 2=681156, 3=3040632, 4=3799908\} [[0, 1, 2, 5, 8, 14], [0, 3, 4, 66, 72, 95], [0, 6, 12, 22, 71, 90], [0, 9, 16, 85, 98, 103], [0, 10, 63, 77, 84, 100]]
\item \{1=39312, 2=561708, 3=3039624, 4=3919356\} [[0, 1, 2, 5, 8, 14], [0, 3, 10, 55, 101, 119], [0, 4, 33, 35, 50, 91], [0, 7, 34, 60, 100, 125], [0, 13, 40, 79, 98, 103]]
\item \{1=39312, 2=567000, 3=2998800, 4=3954888\} [[0, 1, 2, 5, 8, 14], [0, 3, 4, 59, 99, 115], [0, 6, 37, 52, 56, 114], [0, 7, 23, 60, 66, 95], [0, 12, 21, 62, 101, 120]]
\item \{1=39312, 2=575316, 3=2926728, 4=4018644\} [[0, 1, 2, 5, 8, 14], [0, 3, 16, 33, 96, 118], [0, 4, 31, 37, 68, 87], [0, 6, 7, 34, 115, 124], [0, 9, 43, 71, 73, 119]]
\item \{1=39312, 2=588168, 3=3012912, 4=3919608\} [[0, 1, 2, 5, 8, 14], [0, 3, 7, 59, 72, 109], [0, 4, 35, 41, 76, 88], [0, 6, 34, 45, 60, 77], [0, 9, 24, 78, 98, 125]]
\item \{1=39312, 2=588924, 3=3001320, 4=3930444\} [[0, 1, 2, 5, 8, 14], [0, 3, 7, 48, 50, 104], [0, 4, 46, 59, 103, 121], [0, 6, 29, 65, 66, 67], [0, 12, 27, 37, 70, 112]]
\item \{1=39312, 2=592704, 3=2999808, 4=3928176\} [[0, 1, 2, 5, 8, 14], [0, 3, 7, 46, 60, 113], [0, 4, 42, 82, 90, 104], [0, 9, 21, 50, 65, 111], [0, 10, 95, 96, 99, 118]]
\item \{1=39312, 2=594972, 3=2970072, 4=3955644\} [[0, 1, 2, 5, 8, 14], [0, 3, 25, 47, 120, 124], [0, 4, 28, 33, 69, 71], [0, 7, 9, 21, 64, 76], [0, 15, 54, 78, 115, 122]]
\item \{1=39312, 2=599508, 3=2971080, 4=3950100\} [[0, 1, 2, 5, 8, 14], [0, 3, 4, 40, 47, 82], [0, 6, 19, 78, 93, 102], [0, 7, 87, 99, 116, 117], [0, 10, 35, 59, 64, 106]]
\item \{1=39312, 2=609336, 3=3000816, 4=3910536\} [[0, 1, 2, 5, 8, 14], [0, 3, 21, 22, 33, 96], [0, 4, 37, 72, 88, 89], [0, 6, 31, 60, 114, 116], [0, 15, 43, 64, 67, 106]]
\item \{1=39312, 2=613116, 3=3006360, 4=3901212\} [[0, 1, 2, 5, 8, 14], [0, 3, 7, 58, 87, 100], [0, 4, 34, 35, 102, 113], [0, 6, 51, 65, 94, 114], [0, 9, 62, 78, 109, 125]]
\item \{1=39312, 2=613116, 3=3033576, 4=3873996\} [[0, 1, 2, 5, 8, 14], [0, 3, 4, 58, 84, 113], [0, 7, 29, 44, 93, 105], [0, 9, 50, 60, 63, 124], [0, 10, 22, 49, 96, 109]]
\item \{1=39312, 2=617652, 3=3017448, 4=3885588\} [[0, 1, 2, 5, 8, 14], [0, 3, 21, 38, 66, 83], [0, 6, 23, 59, 75, 86], [0, 10, 15, 56, 64, 105], [0, 13, 46, 79, 99, 118], [0, 20, 22, 95, 97, 124], [0, 31, 51, 91, 114, 120]]
\item \{1=39312, 2=618408, 3=2981664, 4=3920616\} [[0, 1, 2, 5, 8, 14], [0, 3, 4, 26, 105, 114], [0, 6, 24, 37, 104, 110], [0, 9, 36, 59, 78, 91], [0, 10, 19, 28, 82, 115], [0, 15, 18, 74, 95, 113], [0, 16, 54, 69, 94, 99]]
\item \{1=39312, 2=619920, 3=3038112, 4=3862656\} [[0, 1, 2, 5, 8, 14], [0, 3, 21, 47, 51, 107], [0, 4, 36, 72, 119, 123], [0, 6, 39, 49, 70, 110], [0, 16, 40, 41, 63, 79]]
\item \{1=39312, 2=630504, 3=3030048, 4=3860136\} [[0, 1, 2, 5, 8, 14], [0, 3, 10, 31, 95, 121], [0, 4, 68, 70, 79, 108], [0, 12, 19, 66, 88, 90], [0, 18, 22, 93, 94, 118]]
\item \{1=39312, 2=632016, 3=2991744, 4=3896928\} [[0, 1, 2, 5, 8, 14], [0, 3, 4, 26, 101, 106], [0, 6, 63, 69, 108, 116], [0, 7, 42, 71, 81, 123], [0, 16, 21, 61, 100, 104]]
\item \{1=39312, 2=635040, 3=2917152, 4=3968496\} [[0, 1, 2, 5, 8, 14], [0, 3, 21, 22, 47, 112], [0, 4, 35, 51, 76, 82], [0, 6, 16, 83, 85, 86], [0, 13, 52, 64, 96, 120]]
\item \{1=39312, 2=635796, 3=2967048, 4=3917844\} [[0, 1, 2, 5, 8, 14], [0, 3, 7, 30, 50, 51], [0, 4, 40, 92, 93, 121], [0, 12, 22, 44, 46, 67], [0, 13, 32, 52, 71, 90]]
\item \{1=39312, 2=646380, 3=3028536, 4=3845772\} [[0, 1, 2, 5, 8, 14], [0, 3, 7, 10, 68, 93], [0, 4, 9, 48, 58, 92], [0, 12, 32, 52, 104, 119], [0, 13, 64, 65, 102, 121]]
\item \{1=39312, 2=653184, 3=2962512, 4=3904992\} [[0, 1, 2, 5, 8, 14], [0, 3, 4, 34, 47, 121], [0, 6, 10, 70, 84, 86], [0, 7, 16, 36, 114, 116], [0, 12, 62, 66, 77, 90]]
\item \{1=39312, 2=659232, 3=2980656, 4=3880800\} [[0, 1, 2, 5, 8, 14], [0, 3, 4, 90, 102, 107], [0, 6, 72, 96, 108, 123], [0, 7, 35, 91, 111, 116], [0, 16, 41, 52, 85, 105]]
\item \{1=39312, 2=668304, 3=2955456, 4=3896928\} [[0, 1, 2, 5, 8, 14], [0, 3, 7, 10, 95, 96], [0, 4, 17, 49, 76, 78], [0, 9, 44, 50, 65, 112], [0, 12, 24, 30, 98, 107]]
\item \{1=39312, 2=676620, 3=2995272, 4=3848796\} [[0, 1, 2, 5, 8, 14], [0, 3, 7, 16, 34, 109], [0, 4, 30, 40, 78, 102], [0, 6, 58, 70, 100, 125], [0, 10, 62, 75, 92, 112]]
\item \{1=40320, 2=585900, 3=3003336, 4=3930444\} [[0, 1, 3, 6, 11, 17], [0, 2, 4, 16, 43, 99], [0, 5, 34, 37, 98, 109], [0, 9, 19, 36, 49, 81], [0, 10, 42, 106, 110, 119], [0, 12, 15, 35, 50, 64], [0, 13, 41, 52, 93, 105]]
\item \{1=40320, 2=585900, 3=3011400, 4=3922380\} [[0, 1, 2, 5, 8, 14], [0, 3, 10, 25, 31, 102], [0, 4, 72, 77, 82, 107], [0, 7, 43, 49, 68, 94], [0, 13, 88, 101, 103, 121]]
\item \{1=40320, 2=593460, 3=3041640, 4=3884580\} [[0, 1, 2, 5, 8, 14], [0, 3, 7, 48, 68, 103], [0, 4, 70, 81, 106, 111], [0, 6, 55, 94, 117, 124], [0, 10, 46, 56, 88, 116], [0, 15, 35, 63, 76, 118], [0, 27, 39, 73, 93, 119]]
\item \{1=40320, 2=594216, 3=2989728, 4=3935736\} [[0, 1, 2, 5, 8, 14], [0, 3, 16, 36, 73, 103], [0, 4, 17, 64, 95, 119], [0, 7, 15, 81, 86, 100], [0, 12, 27, 32, 104, 105]]
\item \{1=40320, 2=602532, 3=2980152, 4=3936996\} [[0, 1, 2, 5, 8, 14], [0, 3, 4, 16, 64, 72], [0, 7, 40, 78, 93, 110], [0, 10, 19, 33, 49, 123], [0, 13, 32, 87, 95, 98]]
\item \{1=40320, 2=607824, 3=2973600, 4=3938256\} [[0, 1, 2, 5, 8, 14], [0, 3, 16, 66, 75, 103], [0, 4, 17, 73, 96, 115], [0, 9, 27, 50, 74, 121], [0, 10, 22, 34, 79, 92]]
\item \{1=40320, 2=610848, 3=2994768, 4=3914064\} [[0, 1, 2, 5, 8, 14], [0, 3, 4, 36, 70, 95], [0, 6, 13, 33, 58, 102], [0, 7, 43, 54, 56, 85], [0, 16, 61, 68, 83, 125]]
\item \{1=40320, 2=616896, 3=3011904, 4=3890880\} [[0, 1, 2, 5, 8, 14], [0, 3, 7, 55, 106, 107], [0, 4, 29, 65, 87, 99], [0, 6, 39, 79, 108, 116], [0, 10, 60, 84, 117, 122]]
\item \{1=40320, 2=618408, 3=2970576, 4=3930696\} [[0, 1, 2, 5, 8, 14], [0, 3, 10, 25, 67, 85], [0, 6, 36, 43, 94, 114], [0, 7, 37, 69, 107, 125], [0, 13, 83, 96, 108, 120]]
\item \{1=40320, 2=622944, 3=2975616, 4=3921120\} [[0, 1, 2, 5, 8, 14], [0, 3, 10, 21, 57, 117], [0, 4, 27, 46, 55, 121], [0, 6, 26, 49, 53, 83], [0, 9, 50, 64, 66, 79]]
\item \{1=40320, 2=624456, 3=2972592, 4=3922632\} [[0, 1, 2, 5, 8, 14], [0, 3, 16, 30, 115, 118], [0, 4, 51, 70, 77, 124], [0, 6, 26, 57, 91, 92], [0, 7, 44, 73, 108, 112]]
\item \{1=40320, 2=624456, 3=3011904, 4=3883320\} [[0, 1, 2, 5, 8, 14], [0, 3, 10, 22, 57, 83], [0, 4, 34, 64, 118, 125], [0, 6, 27, 48, 96, 109], [0, 16, 33, 41, 72, 92]]
\item \{1=40320, 2=626724, 3=2972088, 4=3920868\} [[0, 1, 2, 5, 8, 14], [0, 3, 21, 33, 69, 124], [0, 4, 22, 91, 107, 120], [0, 6, 16, 35, 72, 89], [0, 9, 10, 75, 83, 88], [0, 19, 23, 67, 90, 121], [0, 31, 34, 74, 95, 122]]
\item \{1=40320, 2=628236, 3=2951928, 4=3939516\} [[0, 1, 2, 5, 8, 14], [0, 3, 16, 30, 84, 121], [0, 4, 59, 98, 110, 124], [0, 6, 34, 60, 107, 112], [0, 7, 50, 70, 97, 99]]
\item \{1=40320, 2=629748, 3=2982168, 4=3907764\} [[0, 1, 2, 5, 8, 14], [0, 3, 7, 35, 49, 107], [0, 4, 21, 22, 82, 121], [0, 6, 25, 55, 79, 123], [0, 12, 31, 41, 44, 54]]
\item \{1=40320, 2=632016, 3=3027024, 4=3860640\} [[0, 1, 2, 5, 8, 14], [0, 3, 7, 37, 67, 82], [0, 4, 34, 75, 97, 121], [0, 6, 30, 66, 69, 76], [0, 12, 24, 95, 96, 113]]
\item \{1=40320, 2=634284, 3=3005352, 4=3880044\} [[0, 1, 2, 5, 8, 14], [0, 3, 16, 22, 42, 108], [0, 4, 27, 61, 66, 105], [0, 6, 29, 32, 81, 95], [0, 7, 10, 46, 90, 102]]
\item \{1=40320, 2=634284, 3=3010392, 4=3875004\} [[0, 1, 2, 5, 8, 14], [0, 3, 16, 32, 63, 123], [0, 6, 42, 84, 85, 88], [0, 7, 40, 45, 67, 124], [0, 12, 22, 61, 71, 118], [0, 18, 55, 68, 75, 98], [0, 24, 35, 47, 86, 108]]
\item \{1=40320, 2=644868, 3=3018456, 4=3856356\} [[0, 1, 2, 5, 8, 14], [0, 3, 7, 16, 68, 104], [0, 4, 35, 77, 106, 119], [0, 6, 25, 31, 60, 92], [0, 12, 15, 93, 112, 114]]
\item \{1=40320, 2=656208, 3=2971584, 4=3891888\} [[0, 1, 2, 5, 8, 14], [0, 3, 10, 46, 51, 91], [0, 4, 38, 43, 87, 118], [0, 6, 12, 19, 102, 108], [0, 13, 15, 55, 119, 123]]
\item \{1=40320, 2=664524, 3=2974104, 4=3881052\} [[0, 1, 2, 5, 8, 14], [0, 3, 4, 66, 88, 90], [0, 6, 31, 53, 77, 92], [0, 7, 34, 48, 105, 113], [0, 13, 20, 32, 63, 96]]
\item \{1=41328, 2=572292, 3=2996280, 4=3950100\} [[0, 1, 2, 5, 8, 14], [0, 3, 4, 16, 34, 73], [0, 6, 69, 81, 94, 114], [0, 7, 22, 42, 44, 72], [0, 12, 66, 100, 101, 110]]
\item \{1=41328, 2=583632, 3=3043152, 4=3891888\} [[0, 1, 2, 5, 8, 14], [0, 3, 4, 23, 48, 102], [0, 7, 42, 71, 95, 113], [0, 9, 26, 39, 55, 103], [0, 13, 30, 72, 104, 107]]
\item \{1=41328, 2=584388, 3=2995272, 4=3939012\} [[0, 1, 2, 5, 8, 14], [0, 3, 4, 63, 66, 103], [0, 6, 13, 33, 78, 86], [0, 7, 40, 45, 67, 124], [0, 9, 26, 36, 38, 57], [0, 10, 21, 27, 74, 77], [0, 19, 39, 73, 84, 113], [0, 20, 22, 50, 52, 98], [0, 32, 76, 82, 104, 118]]
\item \{1=41328, 2=585900, 3=2952936, 4=3979836\} [[0, 1, 2, 5, 8, 14], [0, 3, 4, 39, 67, 76], [0, 6, 29, 34, 53, 110], [0, 10, 30, 33, 58, 125], [0, 12, 16, 52, 68, 121]]
\item \{1=41328, 2=587412, 3=3010392, 4=3920868\} [[0, 1, 2, 5, 8, 14], [0, 3, 10, 22, 63, 107], [0, 4, 19, 97, 103, 117], [0, 13, 47, 61, 67, 110], [0, 27, 66, 76, 91, 98]]
\item \{1=41328, 2=597240, 3=3010896, 4=3910536\} [[0, 1, 2, 5, 8, 14], [0, 3, 7, 67, 91, 107], [0, 6, 63, 73, 104, 122], [0, 10, 40, 76, 77, 87], [0, 18, 32, 54, 110, 125]]
\item \{1=41328, 2=599508, 3=2968056, 4=3951108\} [[0, 1, 2, 5, 8, 14], [0, 3, 4, 48, 77, 96], [0, 6, 44, 50, 88, 92], [0, 7, 73, 81, 95, 106], [0, 9, 46, 53, 90, 93]]
\item \{1=41328, 2=604800, 3=3004848, 4=3909024\} [[0, 1, 2, 5, 8, 14], [0, 3, 4, 37, 45, 121], [0, 6, 38, 43, 71, 125], [0, 7, 23, 29, 84, 122], [0, 10, 59, 83, 94, 107]]
\item \{1=41328, 2=611604, 3=2953944, 4=3953124\} [[0, 1, 2, 5, 8, 14], [0, 3, 10, 63, 114, 118], [0, 4, 50, 51, 87, 89], [0, 6, 24, 45, 94, 112], [0, 9, 38, 90, 95, 104]]
\item \{1=41328, 2=621432, 3=3032064, 4=3865176\} [[0, 1, 2, 5, 8, 14], [0, 3, 7, 36, 37, 119], [0, 4, 61, 104, 113, 120], [0, 6, 27, 38, 83, 85], [0, 10, 42, 68, 73, 116]]
\item \{1=41328, 2=622188, 3=2965032, 4=3931452\} [[0, 1, 2, 5, 8, 14], [0, 3, 16, 93, 111, 120], [0, 4, 15, 29, 51, 99], [0, 6, 32, 35, 54, 58], [0, 7, 52, 71, 100, 109]]
\item \{1=41328, 2=628992, 3=3012912, 4=3876768\} [[0, 1, 2, 5, 8, 14], [0, 3, 7, 38, 59, 96], [0, 4, 9, 48, 49, 62], [0, 6, 21, 35, 75, 122], [0, 12, 31, 47, 69, 86]]
\item \{1=41328, 2=630504, 3=2968560, 4=3919608\} [[0, 1, 2, 5, 8, 14], [0, 3, 4, 19, 47, 82], [0, 6, 32, 39, 62, 94], [0, 9, 27, 41, 112, 123], [0, 10, 42, 64, 107, 120]]
\item \{1=41328, 2=632772, 3=2976120, 4=3909780\} [[0, 1, 2, 5, 8, 14], [0, 3, 4, 26, 58, 112], [0, 6, 57, 87, 107, 115], [0, 7, 50, 69, 113, 124], [0, 19, 22, 32, 88, 103]]
\item \{1=41328, 2=635040, 3=2968560, 4=3915072\} [[0, 1, 2, 5, 8, 14], [0, 3, 10, 54, 78, 100], [0, 4, 30, 49, 101, 105], [0, 6, 65, 74, 81, 84], [0, 13, 52, 53, 103, 119]]
\item \{1=41328, 2=635040, 3=2994768, 4=3888864\} [[0, 1, 2, 5, 8, 14], [0, 3, 7, 38, 92, 124], [0, 4, 28, 56, 104, 109], [0, 6, 26, 74, 88, 98], [0, 12, 46, 65, 66, 102]]
\item \{1=41328, 2=638064, 3=3001824, 4=3878784\} [[0, 1, 2, 5, 8, 14], [0, 3, 4, 34, 81, 92], [0, 6, 57, 70, 89, 91], [0, 9, 55, 88, 94, 124], [0, 10, 26, 66, 73, 86]]
\item \{1=41328, 2=638820, 3=3039624, 4=3840228\} [[0, 1, 2, 5, 8, 14], [0, 3, 10, 83, 100, 109], [0, 6, 44, 81, 96, 118], [0, 7, 69, 92, 99, 119], [0, 9, 13, 41, 66, 107]]
\item \{1=41328, 2=639576, 3=2989728, 4=3889368\} [[0, 1, 2, 5, 8, 14], [0, 3, 4, 34, 77, 99], [0, 6, 16, 75, 85, 108], [0, 7, 53, 64, 69, 102], [0, 13, 32, 94, 96, 122]]
\item \{1=41328, 2=641088, 3=2988720, 4=3888864\} [[0, 1, 2, 5, 8, 14], [0, 3, 4, 49, 86, 109], [0, 7, 34, 44, 95, 124], [0, 9, 38, 63, 69, 71], [0, 13, 23, 30, 78, 108]]
\item \{1=41328, 2=644112, 3=3009888, 4=3864672\} [[0, 1, 2, 5, 8, 14], [0, 3, 4, 36, 65, 66], [0, 6, 26, 38, 62, 120], [0, 7, 72, 82, 92, 123], [0, 12, 53, 73, 79, 103]]
\item \{1=41328, 2=648648, 3=2999808, 4=3870216\} [[0, 1, 2, 5, 8, 14], [0, 3, 7, 46, 78, 119], [0, 4, 41, 52, 103, 109], [0, 6, 36, 39, 40, 93], [0, 12, 45, 67, 92, 102]]
\item \{1=41328, 2=650916, 3=2978136, 4=3889620\} [[0, 1, 2, 5, 8, 14], [0, 3, 21, 55, 63, 119], [0, 4, 31, 81, 103, 122], [0, 7, 28, 29, 40, 56], [0, 13, 37, 87, 99, 118]]
\item \{1=41328, 2=655452, 3=2973096, 4=3890124\} [[0, 1, 2, 5, 8, 14], [0, 3, 4, 59, 60, 114], [0, 6, 41, 82, 88, 123], [0, 7, 23, 66, 71, 111], [0, 9, 21, 52, 57, 83]]
\item \{1=41328, 2=656964, 3=2981160, 4=3880548\} [[0, 1, 2, 5, 8, 14], [0, 3, 7, 16, 23, 58], [0, 4, 15, 86, 92, 93], [0, 9, 19, 87, 88, 124], [0, 12, 22, 44, 73, 101]]
\item \{1=41328, 2=657720, 3=2996784, 4=3864168\} [[0, 1, 2, 5, 8, 14], [0, 3, 10, 57, 92, 115], [0, 4, 18, 31, 89, 107], [0, 9, 12, 44, 104, 114], [0, 13, 16, 28, 75, 78]]
\item \{1=41328, 2=691740, 3=2968056, 4=3858876\} [[0, 1, 2, 5, 8, 14], [0, 3, 4, 59, 75, 113], [0, 6, 29, 42, 70, 103], [0, 7, 16, 67, 87, 106], [0, 10, 35, 62, 98, 105]]
\item \{1=42336, 2=588168, 3=2995776, 4=3933720\} [[0, 1, 2, 5, 8, 14], [0, 3, 16, 75, 112, 118], [0, 4, 9, 33, 35, 101], [0, 10, 61, 73, 85, 117], [0, 12, 13, 79, 100, 110]]
\item \{1=42336, 2=588924, 3=3003336, 4=3925404\} [[0, 1, 2, 5, 8, 14], [0, 3, 10, 25, 68, 119], [0, 4, 36, 61, 99, 122], [0, 6, 54, 85, 114, 120], [0, 7, 20, 27, 29, 79]]
\item \{1=42336, 2=590436, 3=2949912, 4=3977316\} [[0, 1, 2, 5, 8, 14], [0, 3, 21, 104, 117, 118], [0, 6, 27, 31, 74, 89], [0, 9, 38, 87, 91, 109], [0, 13, 22, 83, 97, 115]]
\item \{1=42336, 2=595728, 3=3025008, 4=3896928\} [[0, 1, 2, 5, 8, 14], [0, 3, 4, 39, 43, 63], [0, 6, 33, 62, 113, 119], [0, 9, 41, 85, 87, 107], [0, 10, 50, 64, 74, 79]]
\item \{1=42336, 2=596484, 3=2985192, 4=3935988\} [[0, 1, 2, 5, 8, 14], [0, 3, 16, 75, 81, 99], [0, 4, 41, 68, 73, 96], [0, 6, 33, 34, 61, 116], [0, 15, 18, 59, 63, 121]]
\item \{1=42336, 2=597240, 3=2969568, 4=3950856\} [[0, 1, 3, 6, 11, 17], [0, 2, 4, 32, 67, 112], [0, 7, 16, 46, 62, 69], [0, 9, 57, 58, 60, 82], [0, 10, 28, 66, 76, 108]]
\item \{1=42336, 2=597240, 3=3009888, 4=3910536\} [[0, 1, 2, 5, 8, 14], [0, 3, 4, 40, 96, 116], [0, 6, 48, 58, 99, 115], [0, 7, 31, 38, 59, 118], [0, 12, 19, 66, 77, 101]]
\item \{1=42336, 2=598752, 3=3022992, 4=3895920\} [[0, 1, 2, 5, 8, 14], [0, 3, 16, 45, 57, 111], [0, 4, 58, 60, 85, 110], [0, 6, 29, 31, 41, 99], [0, 9, 53, 66, 77, 101]]
\item \{1=42336, 2=602532, 3=2990232, 4=3924900\} [[0, 1, 2, 5, 8, 14], [0, 3, 7, 52, 63, 123], [0, 4, 42, 62, 76, 122], [0, 6, 56, 59, 81, 91], [0, 15, 37, 50, 70, 94]]
\item \{1=42336, 2=603288, 3=3015936, 4=3898440\} [[0, 1, 2, 5, 8, 14], [0, 3, 21, 42, 81, 107], [0, 4, 30, 34, 50, 119], [0, 6, 49, 52, 63, 90], [0, 7, 56, 62, 85, 97]]
\item \{1=42336, 2=603288, 3=3044160, 4=3870216\} [[0, 1, 2, 5, 8, 14], [0, 3, 10, 29, 82, 107], [0, 4, 34, 113, 114, 119], [0, 9, 12, 44, 54, 112], [0, 13, 37, 79, 84, 87]]
\item \{1=42336, 2=611604, 3=3018456, 4=3887604\} [[0, 1, 2, 5, 8, 14], [0, 3, 4, 10, 103, 107], [0, 6, 13, 66, 102, 120], [0, 7, 65, 70, 75, 81], [0, 9, 26, 44, 64, 104]]
\item \{1=42336, 2=623700, 3=3008376, 4=3885588\} [[0, 1, 2, 5, 8, 14], [0, 3, 16, 46, 103, 125], [0, 4, 6, 43, 89, 92], [0, 10, 54, 60, 91, 94], [0, 15, 44, 64, 74, 118], [0, 18, 55, 68, 75, 98], [0, 20, 22, 32, 34, 81]]
\item \{1=42336, 2=640332, 3=3011400, 4=3865932\} [[0, 1, 2, 5, 8, 14], [0, 3, 4, 77, 101, 119], [0, 6, 19, 53, 91, 100], [0, 7, 9, 34, 54, 87], [0, 15, 81, 88, 93, 125]]
\item \{1=42336, 2=642600, 3=2980656, 4=3894408\} [[0, 1, 2, 5, 8, 14], [0, 3, 4, 37, 50, 96], [0, 6, 33, 39, 114, 123], [0, 7, 66, 87, 95, 101], [0, 9, 29, 53, 76, 115]]
\item \{1=42336, 2=644112, 3=2993760, 4=3879792\} [[0, 1, 2, 5, 8, 14], [0, 3, 16, 36, 39, 91], [0, 4, 27, 34, 66, 79], [0, 7, 31, 51, 59, 124], [0, 10, 29, 50, 111, 118]]
\item \{1=42336, 2=645624, 3=2963520, 4=3908520\} [[0, 1, 2, 5, 8, 14], [0, 3, 16, 32, 63, 117], [0, 6, 19, 52, 53, 94], [0, 7, 30, 40, 90, 119], [0, 9, 21, 26, 73, 77], [0, 10, 34, 56, 67, 113], [0, 20, 22, 68, 70, 112]]
\item \{1=42336, 2=650916, 3=2996280, 4=3870468\} [[0, 1, 2, 5, 8, 14], [0, 3, 4, 16, 22, 60], [0, 6, 27, 68, 88, 99], [0, 7, 54, 81, 94, 116], [0, 12, 21, 71, 85, 124], [0, 24, 31, 105, 108, 113], [0, 35, 57, 64, 115, 119]]
\item \{1=42336, 2=661500, 3=3038616, 4=3817548\} [[0, 1, 2, 5, 8, 14], [0, 3, 4, 65, 73, 123], [0, 7, 49, 68, 78, 110], [0, 12, 38, 60, 75, 83], [0, 13, 30, 33, 47, 87]]
\item \{1=43344, 2=578340, 3=2998296, 4=3940020\} [[0, 1, 2, 5, 8, 14], [0, 3, 4, 57, 88, 110], [0, 6, 44, 67, 86, 92], [0, 7, 35, 84, 91, 117], [0, 9, 16, 55, 90, 104]]
\item \{1=43344, 2=579852, 3=2977128, 4=3959676\} [[0, 1, 2, 5, 8, 14], [0, 3, 7, 23, 52, 100], [0, 4, 28, 39, 108, 123], [0, 6, 13, 45, 70, 83], [0, 12, 37, 92, 102, 124]]
\item \{1=43344, 2=591192, 3=2979648, 4=3945816\} [[0, 1, 2, 5, 8, 14], [0, 3, 4, 29, 48, 90], [0, 6, 43, 58, 70, 93], [0, 9, 41, 66, 77, 125], [0, 10, 42, 56, 75, 118], [0, 15, 33, 91, 100, 120], [0, 19, 46, 57, 101, 119]]
\item \{1=43344, 2=591192, 3=3004848, 4=3920616\} [[0, 1, 2, 5, 8, 14], [0, 3, 4, 55, 65, 91], [0, 6, 12, 53, 66, 102], [0, 9, 29, 79, 83, 115], [0, 10, 19, 37, 58, 123]]
\item \{1=43344, 2=601776, 3=2987712, 4=3927168\} [[0, 1, 2, 5, 8, 14], [0, 3, 10, 25, 99, 117], [0, 4, 70, 77, 79, 90], [0, 6, 37, 62, 84, 92], [0, 9, 35, 76, 82, 88]]
\item \{1=43344, 2=604044, 3=3013416, 4=3899196\} [[0, 1, 3, 6, 11, 17], [0, 2, 4, 19, 49, 70], [0, 5, 37, 42, 84, 91], [0, 12, 32, 38, 95, 97], [0, 16, 20, 57, 110, 117]]
\item \{1=43344, 2=625212, 3=3049704, 4=3841740\} [[0, 1, 3, 6, 11, 17], [0, 2, 4, 19, 23, 38], [0, 5, 36, 49, 75, 94], [0, 13, 60, 84, 91, 103], [0, 16, 34, 57, 109, 113]]
\item \{1=43344, 2=637308, 3=2991240, 4=3888108\} [[0, 1, 2, 5, 8, 14], [0, 3, 7, 64, 69, 94], [0, 4, 15, 65, 76, 91], [0, 6, 10, 71, 116, 125], [0, 12, 43, 44, 73, 121]]
\item \{1=43344, 2=638064, 3=3008880, 4=3869712\} [[0, 1, 2, 5, 8, 14], [0, 3, 10, 37, 48, 58], [0, 4, 19, 50, 63, 125], [0, 12, 61, 86, 92, 104], [0, 15, 25, 79, 83, 87]]
\item \{1=43344, 2=643356, 3=3008376, 4=3864924\} [[0, 1, 3, 6, 11, 17], [0, 2, 4, 26, 71, 81], [0, 5, 37, 50, 79, 92], [0, 7, 39, 54, 61, 97], [0, 10, 15, 56, 64, 105], [0, 16, 18, 20, 22, 63], [0, 32, 76, 82, 104, 118]]
\item \{1=43344, 2=647136, 3=2982672, 4=3886848\} [[0, 1, 2, 5, 8, 14], [0, 3, 4, 73, 76, 116], [0, 6, 42, 64, 68, 69], [0, 9, 12, 32, 61, 98], [0, 13, 55, 65, 113, 117]]
\item \{1=43344, 2=647892, 3=2985192, 4=3883572\} [[0, 1, 2, 5, 8, 14], [0, 3, 7, 29, 34, 122], [0, 4, 20, 67, 112, 115], [0, 6, 25, 53, 68, 103], [0, 23, 47, 60, 76, 106]]
\item \{1=43344, 2=653184, 3=2980656, 4=3882816\} [[0, 1, 2, 5, 8, 14], [0, 3, 7, 16, 72, 84], [0, 4, 35, 38, 91, 118], [0, 6, 25, 46, 52, 68], [0, 13, 32, 102, 104, 113]]
\item \{1=43344, 2=654696, 3=3028032, 4=3833928\} [[0, 1, 2, 5, 8, 14], [0, 3, 16, 81, 117, 123], [0, 4, 54, 85, 86, 87], [0, 6, 19, 39, 53, 63], [0, 7, 38, 60, 93, 111]]
\item \{1=43344, 2=655452, 3=2989224, 4=3871980\} [[0, 1, 2, 5, 8, 14], [0, 3, 7, 36, 47, 120], [0, 4, 70, 72, 75, 97], [0, 6, 43, 46, 82, 114], [0, 12, 23, 29, 83, 117]]
\item \{1=43344, 2=657720, 3=3001824, 4=3857112\} [[0, 1, 2, 5, 8, 14], [0, 3, 16, 37, 101, 108], [0, 4, 20, 88, 89, 114], [0, 6, 57, 62, 110, 115], [0, 12, 38, 79, 84, 95]]
\item \{1=43344, 2=662256, 3=2995776, 4=3858624\} [[0, 1, 2, 5, 8, 14], [0, 3, 16, 30, 69, 122], [0, 4, 47, 77, 116, 125], [0, 6, 13, 37, 95, 97], [0, 9, 33, 35, 102, 104]]
\item \{1=44352, 2=568512, 3=2984688, 4=3962448\} [[0, 1, 2, 5, 8, 14], [0, 3, 16, 45, 84, 103], [0, 4, 19, 36, 38, 120], [0, 6, 44, 58, 115, 116], [0, 12, 39, 51, 77, 85]]
\item \{1=44352, 2=570780, 3=3001320, 4=3943548\} [[0, 1, 2, 5, 8, 14], [0, 3, 21, 45, 96, 120], [0, 4, 36, 61, 103, 125], [0, 6, 27, 33, 82, 85], [0, 13, 37, 65, 71, 112]]
\item \{1=44352, 2=605556, 3=3039624, 4=3870468\} [[0, 1, 2, 5, 8, 14], [0, 3, 7, 22, 67, 91], [0, 4, 38, 68, 76, 119], [0, 10, 25, 52, 96, 99], [0, 12, 48, 84, 103, 109]]
\item \{1=44352, 2=624456, 3=3009888, 4=3881304\} [[0, 1, 2, 5, 8, 14], [0, 3, 7, 23, 85, 87], [0, 4, 27, 42, 86, 107], [0, 6, 59, 72, 99, 119], [0, 9, 46, 77, 108, 117]]
\item \{1=44352, 2=626724, 3=3003336, 4=3885588\} [[0, 1, 2, 5, 8, 14], [0, 3, 4, 56, 84, 113], [0, 6, 43, 88, 94, 124], [0, 7, 39, 41, 42, 76], [0, 9, 50, 63, 103, 107]]
\item \{1=44352, 2=628236, 3=3019464, 4=3867948\} [[0, 1, 2, 5, 8, 14], [0, 3, 4, 29, 41, 69], [0, 6, 58, 62, 74, 90], [0, 7, 22, 40, 65, 93], [0, 18, 43, 53, 77, 98]]
\item \{1=44352, 2=635796, 3=2971080, 4=3908772\} [[0, 1, 3, 6, 11, 17], [0, 2, 4, 23, 34, 108], [0, 5, 27, 49, 86, 118], [0, 7, 28, 40, 48, 102], [0, 10, 20, 57, 101, 116]]
\item \{1=44352, 2=638064, 3=3011904, 4=3865680\} [[0, 1, 2, 5, 8, 14], [0, 3, 10, 75, 106, 122], [0, 4, 19, 26, 61, 112], [0, 9, 35, 76, 97, 113], [0, 16, 54, 69, 94, 99], [0, 18, 37, 81, 103, 109], [0, 33, 59, 92, 115, 125]]
\item \{1=44352, 2=656964, 3=3012408, 4=3846276\} [[0, 1, 2, 5, 8, 14], [0, 3, 4, 73, 81, 103], [0, 6, 60, 77, 96, 108], [0, 7, 55, 59, 71, 78], [0, 12, 15, 54, 112, 122]]
\item \{1=44352, 2=657720, 3=2984688, 4=3873240\} [[0, 1, 2, 5, 8, 14], [0, 3, 10, 32, 116, 118], [0, 4, 37, 71, 81, 100], [0, 6, 33, 46, 99, 125], [0, 13, 65, 70, 101, 102]]
\item \{1=44352, 2=669060, 3=2965032, 4=3881556\} [[0, 1, 2, 5, 8, 14], [0, 3, 16, 37, 60, 96], [0, 4, 18, 54, 57, 103], [0, 7, 46, 50, 73, 81], [0, 9, 32, 35, 86, 97]]
\item \{1=44352, 2=682668, 3=3029544, 4=3803436\} [[0, 1, 2, 5, 8, 14], [0, 3, 21, 42, 69, 114], [0, 4, 27, 36, 66, 78], [0, 6, 56, 75, 112, 123], [0, 12, 73, 93, 117, 121]]
\item \{1=45360, 2=560952, 3=2973600, 4=3980088\} [[0, 1, 2, 5, 8, 14], [0, 3, 4, 36, 38, 116], [0, 6, 70, 76, 111, 122], [0, 7, 30, 56, 60, 119], [0, 13, 53, 73, 110, 117]]
\item \{1=45360, 2=566244, 3=3009384, 4=3939012\} [[0, 1, 2, 5, 8, 14], [0, 3, 7, 54, 59, 96], [0, 4, 18, 86, 88, 123], [0, 9, 39, 43, 63, 100], [0, 13, 24, 30, 53, 125]]
\item \{1=45360, 2=597996, 3=2976120, 4=3940524\} [[0, 1, 2, 5, 8, 14], [0, 3, 4, 30, 38, 54], [0, 6, 42, 70, 100, 116], [0, 7, 48, 52, 56, 119], [0, 9, 55, 67, 93, 102]]
\item \{1=45360, 2=602532, 3=2970072, 4=3942036\} [[0, 1, 2, 5, 8, 14], [0, 3, 7, 55, 98, 112], [0, 6, 27, 29, 32, 121], [0, 10, 21, 42, 84, 120], [0, 12, 30, 52, 68, 115]]
\item \{1=45360, 2=614628, 3=2984184, 4=3915828\} [[0, 1, 2, 5, 8, 14], [0, 3, 21, 33, 60, 82], [0, 4, 23, 37, 86, 119], [0, 6, 22, 70, 72, 89], [0, 9, 10, 26, 83, 88]]
\item \{1=45360, 2=619164, 3=2992248, 4=3903228\} [[0, 1, 2, 5, 8, 14], [0, 3, 16, 30, 50, 113], [0, 4, 19, 41, 79, 118], [0, 6, 21, 92, 98, 124], [0, 9, 36, 105, 114, 123], [0, 18, 86, 93, 107, 125], [0, 24, 69, 81, 108, 115]]
\item \{1=45360, 2=628236, 3=2938824, 4=3947580\} [[0, 1, 2, 5, 8, 14], [0, 3, 4, 19, 90, 107], [0, 6, 41, 91, 102, 119], [0, 10, 49, 82, 115, 122], [0, 12, 15, 53, 96, 112]]
\item \{1=45360, 2=635040, 3=3054240, 4=3825360\} [[0, 1, 2, 5, 8, 14], [0, 3, 4, 56, 82, 121], [0, 7, 62, 69, 95, 109], [0, 9, 19, 36, 49, 81], [0, 18, 47, 78, 114, 116], [0, 20, 22, 29, 31, 33], [0, 21, 61, 75, 106, 110]]
\item \{1=45360, 2=637308, 3=3015432, 4=3861900\} [[0, 1, 2, 5, 8, 14], [0, 3, 4, 50, 57, 103], [0, 6, 26, 30, 86, 87], [0, 7, 35, 58, 73, 106], [0, 10, 46, 56, 88, 116], [0, 21, 23, 43, 77, 91], [0, 32, 55, 63, 75, 105]]
\item \{1=45360, 2=654696, 3=2964528, 4=3895416\} [[0, 1, 2, 5, 8, 14], [0, 3, 10, 45, 111, 121], [0, 4, 34, 40, 48, 109], [0, 6, 21, 70, 107, 122], [0, 19, 67, 87, 112, 114]]
\item \{1=45360, 2=676620, 3=3021480, 4=3816540\} [[0, 1, 2, 5, 8, 14], [0, 3, 4, 35, 77, 121], [0, 6, 13, 39, 47, 52], [0, 7, 43, 75, 94, 99], [0, 10, 33, 40, 108, 120]]
\item \{1=46368, 2=591948, 3=2995272, 4=3926412\} [[0, 1, 2, 5, 8, 14], [0, 3, 7, 23, 85, 103], [0, 4, 49, 61, 108, 118], [0, 6, 36, 37, 70, 96], [0, 13, 75, 81, 115, 119]]
\item \{1=46368, 2=592704, 3=2996784, 4=3924144\} [[0, 1, 2, 5, 8, 14], [0, 3, 7, 22, 90, 95], [0, 4, 56, 69, 82, 119], [0, 12, 47, 83, 84, 107], [0, 16, 19, 61, 67, 116]]
\item \{1=46368, 2=597996, 3=3037608, 4=3878028\} [[0, 1, 2, 5, 8, 14], [0, 3, 4, 36, 60, 98], [0, 6, 24, 34, 92, 111], [0, 7, 52, 102, 119, 120], [0, 16, 28, 71, 99, 110]]
\item \{1=46368, 2=609336, 3=3013920, 4=3890376\} [[0, 1, 2, 5, 8, 14], [0, 3, 7, 41, 51, 91], [0, 4, 23, 68, 93, 124], [0, 9, 24, 67, 92, 102], [0, 12, 50, 60, 103, 125]]
\item \{1=46368, 2=619920, 3=2994768, 4=3898944\} [[0, 1, 2, 5, 8, 14], [0, 3, 4, 21, 63, 65], [0, 7, 47, 69, 82, 113], [0, 10, 18, 60, 81, 84], [0, 15, 29, 53, 75, 85]]
\item \{1=46368, 2=632772, 3=2930760, 4=3950100\} [[0, 1, 2, 5, 8, 14], [0, 3, 16, 104, 114, 115], [0, 4, 21, 26, 71, 75], [0, 6, 25, 56, 60, 103], [0, 7, 38, 42, 67, 116]]
\item \{1=46368, 2=634284, 3=2978136, 4=3901212\} [[0, 1, 2, 5, 8, 14], [0, 3, 21, 30, 82, 119], [0, 4, 26, 47, 70, 125], [0, 6, 25, 37, 64, 74], [0, 10, 13, 66, 90, 124]]
\item \{1=46368, 2=640332, 3=2926728, 4=3946572\} [[0, 1, 2, 5, 8, 14], [0, 3, 10, 24, 29, 112], [0, 6, 34, 69, 79, 81], [0, 7, 15, 106, 110, 121], [0, 9, 64, 76, 90, 120], [0, 20, 22, 50, 52, 98], [0, 33, 59, 92, 115, 125]]
\item \{1=46368, 2=641844, 3=2991240, 4=3880548\} [[0, 1, 2, 5, 8, 14], [0, 3, 4, 39, 69, 76], [0, 6, 22, 33, 55, 108], [0, 7, 26, 30, 44, 100], [0, 9, 13, 74, 102, 122]]
\item \{1=47376, 2=576072, 3=3029040, 4=3907512\} [[0, 1, 2, 5, 8, 14], [0, 3, 21, 78, 100, 114], [0, 6, 13, 64, 73, 104], [0, 9, 25, 72, 96, 122], [0, 12, 102, 106, 113, 120]]
\item \{1=47376, 2=597996, 3=2980152, 4=3934476\} [[0, 1, 2, 5, 8, 14], [0, 3, 4, 50, 100, 104], [0, 6, 22, 40, 48, 75], [0, 10, 37, 58, 69, 94], [0, 18, 39, 59, 63, 118]]
\item \{1=47376, 2=607068, 3=2982168, 4=3923388\} [[0, 1, 2, 5, 8, 14], [0, 3, 4, 51, 77, 121], [0, 6, 23, 27, 72, 91], [0, 7, 22, 35, 59, 94], [0, 10, 19, 30, 76, 99]]
\item \{1=47376, 2=617652, 3=2995272, 4=3899700\} [[0, 1, 2, 5, 8, 14], [0, 3, 7, 31, 66, 67], [0, 4, 35, 72, 96, 115], [0, 6, 36, 70, 99, 109], [0, 15, 39, 77, 79, 107]]
\item \{1=47376, 2=619920, 3=2965536, 4=3927168\} [[0, 1, 2, 5, 8, 14], [0, 3, 4, 23, 59, 116], [0, 6, 44, 54, 64, 73], [0, 9, 18, 32, 69, 119], [0, 12, 62, 66, 81, 120]]
\item \{1=47376, 2=621432, 3=3011904, 4=3879288\} [[0, 1, 2, 5, 8, 14], [0, 3, 10, 37, 116, 124], [0, 4, 27, 28, 33, 78], [0, 7, 49, 65, 81, 95], [0, 15, 41, 55, 66, 99]]
\item \{1=47376, 2=628992, 3=2982672, 4=3900960\} [[0, 1, 2, 5, 8, 14], [0, 3, 29, 51, 98, 99], [0, 4, 66, 71, 113, 125], [0, 6, 40, 60, 70, 76], [0, 7, 15, 106, 110, 121], [0, 10, 26, 56, 58, 107], [0, 18, 21, 64, 77, 116]]
\item \{1=47376, 2=632016, 3=2976624, 4=3903984\} [[0, 1, 2, 5, 8, 14], [0, 3, 21, 30, 39, 92], [0, 4, 29, 33, 55, 108], [0, 9, 50, 66, 87, 109], [0, 13, 20, 48, 97, 116]]
\item \{1=47376, 2=635040, 3=2978640, 4=3898944\} [[0, 1, 2, 5, 8, 14], [0, 3, 10, 69, 93, 106], [0, 4, 19, 22, 102, 103], [0, 6, 30, 55, 85, 115], [0, 15, 21, 74, 98, 124]]
\item \{1=47376, 2=646380, 3=2995272, 4=3870972\} [[0, 1, 2, 5, 8, 14], [0, 3, 4, 41, 98, 110], [0, 6, 45, 54, 104, 121], [0, 9, 32, 90, 92, 122], [0, 12, 44, 52, 86, 112]]
\item \{1=47376, 2=682668, 3=2982168, 4=3847788\} [[0, 1, 2, 5, 8, 14], [0, 3, 7, 51, 122, 123], [0, 4, 29, 35, 40, 43], [0, 12, 27, 30, 54, 59], [0, 13, 52, 58, 64, 125]]
\item \{1=48384, 2=601776, 3=3018960, 4=3890880\} [[0, 1, 2, 5, 8, 14], [0, 3, 4, 59, 60, 102], [0, 6, 13, 16, 36, 112], [0, 7, 27, 34, 77, 90], [0, 10, 83, 87, 101, 122]]
\item \{1=48384, 2=610848, 3=3018960, 4=3881808\} [[0, 1, 2, 5, 8, 14], [0, 3, 4, 100, 110, 115], [0, 6, 39, 55, 59, 91], [0, 7, 34, 49, 56, 61], [0, 15, 41, 44, 66, 109]]
\item \{1=48384, 2=622944, 3=3004848, 4=3883824\} [[0, 1, 2, 5, 8, 14], [0, 3, 4, 18, 65, 88], [0, 9, 36, 61, 74, 93], [0, 10, 15, 56, 64, 105], [0, 12, 30, 55, 113, 115], [0, 20, 23, 77, 78, 90], [0, 26, 41, 82, 87, 100]]
\item \{1=48384, 2=624456, 3=3035088, 4=3852072\} [[0, 1, 2, 5, 8, 14], [0, 3, 10, 25, 30, 60], [0, 4, 35, 63, 82, 88], [0, 6, 41, 49, 54, 118], [0, 7, 21, 47, 77, 103], [0, 9, 36, 68, 87, 99], [0, 20, 22, 32, 34, 81]]
\item \{1=49392, 2=624456, 3=2992752, 4=3893400\} [[0, 1, 2, 5, 8, 14], [0, 3, 7, 25, 60, 116], [0, 4, 9, 68, 82, 107], [0, 10, 56, 63, 104, 123], [0, 15, 53, 58, 72, 125], [0, 16, 18, 84, 86, 98], [0, 24, 69, 81, 108, 115]]
\item \{1=49392, 2=641088, 3=2973600, 4=3895920\} [[0, 1, 2, 5, 8, 14], [0, 3, 4, 69, 76, 89], [0, 6, 10, 22, 29, 122], [0, 9, 58, 64, 104, 107], [0, 13, 73, 81, 96, 120]]
\item \{1=49392, 2=647136, 3=2956464, 4=3907008\} [[0, 1, 2, 5, 8, 14], [0, 3, 4, 73, 81, 114], [0, 6, 13, 48, 71, 102], [0, 7, 29, 91, 93, 103], [0, 9, 24, 58, 100, 106]]
\item \{1=49392, 2=656964, 3=3002328, 4=3851316\} [[0, 1, 2, 5, 8, 14], [0, 3, 10, 40, 77, 84], [0, 4, 48, 63, 70, 89], [0, 6, 45, 58, 94, 113], [0, 16, 46, 66, 91, 100]]
\item \{1=50400, 2=592704, 3=2996784, 4=3920112\} [[0, 1, 2, 5, 8, 14], [0, 3, 4, 39, 49, 64], [0, 6, 12, 60, 63, 104], [0, 9, 62, 74, 75, 98], [0, 10, 27, 86, 106, 122]]
\item \{1=50400, 2=612360, 3=3029040, 4=3868200\} [[0, 1, 2, 5, 8, 14], [0, 3, 4, 30, 54, 81], [0, 6, 34, 44, 69, 119], [0, 9, 21, 49, 90, 96], [0, 10, 29, 60, 66, 92]]
\item \{1=50400, 2=642600, 3=2927232, 4=3939768\} [[0, 1, 2, 5, 8, 14], [0, 3, 4, 19, 58, 99], [0, 6, 60, 76, 82, 111], [0, 9, 35, 42, 57, 124], [0, 13, 65, 70, 93, 122]]
\item \{1=50400, 2=672084, 3=2996280, 4=3841236\} [[0, 1, 2, 5, 8, 14], [0, 3, 16, 46, 51, 100], [0, 4, 6, 21, 31, 60], [0, 7, 26, 30, 90, 101], [0, 12, 13, 53, 93, 108]]
\item \{1=51408, 2=613872, 3=3020976, 4=3873744\} [[0, 1, 2, 5, 8, 14], [0, 3, 7, 22, 35, 72], [0, 4, 54, 73, 98, 120], [0, 10, 30, 40, 52, 104], [0, 15, 61, 67, 68, 113]]
\item \{1=51408, 2=678888, 3=2988720, 4=3840984\} [[0, 1, 2, 5, 8, 14], [0, 3, 7, 49, 50, 85], [0, 4, 20, 61, 101, 110], [0, 6, 44, 48, 103, 109], [0, 15, 19, 57, 65, 70]]
\item \{1=52416, 2=591192, 3=3040128, 4=3876264\} [[0, 1, 2, 5, 8, 14], [0, 3, 4, 26, 75, 86], [0, 6, 32, 64, 103, 107], [0, 7, 16, 60, 70, 124], [0, 9, 24, 92, 102, 112]]
\item \{1=52416, 2=598752, 3=3049200, 4=3859632\} [[0, 1, 2, 5, 8, 14], [0, 3, 7, 23, 68, 102], [0, 4, 20, 45, 81, 92], [0, 9, 46, 61, 75, 120], [0, 18, 93, 106, 119, 123]]
\item \{1=52416, 2=615384, 3=3039120, 4=3853080\} [[0, 1, 2, 5, 8, 14], [0, 3, 4, 16, 52, 62], [0, 6, 36, 55, 77, 116], [0, 7, 58, 72, 85, 98], [0, 10, 15, 90, 101, 108], [0, 12, 38, 54, 94, 120], [0, 19, 20, 22, 64, 66]]
\item \{1=52416, 2=626724, 3=2998296, 4=3882564\} [[0, 1, 2, 5, 8, 14], [0, 3, 4, 21, 39, 85], [0, 6, 84, 87, 111, 118], [0, 7, 26, 53, 82, 115], [0, 9, 43, 66, 79, 114]]
\item \{1=53424, 2=608580, 3=2957976, 4=3940020\} [[0, 1, 2, 5, 8, 14], [0, 3, 7, 49, 101, 123], [0, 4, 41, 88, 102, 103], [0, 6, 44, 58, 77, 90], [0, 12, 15, 37, 39, 84]]
\item \{1=53424, 2=615384, 3=2982672, 4=3908520\} [[0, 1, 2, 5, 8, 14], [0, 3, 4, 16, 23, 108], [0, 7, 44, 102, 117, 119], [0, 9, 12, 19, 104, 118], [0, 13, 46, 78, 92, 123]]
\item \{1=53424, 2=642600, 3=3056256, 4=3807720\} [[0, 1, 2, 5, 8, 14], [0, 3, 16, 42, 75, 120], [0, 4, 22, 33, 68, 115], [0, 7, 23, 35, 36, 73], [0, 12, 59, 63, 83, 107]]
\item \{1=53424, 2=651672, 3=3018960, 4=3835944\} [[0, 1, 2, 5, 8, 14], [0, 3, 4, 30, 85, 88], [0, 6, 67, 83, 87, 89], [0, 9, 16, 29, 114, 117], [0, 20, 42, 50, 53, 91]]
\item \{1=54432, 2=615384, 3=2980656, 4=3909528\} [[0, 1, 2, 5, 8, 14], [0, 3, 4, 31, 36, 95], [0, 6, 37, 67, 78, 111], [0, 7, 21, 44, 110, 117], [0, 10, 19, 87, 93, 99]]
\item \{1=54432, 2=655452, 3=3015432, 4=3834684\} [[0, 1, 2, 5, 8, 14], [0, 3, 10, 33, 68, 119], [0, 4, 26, 28, 58, 62], [0, 6, 22, 64, 106, 121], [0, 19, 47, 49, 69, 79]]
\item \{1=54432, 2=666792, 3=2936304, 4=3902472\} [[0, 1, 2, 5, 8, 14], [0, 3, 4, 47, 52, 55], [0, 7, 28, 57, 84, 119], [0, 9, 59, 69, 71, 94], [0, 10, 42, 45, 77, 116], [0, 20, 22, 29, 31, 33], [0, 24, 60, 72, 108, 109]]
\item \{1=57456, 2=707616, 3=2950416, 4=3844512\} [[0, 1, 2, 5, 8, 14], [0, 3, 4, 35, 81, 117], [0, 6, 68, 74, 102, 125], [0, 7, 10, 55, 91, 96], [0, 15, 32, 38, 78, 92]]
\item \{1=59472, 2=669816, 3=3013920, 4=3816792\} [[0, 1, 2, 5, 8, 14], [0, 3, 10, 25, 68, 77], [0, 4, 34, 40, 96, 114], [0, 6, 45, 79, 81, 93], [0, 24, 51, 84, 98, 116]]
\item \{1=59472, 2=694764, 3=2973096, 4=3832668\} [[0, 1, 2, 5, 8, 14], [0, 3, 10, 21, 85, 111], [0, 4, 33, 59, 88, 114], [0, 7, 23, 46, 53, 117], [0, 13, 37, 98, 110, 122]]
\item \{1=62496, 2=658476, 3=3043656, 4=3795372\} [[0, 1, 2, 5, 8, 14], [0, 3, 7, 25, 96, 99], [0, 4, 55, 67, 108, 117], [0, 9, 60, 91, 101, 119], [0, 10, 35, 56, 68, 114], [0, 16, 51, 106, 110, 122], [0, 19, 58, 72, 102, 112]]
\item \{1=63504, 2=706860, 3=2967048, 4=3822588\} [[0, 1, 2, 5, 8, 14], [0, 3, 4, 48, 58, 101], [0, 6, 33, 40, 78, 120], [0, 7, 22, 63, 96, 109], [0, 9, 19, 36, 49, 81], [0, 10, 42, 56, 75, 118], [0, 15, 57, 83, 104, 119]]
\item \{0=1260, 1=18144, 2=577584, 3=2981664, 4=3981348\} [[0, 1, 2, 5, 8, 14], [0, 3, 7, 16, 25, 102], [0, 4, 30, 47, 91, 93], [0, 9, 32, 59, 62, 87], [0, 15, 66, 69, 85, 125]]
\item \{0=1260, 1=19152, 2=601020, 3=2994264, 4=3944304\} [[0, 1, 2, 5, 8, 14], [0, 3, 4, 31, 69, 75], [0, 6, 24, 64, 94, 111], [0, 7, 44, 48, 50, 123], [0, 9, 36, 79, 92, 109], [0, 15, 18, 74, 95, 113], [0, 19, 23, 67, 90, 121]]
\item \{0=1260, 1=20160, 2=593460, 3=3002328, 4=3942792\} [[0, 1, 2, 5, 8, 14], [0, 3, 4, 35, 86, 99], [0, 6, 40, 87, 101, 121], [0, 7, 63, 84, 112, 122], [0, 12, 54, 62, 67, 117]]
\item \{0=1260, 1=20160, 2=624456, 3=3009888, 4=3904236\} [[0, 1, 3, 6, 11, 17], [0, 2, 4, 20, 42, 114], [0, 5, 10, 41, 76, 112], [0, 15, 51, 64, 91, 125], [0, 16, 29, 57, 82, 109]]
\item \{0=1260, 1=22176, 2=528444, 3=3000312, 4=4007808\} [[0, 1, 2, 5, 8, 14], [0, 3, 16, 50, 74, 124], [0, 4, 36, 56, 85, 87], [0, 6, 61, 91, 96, 117], [0, 15, 19, 54, 94, 98], [0, 20, 39, 44, 73, 109], [0, 21, 33, 77, 81, 123]]
\item \{0=1260, 1=22176, 2=563976, 3=3012912, 4=3959676\} [[0, 1, 2, 5, 8, 14], [0, 3, 10, 49, 67, 120], [0, 6, 41, 88, 89, 90], [0, 9, 33, 63, 71, 83], [0, 16, 35, 68, 112, 122]]
\item \{0=1260, 1=22176, 2=565488, 3=2975616, 4=3995460\} [[0, 1, 2, 5, 8, 14], [0, 3, 21, 51, 73, 91], [0, 6, 40, 41, 60, 65], [0, 7, 9, 53, 79, 82], [0, 12, 15, 72, 105, 106]]
\item \{0=1260, 1=22176, 2=601020, 3=2969064, 4=3966480\} [[0, 1, 3, 6, 11, 17], [0, 2, 4, 20, 47, 60], [0, 5, 70, 97, 104, 112], [0, 7, 62, 63, 65, 109], [0, 10, 30, 85, 98, 125]]
\item \{0=1260, 1=23184, 2=560952, 3=2971584, 4=4003020\} [[0, 1, 2, 5, 8, 14], [0, 3, 16, 29, 33, 69], [0, 4, 67, 71, 92, 116], [0, 6, 54, 55, 96, 108], [0, 7, 10, 34, 93, 119]]
\item \{0=1260, 1=23184, 2=573048, 3=3062304, 4=3900204\} [[0, 1, 2, 5, 8, 14], [0, 3, 7, 16, 40, 54], [0, 4, 29, 64, 69, 91], [0, 6, 21, 43, 58, 62], [0, 10, 26, 81, 97, 114]]
\item \{0=1260, 1=23184, 2=602532, 3=2985192, 4=3947832\} [[0, 1, 2, 5, 8, 14], [0, 3, 10, 37, 72, 116], [0, 4, 34, 59, 67, 71], [0, 6, 29, 70, 118, 125], [0, 13, 63, 64, 114, 123]]
\item \{0=1260, 1=23184, 2=604800, 3=3020976, 4=3909780\} [[0, 1, 2, 5, 8, 14], [0, 3, 4, 57, 64, 100], [0, 6, 42, 44, 90, 112], [0, 9, 12, 74, 95, 107], [0, 10, 15, 49, 76, 108], [0, 28, 32, 69, 106, 110], [0, 39, 63, 73, 88, 103]]
\item \{0=1260, 1=24192, 2=543564, 3=2989224, 4=4001760\} [[0, 1, 2, 5, 8, 14], [0, 3, 4, 35, 39, 60], [0, 7, 36, 100, 102, 117], [0, 10, 13, 40, 64, 113], [0, 28, 49, 52, 74, 95], [0, 32, 76, 82, 104, 118], [0, 33, 37, 47, 109, 116]]
\item \{0=1260, 1=24192, 2=613116, 3=2982168, 4=3939264\} [[0, 1, 2, 5, 8, 14], [0, 3, 21, 33, 69, 111], [0, 4, 47, 49, 115, 121], [0, 6, 36, 40, 67, 113], [0, 10, 75, 112, 119, 120]]
\item \{0=1260, 1=25200, 2=543564, 3=3026520, 4=3963456\} [[0, 1, 2, 5, 8, 14], [0, 3, 10, 51, 68, 103], [0, 4, 54, 55, 67, 108], [0, 6, 27, 111, 116, 119], [0, 7, 46, 92, 102, 117]]
\item \{0=1260, 1=25200, 2=570024, 3=2990736, 4=3972780\} [[0, 1, 2, 5, 8, 14], [0, 3, 16, 32, 33, 102], [0, 6, 19, 22, 47, 107], [0, 7, 35, 54, 111, 120], [0, 23, 24, 46, 95, 124]]
\item \{0=1260, 1=25200, 2=585144, 3=2994768, 4=3953628\} [[0, 1, 2, 5, 8, 14], [0, 3, 4, 31, 43, 89], [0, 6, 57, 81, 104, 113], [0, 9, 12, 47, 94, 106], [0, 15, 33, 91, 100, 120], [0, 24, 52, 64, 102, 108], [0, 28, 39, 53, 73, 115]]
\item \{0=1260, 1=25200, 2=601020, 3=3012408, 4=3920112\} [[0, 1, 2, 5, 8, 14], [0, 3, 10, 66, 96, 123], [0, 4, 70, 77, 90, 93], [0, 6, 36, 43, 63, 100], [0, 13, 31, 32, 58, 92]]
\item \{0=1260, 1=25200, 2=627480, 3=3012912, 4=3893148\} [[0, 1, 2, 5, 8, 14], [0, 3, 7, 10, 60, 63], [0, 4, 22, 66, 95, 111], [0, 6, 44, 50, 94, 121], [0, 9, 13, 38, 48, 67]]
\item \{0=1260, 1=26208, 2=502740, 3=3020472, 4=4009320\} [[0, 1, 2, 5, 8, 14], [0, 3, 7, 66, 81, 121], [0, 4, 18, 43, 46, 79], [0, 9, 29, 39, 88, 124], [0, 15, 55, 75, 87, 114], [0, 23, 24, 97, 107, 108], [0, 31, 37, 69, 83, 109]]
\item \{0=1260, 1=26208, 2=583632, 3=2955456, 4=3993444\} [[0, 1, 2, 5, 8, 14], [0, 3, 21, 25, 63, 119], [0, 4, 51, 95, 112, 114], [0, 7, 28, 61, 66, 86], [0, 9, 40, 44, 69, 85]]
\item \{0=1260, 1=26208, 2=589680, 3=2991744, 4=3951108\} [[0, 1, 2, 5, 8, 14], [0, 3, 4, 40, 46, 101], [0, 6, 19, 50, 52, 112], [0, 7, 28, 58, 70, 118], [0, 10, 41, 49, 64, 94]]
\item \{0=1260, 1=26208, 2=598752, 3=3039120, 4=3894660\} [[0, 1, 2, 5, 8, 14], [0, 3, 10, 29, 75, 119], [0, 4, 9, 30, 68, 116], [0, 13, 53, 65, 106, 107], [0, 15, 27, 70, 84, 99]]
\item \{0=1260, 1=27216, 2=569268, 3=2966040, 4=3996216\} [[0, 1, 2, 5, 8, 14], [0, 3, 7, 30, 48, 102], [0, 6, 50, 57, 101, 123], [0, 9, 32, 36, 51, 64], [0, 10, 52, 56, 85, 122], [0, 12, 46, 67, 69, 116], [0, 39, 42, 74, 95, 124]]
\item \{0=1260, 1=27216, 2=570024, 3=2956464, 4=4005036\} [[0, 1, 2, 5, 8, 14], [0, 3, 16, 21, 45, 120], [0, 4, 33, 67, 70, 76], [0, 9, 49, 85, 113, 115], [0, 20, 25, 27, 38, 102]]
\item \{0=1260, 1=27216, 2=591192, 3=2996784, 4=3943548\} [[0, 1, 2, 5, 8, 14], [0, 3, 4, 29, 94, 104], [0, 6, 26, 48, 107, 123], [0, 7, 40, 78, 108, 113], [0, 9, 41, 53, 76, 79]]
\item \{0=1260, 1=27216, 2=594216, 3=2978640, 4=3958668\} [[0, 1, 2, 5, 8, 14], [0, 3, 4, 41, 52, 82], [0, 6, 26, 30, 107, 111], [0, 7, 15, 39, 63, 109], [0, 13, 24, 79, 116, 124]]
\item \{0=1260, 1=27216, 2=604044, 3=3000312, 4=3927168\} [[0, 1, 2, 5, 8, 14], [0, 3, 10, 32, 36, 76], [0, 4, 16, 77, 93, 123], [0, 6, 13, 81, 94, 117], [0, 15, 60, 72, 85, 91]]
\item \{0=1260, 1=27216, 2=606312, 3=2950416, 4=3974796\} [[0, 1, 2, 5, 8, 14], [0, 3, 4, 45, 60, 90], [0, 6, 57, 96, 101, 109], [0, 9, 56, 87, 103, 112], [0, 12, 68, 72, 106, 114]]
\item \{0=1260, 1=27216, 2=628992, 3=3027024, 4=3875508\} [[0, 1, 2, 5, 8, 14], [0, 3, 4, 29, 75, 79], [0, 6, 82, 96, 107, 115], [0, 7, 45, 53, 65, 119], [0, 12, 37, 47, 69, 110]]
\item \{0=1260, 1=28224, 2=554148, 3=3003336, 4=3973032\} [[0, 1, 2, 5, 8, 14], [0, 3, 16, 45, 96, 99], [0, 4, 51, 64, 67, 87], [0, 6, 55, 116, 117, 122], [0, 7, 27, 59, 107, 111]]
\item \{0=1260, 1=28224, 2=561708, 3=2953944, 4=4014864\} [[0, 1, 2, 5, 8, 14], [0, 3, 10, 51, 57, 92], [0, 4, 9, 17, 61, 67], [0, 7, 27, 42, 48, 107], [0, 13, 24, 63, 65, 116]]
\item \{0=1260, 1=28224, 2=594216, 3=2969568, 4=3966732\} [[0, 1, 2, 5, 8, 14], [0, 3, 4, 75, 88, 106], [0, 6, 57, 60, 90, 100], [0, 7, 62, 65, 117, 123], [0, 12, 29, 79, 102, 121]]
\item \{0=1260, 1=28224, 2=603288, 3=3003840, 4=3923388\} [[0, 1, 2, 5, 8, 14], [0, 3, 4, 45, 107, 124], [0, 6, 55, 69, 88, 101], [0, 9, 18, 61, 102, 108], [0, 10, 66, 96, 99, 115]]
\item \{0=1260, 1=28224, 2=605556, 3=2983176, 4=3941784\} [[0, 1, 2, 5, 8, 14], [0, 3, 16, 33, 84, 101], [0, 6, 44, 55, 96, 120], [0, 7, 34, 49, 103, 118], [0, 13, 53, 66, 72, 114]]
\item \{0=1260, 1=28224, 2=605556, 3=3002328, 4=3922632\} [[0, 1, 2, 5, 8, 14], [0, 3, 4, 68, 106, 119], [0, 6, 46, 86, 96, 99], [0, 7, 28, 43, 58, 87], [0, 10, 15, 41, 103, 124]]
\item \{0=1260, 1=28224, 2=619920, 3=3020976, 4=3889620\} [[0, 1, 2, 5, 8, 14], [0, 3, 4, 49, 75, 92], [0, 6, 38, 54, 61, 68], [0, 7, 10, 27, 85, 120], [0, 13, 78, 105, 112, 123]]
\item \{0=1260, 1=28224, 2=629748, 3=2993256, 4=3907512\} [[0, 1, 2, 5, 8, 14], [0, 3, 10, 36, 46, 87], [0, 4, 16, 92, 111, 122], [0, 6, 19, 62, 95, 121], [0, 7, 48, 88, 101, 102]]
\item \{0=1260, 1=28224, 2=632772, 3=2948904, 4=3948840\} [[0, 1, 2, 5, 8, 14], [0, 3, 21, 30, 113, 123], [0, 6, 19, 23, 43, 57], [0, 7, 52, 63, 64, 88], [0, 9, 47, 76, 91, 105]]
\item \{0=1260, 1=29232, 2=550368, 3=2976624, 4=4002516\} [[0, 1, 2, 5, 8, 14], [0, 3, 10, 72, 91, 115], [0, 4, 18, 38, 43, 100], [0, 7, 30, 60, 64, 111], [0, 9, 47, 98, 101, 123]]
\item \{0=1260, 1=29232, 2=568512, 3=3018960, 4=3942036\} [[0, 1, 3, 6, 11, 17], [0, 2, 4, 19, 40, 78], [0, 5, 27, 42, 48, 119], [0, 10, 35, 64, 108, 123], [0, 13, 41, 46, 69, 72]]
\item \{0=1260, 1=29232, 2=587412, 3=2974104, 4=3967992\} [[0, 1, 2, 5, 8, 14], [0, 3, 4, 69, 77, 122], [0, 6, 43, 73, 79, 83], [0, 7, 20, 22, 46, 48], [0, 9, 27, 36, 58, 90], [0, 10, 28, 34, 40, 93], [0, 18, 55, 68, 75, 98]]
\item \{0=1260, 1=29232, 2=612360, 3=2990736, 4=3926412\} [[0, 1, 2, 5, 8, 14], [0, 3, 7, 38, 52, 63], [0, 4, 16, 28, 55, 92], [0, 6, 13, 57, 78, 124], [0, 9, 27, 64, 98, 120]]
\item \{0=1260, 1=29232, 2=617652, 3=2954952, 4=3956904\} [[0, 1, 2, 5, 8, 14], [0, 3, 10, 25, 51, 69], [0, 6, 30, 53, 68, 90], [0, 13, 20, 22, 55, 57], [0, 15, 27, 74, 86, 124], [0, 18, 76, 114, 115, 118], [0, 21, 33, 77, 81, 123]]
\item \{0=1260, 1=29232, 2=625968, 3=3010896, 4=3892644\} [[0, 1, 2, 5, 8, 14], [0, 3, 4, 37, 68, 124], [0, 6, 33, 61, 108, 116], [0, 7, 9, 43, 71, 73], [0, 13, 22, 52, 65, 94]]
\item \{0=1260, 1=29232, 2=666792, 3=2993760, 4=3868956\} [[0, 1, 2, 5, 8, 14], [0, 3, 4, 60, 107, 122], [0, 6, 45, 56, 101, 103], [0, 12, 41, 48, 105, 123], [0, 16, 28, 35, 40, 68]]
\item \{0=1260, 1=29232, 2=678132, 3=3032568, 4=3818808\} [[0, 1, 2, 5, 8, 14], [0, 3, 16, 51, 115, 116], [0, 4, 59, 77, 92, 108], [0, 6, 61, 66, 68, 99], [0, 7, 30, 91, 113, 120], [0, 9, 19, 36, 49, 81], [0, 12, 15, 39, 75, 85]]
\item \{0=1260, 1=30240, 2=566244, 3=3023496, 4=3938760\} [[0, 1, 2, 5, 8, 14], [0, 3, 4, 48, 92, 106], [0, 6, 55, 65, 79, 124], [0, 7, 35, 38, 107, 125], [0, 15, 29, 57, 96, 109]]
\item \{0=1260, 1=30240, 2=601776, 3=3013920, 4=3912804\} [[0, 1, 2, 5, 8, 14], [0, 3, 7, 55, 103, 108], [0, 4, 17, 26, 27, 112], [0, 9, 21, 60, 93, 107], [0, 10, 28, 71, 72, 110], [0, 23, 68, 69, 90, 122], [0, 35, 57, 64, 115, 119]]
\item \{0=1260, 1=30240, 2=610092, 3=2958984, 4=3959424\} [[0, 1, 2, 5, 8, 14], [0, 3, 16, 36, 54, 66], [0, 6, 12, 44, 85, 104], [0, 10, 26, 29, 74, 84], [0, 13, 38, 58, 101, 120]]
\item \{0=1260, 1=30240, 2=615384, 3=2961504, 4=3951612\} [[0, 1, 2, 5, 8, 14], [0, 3, 4, 27, 90, 123], [0, 6, 54, 57, 78, 104], [0, 9, 63, 68, 95, 105], [0, 15, 18, 38, 58, 109]]
\item \{0=1260, 1=30240, 2=622944, 3=2967552, 4=3938004\} [[0, 1, 2, 5, 8, 14], [0, 3, 4, 37, 60, 105], [0, 6, 29, 40, 75, 92], [0, 9, 46, 62, 103, 119], [0, 10, 26, 63, 87, 112]]
\item \{0=1260, 1=31248, 2=525420, 3=3012408, 4=3989664\} [[0, 1, 2, 5, 8, 14], [0, 3, 21, 25, 65, 108], [0, 4, 6, 36, 104, 123], [0, 7, 10, 72, 88, 109], [0, 15, 67, 73, 93, 118]]
\item \{0=1260, 1=31248, 2=529200, 3=3025008, 4=3973284\} [[0, 1, 2, 5, 8, 14], [0, 3, 10, 36, 50, 91], [0, 4, 60, 98, 101, 119], [0, 6, 13, 23, 33, 121], [0, 7, 44, 96, 106, 114]]
\item \{0=1260, 1=31248, 2=581364, 3=3034584, 4=3911544\} [[0, 1, 2, 5, 8, 14], [0, 3, 10, 37, 91, 114], [0, 4, 55, 86, 97, 122], [0, 6, 33, 54, 71, 85], [0, 7, 36, 70, 99, 118]]
\item \{0=1260, 1=31248, 2=606312, 3=2994768, 4=3926412\} [[0, 1, 2, 5, 8, 14], [0, 3, 4, 34, 43, 120], [0, 6, 13, 40, 63, 82], [0, 12, 23, 38, 113, 115], [0, 16, 59, 72, 105, 108]]
\item \{0=1260, 1=31248, 2=628236, 3=2991240, 4=3908016\} [[0, 1, 2, 5, 8, 14], [0, 3, 4, 22, 76, 114], [0, 6, 79, 84, 87, 90], [0, 7, 26, 53, 55, 94], [0, 12, 39, 63, 98, 119]]
\item \{0=1260, 1=31248, 2=632016, 3=2972592, 4=3922884\} [[0, 1, 2, 5, 8, 14], [0, 3, 7, 10, 31, 102], [0, 4, 30, 64, 112, 125], [0, 9, 24, 68, 75, 107], [0, 12, 50, 88, 94, 118]]
\item \{0=1260, 1=31248, 2=635040, 3=3015936, 4=3876516\} [[0, 1, 2, 5, 8, 14], [0, 3, 16, 60, 74, 96], [0, 4, 9, 56, 71, 87], [0, 6, 41, 97, 113, 124], [0, 15, 67, 77, 91, 98]]
\item \{0=1260, 1=31248, 2=650160, 3=3015936, 4=3861396\} [[0, 1, 2, 5, 8, 14], [0, 3, 7, 46, 81, 122], [0, 4, 17, 20, 79, 99], [0, 9, 25, 27, 38, 125], [0, 10, 51, 82, 103, 108]]
\item \{0=1260, 1=32256, 2=588168, 3=2954448, 4=3983868\} [[0, 1, 2, 5, 8, 14], [0, 3, 4, 16, 36, 76], [0, 6, 42, 50, 57, 109], [0, 7, 28, 78, 94, 114], [0, 13, 53, 95, 96, 122]]
\item \{0=1260, 1=32256, 2=591948, 3=2999304, 4=3935232\} [[0, 1, 2, 5, 8, 14], [0, 3, 4, 35, 108, 115], [0, 6, 23, 51, 58, 60], [0, 9, 37, 43, 59, 92], [0, 12, 27, 77, 98, 116], [0, 20, 22, 103, 105, 125], [0, 32, 76, 82, 104, 118]]
\item \{0=1260, 1=32256, 2=595728, 3=2976624, 4=3954132\} [[0, 1, 2, 5, 8, 14], [0, 3, 4, 41, 57, 69], [0, 6, 30, 99, 104, 109], [0, 7, 23, 46, 94, 100], [0, 13, 18, 45, 58, 61]]
\item \{0=1260, 1=32256, 2=595728, 3=2978640, 4=3952116\} [[0, 1, 2, 5, 8, 14], [0, 3, 4, 37, 68, 75], [0, 6, 54, 71, 106, 122], [0, 7, 9, 29, 35, 70], [0, 16, 61, 76, 84, 92]]
\item \{0=1260, 1=32256, 2=598752, 3=2962512, 4=3965220\} [[0, 1, 2, 5, 8, 14], [0, 3, 21, 55, 63, 113], [0, 4, 40, 45, 59, 82], [0, 7, 29, 51, 56, 117], [0, 13, 23, 58, 90, 118], [0, 16, 37, 52, 109, 112], [0, 34, 39, 50, 73, 121]]
\item \{0=1260, 1=32256, 2=601776, 3=2964528, 4=3960180\} [[0, 1, 2, 5, 8, 14], [0, 3, 10, 67, 102, 112], [0, 4, 79, 85, 89, 107], [0, 6, 45, 46, 103, 113], [0, 12, 38, 82, 88, 124]]
\item \{0=1260, 1=32256, 2=601776, 3=2990736, 4=3933972\} [[0, 1, 2, 5, 8, 14], [0, 3, 7, 10, 34, 36], [0, 4, 47, 61, 97, 108], [0, 6, 21, 74, 85, 96], [0, 16, 54, 69, 94, 99], [0, 18, 76, 114, 115, 118], [0, 29, 67, 87, 91, 120]]
\item \{0=1260, 1=32256, 2=604800, 3=2997792, 4=3923892\} [[0, 1, 2, 5, 8, 14], [0, 3, 7, 39, 64, 111], [0, 4, 17, 63, 68, 70], [0, 9, 46, 67, 83, 94], [0, 10, 51, 59, 91, 110]]
\item \{0=1260, 1=32256, 2=612360, 3=3000816, 4=3913308\} [[0, 1, 2, 5, 8, 14], [0, 3, 4, 29, 72, 114], [0, 6, 23, 53, 79, 106], [0, 7, 57, 73, 95, 109], [0, 9, 26, 104, 112, 115]]
\item \{0=1260, 1=32256, 2=617652, 3=2987208, 4=3921624\} [[0, 1, 2, 5, 8, 14], [0, 3, 7, 46, 81, 99], [0, 4, 17, 23, 50, 96], [0, 9, 38, 41, 61, 119], [0, 10, 56, 69, 101, 125], [0, 13, 39, 40, 102, 114], [0, 24, 30, 104, 108, 112]]
\item \{0=1260, 1=32256, 2=624456, 3=2993760, 4=3908268\} [[0, 1, 2, 5, 8, 14], [0, 3, 16, 32, 91, 111], [0, 6, 55, 63, 69, 104], [0, 7, 44, 64, 103, 117], [0, 9, 36, 68, 87, 99], [0, 19, 46, 57, 101, 119], [0, 21, 39, 60, 77, 106]]
\item \{0=1260, 1=32256, 2=628992, 3=2953440, 4=3944052\} [[0, 1, 2, 5, 8, 14], [0, 3, 4, 63, 90, 92], [0, 6, 10, 32, 64, 99], [0, 12, 35, 61, 81, 122], [0, 16, 42, 58, 105, 108]]
\item \{0=1260, 1=32256, 2=633528, 3=3014928, 4=3878028\} [[0, 1, 2, 5, 8, 14], [0, 3, 4, 55, 60, 66], [0, 6, 16, 70, 96, 118], [0, 7, 10, 54, 62, 85], [0, 9, 86, 92, 99, 119]]
\item \{0=1260, 1=32256, 2=669060, 3=3007368, 4=3850056\} [[0, 1, 2, 5, 8, 14], [0, 3, 7, 49, 73, 115], [0, 4, 26, 38, 71, 79], [0, 13, 15, 57, 113, 124], [0, 20, 31, 93, 97, 121]]
\item \{0=1260, 1=33264, 2=588168, 3=3007872, 4=3929436\} [[0, 1, 2, 5, 8, 14], [0, 3, 10, 21, 78, 113], [0, 4, 29, 76, 89, 91], [0, 6, 25, 34, 84, 108], [0, 12, 35, 63, 82, 122]]
\item \{0=1260, 1=33264, 2=593460, 3=2961000, 4=3971016\} [[0, 1, 2, 5, 8, 14], [0, 3, 7, 57, 92, 115], [0, 4, 34, 47, 62, 116], [0, 6, 22, 64, 71, 87], [0, 10, 21, 25, 79, 98]]
\item \{0=1260, 1=33264, 2=600264, 3=2959488, 4=3965724\} [[0, 1, 2, 5, 8, 14], [0, 3, 4, 49, 52, 115], [0, 6, 13, 23, 33, 37], [0, 7, 21, 58, 69, 73], [0, 9, 19, 56, 93, 108]]
\item \{0=1260, 1=33264, 2=606312, 3=2954448, 4=3964716\} [[0, 1, 2, 5, 8, 14], [0, 3, 4, 70, 72, 120], [0, 6, 57, 64, 90, 94], [0, 7, 27, 34, 99, 108], [0, 13, 79, 109, 112, 114]]
\item \{0=1260, 1=33264, 2=608580, 3=2949912, 4=3966984\} [[0, 1, 2, 5, 8, 14], [0, 3, 4, 37, 104, 119], [0, 6, 48, 54, 68, 106], [0, 7, 36, 41, 85, 92], [0, 10, 50, 61, 69, 111]]
\item \{0=1260, 1=33264, 2=616896, 3=2964528, 4=3944052\} [[0, 1, 2, 5, 8, 14], [0, 3, 16, 45, 81, 117], [0, 4, 29, 40, 89, 104], [0, 6, 53, 98, 100, 125], [0, 21, 23, 50, 107, 115]]
\item \{0=1260, 1=33264, 2=616896, 3=3006864, 4=3901716\} [[0, 1, 2, 5, 8, 14], [0, 3, 4, 46, 114, 117], [0, 6, 12, 61, 81, 110], [0, 9, 37, 65, 69, 119], [0, 15, 38, 79, 94, 124]]
\item \{0=1260, 1=33264, 2=632016, 3=2933280, 4=3960180\} [[0, 1, 2, 5, 8, 14], [0, 3, 10, 25, 48, 60], [0, 4, 15, 19, 47, 107], [0, 9, 16, 51, 67, 104], [0, 18, 24, 46, 116, 125]]
\item \{0=1260, 1=33264, 2=661500, 3=2984184, 4=3879792\} [[0, 1, 3, 6, 11, 17], [0, 2, 4, 26, 41, 122], [0, 5, 37, 52, 114, 118], [0, 7, 63, 70, 105, 124], [0, 9, 85, 90, 97, 113]]
\item \{0=1260, 1=34272, 2=557928, 3=2988720, 4=3977820\} [[0, 1, 2, 5, 8, 14], [0, 3, 7, 99, 104, 122], [0, 4, 18, 39, 85, 115], [0, 10, 50, 87, 98, 113], [0, 16, 28, 29, 100, 102]]
\item \{0=1260, 1=34272, 2=560952, 3=2989728, 4=3973788\} [[0, 1, 2, 5, 8, 14], [0, 3, 16, 30, 42, 117], [0, 4, 22, 41, 95, 123], [0, 9, 40, 88, 98, 118], [0, 21, 39, 102, 104, 107]]
\item \{0=1260, 1=34272, 2=593460, 3=2983176, 4=3947832\} [[0, 1, 2, 5, 8, 14], [0, 3, 4, 66, 94, 123], [0, 6, 50, 53, 73, 111], [0, 7, 45, 64, 96, 103], [0, 10, 19, 61, 91, 95]]
\item \{0=1260, 1=34272, 2=600264, 3=3000816, 4=3923388\} [[0, 1, 2, 5, 8, 14], [0, 3, 7, 29, 108, 112], [0, 4, 39, 65, 93, 97], [0, 6, 10, 41, 76, 115], [0, 20, 23, 69, 85, 92]]
\item \{0=1260, 1=34272, 2=612360, 3=3018960, 4=3893148\} [[0, 1, 2, 5, 8, 14], [0, 3, 29, 33, 82, 122], [0, 4, 23, 57, 89, 103], [0, 6, 35, 50, 98, 116], [0, 13, 24, 45, 46, 88]]
\item \{0=1260, 1=34272, 2=626724, 3=3007368, 4=3890376\} [[0, 1, 2, 5, 8, 14], [0, 3, 4, 86, 91, 124], [0, 6, 43, 59, 79, 123], [0, 7, 23, 49, 90, 112], [0, 10, 52, 56, 85, 122], [0, 12, 28, 61, 76, 104], [0, 32, 55, 63, 75, 105]]
\item \{0=1260, 1=34272, 2=631260, 3=2993256, 4=3899952\} [[0, 1, 2, 5, 8, 14], [0, 3, 7, 10, 25, 109], [0, 4, 35, 90, 93, 114], [0, 9, 29, 76, 88, 98], [0, 13, 22, 70, 96, 122]]
\item \{0=1260, 1=34272, 2=649404, 3=2991240, 4=3883824\} [[0, 1, 2, 5, 8, 14], [0, 3, 4, 60, 101, 117], [0, 6, 38, 48, 54, 58], [0, 7, 39, 96, 109, 112], [0, 10, 28, 56, 70, 102], [0, 13, 20, 93, 107, 120], [0, 19, 31, 106, 110, 125]]
\item \{0=1260, 1=34272, 2=659988, 3=2994264, 4=3870216\} [[0, 1, 2, 5, 8, 14], [0, 3, 4, 54, 119, 125], [0, 7, 48, 56, 64, 104], [0, 9, 13, 41, 83, 88], [0, 12, 27, 70, 102, 121]]
\item \{0=1260, 1=34272, 2=668304, 3=2986704, 4=3869460\} [[0, 1, 2, 5, 8, 14], [0, 3, 4, 52, 59, 73], [0, 6, 12, 51, 84, 106], [0, 9, 42, 62, 102, 122], [0, 15, 44, 50, 90, 120]]
\item \{0=1260, 1=34272, 2=681156, 3=2987208, 4=3856104\} [[0, 1, 2, 5, 8, 14], [0, 3, 4, 59, 91, 125], [0, 6, 12, 27, 50, 63], [0, 9, 32, 43, 54, 85], [0, 13, 38, 52, 83, 118]]
\item \{0=1260, 1=35280, 2=555660, 3=2991240, 4=3976560\} [[0, 1, 2, 5, 8, 14], [0, 3, 4, 35, 65, 66], [0, 6, 33, 47, 95, 104], [0, 7, 21, 26, 109, 118], [0, 10, 41, 45, 61, 116]]
\item \{0=1260, 1=35280, 2=558684, 3=3045672, 4=3919104\} [[0, 1, 2, 5, 8, 14], [0, 3, 7, 30, 47, 99], [0, 4, 10, 17, 68, 84], [0, 9, 42, 66, 71, 94], [0, 16, 46, 102, 112, 113]]
\item \{0=1260, 1=35280, 2=561708, 3=2975112, 4=3986640\} [[0, 1, 3, 6, 11, 17], [0, 2, 4, 20, 34, 88], [0, 5, 46, 47, 96, 107], [0, 7, 10, 41, 76, 106], [0, 19, 60, 75, 94, 100]]
\item \{0=1260, 1=35280, 2=567756, 3=2992248, 4=3963456\} [[0, 1, 2, 5, 8, 14], [0, 3, 4, 63, 96, 108], [0, 6, 10, 37, 99, 101], [0, 7, 59, 64, 71, 87], [0, 13, 18, 90, 100, 120]]
\item \{0=1260, 1=35280, 2=573804, 3=2954952, 4=3994704\} [[0, 1, 2, 5, 8, 14], [0, 3, 4, 16, 40, 43], [0, 6, 55, 94, 113, 120], [0, 7, 22, 31, 88, 109], [0, 10, 42, 67, 73, 87]]
\item \{0=1260, 1=35280, 2=581364, 3=3026520, 4=3915576\} [[0, 1, 2, 5, 8, 14], [0, 3, 10, 49, 107, 112], [0, 4, 29, 76, 102, 113], [0, 6, 26, 27, 55, 104], [0, 12, 43, 45, 94, 106]]
\item \{0=1260, 1=35280, 2=588924, 3=2993256, 4=3941280\} [[0, 1, 2, 5, 8, 14], [0, 3, 16, 32, 63, 100], [0, 6, 7, 36, 44, 111], [0, 13, 22, 37, 73, 109], [0, 15, 40, 72, 86, 119], [0, 18, 21, 64, 77, 116], [0, 23, 50, 51, 90, 113]]
\item \{0=1260, 1=35280, 2=598752, 3=2967552, 4=3957156\} [[0, 1, 3, 6, 11, 17], [0, 2, 4, 43, 65, 72], [0, 5, 20, 98, 102, 122], [0, 10, 46, 86, 99, 124], [0, 12, 19, 37, 97, 103]]
\item \{0=1260, 1=35280, 2=609336, 3=3007872, 4=3906252\} [[0, 1, 2, 5, 8, 14], [0, 3, 10, 75, 78, 99], [0, 4, 21, 23, 52, 105], [0, 6, 26, 29, 66, 117], [0, 7, 35, 41, 57, 113]]
\item \{0=1260, 1=35280, 2=609336, 3=3009888, 4=3904236\} [[0, 1, 2, 5, 8, 14], [0, 3, 4, 91, 102, 108], [0, 6, 31, 51, 59, 60], [0, 9, 21, 53, 58, 82], [0, 12, 19, 43, 75, 121]]
\item \{0=1260, 1=35280, 2=616140, 3=2994264, 4=3913056\} [[0, 1, 2, 5, 8, 14], [0, 3, 16, 55, 63, 124], [0, 4, 59, 62, 104, 110], [0, 7, 30, 78, 97, 99], [0, 9, 39, 54, 112, 117]]
\item \{0=1260, 1=35280, 2=616140, 3=3052728, 4=3854592\} [[0, 1, 2, 5, 8, 14], [0, 3, 16, 55, 74, 115], [0, 4, 19, 39, 116, 125], [0, 6, 36, 37, 61, 96], [0, 10, 13, 20, 52, 105]]
\item \{0=1260, 1=35280, 2=635796, 3=2999304, 4=3888360\} [[0, 1, 3, 6, 11, 17], [0, 2, 4, 23, 67, 112], [0, 5, 40, 63, 83, 86], [0, 7, 16, 66, 68, 105], [0, 9, 38, 51, 82, 114]]
\item \{0=1260, 1=36288, 2=604800, 3=2991744, 4=3925908\} [[0, 1, 2, 5, 8, 14], [0, 3, 25, 29, 83, 104], [0, 4, 94, 98, 101, 106], [0, 6, 38, 60, 76, 82], [0, 13, 20, 64, 79, 92]]
\item \{0=1260, 1=36288, 2=614628, 3=3007368, 4=3900456\} [[0, 1, 2, 5, 8, 14], [0, 3, 7, 59, 81, 91], [0, 4, 18, 20, 43, 85], [0, 9, 35, 66, 87, 94], [0, 12, 21, 61, 90, 102]]
\item \{0=1260, 1=36288, 2=622944, 3=3027024, 4=3872484\} [[0, 1, 2, 5, 8, 14], [0, 3, 10, 67, 91, 101], [0, 4, 18, 27, 59, 110], [0, 7, 22, 23, 41, 87], [0, 15, 19, 54, 94, 98], [0, 32, 55, 63, 75, 105], [0, 35, 37, 48, 99, 109]]
\item \{0=1260, 1=36288, 2=623700, 3=3014424, 4=3884328\} [[0, 1, 2, 5, 8, 14], [0, 3, 10, 91, 101, 119], [0, 4, 77, 102, 103, 111], [0, 15, 51, 72, 94, 104], [0, 18, 32, 79, 87, 116]]
\item \{0=1260, 1=36288, 2=632772, 3=3024504, 4=3865176\} [[0, 1, 2, 5, 8, 14], [0, 3, 4, 10, 79, 99], [0, 6, 26, 38, 89, 124], [0, 9, 58, 69, 86, 108], [0, 15, 35, 63, 76, 118], [0, 16, 57, 105, 119, 121], [0, 18, 21, 64, 77, 116]]
\item \{0=1260, 1=36288, 2=635040, 3=2994768, 4=3892644\} [[0, 1, 2, 5, 8, 14], [0, 3, 16, 54, 81, 106], [0, 4, 60, 70, 96, 112], [0, 6, 56, 101, 102, 121], [0, 9, 38, 39, 73, 111], [0, 23, 46, 90, 103, 104], [0, 24, 35, 47, 86, 108]]
\item \{0=1260, 1=36288, 2=645624, 3=2951424, 4=3925404\} [[0, 1, 2, 5, 8, 14], [0, 3, 10, 21, 33, 119], [0, 4, 22, 29, 49, 84], [0, 6, 23, 52, 58, 59], [0, 7, 63, 91, 104, 122]]
\item \{0=1260, 1=37296, 2=601776, 3=3018960, 4=3900708\} [[0, 1, 2, 5, 8, 14], [0, 3, 7, 31, 78, 125], [0, 4, 38, 52, 90, 102], [0, 6, 13, 37, 71, 77], [0, 9, 42, 48, 108, 111]]
\item \{0=1260, 1=37296, 2=622944, 3=2968560, 4=3929940\} [[0, 1, 2, 5, 8, 14], [0, 3, 10, 48, 60, 110], [0, 4, 38, 111, 112, 125], [0, 7, 53, 102, 108, 120], [0, 9, 25, 71, 107, 123]]
\item \{0=1260, 1=37296, 2=635796, 3=2980152, 4=3905496\} [[0, 1, 2, 5, 8, 14], [0, 3, 4, 73, 85, 123], [0, 6, 47, 61, 69, 88], [0, 9, 29, 102, 117, 124], [0, 12, 31, 38, 77, 92]]
\item \{0=1260, 1=37296, 2=635796, 3=3020472, 4=3865176\} [[0, 1, 2, 5, 8, 14], [0, 3, 7, 77, 102, 106], [0, 4, 36, 38, 87, 119], [0, 10, 18, 76, 88, 95], [0, 19, 20, 81, 103, 124]]
\item \{0=1260, 1=37296, 2=641088, 3=3009888, 4=3870468\} [[0, 1, 2, 5, 8, 14], [0, 3, 29, 39, 47, 106], [0, 6, 31, 35, 71, 113], [0, 7, 16, 34, 72, 105], [0, 9, 25, 48, 68, 114]]
\item \{0=1260, 1=37296, 2=642600, 3=3022992, 4=3855852\} [[0, 1, 2, 5, 8, 14], [0, 3, 4, 114, 117, 119], [0, 7, 37, 64, 109, 123], [0, 9, 27, 54, 76, 118], [0, 10, 45, 56, 87, 115], [0, 12, 21, 35, 62, 113], [0, 26, 48, 63, 67, 100]]
\item \{0=1260, 1=37296, 2=651672, 3=2958480, 4=3911292\} [[0, 1, 3, 6, 11, 17], [0, 2, 4, 37, 101, 108], [0, 5, 20, 61, 105, 118], [0, 7, 53, 56, 97, 125], [0, 12, 43, 47, 67, 93]]
\item \{0=1260, 1=37296, 2=661500, 3=2980152, 4=3879792\} [[0, 1, 2, 5, 8, 14], [0, 3, 4, 45, 70, 103], [0, 6, 33, 48, 104, 118], [0, 7, 9, 46, 73, 124], [0, 10, 58, 94, 108, 121]]
\item \{0=1260, 1=38304, 2=563220, 3=3011400, 4=3945816\} [[0, 1, 3, 6, 11, 17], [0, 2, 4, 40, 47, 69], [0, 5, 59, 66, 86, 87], [0, 12, 16, 45, 61, 119], [0, 19, 58, 72, 102, 112], [0, 20, 22, 103, 105, 125], [0, 24, 54, 55, 75, 97]]
\item \{0=1260, 1=38304, 2=579096, 3=3028032, 4=3913308\} [[0, 1, 2, 5, 8, 14], [0, 3, 4, 52, 106, 119], [0, 6, 67, 89, 96, 99], [0, 10, 46, 60, 103, 105], [0, 13, 23, 65, 93, 114]]
\item \{0=1260, 1=38304, 2=604044, 3=3004344, 4=3912048\} [[0, 1, 2, 5, 8, 14], [0, 3, 7, 31, 49, 50], [0, 6, 44, 46, 86, 112], [0, 12, 59, 69, 95, 96], [0, 18, 38, 41, 82, 85], [0, 20, 22, 53, 99, 101], [0, 28, 51, 76, 103, 118]]
\item \{0=1260, 1=38304, 2=607068, 3=3021480, 4=3891888\} [[0, 1, 2, 5, 8, 14], [0, 3, 4, 10, 51, 76], [0, 6, 21, 41, 65, 123], [0, 7, 62, 75, 87, 107], [0, 9, 58, 109, 113, 125]]
\item \{0=1260, 1=38304, 2=637308, 3=2973096, 4=3910032\} [[0, 1, 2, 5, 8, 14], [0, 3, 4, 94, 98, 104], [0, 6, 12, 53, 59, 120], [0, 9, 36, 105, 114, 123], [0, 10, 29, 56, 71, 103], [0, 13, 30, 81, 96, 110], [0, 20, 39, 44, 73, 109]]
\item \{0=1260, 1=38304, 2=640332, 3=2981160, 4=3898944\} [[0, 1, 3, 6, 11, 17], [0, 2, 4, 8, 10, 21], [0, 5, 23, 93, 94, 125], [0, 7, 36, 38, 69, 119], [0, 13, 22, 46, 103, 108], [0, 15, 40, 71, 112, 124], [0, 35, 39, 73, 101, 122]]
\item \{0=1260, 1=38304, 2=654696, 3=3011904, 4=3853836\} [[0, 1, 2, 5, 8, 14], [0, 3, 7, 19, 38, 48], [0, 6, 65, 66, 87, 117], [0, 9, 55, 67, 95, 112], [0, 21, 49, 73, 115, 119]]
\item \{0=1260, 1=38304, 2=666792, 3=2941344, 4=3912300\} [[0, 1, 2, 5, 8, 14], [0, 3, 7, 10, 67, 68], [0, 4, 28, 49, 54, 62], [0, 6, 26, 40, 60, 95], [0, 16, 20, 88, 93, 118]]
\item \{0=1260, 1=38304, 2=697032, 3=2979648, 4=3843756\} [[0, 1, 2, 5, 8, 14], [0, 3, 21, 60, 84, 101], [0, 6, 82, 88, 112, 124], [0, 7, 41, 42, 53, 99], [0, 9, 50, 57, 92, 97]]
\item \{0=1260, 1=39312, 2=591192, 3=3005856, 4=3922380\} [[0, 1, 2, 5, 8, 14], [0, 3, 10, 50, 51, 96], [0, 4, 21, 22, 33, 70], [0, 7, 9, 38, 58, 120], [0, 12, 46, 67, 82, 93]]
\item \{0=1260, 1=39312, 2=596484, 3=2997288, 4=3925656\} [[0, 1, 2, 5, 8, 14], [0, 3, 4, 43, 48, 90], [0, 6, 12, 22, 32, 77], [0, 9, 68, 71, 91, 103], [0, 24, 51, 72, 84, 110]]
\item \{0=1260, 1=39312, 2=602532, 3=3017448, 4=3899448\} [[0, 1, 3, 6, 11, 17], [0, 2, 4, 16, 34, 75], [0, 5, 24, 35, 113, 115], [0, 9, 29, 79, 114, 121], [0, 23, 30, 76, 97, 106]]
\item \{0=1260, 1=39312, 2=611604, 3=3006360, 4=3901464\} [[0, 1, 2, 5, 8, 14], [0, 3, 4, 72, 86, 98], [0, 6, 16, 42, 64, 110], [0, 10, 33, 35, 62, 91], [0, 12, 37, 61, 81, 92]]
\item \{0=1260, 1=39312, 2=616896, 3=2927232, 4=3975300\} [[0, 1, 2, 5, 8, 14], [0, 3, 10, 45, 57, 120], [0, 4, 16, 17, 59, 118], [0, 7, 58, 87, 111, 112], [0, 13, 22, 78, 99, 105]]
\item \{0=1260, 1=39312, 2=616896, 3=3042144, 4=3860388\} [[0, 1, 2, 5, 8, 14], [0, 3, 16, 57, 78, 111], [0, 4, 31, 75, 99, 107], [0, 6, 19, 70, 101, 110], [0, 9, 54, 84, 117, 121]]
\item \{0=1260, 1=39312, 2=628236, 3=2981160, 4=3910032\} [[0, 1, 2, 5, 8, 14], [0, 3, 4, 21, 55, 114], [0, 6, 48, 53, 62, 85], [0, 9, 57, 99, 103, 108], [0, 10, 42, 64, 66, 113]]
\item \{0=1260, 1=39312, 2=647136, 3=2942352, 4=3929940\} [[0, 1, 2, 5, 8, 14], [0, 3, 21, 33, 46, 114], [0, 4, 60, 72, 103, 122], [0, 7, 44, 48, 51, 112], [0, 12, 47, 62, 93, 119]]
\item \{0=1260, 1=39312, 2=650160, 3=2975616, 4=3893652\} [[0, 1, 2, 5, 8, 14], [0, 3, 4, 43, 65, 124], [0, 6, 13, 35, 38, 82], [0, 9, 87, 95, 106, 113], [0, 19, 30, 46, 49, 108]]
\item \{0=1260, 1=40320, 2=641844, 3=3026520, 4=3850056\} [[0, 1, 2, 5, 8, 14], [0, 3, 16, 32, 36, 120], [0, 6, 7, 21, 88, 101], [0, 10, 19, 30, 58, 99], [0, 12, 37, 47, 79, 92]]
\item \{0=1260, 1=40320, 2=652428, 3=2985192, 4=3880800\} [[0, 1, 2, 5, 8, 14], [0, 3, 21, 47, 57, 96], [0, 4, 18, 27, 61, 116], [0, 9, 10, 70, 90, 114], [0, 16, 34, 91, 120, 125], [0, 19, 20, 22, 64, 66], [0, 24, 31, 105, 108, 113]]
\item \{0=1260, 1=40320, 2=652428, 3=3001320, 4=3864672\} [[0, 1, 2, 5, 8, 14], [0, 3, 4, 27, 39, 117], [0, 6, 45, 72, 85, 122], [0, 7, 54, 76, 95, 114], [0, 10, 26, 63, 86, 113]]
\item \{0=1260, 1=40320, 2=667548, 3=2987208, 4=3863664\} [[0, 1, 2, 5, 8, 14], [0, 3, 10, 49, 67, 121], [0, 6, 21, 36, 79, 107], [0, 7, 35, 66, 78, 98], [0, 15, 44, 59, 119, 123]]
\item \{0=1260, 1=41328, 2=604800, 3=2995776, 4=3916836\} [[0, 1, 2, 5, 8, 14], [0, 3, 16, 21, 60, 74], [0, 4, 66, 89, 111, 121], [0, 9, 18, 35, 101, 123], [0, 15, 24, 34, 75, 90]]
\item \{0=1260, 1=41328, 2=611604, 3=3023496, 4=3882312\} [[0, 1, 2, 5, 8, 14], [0, 3, 25, 45, 69, 100], [0, 4, 54, 72, 94, 110], [0, 6, 24, 31, 62, 103], [0, 7, 55, 75, 82, 101], [0, 9, 27, 51, 81, 114], [0, 37, 40, 78, 92, 109]]
\item \{0=1260, 1=41328, 2=627480, 3=3013920, 4=3876012\} [[0, 1, 2, 5, 8, 14], [0, 3, 4, 55, 81, 117], [0, 7, 40, 73, 100, 116], [0, 9, 13, 57, 78, 102], [0, 12, 16, 46, 61, 92]]
\item \{0=1260, 1=41328, 2=632016, 3=2946384, 4=3939012\} [[0, 1, 3, 6, 11, 17], [0, 2, 4, 28, 93, 108], [0, 5, 32, 60, 92, 104], [0, 7, 53, 96, 112, 113], [0, 10, 12, 75, 102, 121]]
\item \{0=1260, 1=41328, 2=649404, 3=2956968, 4=3911040\} [[0, 1, 2, 5, 8, 14], [0, 3, 7, 36, 81, 107], [0, 4, 23, 56, 85, 103], [0, 6, 40, 52, 117, 120], [0, 12, 30, 37, 77, 119]]
\item \{0=1260, 1=41328, 2=651672, 3=2967552, 4=3898188\} [[0, 1, 2, 5, 8, 14], [0, 3, 7, 60, 106, 114], [0, 4, 16, 17, 67, 86], [0, 9, 40, 71, 74, 87], [0, 10, 28, 56, 70, 102], [0, 21, 57, 58, 72, 97], [0, 35, 39, 73, 101, 122]]
\item \{0=1260, 1=41328, 2=678888, 3=3013920, 4=3824604\} [[0, 1, 2, 5, 8, 14], [0, 3, 7, 35, 110, 119], [0, 4, 22, 29, 79, 99], [0, 6, 44, 50, 66, 117], [0, 9, 21, 26, 73, 77], [0, 16, 34, 91, 120, 125], [0, 19, 58, 72, 102, 112]]
\item \{0=1260, 1=41328, 2=682668, 3=2996280, 4=3838464\} [[0, 1, 2, 5, 8, 14], [0, 3, 33, 51, 72, 121], [0, 6, 16, 79, 109, 123], [0, 9, 36, 50, 69, 82], [0, 10, 20, 56, 61, 94], [0, 12, 32, 62, 98, 125], [0, 23, 24, 97, 107, 108]]
\item \{0=1260, 1=42336, 2=582120, 3=3003840, 4=3930444\} [[0, 1, 2, 5, 8, 14], [0, 3, 4, 18, 56, 100], [0, 9, 16, 50, 72, 109], [0, 12, 61, 98, 110, 114], [0, 15, 57, 83, 104, 119], [0, 20, 22, 53, 99, 101], [0, 26, 39, 41, 73, 118]]
\item \{0=1260, 1=42336, 2=591192, 3=3065328, 4=3859884\} [[0, 1, 2, 5, 8, 14], [0, 3, 10, 42, 69, 84], [0, 4, 23, 45, 99, 109], [0, 6, 29, 33, 117, 124], [0, 9, 36, 53, 85, 112], [0, 19, 46, 57, 101, 119], [0, 21, 41, 44, 74, 95]]
\item \{0=1260, 1=42336, 2=600264, 3=2977632, 4=3938508\} [[0, 1, 2, 5, 8, 14], [0, 3, 7, 29, 49, 71], [0, 4, 16, 76, 81, 101], [0, 12, 50, 70, 112, 117], [0, 13, 84, 90, 91, 108]]
\item \{0=1260, 1=42336, 2=625212, 3=2965032, 4=3926160\} [[0, 1, 2, 5, 8, 14], [0, 3, 4, 41, 48, 57], [0, 6, 26, 61, 77, 116], [0, 9, 24, 51, 90, 122], [0, 15, 43, 59, 67, 121]]
\item \{0=1260, 1=42336, 2=629748, 3=3011400, 4=3875256\} [[0, 1, 2, 5, 8, 14], [0, 3, 4, 29, 39, 124], [0, 6, 36, 44, 46, 121], [0, 7, 75, 95, 109, 114], [0, 10, 20, 49, 101, 110]]
\item \{0=1260, 1=42336, 2=636552, 3=2939328, 4=3940524\} [[0, 1, 2, 5, 8, 14], [0, 3, 21, 48, 96, 114], [0, 4, 33, 38, 70, 124], [0, 7, 22, 62, 63, 111], [0, 9, 36, 103, 117, 122], [0, 10, 51, 56, 84, 121], [0, 27, 37, 39, 91, 109]]
\item \{0=1260, 1=42336, 2=654696, 3=3003840, 4=3857868\} [[0, 1, 2, 5, 8, 14], [0, 3, 16, 22, 54, 66], [0, 4, 55, 86, 87, 124], [0, 6, 46, 88, 98, 113], [0, 7, 34, 57, 67, 110]]
\item \{0=1260, 1=42336, 2=672084, 3=2978136, 4=3866184\} [[0, 1, 3, 6, 11, 17], [0, 2, 4, 19, 43, 72], [0, 5, 36, 48, 112, 123], [0, 10, 20, 66, 69, 85], [0, 15, 27, 65, 88, 101]]
\item \{0=1260, 1=43344, 2=585144, 3=3046176, 4=3884076\} [[0, 1, 2, 5, 8, 14], [0, 3, 16, 66, 104, 121], [0, 4, 17, 60, 96, 115], [0, 7, 63, 67, 92, 116], [0, 9, 38, 39, 73, 111], [0, 15, 33, 91, 100, 120], [0, 20, 22, 50, 52, 98]]
\item \{0=1260, 1=43344, 2=589680, 3=3040128, 4=3885588\} [[0, 1, 2, 5, 8, 14], [0, 3, 4, 19, 81, 113], [0, 6, 48, 70, 102, 111], [0, 10, 41, 83, 87, 107], [0, 16, 35, 112, 122, 123]]
\item \{0=1260, 1=43344, 2=630504, 3=2982672, 4=3902220\} [[0, 1, 2, 5, 8, 14], [0, 3, 7, 16, 25, 121], [0, 4, 23, 71, 72, 99], [0, 10, 49, 60, 108, 115], [0, 13, 20, 38, 70, 91]]
\item \{0=1260, 1=43344, 2=639576, 3=2976624, 4=3899196\} [[0, 1, 2, 5, 8, 14], [0, 3, 7, 25, 73, 123], [0, 4, 16, 61, 94, 121], [0, 9, 18, 35, 57, 62], [0, 20, 22, 32, 34, 81], [0, 21, 41, 44, 74, 95], [0, 23, 50, 51, 90, 113]]
\item \{0=1260, 1=43344, 2=643356, 3=2971080, 4=3900960\} [[0, 1, 2, 5, 8, 14], [0, 3, 4, 43, 100, 109], [0, 6, 48, 64, 91, 118], [0, 7, 28, 56, 88, 123], [0, 9, 16, 35, 40, 79]]
\item \{0=1260, 1=43344, 2=648648, 3=3014928, 4=3851820\} [[0, 1, 2, 5, 8, 14], [0, 3, 10, 74, 116, 124], [0, 4, 18, 78, 84, 86], [0, 7, 9, 76, 91, 100], [0, 12, 22, 66, 67, 109]]
\item \{0=1260, 1=43344, 2=654696, 3=2972592, 4=3888108\} [[0, 1, 2, 5, 8, 14], [0, 3, 7, 23, 35, 88], [0, 4, 55, 81, 97, 115], [0, 9, 52, 58, 91, 118], [0, 18, 22, 37, 41, 98]]
\item \{0=1260, 1=43344, 2=662256, 3=2978640, 4=3874500\} [[0, 1, 2, 5, 8, 14], [0, 3, 7, 56, 60, 76], [0, 4, 20, 70, 82, 111], [0, 6, 25, 50, 64, 91], [0, 13, 46, 67, 86, 113]]
\item \{0=1260, 1=44352, 2=624456, 3=3028032, 4=3861900\} [[0, 1, 2, 5, 8, 14], [0, 3, 7, 25, 63, 122], [0, 4, 35, 37, 49, 117], [0, 9, 36, 52, 65, 84], [0, 10, 27, 56, 59, 108], [0, 13, 20, 47, 101, 112], [0, 39, 43, 54, 74, 94]]
\item \{0=1260, 1=44352, 2=666036, 3=3017448, 4=3830904\} [[0, 1, 2, 5, 8, 14], [0, 3, 16, 22, 32, 93], [0, 6, 37, 43, 47, 79], [0, 7, 15, 106, 110, 121], [0, 12, 55, 59, 75, 90], [0, 21, 25, 49, 73, 123], [0, 24, 52, 64, 102, 108]]
\item \{0=1260, 1=45360, 2=575316, 3=2953944, 4=3984120\} [[0, 1, 2, 5, 8, 14], [0, 3, 4, 99, 105, 121], [0, 6, 76, 98, 116, 122], [0, 9, 38, 67, 82, 124], [0, 13, 37, 70, 90, 107]]
\item \{0=1260, 1=45360, 2=624456, 3=2984688, 4=3904236\} [[0, 1, 2, 5, 8, 14], [0, 3, 21, 42, 91, 98], [0, 4, 75, 76, 87, 101], [0, 9, 18, 58, 93, 121], [0, 10, 45, 81, 106, 123]]
\item \{0=1260, 1=45360, 2=632016, 3=2999808, 4=3881556\} [[0, 1, 2, 5, 8, 14], [0, 3, 4, 39, 103, 109], [0, 6, 30, 32, 47, 92], [0, 7, 31, 86, 87, 111], [0, 10, 46, 66, 95, 125]]
\item \{0=1260, 1=45360, 2=640332, 3=2999304, 4=3873744\} [[0, 1, 2, 5, 8, 14], [0, 3, 16, 69, 81, 125], [0, 4, 19, 30, 78, 106], [0, 6, 25, 29, 32, 111], [0, 10, 42, 73, 104, 120]]
\item \{0=1260, 1=45360, 2=650160, 3=3049200, 4=3814020\} [[0, 1, 3, 6, 11, 17], [0, 2, 7, 52, 68, 78], [0, 4, 27, 83, 101, 116], [0, 9, 16, 34, 102, 118], [0, 12, 55, 59, 75, 90], [0, 20, 22, 95, 97, 124], [0, 21, 30, 51, 77, 98]]
\item \{0=1260, 1=47376, 2=608580, 3=3015432, 4=3887352\} [[0, 1, 2, 5, 8, 14], [0, 3, 21, 36, 117, 123], [0, 4, 82, 100, 104, 120], [0, 6, 35, 49, 63, 88], [0, 7, 27, 28, 41, 78]]
\item \{0=1260, 1=47376, 2=627480, 3=2999808, 4=3884076\} [[0, 1, 2, 5, 8, 14], [0, 3, 4, 16, 19, 118], [0, 6, 22, 57, 77, 81], [0, 9, 44, 58, 72, 124], [0, 10, 30, 52, 68, 102], [0, 24, 31, 105, 108, 113], [0, 34, 39, 50, 73, 121]]
\item \{0=1260, 1=47376, 2=634284, 3=2971080, 4=3906000\} [[0, 1, 2, 5, 8, 14], [0, 3, 4, 21, 27, 84], [0, 6, 22, 53, 63, 111], [0, 7, 31, 57, 59, 107], [0, 9, 42, 66, 71, 90]]
\item \{0=1260, 1=48384, 2=622188, 3=3019464, 4=3868704\} [[0, 1, 2, 5, 8, 14], [0, 3, 4, 37, 67, 71], [0, 6, 51, 94, 123, 124], [0, 7, 10, 45, 61, 76], [0, 12, 30, 59, 79, 95]]
\item \{0=1260, 1=48384, 2=639576, 3=2984688, 4=3886092\} [[0, 1, 2, 5, 8, 14], [0, 3, 10, 63, 67, 95], [0, 4, 36, 42, 103, 107], [0, 6, 35, 44, 94, 113], [0, 12, 51, 77, 105, 119]]
\item \{0=1260, 1=48384, 2=641088, 3=2996784, 4=3872484\} [[0, 1, 2, 5, 8, 14], [0, 3, 16, 45, 67, 100], [0, 4, 61, 73, 88, 106], [0, 6, 25, 29, 74, 114], [0, 10, 59, 104, 107, 115]]
\item \{0=1260, 1=48384, 2=656964, 3=2987208, 4=3866184\} [[0, 1, 2, 5, 8, 14], [0, 3, 4, 56, 88, 117], [0, 6, 38, 43, 85, 125], [0, 7, 21, 40, 73, 120], [0, 9, 42, 64, 76, 99]]
\item \{0=1260, 1=48384, 2=710640, 3=2983680, 4=3816036\} [[0, 1, 2, 5, 8, 14], [0, 3, 4, 57, 107, 117], [0, 6, 31, 43, 84, 123], [0, 7, 29, 60, 65, 115], [0, 12, 30, 38, 52, 122]]
\item \{0=1260, 1=49392, 2=647892, 3=3001320, 4=3860136\} [[0, 1, 2, 5, 8, 14], [0, 3, 10, 36, 83, 124], [0, 6, 42, 46, 55, 99], [0, 12, 38, 60, 100, 125], [0, 13, 18, 48, 90, 108]]
\item \{0=1260, 1=49392, 2=648648, 3=3016944, 4=3843756\} [[0, 1, 2, 5, 8, 14], [0, 3, 16, 66, 87, 103], [0, 4, 21, 71, 91, 105], [0, 6, 39, 41, 107, 109], [0, 7, 40, 45, 67, 124], [0, 18, 76, 114, 115, 118], [0, 24, 30, 104, 108, 112]]
\item \{0=1260, 1=50400, 2=605556, 3=3022488, 4=3880296\} [[0, 1, 2, 5, 8, 14], [0, 3, 25, 29, 47, 72], [0, 7, 45, 93, 98, 104], [0, 9, 38, 60, 100, 106], [0, 10, 69, 85, 96, 116]]
\item \{0=1260, 1=50400, 2=611604, 3=2977128, 4=3919608\} [[0, 1, 2, 5, 8, 14], [0, 3, 4, 47, 57, 94], [0, 6, 39, 59, 62, 112], [0, 7, 29, 42, 49, 71], [0, 12, 45, 51, 53, 75]]
\item \{0=1260, 1=51408, 2=616140, 3=2965032, 4=3926160\} [[0, 1, 2, 5, 8, 14], [0, 3, 4, 82, 106, 114], [0, 6, 26, 38, 52, 59], [0, 12, 23, 43, 100, 124], [0, 21, 31, 51, 79, 83]]
\item \{0=1260, 1=51408, 2=646380, 3=2982168, 4=3878784\} [[0, 1, 2, 5, 8, 14], [0, 3, 7, 38, 57, 73], [0, 4, 36, 55, 88, 102], [0, 6, 35, 62, 86, 106], [0, 15, 26, 27, 59, 84]]
\item \{0=1260, 1=53424, 2=622944, 3=3048192, 4=3834180\} [[0, 1, 2, 5, 8, 14], [0, 3, 4, 41, 92, 104], [0, 6, 26, 65, 68, 125], [0, 7, 9, 39, 44, 69], [0, 13, 40, 50, 64, 121]]
\item \{0=1260, 1=53424, 2=625212, 3=3070872, 4=3809232\} [[0, 1, 2, 5, 8, 14], [0, 3, 4, 19, 66, 79], [0, 6, 40, 47, 82, 101], [0, 9, 44, 77, 85, 99], [0, 15, 29, 69, 93, 107], [0, 16, 18, 20, 22, 63], [0, 24, 43, 55, 94, 108]]
\item \{0=1260, 1=54432, 2=580608, 3=2964528, 4=3959172\} [[0, 1, 2, 5, 8, 14], [0, 3, 7, 23, 87, 99], [0, 4, 39, 45, 71, 95], [0, 6, 62, 63, 64, 121], [0, 9, 55, 88, 106, 111]]
\item \{0=2520, 1=23184, 2=613116, 3=2994264, 4=3926916\} [[0, 1, 3, 6, 11, 17], [0, 2, 4, 8, 10, 21], [0, 5, 24, 50, 52, 115], [0, 9, 38, 63, 78, 99], [0, 12, 13, 48, 55, 83], [0, 23, 35, 85, 90, 125], [0, 29, 58, 66, 72, 105]]
\item \{0=2520, 1=24192, 2=601776, 3=3004848, 4=3926664\} [[0, 1, 2, 5, 8, 14], [0, 3, 10, 68, 75, 103], [0, 4, 33, 43, 49, 92], [0, 7, 54, 81, 94, 116], [0, 12, 32, 37, 44, 98], [0, 23, 24, 97, 107, 108], [0, 29, 70, 84, 106, 110]]
\item \{0=2520, 1=25200, 2=557172, 3=3001320, 4=3973788\} [[0, 1, 2, 5, 8, 14], [0, 3, 10, 46, 54, 119], [0, 4, 23, 38, 82, 125], [0, 7, 9, 40, 88, 123], [0, 12, 27, 32, 77, 84]]
\item \{0=2520, 1=26208, 2=591948, 3=3023496, 4=3915828\} [[0, 1, 2, 5, 8, 14], [0, 3, 10, 48, 54, 119], [0, 4, 38, 69, 78, 117], [0, 9, 12, 72, 100, 125], [0, 13, 15, 32, 65, 106]]
\item \{0=2520, 1=28224, 2=602532, 3=3027528, 4=3899196\} [[0, 1, 2, 5, 8, 14], [0, 3, 25, 84, 87, 100], [0, 4, 6, 29, 31, 114], [0, 7, 45, 88, 99, 118], [0, 9, 38, 58, 77, 104]]
\item \{0=2520, 1=29232, 2=601020, 3=3015432, 4=3911796\} [[0, 1, 2, 5, 8, 14], [0, 3, 4, 23, 65, 104], [0, 6, 26, 57, 63, 84], [0, 9, 39, 70, 74, 76], [0, 13, 33, 38, 78, 123]]
\item \{0=2520, 1=29232, 2=624456, 3=2996784, 4=3907008\} [[0, 1, 2, 5, 8, 14], [0, 3, 4, 23, 46, 96], [0, 6, 16, 55, 100, 124], [0, 7, 9, 35, 48, 125], [0, 12, 24, 52, 85, 95]]
\item \{0=2520, 1=31248, 2=601776, 3=3042144, 4=3882312\} [[0, 1, 2, 5, 8, 14], [0, 3, 4, 30, 65, 76], [0, 6, 53, 93, 106, 111], [0, 7, 27, 68, 79, 112], [0, 10, 54, 61, 90, 104]]
\item \{0=2520, 1=31248, 2=656964, 3=2997288, 4=3871980\} [[0, 1, 2, 5, 8, 14], [0, 3, 4, 39, 88, 118], [0, 6, 44, 50, 92, 106], [0, 7, 20, 42, 62, 68], [0, 9, 18, 77, 98, 104]]
\item \{0=2520, 1=32256, 2=585900, 3=3011400, 4=3927924\} [[0, 1, 2, 5, 8, 14], [0, 3, 10, 49, 64, 103], [0, 4, 6, 37, 52, 59], [0, 7, 31, 66, 67, 69], [0, 26, 48, 63, 88, 100]]
\item \{0=2520, 1=32256, 2=600264, 3=2995776, 4=3929184\} [[0, 1, 2, 5, 8, 14], [0, 3, 4, 41, 49, 95], [0, 6, 30, 34, 89, 118], [0, 16, 69, 73, 96, 108], [0, 18, 35, 37, 84, 94]]
\item \{0=2520, 1=32256, 2=635040, 3=2991744, 4=3898440\} [[0, 1, 2, 5, 8, 14], [0, 3, 4, 29, 88, 112], [0, 6, 38, 40, 54, 59], [0, 7, 31, 53, 77, 84], [0, 13, 28, 33, 75, 95]]
\item \{0=2520, 1=32256, 2=639576, 3=3006864, 4=3878784\} [[0, 1, 2, 5, 8, 14], [0, 3, 10, 25, 78, 101], [0, 4, 28, 31, 37, 51], [0, 6, 19, 65, 72, 100], [0, 20, 26, 57, 81, 123]]
\item \{0=2520, 1=33264, 2=600264, 3=3008880, 4=3915072\} [[0, 1, 2, 5, 8, 14], [0, 3, 7, 37, 60, 66], [0, 4, 30, 54, 56, 76], [0, 6, 39, 84, 85, 120], [0, 12, 15, 79, 82, 95]]
\item \{0=2520, 1=34272, 2=566244, 3=3007368, 4=3949596\} [[0, 1, 2, 5, 8, 14], [0, 3, 4, 16, 57, 64], [0, 6, 31, 43, 67, 93], [0, 7, 27, 37, 102, 111], [0, 15, 51, 66, 70, 125]]
\item \{0=2520, 1=34272, 2=620676, 3=3023496, 4=3879036\} [[0, 1, 2, 5, 8, 14], [0, 3, 4, 95, 121, 124], [0, 6, 22, 24, 53, 77], [0, 9, 18, 81, 84, 94], [0, 10, 37, 83, 98, 118]]
\item \{0=2520, 1=34272, 2=653184, 3=3002832, 4=3867192\} [[0, 1, 2, 5, 8, 14], [0, 3, 7, 16, 78, 123], [0, 4, 58, 98, 108, 110], [0, 6, 13, 40, 46, 75], [0, 9, 10, 34, 79, 103]]
\item \{0=2520, 1=35280, 2=573804, 3=3020472, 4=3927924\} [[0, 1, 2, 5, 8, 14], [0, 3, 4, 21, 96, 98], [0, 6, 37, 45, 67, 121], [0, 7, 27, 50, 61, 118], [0, 16, 40, 113, 117, 124], [0, 24, 30, 104, 108, 112], [0, 35, 57, 64, 115, 119]]
\item \{0=2520, 1=36288, 2=589680, 3=2984688, 4=3946824\} [[0, 1, 2, 5, 8, 14], [0, 3, 16, 37, 73, 113], [0, 4, 18, 85, 86, 111], [0, 7, 22, 79, 87, 91], [0, 9, 27, 38, 66, 124]]
\item \{0=2520, 1=37296, 2=557172, 3=2995272, 4=3967740\} [[0, 1, 2, 5, 8, 14], [0, 3, 29, 30, 112, 120], [0, 4, 22, 42, 47, 73], [0, 6, 34, 69, 96, 101], [0, 10, 27, 66, 98, 124]]
\item \{0=2520, 1=37296, 2=584388, 3=2998296, 4=3937500\} [[0, 1, 2, 5, 8, 14], [0, 3, 21, 25, 46, 108], [0, 4, 48, 93, 95, 114], [0, 7, 58, 66, 69, 110], [0, 9, 37, 49, 60, 103]]
\item \{0=2520, 1=37296, 2=661500, 3=3030552, 4=3828132\} [[0, 1, 3, 6, 11, 17], [0, 2, 4, 16, 34, 99], [0, 5, 47, 63, 86, 122], [0, 9, 27, 49, 78, 105], [0, 15, 41, 54, 79, 112]]
\item \{0=2520, 1=38304, 2=624456, 3=2972592, 4=3922128\} [[0, 1, 2, 5, 8, 14], [0, 3, 7, 30, 46, 96], [0, 4, 10, 31, 94, 108], [0, 9, 40, 70, 111, 114], [0, 18, 19, 22, 47, 88]]
\item \{0=2520, 1=39312, 2=624456, 3=2994768, 4=3898944\} [[0, 1, 2, 5, 8, 14], [0, 3, 10, 21, 33, 93], [0, 4, 47, 67, 87, 111], [0, 6, 48, 84, 103, 108], [0, 7, 53, 63, 86, 97]]
\item \{0=2520, 1=39312, 2=633528, 3=2984688, 4=3899952\} [[0, 1, 2, 5, 8, 14], [0, 3, 7, 65, 91, 118], [0, 4, 47, 70, 71, 109], [0, 6, 62, 72, 94, 119], [0, 15, 60, 67, 99, 115]]
\item \{0=2520, 1=39312, 2=654696, 3=2942352, 4=3921120\} [[0, 1, 2, 5, 8, 14], [0, 3, 4, 49, 71, 123], [0, 6, 39, 63, 66, 121], [0, 9, 26, 65, 95, 118], [0, 13, 18, 24, 78, 115]]
\item \{0=2520, 1=40320, 2=606312, 3=2945376, 4=3965472\} [[0, 1, 2, 5, 8, 14], [0, 3, 16, 45, 87, 117], [0, 4, 39, 41, 76, 89], [0, 6, 32, 54, 55, 107], [0, 20, 51, 81, 99, 110]]
\item \{0=2520, 1=41328, 2=575316, 3=2980152, 4=3960684\} [[0, 1, 2, 5, 8, 14], [0, 3, 4, 26, 50, 117], [0, 6, 24, 56, 101, 119], [0, 7, 30, 53, 86, 106], [0, 9, 21, 60, 93, 107], [0, 19, 58, 72, 102, 112], [0, 20, 22, 83, 85, 87]]
\item \{0=2520, 1=41328, 2=653940, 3=3031560, 4=3830652\} [[0, 1, 2, 5, 8, 14], [0, 3, 16, 45, 67, 125], [0, 6, 21, 64, 92, 106], [0, 7, 42, 55, 95, 97], [0, 10, 35, 71, 72, 115]]
\item \{0=3780, 1=34272, 2=606312, 3=3033072, 4=3882564\} [[0, 1, 2, 5, 8, 14], [0, 3, 4, 21, 117, 122], [0, 6, 26, 52, 81, 93], [0, 10, 68, 82, 108, 111], [0, 13, 37, 50, 79, 118]]
\item \{0=3780, 1=42336, 2=611604, 3=2981160, 4=3921120\} [[0, 1, 2, 5, 8, 14], [0, 3, 7, 23, 40, 103], [0, 6, 13, 38, 54, 119], [0, 16, 42, 43, 98, 118], [0, 18, 44, 63, 86, 106]]
\end{enumerate}
\end{example}

\begin{example} {\tt SmallGroup(126,12) = C21 x S3}
\begin{enumerate}
\item \{1=20160, 2=536004, 3=2952936, 4=4050900\} [[0, 1, 4, 7, 13, 19], [0, 2, 6, 46, 76, 84], [0, 3, 23, 52, 105, 106], [0, 5, 36, 39, 68, 80], [0, 11, 12, 29, 77, 120]]
\item \{1=21168, 2=594216, 3=3030048, 4=3914568\} [[0, 1, 4, 7, 13, 19], [0, 2, 12, 26, 73, 112], [0, 5, 69, 81, 92, 107], [0, 6, 42, 43, 50, 74], [0, 10, 39, 71, 105, 124]]
\item \{1=22176, 2=594216, 3=2990736, 4=3952872\} [[0, 1, 4, 7, 13, 19], [0, 2, 3, 23, 75, 98], [0, 6, 22, 26, 92, 122], [0, 9, 53, 59, 61, 102], [0, 25, 32, 49, 56, 89]]
\item \{1=23184, 2=605556, 3=3076920, 4=3854340\} [[0, 1, 4, 7, 13, 19], [0, 2, 3, 28, 55, 69], [0, 6, 56, 58, 91, 92], [0, 12, 22, 114, 123, 124], [0, 21, 43, 64, 81, 121]]
\item \{1=25200, 2=617652, 3=2992248, 4=3924900\} [[0, 1, 4, 7, 13, 19], [0, 2, 6, 20, 21, 43], [0, 3, 32, 55, 66, 94], [0, 11, 23, 85, 116, 118], [0, 16, 42, 44, 100, 105]]
\item \{1=27216, 2=605556, 3=2943864, 4=3983364\} [[0, 1, 4, 7, 13, 19], [0, 2, 6, 43, 54, 125], [0, 3, 60, 77, 98, 99], [0, 9, 56, 86, 94, 120], [0, 16, 18, 23, 36, 116]]
\item \{1=27216, 2=611604, 3=3025512, 4=3895668\} [[0, 1, 4, 7, 13, 19], [0, 2, 3, 61, 63, 116], [0, 5, 47, 56, 73, 97], [0, 6, 10, 69, 82, 106], [0, 11, 27, 75, 117, 124]]
\item \{1=28224, 2=573048, 3=2996784, 4=3961944\} [[0, 1, 4, 7, 13, 19], [0, 2, 6, 33, 55, 85], [0, 3, 41, 48, 74, 75], [0, 5, 29, 50, 97, 100], [0, 12, 25, 34, 47, 72]]
\item \{1=29232, 2=557172, 3=3001320, 4=3972276\} [[0, 1, 4, 7, 13, 19], [0, 2, 6, 23, 64, 66], [0, 3, 25, 36, 65, 96], [0, 5, 62, 73, 105, 109], [0, 17, 46, 58, 107, 110]]
\item \{1=29232, 2=632016, 3=3005856, 4=3892896\} [[0, 1, 4, 7, 13, 19], [0, 2, 3, 43, 84, 92], [0, 5, 9, 34, 35, 60], [0, 10, 18, 22, 58, 125], [0, 15, 37, 51, 82, 102]]
\item \{1=30240, 2=550368, 3=2981664, 4=3997728\} [[0, 1, 4, 7, 13, 19], [0, 2, 6, 34, 82, 124], [0, 3, 15, 54, 56, 119], [0, 9, 53, 59, 95, 116], [0, 17, 46, 66, 77, 114]]
\item \{1=30240, 2=558684, 3=2989224, 4=3981852\} [[0, 1, 4, 7, 13, 19], [0, 2, 6, 96, 114, 122], [0, 5, 17, 55, 67, 84], [0, 9, 24, 62, 66, 73], [0, 18, 26, 48, 64, 70]]
\item \{1=30240, 2=615384, 3=3027024, 4=3887352\} [[0, 1, 4, 7, 13, 19], [0, 2, 6, 31, 63, 96], [0, 3, 47, 64, 69, 103], [0, 5, 11, 43, 49, 76], [0, 15, 24, 54, 95, 124]]
\item \{1=30240, 2=619920, 3=3014928, 4=3894912\} [[0, 1, 4, 7, 13, 19], [0, 2, 6, 23, 84, 113], [0, 3, 68, 85, 98, 123], [0, 5, 47, 71, 105, 112], [0, 9, 34, 43, 66, 94]]
\item \{1=30240, 2=628992, 3=2987712, 4=3913056\} [[0, 1, 4, 7, 13, 19], [0, 2, 6, 20, 58, 69], [0, 3, 21, 77, 111, 123], [0, 10, 31, 33, 79, 99], [0, 11, 62, 84, 89, 119]]
\item \{1=30240, 2=666036, 3=2986200, 4=3877524\} [[0, 1, 4, 7, 13, 19], [0, 2, 3, 30, 64, 111], [0, 5, 39, 56, 95, 120], [0, 6, 43, 68, 92, 107], [0, 9, 18, 70, 91, 123]]
\item \{1=31248, 2=628992, 3=3014928, 4=3884832\} [[0, 1, 4, 7, 13, 19], [0, 2, 6, 76, 111, 122], [0, 3, 16, 47, 57, 97], [0, 5, 45, 66, 77, 94], [0, 9, 11, 59, 69, 73]]
\item \{1=32256, 2=616140, 3=2977128, 4=3934476\} [[0, 1, 4, 7, 13, 19], [0, 2, 6, 30, 31, 82], [0, 3, 43, 59, 79, 115], [0, 5, 42, 63, 74, 77], [0, 12, 16, 52, 75, 100]]
\item \{1=32256, 2=619164, 3=2943864, 4=3964716\} [[0, 1, 4, 7, 13, 19], [0, 2, 3, 26, 48, 82], [0, 5, 21, 52, 68, 87], [0, 11, 32, 89, 94, 101], [0, 22, 42, 50, 92, 120]]
\item \{1=32256, 2=621432, 3=3002832, 4=3903480\} [[0, 1, 4, 7, 13, 19], [0, 2, 3, 37, 87, 94], [0, 5, 11, 49, 86, 105], [0, 9, 32, 51, 85, 90], [0, 18, 48, 88, 103, 118]]
\item \{1=32256, 2=664524, 3=3037608, 4=3825612\} [[0, 1, 4, 7, 13, 19], [0, 2, 6, 40, 94, 106], [0, 3, 15, 65, 75, 113], [0, 10, 39, 46, 55, 116], [0, 11, 12, 70, 122, 124]]
\item \{1=33264, 2=624456, 3=2948400, 4=3953880\} [[0, 1, 4, 7, 13, 19], [0, 2, 6, 26, 63, 64], [0, 3, 33, 59, 79, 90], [0, 5, 30, 73, 77, 80], [0, 9, 40, 42, 85, 100]]
\item \{1=33264, 2=637308, 3=3001320, 4=3888108\} [[0, 1, 4, 7, 13, 19], [0, 2, 6, 20, 104, 122], [0, 3, 26, 63, 95, 98], [0, 10, 55, 56, 99, 107], [0, 15, 34, 69, 94, 113]]
\item \{1=34272, 2=588168, 3=2987712, 4=3949848\} [[0, 1, 4, 7, 13, 19], [0, 2, 6, 20, 81, 94], [0, 3, 43, 88, 116, 120], [0, 11, 21, 27, 68, 75], [0, 16, 61, 67, 90, 117]]
\item \{1=34272, 2=651672, 3=2978640, 4=3895416\} [[0, 1, 4, 7, 13, 19], [0, 2, 6, 20, 82, 114], [0, 3, 47, 49, 81, 83], [0, 15, 24, 37, 54, 122], [0, 17, 23, 43, 72, 96]]
\item \{1=35280, 2=568512, 3=2999808, 4=3956400\} [[0, 1, 4, 7, 13, 19], [0, 2, 6, 61, 100, 109], [0, 3, 33, 46, 50, 75], [0, 5, 41, 44, 97, 123], [0, 9, 38, 49, 76, 119]]
\item \{1=35280, 2=628992, 3=2999808, 4=3895920\} [[0, 1, 4, 7, 13, 19], [0, 2, 3, 54, 61, 112], [0, 5, 41, 69, 89, 121], [0, 6, 38, 53, 92, 107], [0, 9, 18, 65, 100, 119]]
\item \{1=36288, 2=626724, 3=2980152, 4=3916836\} [[0, 1, 4, 7, 13, 19], [0, 2, 6, 26, 61, 73], [0, 3, 24, 30, 44, 81], [0, 5, 57, 80, 94, 121], [0, 16, 23, 69, 76, 82]]
\item \{1=36288, 2=678888, 3=2996784, 4=3848040\} [[0, 1, 4, 7, 13, 19], [0, 2, 6, 30, 58, 79], [0, 3, 38, 41, 60, 83], [0, 5, 51, 54, 95, 105], [0, 12, 62, 66, 100, 110]]
\item \{1=38304, 2=584388, 3=3031560, 4=3905748\} [[0, 1, 4, 7, 13, 19], [0, 2, 3, 61, 101, 124], [0, 5, 27, 30, 53, 75], [0, 9, 29, 58, 85, 90], [0, 12, 24, 25, 66, 122]]
\item \{1=38304, 2=607824, 3=2981664, 4=3932208\} [[0, 1, 4, 7, 13, 19], [0, 2, 6, 63, 64, 121], [0, 3, 21, 26, 52, 58], [0, 5, 74, 77, 84, 99], [0, 10, 39, 44, 83, 104]]
\item \{1=38304, 2=638064, 3=3011904, 4=3871728\} [[0, 1, 4, 7, 13, 19], [0, 2, 6, 23, 79, 90], [0, 3, 67, 92, 122, 123], [0, 9, 33, 57, 89, 97], [0, 11, 31, 55, 66, 77]]
\item \{1=39312, 2=584388, 3=3010392, 4=3925908\} [[0, 1, 4, 7, 13, 19], [0, 2, 3, 43, 45, 82], [0, 5, 35, 42, 68, 94], [0, 6, 24, 37, 62, 110], [0, 10, 36, 46, 70, 92]]
\item \{1=39312, 2=596484, 3=2998296, 4=3925908\} [[0, 1, 4, 7, 13, 19], [0, 2, 3, 28, 69, 118], [0, 6, 21, 73, 108, 124], [0, 9, 18, 79, 83, 85], [0, 10, 26, 49, 90, 114]]
\item \{1=40320, 2=576072, 3=3008880, 4=3934728\} [[0, 1, 4, 7, 13, 19], [0, 2, 3, 40, 84, 107], [0, 5, 35, 36, 102, 106], [0, 6, 33, 50, 89, 120], [0, 10, 100, 108, 111, 112]]
\item \{1=40320, 2=579096, 3=3008880, 4=3931704\} [[0, 1, 4, 7, 13, 19], [0, 2, 6, 17, 91, 104], [0, 5, 50, 66, 83, 114], [0, 10, 20, 26, 51, 115], [0, 12, 16, 52, 58, 74]]
\item \{1=40320, 2=588924, 3=3055752, 4=3875004\} [[0, 1, 4, 7, 13, 19], [0, 2, 3, 21, 69, 108], [0, 5, 37, 43, 79, 83], [0, 9, 29, 55, 88, 117], [0, 18, 70, 71, 91, 123]]
\item \{1=40320, 2=623700, 3=3001320, 4=3894660\} [[0, 1, 4, 7, 13, 19], [0, 2, 6, 64, 98, 104], [0, 3, 26, 68, 74, 119], [0, 9, 59, 66, 69, 80], [0, 10, 34, 65, 70, 90]]
\item \{1=40320, 2=638064, 3=2990736, 4=3890880\} [[0, 1, 4, 7, 13, 19], [0, 2, 3, 30, 37, 102], [0, 5, 11, 50, 63, 110], [0, 6, 31, 51, 108, 117], [0, 9, 18, 61, 90, 118]]
\item \{1=40320, 2=703836, 3=3001320, 4=3814524\} [[0, 1, 4, 7, 13, 19], [0, 2, 3, 23, 87, 94], [0, 5, 33, 34, 73, 116], [0, 9, 35, 64, 70, 122], [0, 11, 86, 98, 105, 118]]
\item \{1=41328, 2=602532, 3=2986200, 4=3929940\} [[0, 1, 4, 7, 13, 19], [0, 2, 3, 26, 77, 112], [0, 5, 38, 58, 84, 96], [0, 6, 25, 73, 90, 125], [0, 9, 53, 74, 104, 116]]
\item \{1=41328, 2=628992, 3=3051216, 4=3838464\} [[0, 1, 4, 7, 13, 19], [0, 2, 6, 33, 49, 116], [0, 3, 85, 94, 111, 119], [0, 5, 21, 54, 69, 70], [0, 18, 55, 79, 115, 118]]
\item \{1=41328, 2=638064, 3=3057264, 4=3823344\} [[0, 1, 4, 7, 13, 19], [0, 2, 3, 43, 74, 98], [0, 5, 23, 86, 96, 120], [0, 6, 24, 77, 108, 125], [0, 11, 15, 32, 75, 113]]
\item \{1=42336, 2=616140, 3=3007368, 4=3894156\} [[0, 1, 4, 7, 13, 19], [0, 2, 3, 61, 87, 118], [0, 5, 25, 26, 47, 106], [0, 10, 49, 73, 85, 90], [0, 12, 51, 74, 100, 115]]
\item \{1=42336, 2=619920, 3=2996784, 4=3900960\} [[0, 1, 4, 7, 13, 19], [0, 2, 12, 55, 91, 112], [0, 3, 33, 56, 66, 80], [0, 5, 23, 73, 102, 118], [0, 6, 21, 53, 74, 85]]
\item \{1=42336, 2=621432, 3=3002832, 4=3893400\} [[0, 1, 4, 7, 13, 19], [0, 2, 6, 31, 43, 91], [0, 3, 15, 39, 47, 115], [0, 9, 25, 40, 98, 117], [0, 11, 29, 55, 96, 108]]
\item \{1=42336, 2=666036, 3=3058776, 4=3792852\} [[0, 1, 4, 7, 13, 19], [0, 2, 6, 37, 45, 119], [0, 3, 15, 81, 85, 109], [0, 9, 27, 90, 92, 125], [0, 11, 31, 68, 69, 123]]
\item \{1=43344, 2=571536, 3=3020976, 4=3924144\} [[0, 1, 4, 7, 13, 19], [0, 2, 6, 34, 54, 95], [0, 3, 49, 99, 118, 123], [0, 9, 68, 75, 79, 120], [0, 10, 29, 92, 104, 119]]
\item \{1=43344, 2=601776, 3=2999808, 4=3915072\} [[0, 1, 4, 7, 13, 19], [0, 2, 3, 33, 58, 105], [0, 5, 36, 48, 65, 85], [0, 6, 32, 59, 95, 100], [0, 9, 21, 27, 49, 90]]
\item \{1=44352, 2=646380, 3=3049704, 4=3819564\} [[0, 1, 4, 7, 13, 19], [0, 2, 3, 48, 104, 118], [0, 5, 17, 85, 92, 111], [0, 6, 46, 53, 58, 117], [0, 9, 38, 77, 100, 110]]
\item \{1=45360, 2=641844, 3=2998296, 4=3874500\} [[0, 1, 4, 7, 13, 19], [0, 2, 3, 26, 92, 112], [0, 5, 29, 45, 56, 80], [0, 9, 42, 52, 73, 113], [0, 11, 51, 59, 96, 100]]
\item \{1=46368, 2=600264, 3=3054240, 4=3859128\} [[0, 1, 4, 7, 13, 19], [0, 2, 3, 15, 67, 120], [0, 9, 27, 48, 64, 95], [0, 11, 77, 83, 90, 94], [0, 16, 43, 65, 75, 85]]
\item \{1=46368, 2=645624, 3=2996784, 4=3871224\} [[0, 1, 4, 7, 13, 19], [0, 2, 6, 23, 86, 118], [0, 3, 30, 65, 97, 108], [0, 5, 11, 40, 57, 96], [0, 18, 21, 64, 66, 73]]
\item \{1=46368, 2=659232, 3=2978640, 4=3875760\} [[0, 1, 4, 7, 13, 19], [0, 2, 6, 79, 98, 109], [0, 3, 24, 43, 99, 105], [0, 5, 48, 62, 104, 119], [0, 10, 31, 55, 91, 93]]
\item \{1=47376, 2=613872, 3=2987712, 4=3911040\} [[0, 1, 4, 7, 13, 19], [0, 2, 3, 45, 91, 111], [0, 5, 21, 59, 61, 63], [0, 6, 25, 68, 90, 109], [0, 9, 27, 32, 49, 116]]
\item \{1=49392, 2=630504, 3=3060288, 4=3819816\} [[0, 1, 4, 7, 13, 19], [0, 2, 3, 23, 61, 101], [0, 5, 35, 60, 72, 115], [0, 15, 33, 65, 97, 122], [0, 17, 50, 59, 79, 108]]
\item \{1=49392, 2=690984, 3=2975616, 4=3844008\} [[0, 1, 4, 7, 13, 19], [0, 2, 3, 23, 87, 94], [0, 5, 33, 34, 73, 116], [0, 9, 32, 85, 100, 121], [0, 11, 86, 98, 105, 118]]
\item \{1=53424, 2=711396, 3=2986200, 4=3808980\} [[0, 1, 4, 7, 13, 19], [0, 2, 3, 26, 77, 112], [0, 5, 11, 29, 50, 89], [0, 6, 40, 48, 69, 116], [0, 9, 52, 68, 115, 118]]
\item \{0=1260, 1=28224, 2=579096, 3=3040128, 4=3911292\} [[0, 1, 4, 7, 13, 19], [0, 2, 3, 23, 45, 119], [0, 5, 39, 50, 66, 97], [0, 9, 82, 90, 102, 109], [0, 11, 21, 77, 83, 112]]
\item \{0=1260, 1=28224, 2=595728, 3=3009888, 4=3924900\} [[0, 1, 4, 7, 13, 19], [0, 2, 6, 72, 73, 79], [0, 3, 26, 48, 65, 117], [0, 5, 29, 43, 104, 113], [0, 17, 56, 61, 68, 119]]
\item \{0=1260, 1=30240, 2=644112, 3=3034080, 4=3850308\} [[0, 1, 4, 7, 13, 19], [0, 2, 6, 43, 98, 115], [0, 3, 44, 47, 106, 108], [0, 5, 21, 31, 66, 83], [0, 10, 40, 85, 100, 118]]
\item \{0=1260, 1=31248, 2=572292, 3=2987208, 4=3967992\} [[0, 1, 4, 7, 13, 19], [0, 2, 6, 37, 54, 77], [0, 3, 22, 52, 76, 99], [0, 5, 25, 96, 114, 119], [0, 11, 48, 59, 69, 89]]
\item \{0=1260, 1=31248, 2=597996, 3=3050712, 4=3878784\} [[0, 1, 4, 7, 13, 19], [0, 2, 6, 31, 61, 97], [0, 3, 16, 37, 44, 90], [0, 5, 25, 36, 65, 110], [0, 11, 69, 83, 85, 93]]
\item \{0=1260, 1=32256, 2=623700, 3=2984184, 4=3918600\} [[0, 1, 4, 7, 13, 19], [0, 2, 6, 34, 55, 63], [0, 3, 14, 68, 99, 103], [0, 9, 33, 37, 76, 123], [0, 16, 70, 71, 116, 117]]
\item \{0=1260, 1=34272, 2=599508, 3=2996280, 4=3928680\} [[0, 1, 4, 7, 13, 19], [0, 2, 6, 30, 58, 121], [0, 3, 37, 72, 73, 85], [0, 5, 44, 45, 48, 105], [0, 10, 26, 56, 88, 117]]
\item \{0=1260, 1=34272, 2=612360, 3=3028032, 4=3884076\} [[0, 1, 4, 7, 13, 19], [0, 2, 3, 23, 25, 63], [0, 6, 24, 80, 116, 125], [0, 15, 61, 86, 99, 104], [0, 16, 68, 91, 109, 110]]
\item \{0=1260, 1=36288, 2=586656, 3=2985696, 4=3950100\} [[0, 1, 4, 7, 13, 19], [0, 2, 6, 30, 64, 103], [0, 3, 36, 46, 88, 118], [0, 5, 66, 80, 82, 111], [0, 21, 27, 39, 68, 94]]
\item \{0=1260, 1=37296, 2=607824, 3=2967552, 4=3946068\} [[0, 1, 4, 7, 13, 19], [0, 2, 3, 28, 88, 116], [0, 6, 46, 68, 75, 117], [0, 9, 54, 64, 69, 119], [0, 11, 23, 62, 111, 125]]
\item \{0=1260, 1=41328, 2=644868, 3=3050712, 4=3821832\} [[0, 1, 4, 7, 13, 19], [0, 2, 6, 17, 76, 96], [0, 5, 23, 82, 111, 116], [0, 9, 39, 73, 89, 109], [0, 16, 20, 35, 90, 94]]
\item \{0=1260, 1=42336, 2=632016, 3=2964528, 4=3919860\} [[0, 1, 4, 7, 13, 19], [0, 2, 6, 40, 100, 109], [0, 3, 38, 75, 79, 113], [0, 5, 12, 70, 82, 104], [0, 20, 33, 68, 84, 95]]
\item \{0=1260, 1=43344, 2=640332, 3=2990232, 4=3884832\} [[0, 1, 4, 7, 13, 19], [0, 2, 3, 23, 30, 80], [0, 5, 49, 76, 81, 90], [0, 9, 25, 56, 68, 103], [0, 27, 47, 64, 66, 88]]
\item \{0=1260, 1=47376, 2=604044, 3=3032568, 4=3874752\} [[0, 1, 4, 7, 13, 19], [0, 2, 6, 17, 37, 102], [0, 5, 61, 72, 83, 109], [0, 9, 11, 47, 77, 111], [0, 15, 24, 33, 94, 95]]
\item \{0=1260, 1=48384, 2=605556, 3=3056760, 4=3848040\} [[0, 1, 4, 7, 13, 19], [0, 2, 6, 17, 104, 107], [0, 5, 20, 26, 51, 119], [0, 9, 57, 68, 89, 108], [0, 10, 41, 70, 90, 120]]
\item \{0=2520, 1=39312, 2=599508, 3=3033576, 4=3885084\} [[0, 1, 4, 7, 13, 19], [0, 2, 3, 58, 87, 121], [0, 5, 45, 56, 60, 77], [0, 9, 26, 51, 59, 118], [0, 12, 25, 32, 67, 104]]
\end{enumerate}
\end{example}

\end{document}